\documentclass[10pt,letterpaper]{article}

\usepackage[margin=0.8in]{geometry}
\usepackage{lmodern}
\usepackage{authblk}
\usepackage[expansion=false]{microtype}

\usepackage{amsmath,amssymb,amsthm,mathtools}
\usepackage{booktabs,array}
\usepackage{upquote}
\newcolumntype{L}[1]{>{\raggedright\arraybackslash}p{#1}}

\usepackage{enumitem,graphicx,xcolor}
\usepackage{float}

\definecolor{figureInk}{HTML}{30323D}
\definecolor{figureBridge}{HTML}{7226BA}
\definecolor{figurePath}{HTML}{5C80BC}
\definecolor{figureExtra}{HTML}{AD5D4E}
\definecolor{figureCream}{HTML}{FFE9CE}
\definecolor{figureNavy}{HTML}{004E89}
\definecolor{figureGreen}{HTML}{69B578}
\definecolor{figureRed}{HTML}{A4243B}

\colorlet{figureTwinA}{figureNavy}
\colorlet{figureTwinB}{figureGreen}
\colorlet{faceSimple}{figurePath!18}
\colorlet{faceRepeated}{figureCream}
\colorlet{faceOther}{figureInk!8}
\colorlet{branchA}{figureNavy}
\colorlet{branchB}{figureGreen}
\colorlet{branchC}{figureExtra}
\colorlet{branchD}{figureRed}
\colorlet{branchE}{figureBridge}
\usepackage{tikz}
\usetikzlibrary{arrows.meta,shapes.geometric}

\newcommand{\FigScaleAll}{0.8}
\newcommand{\scaledfigure}[2]{%
  \scalebox{\FigScaleAll}{\scalebox{#1}{#2}}%
}
\newcommand{\scaledtable}[2]{\scaledfigure{#1}{#2}}

\setlist[enumerate]{itemsep=1pt,topsep=3pt}
\setlist[itemize]{itemsep=1pt,topsep=3pt}
\newtheorem{theorem}{Theorem}[section]
\newtheorem{lemma}[theorem]{Lemma}
\newtheorem{proposition}[theorem]{Proposition}
\newtheorem{corollary}[theorem]{Corollary}
\theoremstyle{definition}
\newtheorem{construction}[theorem]{Construction}

\newcommand{\supp}[1]{\textup{[Supp.~#1]}}

\usetikzlibrary{external}

\tikzexternaldisable
\tikzifexternalizing{\PassOptionsToPackage{bookmarks=false}{hyperref}}{}

\usepackage[hidelinks]{hyperref}
\hypersetup{
  pdftitle={Toroidal Obstructions with No K5 Minor},
  pdfauthor={Austin Ulrigg},
  pdfsubject={Toroidal graphs and excluded minors},
  pdfkeywords={Toroidal graph, graph embedding, forbidden minor, graph genus}
}

\title{Toroidal Obstructions with No $K_5$ Minor}
\author{Austin Ulrigg}
\affil{University of Washington, \texttt{austinul@uw.edu}}
\date{}

\begin{document}
\maketitle

\begin{abstract}
Finding the complete set of forbidden minors for the torus remains an open problem.
We determine the complete list of toroidal minor obstructions with no
$K_5$ minor and prove that there are exactly $74$, each with between
$10$ and $14$ vertices. Combined with the known classification of toroidal
obstructions with no $K_{3,3}$ minor, this determines every toroidal
obstruction that excludes one of the forbidden minors
for planarity. The proof uses three-vertex cuts to show that,
apart from $K_{3,7}$, every three-connected obstruction in this class
has two distinct independent triples, each the common neighborhood of
at least two degree-three vertices. We consider the three possible sizes of the intersection of the two triples. In each case, planarity tests reduce the proof to a finite enumeration. From the $74$ minor obstructions we also obtain a complete list of $381$ topological obstructions for the torus with no $K_5$ minor.


\end{abstract}

\section{Introduction}\label{sec:introduction}
In 1930, Kuratowski proved that a graph can be drawn in the
plane without any edge crossings if and only if it contains no
subdivision of the complete graph $K_5$ or the complete bipartite graph $K_{3,3}$~\cite{Kuratowski1930}. In 1937, Wagner showed that the same two graphs are exactly the complete set of forbidden minors for the plane~\cite{Wagner1937}. Both graphs can, however, be embedded on the torus. A natural next question is: Which graphs that cannot be embedded in the plane can still be embedded in the torus?

A graph is \emph{toroidal} if it can be embedded on the torus. A
\emph{minor obstruction}, or forbidden minor, for a surface is a graph that cannot be embedded on that surface while every proper minor of it can. A toroidal obstruction is a minor obstruction for the torus. Robertson and Seymour proved that every fixed surface has only finitely many
minor obstructions~\cite{RobertsonSeymour1990}. Outside of planar graphs, the complete set of forbidden minors is also known for the projective plane. The work of Glover, Huneke and Wang, together with Archdeacon's completeness proof, gives a list of $35$ forbidden minors~\cite{GloverHunekeWang1979,Archdeacon1981}. Despite the results of Robertson and Seymour, outside of the plane and projective plane there is little known about the complete set of forbidden minors for other surfaces.

Myrvold and Woodcock collected $17{,}535$ minor obstructions for the
torus~\cite{MyrvoldWoodcock2018}. Of the $74$ graphs in
Theorem~\ref{thm:classification}, only $\mathfrak T_{74}$ is absent
from their collection~\supp{O4}. We also give $120$ graphs absent
from their larger collection of topological obstructions. Of these,
$22$ are minor obstructions, including $\mathfrak T_{74}$, and the
other $98$ are topological obstructions only~\supp{N1}. All $120$ are
available in the author's database.\footnote{\url{https://austinulrigg.com/toroidal.html}, release of November 11, 2026.}

A \emph{block} of a graph is a maximal connected subgraph with no cut vertex. Battle, Harary, Kodama and
Youngs proved that the orientable genus of a graph is additive over its blocks~\cite{BattleHararyKodamaYoungs1962}. Thus, there is exactly one disconnected toroidal obstruction with no $K_5$ minor and exactly one with a cut vertex, namely $\mathfrak{T}_{25}$ and $\mathfrak{T}_{8}$. Mohar and \v{S}koda classified all $68$ toroidal obstructions of
vertex-connectivity two~\cite[Theorem~8.10 and Corollary~8.11]{MoharSkoda2014}. A plausible approach to completing the classification of the forbidden toroidal minors could be a piecemeal approach as suggested in~\cite{MyrvoldWoodcock2018}.

Recently, Gagarin and Kocay studied
the structure of graphs containing
a $K_5$ subdivision but no $K_{3,3}$ subdivision~\cite{GagarinKocay2002} and later, Gagarin, Myrvold and Chambers proved that there are exactly four toroidal minor obstructions with no $K_{3,3}$ minor, and eleven topological obstructions~\cite{GagarinMyrvoldChambers2009}. Since $K_5$ is the only three-connected nonplanar graph with no $K_{3,3}$ minor~\cite{Hall1943} and $K_5$ is toroidal, each of the four minor obstructions has connectivity at most two. Because $K_{3,3}$ is cubic, excluding a $K_{3,3}$ minor is equivalent to excluding a $K_{3,3}$ subdivision~\cite[Section~1]{GagarinMyrvoldChambers2009}. In this paper, we classify the toroidal obstructions with no $K_5$ minor, which instead form a proper subclass of those with no $K_5$ subdivision. 

\begin{theorem}\label{thm:classification}
A graph that has no $K_5$ minor is a toroidal minor obstruction if and only if it is isomorphic to one of $\mathfrak{T}_{1},\ldots,\mathfrak{T}_{74}$, given in Table~\ref{tab:graph6}.
\end{theorem}

Combining this with the result of Gagarin, Myrvold and Chambers
determines every toroidal obstruction that excludes at least one of
the two forbidden minors for planarity. Any obstruction outside these
two classifications must therefore contain both a $K_5$ minor and a
$K_{3,3}$ minor. 

As noted, Mohar and \v{S}koda's results classify the toroidal obstructions of connectivity two, including those without a $K_5$ minor. Additionally, any nonplanar four-connected graph contains a $K_5$ minor
by Wagner's theorem~\cite{Wagner1937}. Hence, it remains only to consider the three-connected obstructions which is the goal of this paper. For the three-connected case, we use
Halin's theorem~\cite{Halin1967,Halin1968} to obtain three-cuts and study
the components they separate. Apart from $K_{3,7}$, we show every three-connected obstruction with no $K_5$ minor has two distinct
independent triples, each the common neighborhood of at least two
degree-three vertices (Theorem~\ref{thm:two-pages}). We call such a triple a twin triple. Section~\ref{sec:core} proves this theorem and develops lemmas on bridges used
throughout. The classification then splits according to whether the
triples have no common vertex, two common vertices or one, treated in
Sections~\ref{sec:disjoint}, \ref{sec:two-common}
and~\ref{sec:one-common}. Each of Sections~\ref{sec:disjoint}, \ref{sec:two-common} and~\ref{sec:one-common} follows the same steps. We first find the faces of a subgraph containing the two triples and four of the degree-three vertices, and reduce the toroidality of $G$ to planarity tests on the bridges of this subgraph. By minimality, no bridge can then be reduced while its nonplanar tests stay nonplanar (Lemma~\ref{core:relevant}). We finally retain a $K_{3,3}$ subdivision in a nonplanar test, which bounds the size of the bridge, and enumerate the possible forms (Corollary~\ref{common:reduced-forms}). The three sections are independent of one another, with a few small exceptions. Section~\ref{sec:computations} assembles the
proof of Theorem~\ref{thm:classification}. A \emph{topological obstruction} for the torus is a nontoroidal graph
of minimum degree at least three in which every edge deletion is
toroidal. In Section~\ref{sec:topological} we use
Theorem~\ref{thm:classification} to reduce the topological
obstructions with no $K_5$ minor to a search by repeated vertex splitting starting from
$\mathfrak{T}_{1},\ldots,\mathfrak{T}_{74}$, and find that there are
exactly $381$ (Theorem~\ref{thm:topological}). Appendix~\ref{app:computations} lists the supplemental data and describes how each is checked. Table~\ref{tab:cases} lists the cases and the graphs obtained in each.

\ifdefined\thetable\else
  \documentclass[10pt,letterpaper]{article}
  \usepackage[margin=.8in]{geometry}
  \usepackage{lmodern,amssymb,array,booktabs,float}
  \usepackage[draft]{hyperref}
  \newcolumntype{L}[1]{>{\raggedright\arraybackslash}p{#1}}
  \newcommand{\previewlabel}[2]{\newlabel{#1}{{#2}{}{}{}{}}}
  \previewlabel{sec:completion}{7}
  \previewlabel{thm:two-pages}{3.6}
  \previewlabel{thm:disjoint}{4.1}
  \previewlabel{dis:endpoint-only}{4.22}
  \previewlabel{thm:two-common}{5.1}
  \previewlabel{r2:component-count}{5.28}
  \previewlabel{r2:chain-outcome}{5.33}
  \previewlabel{thm:one-common}{6.1}
  \previewlabel{r1:five-outcome}{6.8}
  \previewlabel{r1:missing-a-outcome}{6.13}
  \previewlabel{r1:missing-b-outcome}{6.16}
  \previewlabel{r1:three-outcome}{6.25}
  \bibcite{MoharSkoda2014}{17}
  \newcommand{\finishcasepreview}{\end{document}}
  \pagestyle{empty}
  \begin{document}
\fi

\begin{table}[H]
\centering
\normalsize
\setlength{\tabcolsep}{4pt}
\renewcommand{\arraystretch}{1.03}
\newcommand{\caselink}[2]{\hyperref[#1]{#2}}
\newcommand{\result}[2]{\caselink{#1}{#2~\ref*{#1}}}
\newcommand{\subcase}[4]{\caselink{#1}{#2} & \result{#1}{#3} & #4\\}
\newcommand{\detail}[4]{\leftskip=1em\relax\subcase{#1}{#2}{#3}{#4}}
\begin{tabular}{@{}L{.18\linewidth}L{.36\linewidth}L{.18\linewidth}L{\dimexpr.28\linewidth-6\tabcolsep\relax}@{}}
\toprule
Case & Subcase & Settled by & $i$ with $G\cong\mathfrak T_i$\\
\midrule
\multicolumn{2}{@{}l}{\caselink{sec:completion}{Not two-connected}} &
\result{sec:completion}{Section} & 8, 25\\
\multicolumn{2}{@{}l}{Connectivity two} &
Mohar and \v{S}koda~\cite{MoharSkoda2014} &
3, 6, 12, 13, 27, 28, 31, 36, 42, 44, 46, 47, 50, 60, 62, 64, 65, 71\\
\multicolumn{2}{@{}l}{\caselink{thm:two-pages}{Three-connected, $K_{3,7}$}} &
\result{thm:two-pages}{Theorem} & 1\\
\midrule
\caselink{thm:disjoint}{Disjoint twin triples} &
\subcase{dis:endpoint-only}{Three disjoint edges join the two triples}{Theorem}{5, 48, 58, 68, 69}
& \subcase{thm:disjoint}{No two disjoint twin triples are joined by three disjoint edges}{Theorem}{4, 21, 32, 38, 43, 45, 61, 63, 66, 67}
\addlinespace
\caselink{thm:two-common}{Two common vertices} &
\subcase{thm:two-common}{A component outside the two triples is adjacent to all four vertices in their union}{Theorem}{2, 11, 17, 18, 19, 22, 24, 26, 30, 34, 39, 40, 51, 52, 53, 56, 57, 72, 73, 74}
& \textit{Each component has three neighbors in the triples} & &\\
& \detail{r2:chain-outcome}{A component contains an edge}{Proposition}{15}
& \detail{r2:component-count}{There are only single-vertex components}{Lemma}{9, 10, 16}
\addlinespace
\caselink{thm:one-common}{One common vertex} &
\subcase{r1:five-outcome}{A component outside the two triples is adjacent to all five vertices in their union}{Proposition}{20, 23}
& \textit{Some component has four neighbors in the triples and none has five} & &\\
& \detail{r1:missing-a-outcome}{A component is adjacent to the four vertices other than the common vertex}{Proposition}{29, 54, 55}
& \detail{r1:missing-b-outcome}{Every component with four neighbors is adjacent to the common vertex}{Proposition}{14, 37, 41, 59}
& \subcase{r1:three-outcome}{Every component has three neighbors in the two triples}{Proposition}{7, 33, 35, 49, 70}
\bottomrule
\end{tabular}
\caption{The cases of the classification and the graphs obtained in each.}
\label{tab:cases}
\end{table}

\ifdefined\finishcasepreview\expandafter\finishcasepreview\fi

\section{Preliminaries and basic reductions}\label{sec:preliminaries}
For the convenience of the reader we provide the following preliminary results and definitions. When possible, we follow the terminology of \emph{Graphs on Surfaces} by
Mohar and Thomassen~\cite{MoharThomassen2001}. Graphs are finite,
simple, and undirected, unless stated otherwise. We use $\cup$ for graph
unions, $+$ for adding specified vertices or edges, and $K_n[X]$ for the complete graph on a set $X$ of $n$ vertices. We write $N_G(X)$ for the vertices outside $X$ with a neighbor in $X$. By Youngs's theorem~\cite{Youngs1963}, every embedding of a connected nonplanar graph on the torus is cellular. In this paper, an \emph{obstruction} will mean a toroidal minor obstruction with no
$K_5$ minor. Unless stated otherwise, $G$ denotes a three-connected
obstruction. By minor-minimality, any nontoroidal minor of $G$ is isomorphic to $G$.

For a subgraph $K\subseteq G$, a \emph{$K$-bridge} $B$ is either a
component of $G-V(K)$ together with its neighbors in $K$ and all
edges joining that component to $K$, or an edge of
$E(G)\setminus E(K)$ with both ends in $K$. Its vertices in $K$
form its \emph{attachment set} $A$, and
$\operatorname{Int}(B)=V(B)\setminus A$ is its \emph{interior} vertex set.
When $K$ has no edges, we call its bridges \emph{bridges on $V(K)$}. When we write $B_X$ for a bridge on $X$ with nonempty interior, we specify its component separately. Auxiliary graphs are denoted by $H$, with subscripts when needed. The \textit{exterior} of a $K$-bridge $B$ is obtained from $G$ by deleting the interior vertices of $B$ and all edges of $B$. Every edge of $G$ at an interior vertex of $B$ belongs to $B$, and adding edges between attachments creates no new edge there. Thus we note that, for any vertex set $X$ and any set $E$ of edges between attachments, a set of interior vertices of $B$ that forms a component of one of $G-X$, $B-X$ and $(B+E)-X$ forms a component of all three.

Let $K$ be a connected graph embedded in an orientable surface
$\Sigma$, and let $\mathcal N\subseteq\Sigma$ be a closed regular neighborhood
of $K$. Then $\mathcal N$ is a compact connected orientable surface with
boundary, $K\subseteq\operatorname{int}\mathcal N$. The genus of $\mathcal N$ is that of the
closed surface obtained by gluing a disk to each boundary
component of $\mathcal N$.

\begin{lemma}\label{lem:torus-complement}
Let $\mathcal N$ be a closed regular neighborhood of a connected graph embedded
on the torus, and let $M$ be the closure of its complement.
If $\mathcal N$ has genus one, every component of $M$ is a disk.
If $\mathcal N$ has genus zero, $M$ consists of disks and exactly one annulus
or one torus with one open disk removed.
\end{lemma}\begin{proof}
Let $g,b$ be the genus and number of boundary components of $\mathcal N$,
and let $M_1,\ldots,M_k$ be the components of $M$, with genera $h_i$
and boundary counts $b_i$. Every boundary circle of $\mathcal N$
borders one of the $M_i$, so $\sum_i b_i=b$. Since each $b_i\ge1$,
we have $b\ge k$. The Euler characteristic of a connected orientable surface of genus $h$
with $b$ boundary components is $2-2h-b$. Gluing along the boundary
circles gives the torus, so additivity of Euler~characteristic~gives
\[
0=\chi(\mathcal N)+\sum_{i=1}^k\chi(M_i)
  =(2-2g-b)+\left(2k-2\sum_{i=1}^k h_i-b\right).
\]
Thus $b=1-g+k-\sum_i h_i$. If $g=1$, then $b=k-\sum_i h_i$. If any $h_i>0$, then $b<k$,
a contradiction. Thus every $h_i=0$ and $b=k$, so every $b_i=1$.
Every component of $M$ is therefore a disk. If $g=0$, then $b=1+k-\sum_i h_i$. Since $b\ge k$, we have
$\sum_i h_i\le1$. If one $h_i=1$, then $b=k$ and every $b_i=1$, and
we obtain one torus with an open disk removed, and all other components
are disks. Otherwise every $h_i=0$ and $b=k+1$, so exactly one $b_i=2$
and all the others are one, and we obtain one annulus, and all other
components are disks.
\end{proof}

\subsection{Extending an embedding}\label{sec:extension}

Let $H$ be a graph, $K\subseteq H$ a connected subgraph, and
$\Pi$ an embedding of $K$ on the torus. We say that $\Pi$ extends
to an embedding of $H$ if $H$ has an embedding that agrees with
$\Pi$ on every vertex and edge of $K$. \textit{Assigning} a
$K$-bridge $B$ to a face $F$ of $\Pi$ means choosing that face
as the region in which $B$ is to be drawn outside of its
attachments. Every attachment of $B$ must occur on a
boundary walk of $F$.

We refer to the cyclic list of attachment appearances along a facial walk as its \emph{attachment sequence}, written for example as $abcad$. A vertical bar is used to separate different boundary components, as in $ac\mid ad$.
In an auxiliary drawing with separated boundary appearances,
$a_1,a_2$ denote appearances of the same original vertex $a$.

Suppose $\Pi$ extends to an embedding of $H$. For each $K$-bridge
$B$, remove only the points representing its attachments
from its drawing. The remaining set is connected and disjoint
from $\Pi(K)$, so it lies in one face of $\Pi$. Thus, deciding
whether $\Pi$ extends to an embedding of $H$ amounts to deciding
whether all its bridges can be assigned to faces and drawn
all together at once in their chosen faces without crossings.

The next lemma gives a planarity test for a single bridge in a disk face.

\begin{lemma}
\label{cor:boundary-completions}\label{dis:replacement}
Let $B$ be a bridge with nonempty connected interior.
Suppose its attachments appear around the boundary of a disk
face $F$ in one of the cyclic sequences
\postdisplaypenalty=10000
\[
a,\quad ab,\quad abc,\quad abac,\quad abcd,\quad
abcad,\quad abacad,\quad abcdcb,
\]
up to relabeling and reversal.
Then $B$ can be drawn in $F$ if and only if the graph obtained
from $B$ by joining each pair of distinct consecutive
attachments in the sequence is planar.
\end{lemma}

\begin{proof}
For the sequence $a$, a plane drawing of $B$ with $a$ on the outer face
can be placed in $F$, and a drawing in $F$ is planar. For the other
sequences, let $H_0$ be the graph consisting of the attachments and the
edges joining consecutive attachments. It has a plane embedding with
the listed facial walk, and every facial walk containing all its vertices
has the same cyclic order, up to reversal. If $B\cup H_0$ is planar, the connected interior of $B$ lies
in one face of $H_0$. Its boundary contains every attachment,
so its boundary walk gives the required drawing in $F$.
Conversely, inserting a drawing of $B$ into the corresponding
face of $H_0$ gives a planar drawing of $B\cup H_0$.
\end{proof}

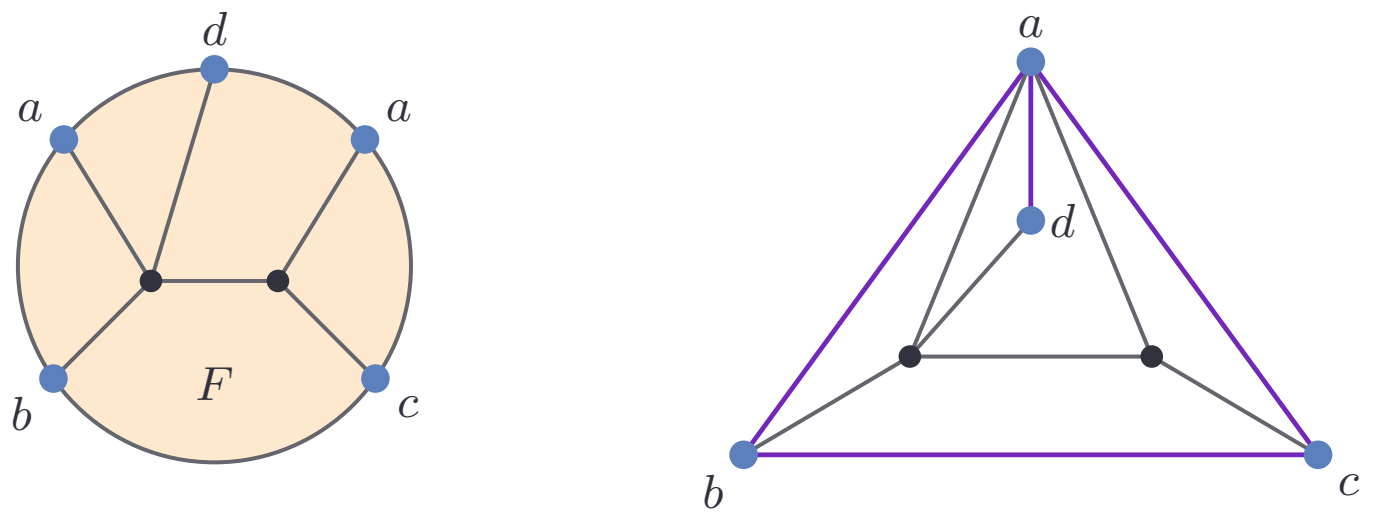
\begin{figure}[htbp]
\centering
\tikzsetnextfilename{boundary-completion-example}%
\begin{tikzpicture}[
  font=\small,text=figureInk,
  vertex/.style={circle,fill=figureInk,inner sep=1.5pt},
  root/.style={circle,fill=figurePath,inner sep=1.9pt},
  edge/.style={draw=figureInk!75,line width=.75pt,
    line cap=round,line join=round},
  added/.style={edge,draw=figureBridge,line width=.9pt},
  every label/.style={inner sep=1.5pt}
]
\begin{scope}
  \fill[faceRepeated] (0,0) circle (1.3);
  \draw[edge] (0,0) circle (1.3);
  \coordinate (a1) at (140:1.3);
  \coordinate (b) at (215:1.3);
  \coordinate (c) at (325:1.3);
  \coordinate (a2) at (40:1.3);
  \coordinate (d) at (90:1.3);
  \coordinate (w1) at (-.42,-.1);
  \coordinate (w2) at (.42,-.1);
  \draw[edge] (w1)--(a1) (w1)--(b) (w1)--(d) (w2)--(c) (w2)--(a2) (w1)--(w2);
  \node[vertex] at (w1) {};
  \node[vertex] at (w2) {};
  \node[root,label=140:{$a$}] at (a1) {};
  \node[root,label=215:{$b$}] at (b) {};
  \node[root,label=325:{$c$}] at (c) {};
  \node[root,label=40:{$a$}] at (a2) {};
  \node[root,label=90:{$d$}] at (d) {};
  \node at (0,-.78) {$F$};
\end{scope}
\begin{scope}[xshift=5.4cm,yshift=-.25cm]
  \coordinate (a) at (0,1.6);
  \coordinate (b) at (-1.9,-1.0);
  \coordinate (c) at (1.9,-1.0);
  \coordinate (d) at (0,.55);
  \coordinate (w1) at (-.8,-.35);
  \coordinate (w2) at (.8,-.35);
  \draw[added] (a)--(b)--(c)--cycle (a)--(d);
  \draw[edge] (w1)--(a) (w1)--(b) (w1)--(d) (w2)--(c) (w2)--(a) (w1)--(w2);
  \node[vertex] at (w1) {};
  \node[vertex] at (w2) {};
  \node[root,label=above:{$a$}] at (a) {};
  \node[root,label=below left:{$b$}] at (b) {};
  \node[root,label=below right:{$c$}] at (c) {};
  \node[root,label={[label distance=-1pt]right:{$d$}}] at (d) {};
\end{scope}
\end{tikzpicture}
\caption{Left, a bridge in a disk face $F$ with attachment sequence $abcad$. Right, the same bridge with the consecutive attachments joined.}
\label{fig:boundary-completion-example}
\end{figure}

Figure~\ref{fig:boundary-completion-example} shows the sequence $abcad$.

Let $\mathcal N$ be a closed regular neighborhood of $\Pi(K)$. For a graph $H$ containing a cycle $C$, let
$\operatorname{Aux}(H,C)$ be $H$ plus a new vertex adjacent precisely
to the vertices of $C$. We consider the
cases where $\mathcal N$ has genus one and genus zero separately. Suppose first that $\mathcal N$ has genus one. By Lemma~\ref{lem:torus-complement}, every face of $\Pi$ is a disk. The following construction adapts the split auxiliary graph of
Mohar and Thomassen~\cite[p.~191]{MoharThomassen2001} to record
the specific chosen boundary appearances.

For each disk face $F$ chosen for at least one bridge $B$, we construct
graphs $H_F$, $C_F$ and $H_F^+$ as follows. For each edge of a bridge $B$ assigned to $F$, and each endpoint
$v$ of that edge belonging to $K$, choose an appearance of $v$
on the boundary walk of $F$. List the distinct chosen appearances
in cyclic order and represent them by distinct new vertices
$x_1,\dots,x_k$. Define $V(H_F)$ to be $\{x_1,\dots,x_k\}$
together with all the interior vertices of the bridges assigned
to $F$.
Every edge of a bridge assigned to $F$ defines one edge of $H_F$.
An edge $uw$ between interior vertices is kept. An edge from an
interior vertex $u$ to an attachment $v$ becomes $ux_i$, where
$x_i$ represents the appearance chosen at $v$. An edge between
attachments becomes $x_ix_j$, where $x_i$ and $x_j$ represent
the appearances chosen for that edge at its two endpoints. Set $m=\max\{k,3\}$. If $k<3$, introduce new vertices
$x_{k+1},\dots,x_m$, distinct from all vertices of $H_F$.
Define $C_F$ to be the cycle $x_1x_2\cdots x_mx_1$.
The vertices $x_1,\dots,x_k$ follow the chosen cyclic boundary
order. Every edge of $C_F$ is a new edge. Set
$H_F^+=\operatorname{Aux}(H_F\cup C_F,C_F)$, retaining any parallel edges.

\begin{lemma}[Adding bridges in a face homeomorphic to a disk]
\label{lem:embedding-extension}
With the notation above, the bridges assigned to a disk face
$F$ can be added to $\Pi$ inside $F$, through their chosen
boundary appearances, if and only if $H_F^+$ is planar.
\end{lemma}
\begin{proof}
By Mohar's disk criterion~\cite[Section~3, first paragraph]{Mohar1994},
applied to $H_F\cup C_F$ with boundary cycle $C_F$, the bridges
assigned to a disk face $F$ can be added to $\Pi$ inside $F$
through their chosen boundary appearances if and only if $H_F^+$
is planar.
\end{proof}

Consequently, when $\mathcal N$ has genus one, $\Pi$ extends to an
embedding of $H$ if and only if each $K$-bridge can be assigned
to a face of $\Pi$ and its boundary appearances chosen so that
every graph $H_F^+$ is planar.

Now suppose that $\mathcal N$ has genus zero. By
Lemma~\ref{lem:torus-complement}, the faces of $\Pi$ are
disks together with exactly one annulus or exactly one torus
with one open disk removed. Suppose the case where some face is a torus with one open disk removed has been excluded, so that there is one annular
face $F$ and all other faces are disks. Represent
$F$ by the interior of a closed annulus $\Sigma_0$ whose boundary
components $\Gamma_1,\Gamma_2$ trace the two boundary walks
of $F$. Assign a collection of $K$-bridges in $H$ to $F$. For each
edge of an assigned bridge and each endpoint $v\in V(K)$,
choose an appearance of $v$ on one of these boundary walks,
with at most two distinct chosen appearances on each walk.
Define $H_{\mathrm{ann}}$ to have one new vertex for each distinct chosen
appearance, together with all the interior vertices of the
assigned bridges. Define its edges by the same endpoint
rules as $H_F$. Appearances of the same attachment
on different boundary walks are represented by distinct
vertices of $H_{\mathrm{ann}}$. For $i=1,2$, let $A_i\subseteq V(H_{\mathrm{ann}})$ consist of the vertices
representing chosen appearances on the boundary walk traced
by $\Gamma_i$. For each $v\in A_i$, let $p_v\in\Gamma_i$
be the boundary point representing that appearance. To define $H_{\mathrm{ann}}^+$, label the vertices of $A_1$ by $a,b$ and those of
$A_2$ by $c,d$, introducing a new vertex for each unused label.
The new vertices are distinct from one another and from
every vertex of $H_{\mathrm{ann}}$. Then $H_{\mathrm{ann}}^+$ has vertex set
$V(H_{\mathrm{ann}})\cup\{a,b,c,d\}$, the edges of $H_{\mathrm{ann}}$, and two new edges
$ab,cd$.

\begin{lemma}[Adding bridges in a face homeomorphic to an annulus]
\label{lem:annulus-embedding}
With the notation above, $H_{\mathrm{ann}}$ has an embedding in $\Sigma_0$ placing
each $v\in A_1\cup A_2$ at $p_v$, with all other vertices and
all edge interiors in $\operatorname{int}\Sigma_0$, if and only if
$H_{\mathrm{ann}}^+$ is planar.
\end{lemma}
\begin{proof}
If $H_{\mathrm{ann}}$ has the required annulus embedding, place any
auxiliary vertices on the corresponding boundary components, cap the
boundary components with disks, and add the edges $ab,cd$ through the
respective disks. This gives a sphere embedding of $H_{\mathrm{ann}}^+$. Conversely, suppose that $H_{\mathrm{ann}}^+$ has a sphere embedding.
Take this embedding and trace two disjoint closed disks around $ab$ and
$cd$, meeting the drawing precisely in these respective edges, with
the endpoints on the corresponding disk boundaries. Remove the disk
interiors, the extra edges, and any auxiliary vertices that were
introduced, and we have an annulus embedding of $H_{\mathrm{ann}}$. There are at most two marked vertices on each boundary component, so
their cyclic order has no further restriction. A homeomorphism
from this annulus to $\Sigma_0$ can therefore place every vertex of
$A_1\cup A_2$ at its prescribed point.
\end{proof}

Therefore, when we assign every $K$-bridge to the annular face or to a disk face,
choosing their boundary appearances, if these choices satisfy
the hypothesis of Lemma~\ref{lem:annulus-embedding},
the bridges can be drawn together in their assigned faces
through those appearances if and only if $H_{\mathrm{ann}}^+$ and every $H_F^+$
for a chosen disk face are planar. As drawings in distinct faces
have disjoint interiors, together with $\Pi(K)$ they give
an embedding of $H$. 

\subsection{Graphs with no \texorpdfstring{$K_5$}{K5} minor}

The \emph{Wagner graph} $V_8$ is obtained from an $8$-cycle by joining
each pair of opposite vertices. It is triangle-free, and toroidal. A \emph{three-cut} is a set of three vertices whose deletion
disconnects the graph. Halin's theorem implies that every three-connected nonplanar graph
with no $K_5$ minor is isomorphic to $V_8$ or has a three-cut
whose deletion leaves at least three components
\cite{Halin1967,Halin1968}, as stated in~\cite[Theorem~7]{KimLeeOum2016}. In particular, such a graph has a three-cut of this kind whenever it contains a triangle.

\subsection{Planarity lemmas}\label{sec:planarity-lemmas}

The next lemma gives conditions on individual bridges
that determine whether they can be drawn together.

\begin{lemma}\label{core:assembly-general}
Let $H$ be a graph and let $A\subseteq V(H)$ have three vertices. Consider the
bridges on $A$.
\begin{enumerate}[label=\textup{(\roman*)}]
\item Three bridges adjacent to all of $A$ make $H$ nonplanar.
\item If exactly two bridges $B_1,B_2$ are adjacent to all of $A$,
then $H$ is planar if and only if $B_1\cup K_3[A]$ and $B_2\cup K_3[A]$ are
planar, every other bridge $B$ with two attachments $\alpha,\beta$
has $B+\alpha\beta$ planar, and every bridge with at most one
attachment is planar.
\item If exactly one bridge $B_1$ is adjacent to all of $A$, let $E^*$
consist of the edges $\alpha\beta$ for which another bridge has
attachment set $\{\alpha,\beta\}$. Then $H$ is planar if and only if
$B_1+E^*$ is planar and the other bridges satisfy the last two
conditions in (ii).
\end{enumerate}
\end{lemma}
\begin{proof}
For (i), contracting the interiors of three bridges adjacent to all of
$A$ gives a $K_{3,3}$ minor. Hence $H$ is~\mbox{not planar}.

For the forward direction of (ii), suppose that $H$ is planar. With
$\{i,j\}=\{1,2\}$, delete all bridges besides $B_i,B_j$, contract the
interior of $B_i$ to $w$, and then replace the three edges at $w$ with
a triangle on $A$. This gives a plane drawing of $B_j\cup K_3[A]$
for both choices of $j$.

If a bridge with at most one attachment were nonplanar, then $H$
would be nonplanar. Such planar bridges can be removed and reinserted
without affecting planarity. If a bridge $B$ with two attachments
$\alpha,\beta$ has $B+\alpha\beta$ nonplanar, then an
$\alpha$--$\beta$ path through the interior of $B_1$ shows that $H$
has $B+\alpha\beta$ as a minor and is also nonplanar. This applies to
(iii) as well.

For (iii), suppose the other bridges satisfy the stated conditions.
Remove the planar bridges with at most one attachment. By the planar
bridge replacement theorem~\cite[Theorem~3.1]{MyrvoldWoodcock2018},
each nontrivial two-attachment bridge can be replaced by the edge
joining its attachments without changing planarity. The resulting
graph is $B_1+E^*$. Since these replacements and the removal of the
other planar bridges do not change planarity, this proves both
directions of (iii).

For the reverse direction of (ii), suppose all the stated conditions
hold. Since $B_1\cup K_3[A]$ and $B_2\cup K_3[A]$ are planar,
Lemma~\ref{cor:boundary-completions}
gives disk drawings of $B_1,B_2$ with $A$ on their boundaries. Gluing
them along their boundary triangles gives a sphere drawing of
$B_1\cup B_2\cup K_3[A]$. Delete the triangle edges not arising from
two-attachment bridges and reverse the replacements above. Reinsert
the planar bridges with at most one attachment. This gives a plane
drawing of $H$.
\end{proof}

In a subdivision of $K_{3,3}$, we call the six degree-three vertices its \emph{corner vertices} and the nine paths corresponding to the edges of $K_{3,3}$ its \emph{branch paths}. The degree-two vertices on these paths are its \emph{interior vertices}. For vertices $u,v$ on a path $P$, write $P[u,v]$ for the subpath with ends $u,v$.

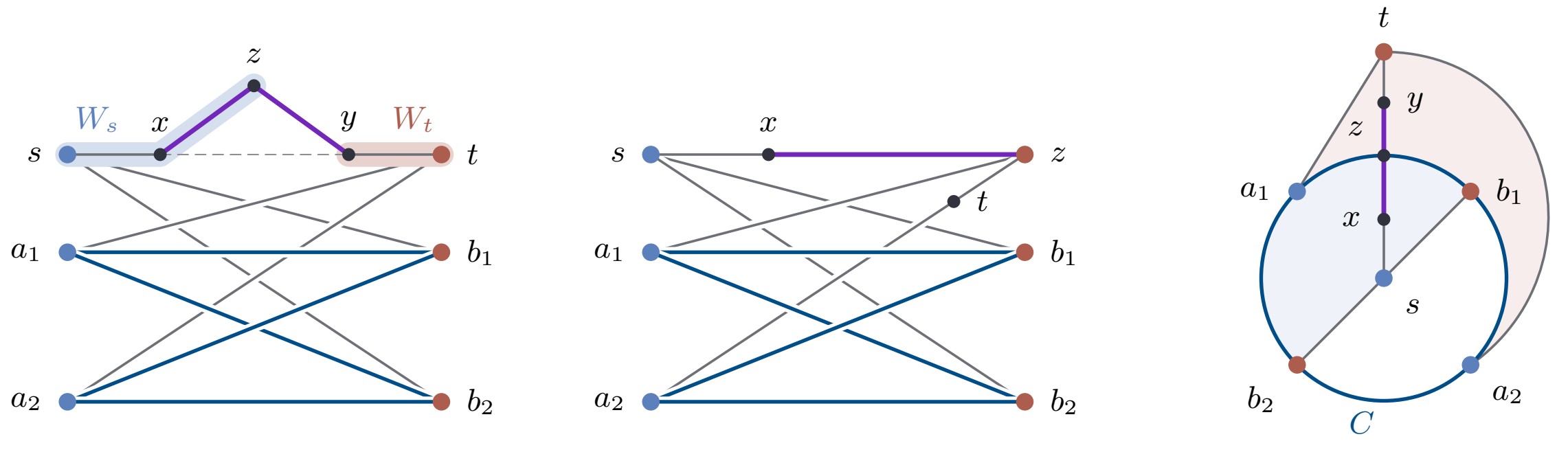
\begin{figure}[htbp]
\centering
\tikzsetnextfilename{position-of-third}%
\scaledfigure{1.2}{%
\begin{tikzpicture}[
  font=\small, line cap=round, line join=round,
  fixed/.style={draw=figureInk!70, line width=.7pt},
  cycle path/.style={draw=figureNavy, line width=1.05pt},
  added/.style={draw=figureBridge, line width=1.35pt},
  over/.style={preaction={draw=white, line width=3pt}},
  vertex/.style={circle, inner sep=0pt, minimum size=3.7pt, fill=figureInk},
  left corner/.style={vertex, minimum size=5pt, fill=figurePath},
  right corner/.style={vertex, minimum size=5pt, fill=figureExtra},
  label distance=1pt
]
\newcommand{\cornerpaths}{%
  \coordinate (s)  at (0,1.24);
  \coordinate (a1) at (0,.26);
  \coordinate (a2) at (0,-1.24);
  \coordinate (r)  at (3.75,1.24);
  \coordinate (b1) at (3.75,.26);
  \coordinate (b2) at (3.75,-1.24);
  \foreach \u/\v in {s/b1,s/b2,a1/r,a2/r}
    \draw[fixed,over] (\u)--(\v);
  \foreach \u/\v in {a1/b1,a1/b2,a2/b1,a2/b2}
    \draw[cycle path,over] (\u)--(\v);
}
\newcommand{\cornerlabels}[1]{%
  \foreach \u/\lab in {s/s,a1/a_1,a2/a_2}
    \node[left corner,label=left:{$\lab$}] at (\u) {};
  \foreach \u/\lab in {r/#1,b1/b_1,b2/b_2}
    \node[right corner,label=right:{$\lab$}] at (\u) {};
}

\begin{scope}[name prefix=L-]
  \cornerpaths
  \coordinate (x) at (.93,1.24);
  \coordinate (y) at (2.82,1.24);
  \coordinate (z) at (1.87,1.93);
  \draw[figurePath!25,line width=7pt] (s)--(x)--(z);
  \draw[figureExtra!28,line width=7pt] (y)--(r);
  \draw[fixed] (s)--(x) (y)--(r);
  \draw[figureInk!55,densely dashed] (x)--(y);
  \draw[added] (x)--(z)--(y);
  \cornerlabels{t}
  \node[vertex,label=above:{$x$}] at (x) {};
  \node[vertex,label=above:{$y$}] at (y) {};
  \node[vertex,label=above:{$z$}] at (z) {};
  \node[text=figurePath]  at (.30,1.59) {$W_s$};
  \node[text=figureExtra] at (3.47,1.59) {$W_t$};
  \node at (1.875,-1.92) {$z\notin V(H)$};
\end{scope}

\begin{scope}[xshift=5.85cm,name prefix=M-]
  \cornerpaths
  \coordinate (x) at (1.18,1.24);
  \path (a2)--(r) coordinate[pos=.81] (t);
  \draw[fixed] (s)--(x);
  \draw[added] (x)--(r);
  \cornerlabels{z}
  \node[vertex,label=above:{$x$}] at (x) {};
  \node[vertex,label=right:{$t$}] at (t) {};
  \node at (1.875,-1.92) {$z$ internal to the $a_1$--$t$ branch};
\end{scope}

\begin{scope}[xshift=13.20cm,name prefix=R-]
  \def\rad{1.23}
  \coordinate (s)  at (0,0);
  \coordinate (t)  at (0,2.27);
  \coordinate (a1) at (135:\rad);
  \coordinate (b1) at (45:\rad);
  \coordinate (a2) at (-45:\rad);
  \coordinate (b2) at (225:\rad);
  \coordinate (x)  at (0,.59);
  \coordinate (z)  at (90:\rad);
  \coordinate (y)  at (0,1.76);
  \def\outerpath{(a1)--(t) .. controls (1.75,2.26) and
        (2.22,.08) .. (a2)}
  \fill[figurePath!10] (b1) arc (45:225:\rad)--(s)--cycle;
  \fill[figureExtra!11] \outerpath arc (-45:135:\rad)--cycle;
  \draw[fixed] \outerpath;
  \draw[fixed] (b1)--(s)--(b2) (s)--(x) (y)--(t);
  \draw[cycle path] (s) circle (\rad);
  \draw[added] (x)--(z)--(y);
  \node[left corner,label=below right:{$s$}] at (s) {};
  \node[right corner,label=above:{$t$}] at (t) {};
  \node[left corner,label=left:{$a_1$}] at (a1) {};
  \node[left corner,label=below right:{$a_2$}] at (a2) {};
  \node[right corner,label=right:{$b_1$}] at (b1) {};
  \node[right corner,label=below left:{$b_2$}] at (b2) {};
  \node[vertex,label=left:{$x$}] at (x) {};
  \node[vertex,label=above left:{$z$}] at (z) {};
  \node[vertex,label=right:{$y$}] at (y) {};
  \node[text=figureNavy] at (-.22,-1.45) {$C$};
  \node at (.12,-1.92) {$z\in V(C)$};
\end{scope}
\end{tikzpicture}%
}
\caption{Lemma~\ref{r2:position-of-third}. The dark blue paths form $C$.
Left, $xzy$ replaces the dashed edge $xy$.
Middle, the displayed $K_{3,3}$ subdivision has corner $z$ in place of $t$,
with $t$ on its $a_2$--$z$ branch.
Right, $xz$ and $yz$ meet $C$ from opposite sides.}
\label{fig:position-of-third}
\end{figure}

\begin{lemma}\label{r2:position-of-third}
Let $H$ be a $K_{3,3}$ subdivision containing the edge $xy$, let $P$ be its
branch path through $xy$, with ends $s,t$, and let $C$ be the cycle formed by
the four branch paths of $H$ avoiding $s$ and $t$. For a vertex $z\notin\{x,y\}$, possibly outside $H$,
the graph $H-xy+xz+yz$ is nonplanar if and only if $z\notin V(C)$.
\end{lemma}
\begin{proof}
Label the corners as in Figure~\ref{fig:position-of-third}, with $x$ before
$y$ on $P$ from $s$ to $t$. If $z\notin V(H)$, the edge $zy$ joins the
two highlighted paths $W_s=P[s,x]\cup xz$ and $W_t=P[y,t]$. If $z$ lies
on $P$, one of $xz,yz$ joins the two components of $P-xy$. In either
case we recover an $s$--$t$ path and hence a $K_{3,3}$ subdivision.

\looseness=-1
If $z$ lies in $H$ but off both $P$ and $C$, we may assume by symmetry
that it lies on the $a_1$--$t$ branch. As in the middle
drawing, use $z$ in place of $t$ as a corner. From $s$, follow $P$ to
$x$ and take $xz$. From $a_1$, follow the $a_1$--$t$ branch to $z$.
From $a_2$, follow the $a_2$--$t$ branch to $t$ and then the
$a_1$--$t$ branch to $z$. The six branches
to $b_1,b_2$ complete the subdivision.

Finally, suppose $z\in V(C)$. Draw $C$ as a circle with the two branches
from $s$ to $C$ inside and the two from $t$ outside. Every vertex of $C$
is on an inside face with $s$ and an outside face with $t$. We can
therefore draw $P[s,x]\cup xz$ in the inside face at $z$ and
$P[y,t]\cup yz$ in the outside face at $z$, giving a plane drawing of
$H-xy+xz+yz$.
\end{proof}

\begin{lemma}\label{r2:planar-splitting}
If $H$ is connected and planar, $ab\in E(H)$, $d\notin\{a,b\}$,
and $H+ad,H+bd$ are nonplanar, then $(H+ad+bd)/ab$ is nonplanar.
\end{lemma}
\begin{proof}
Suppose for the sake of contradiction that $(H+ad+bd)/ab$ is planar, and let $r$
be the vertex obtained by contracting $ab$. Neither $ad$ nor $bd$ is
an edge of $H$, since adding either one makes $H$ nonplanar.

In a plane drawing of $H$, subdivide $au$ once for each common
neighbor $u$ of $a$ and $b$, then contract $ab$ to $r$. Call the
resulting graph $J$. Around $r$, the edges that came from $a$ are
consecutive, as are those that came from $b$. To draw $J+rd$, start
with a plane drawing of $(H+ad+bd)/ab$ and add a path of length two beside $ru$
for each common neighbor $u$.

If $r$ and $d$ do not already share a face in the drawing of $J$
obtained from $H$, apply Lemma~6 of~\cite{ChimaniHlinenyMutzel2012}
to make them share one. Each flipped part contains every edge at $r$,
so each flip preserves or reverses their cyclic order.

Draw $rd$ in the common face, and let $rv$ be the next edge clockwise
from $rd$ at $r$. Undo the contraction inside a small disk around $r$,
using the original drawing there or its reflection. Inside the disk,
draw the edge from $d$ alongside $rv$ to $a$ if $rv$ came from $av$,
and to $b$ if it came from $bv$. Undoing the subdivisions gives a
plane drawing of $H+ad$ or $H+bd$, a contradiction.
\end{proof}

For a cycle $C$, a \emph{tripod} consists of a subdivision of $K_{2,3}$
and three disjoint attachment paths, possibly trivial, joining the vertices
corresponding to the three degree-two vertices of $K_{2,3}$ to $C$.
Each path meets the subdivision
only at its first vertex, and the tripod meets $C$ only at the last vertices of
these paths, its \emph{attachment triple}~\cite{Mohar1994}.

\begin{theorem}[{\cite[Theorem~3.1]{Mohar1994}}]\label{thm:mohar}
Let $C$ be a cycle in a graph $H$. If $H$ has no embedding in a closed
disk with $C$ on its boundary, then $H-E(C)$ contains \textup{(a)} two
disjoint crossing paths whose ends $a,c$ and $b,d$ are distinct
vertices of $C$ in the cyclic order $a,b,c,d$ and whose other vertices
avoid $C$, \textup{(b)} a tripod, or \textup{(c)} a Kuratowski subgraph
in a three-connected component of $\operatorname{Aux}(H,C)$ other than
that containing $C$.
\end{theorem}

\begin{lemma}\label{dis:four-disk-order}
Let $B$ consist of a nonempty connected graph $D$, four distinct
vertices $a,b,c,d$ outside $D$, and edges joining each of these
vertices to $D$. Suppose $B$ is drawn in a disk with $a,b,c,d$
on the boundary in this cyclic order and $D$ in the interior.
Then either $B$ can be redrawn with boundary order $a,b,d,c$,
keeping every vertex of degree one on the original outer face on the outer face,
or $B$ contains two vertex-disjoint paths, one from $a$ to $d$
and one from $b$ to $c$, whose internal vertices lie in $D$.
\end{lemma}

\begin{proof}
Consider the drawing as a plane drawing with outer facial walk $W$.
Each of $a,b,c,d$ occurs once on $W$. Consider the subwalks from $d$ to $a$
and from $b$ to $c$ that avoid the other two vertices. If they
are vertex-disjoint, they contain the required paths. Otherwise, choose a common vertex $x\in V(D)$ and an appearance
of $x$ on each subwalk. An arc in the outer face between these
appearances closes at $x$ to form a simple closed curve meeting
the drawing only at $x$ and separating $a,b$ from $c,d$.
Draw the two parts in disks meeting only at $x$ and reflect
the part containing $c,d$. The new outer boundary walk follows the segment of $W$ through $a,b$ in its
original order and the segment through $c,d$ in reverse. Thus the
boundary order is $a,b,d,c$, and every vertex of $W$ still lies
on the outer face.
\end{proof}

\section{Two distinct common neighborhoods}\label{sec:core}

We call distinct vertices of degree three with the same neighborhood \emph{twins}.
Two twins, their three common neighbors and the six edges between them
form a copy of $K_{2,3}$, which we call their \emph{twin subgraph}.
The main result of this section,
Theorem~\ref{thm:two-pages}, shows that apart from $K_{3,7}$, $G$ has
two distinct independent triples, each the neighborhood of at least two
twins. 

We find twins through three-cuts, using Halin's theorem. For a
three-cut $X$ and a component $D$ of $G-X$, write
$B_X=G[V(D)\cup X]-E(G[X])$ for its bridge in $G$ and call
$H_D=B_X\cup K_3[X]$
the \emph{triangle completion} of $D$. When $G-X$ has at least three
components, $H_D$ is a three-connected proper minor of $G$ by
Lemma~\ref{lem:triangle-completion}. If $H_D$ is planar, then
$D$ is a single vertex with neighborhood $X$ by
Lemmas~\ref{lem:triangle-completion} and~\ref{lem:flat-replacement},
and two such vertices are twins.

\begin{lemma}\label{lem:triangle-completion}
For a three-cut $X$ and a component $D$ of $G-X$, we have $N(D)=X$, and $H_D$ is three-connected. If $G-X$ has at least three components, then $H_D$ is a proper minor of $G$, hence $K_5$-minor-free and toroidal.
\end{lemma}

\begin{proof}
Three-connectivity gives $N(D)=X$. Since $X$ is a clique in
$H_D$, any vertex cut of $H_D$ with at most two vertices would also separate $G$. If two other components besides $D$ exist, contract them onto distinct vertices of $X$
to create a triangle, then delete any vertices outside $V(D)\cup X$.
\end{proof}

\begin{lemma}\label{lem:flat-replacement}
If $H_D$ is planar, then $D$ is a single vertex.
\end{lemma}

\begin{proof}
Since $H_D$ is planar, Lemma~\ref{cor:boundary-completions} gives a disk drawing
of $B_X$ meeting the boundary precisely in $X$.
If $|V(D)|>1$, contracting $D$ to a vertex $v$ gives a toroidal
proper minor. Replacing $v$ and its three incident edges by
that drawing inside a thin closed disk meeting the rest of
the embedding precisely in $X$, and using
the Jordan--Schoenflies theorem~\cite[Section~2.2]{MoharThomassen2001} then yields a toroidal embedding of $G$, a contradiction. 
\end{proof}

\begin{lemma}\label{lem:avoiding-components}
Let $B$ be a planar bridge of $G$ with attachment set $A$, and let $X$ be a
set of three vertices of $B$. Every component of $B-X$ disjoint from $A$ is a
component of $G-X$ with neighborhood $X$. There are at most two such
components, and if there are two, both are single vertices.
\end{lemma}
\begin{proof}
Such a component $D$ consists of interior vertices of $B$. Every edge of $G$
at an interior vertex of $B$ belongs to $B$, so every neighbor of $D$ in $G$
outside $D$ lies in $X$. Thus $D$ is a component of $G-X$, and
three-connectivity gives $N_G(D)=X$. Its bridge on $X$ in $B$ is therefore
adjacent to all of $X$, so Lemma~\ref{core:assembly-general}(i) allows at
most two such components. If there are two, no other bridge on $X$ in $B$ is
adjacent to all of $X$, so part (ii) makes both triangle completions planar,
and Lemma~\ref{lem:flat-replacement} makes both components single vertices.
\end{proof}

\begin{lemma}\label{lem:independent-neighborhood}
Let $p$ be a vertex of degree three in an obstruction. Then (a) the
neighbors of $p$ are independent, and (b) if the obstruction is
three-connected and $p,r$ have twins with $N(p)\ne N(r)$, then $p$ and $r$ are not adjacent.
\end{lemma}
\begin{proof}
(a) follows from~\cite[Lemma~24]{MoharSkoda2017}.
(b) Let $q$ and $z$ be the twins of $p$ and $r$, and suppose
$pr$ is an edge. Then
$N(p)=N(q)=\{r,z,b\}$ and $N(r)=N(z)=\{p,q,d\}$
for some $b,d$. The vertices $p,q,r,z$ are distinct
as $N(p)\ne N(r)$, and $b,d$ lie outside them as twins are not
adjacent. Thus deleting $b,d$ separates $\{p,q,r,z\}$ from the other
vertices, which exist since obstructions have more than seven vertices.
This contradicts three-connectivity.
\end{proof}

The next lemma finds two twins inside $D$. 

\begin{lemma}\label{lem:minimal-side}
Let $Z\subseteq V(G)$. Suppose $G$ has a three-cut $X$ such that $G-X$
has at least three components, one of which, $D$, avoids $Z$ and has
$H_D$ nonplanar, and choose $X$ and $D$ with $|V(D)|$ as small as
possible. Then $G$ has a three-cut $X'$ such that the components of
$G-X'$ are
\begin{enumerate}[label=\textup{(\alph*)}]
\item at least two single vertices in $D$, each with neighborhood $X'$,
and
\item one further component $D'$, which contains $X-X'$ and every
component of $G-X$ other than $D$, and has~$H_{D'}$~nonplanar.
\end{enumerate}
\end{lemma}
\begin{proof}
By Lemma~\ref{lem:triangle-completion}, $H_D$ is three-connected and has no
$K_5$ minor. It is nonplanar and has a triangle on $X$, so Halin's theorem
gives a three-cut $X'$ of $H_D$ leaving at least
three components. Since $H_D-X=D$ is connected, $X'\ne X$. Also, as
$X'\subseteq V(H_D)=V(D)\cup X$ and $|X'|=|X|=3$, $X'$ meets $D$ and
$X-X'$ is nonempty. The vertices of $X-X'$ are a clique in $H_D$ and thus lie in
one component $D_0$ of $H_D-X'$. Every other component $D_i$ lies in
$D-X'$, and there are at least two. Each $D_i$ consists of interior
vertices of the bridge $B_X$, and $H_D$ is obtained from $B_X$ by adding
edges between attachments. Thus $D_i$ is also a component of $G-X'$.

The $D_i$ give at least two components of $G-X'$, and any vertex
of $D_0$ belongs to another. Since $D$ avoids $Z$ and $X'$ meets
$D$, each $D_i\subseteq D-X'$ avoids $Z$ and has fewer vertices
than $D$. If some $H_{D_i}$ were nonplanar, we could have chosen
$X'$ and $D_i$ in place of $X$ and $D$, contrary to the minimality
of $|V(D)|$. Thus every $H_{D_i}$ is planar. By
Lemma~\ref{lem:flat-replacement}, each $D_i$ is a single vertex,
and by Lemma~\ref{lem:triangle-completion}, its neighborhood is $X'$.

Let $D'$ be the subgraph of $G-X'$ induced by all vertices
outside the $D_i$. These are the vertices of $D_0$ and of
every component of $G-X$ other than $D$.
Each component of $G-X$ other than $D$ has neighborhood $X$
by Lemma~\ref{lem:triangle-completion}, so it is joined to
$X-X'\subseteq V(D_0)$. Any added triangle edge used in $D_0$
can be replaced by a path through one of these components.
Thus $D'$ is connected in $G-X'$. Since the $D_i$ are
components of $G-X'$, no edge joins $D'$ to any $D_i$, so
$D'$ is also a component. Finally, $D'$ contains $X-X'$
and at least two components of $G-X$, so $|V(D')|>1$.
Lemma~\ref{lem:flat-replacement} therefore implies that
$H_{D'}$ is nonplanar.
\end{proof}

\begin{theorem}\label{thm:two-pages}
If $G\not\cong K_{3,7}$, then $G$ has two distinct independent triples,
each the neighborhood of at least two degree-three vertices.
\end{theorem}
\begin{proof}
As $G$ is nonplanar and $V_8$ is toroidal, Halin's theorem gives a
three-cut $X$ leaving at least three components. If every
$H_D$ is planar, every component of $G-X$ is a single vertex
with neighborhood $X$ (Lemmas~\ref{lem:triangle-completion}
and~\ref{lem:flat-replacement}), so $X$ is independent by
Lemma~\ref{lem:independent-neighborhood}(a) and $G=K_{3,k}$. By
Ringel's formula~\cite{Ringel1965}, $K_{3,k}$ has genus
$\lceil (k-2)/4\rceil$, so $k\ge7$ and minimality gives $G\cong K_{3,7}$. 

When $G\not\cong K_{3,7}$, some component must have a nonplanar triangle completion. Choose a three-cut $X$ leaving at least three components and a component
$D$ with $H_D$ nonplanar, minimizing $|V(D)|$.
Lemma~\ref{lem:minimal-side}, with $Z=\varnothing$, gives a
three-cut $A_1$, at least two single-vertex components with
neighborhood $A_1$, and one further component $D'$ with
$H_{D'}$ nonplanar.

Let $W$ be the set of all vertices with neighborhood $A_1$.
Every vertex of $W$ is a single-vertex component of $G-A_1$.
Since $H_{D'}$ is nonplanar, $D'$ has more than one vertex,
so $V(D')\cap W=\varnothing$. Now, choose a three-cut $X$ leaving
at least three components and a component $D$ disjoint from
$W$ with $H_D$ nonplanar, again minimizing $|V(D)|$.
Such a choice exists, as $G-A_1$ has at least three components,
and $D'$ is one of them, is disjoint from $W$, and has
nonplanar triangle completion.
Lemma~\ref{lem:minimal-side}, now with $Z=W$, gives a three-cut
$A_2$ and at least two vertices in $D$ with neighborhood $A_2$.
These vertices lie outside $W$, so their neighborhood cannot
be $A_1$. Thus $A_2\ne A_1$. Both triples are independent by
Lemma~\ref{lem:independent-neighborhood}(a).
\end{proof}

\begin{figure}[htbp]
\centering
\begingroup
\scaledfigure{1.1}{\scalebox{0.94}{%
\tikzsetnextfilename{three-cases}%
\begin{tikzpicture}[
  font=\small,text=figureInk,
  vertex/.style={circle,fill=figureInk,inner sep=1.65pt},
  twinA/.style={circle,fill=figureTwinA,inner sep=1.9pt},
  twinB/.style={circle,fill=figureTwinB,inner sep=1.9pt},
  edgeA/.style={draw=figureTwinA,line width=.75pt,
    preaction={draw=white,line width=1.8pt}},
  edgeB/.style={draw=figureTwinB,line width=.75pt,
    preaction={draw=white,line width=1.8pt}},
  bandA/.style={draw=figureTwinA,fill=figureTwinA!9,
    line width=.6pt,rounded corners=4pt},
  bandB/.style={draw=figureTwinB,fill=figureTwinB!9,
    fill opacity=.6,line width=.6pt,rounded corners=4pt},
  rest/.style={draw=figureInk!45,fill=figureInk!3,
    line width=.6pt,rounded corners=5pt},
  heading/.style={font=\normalsize\bfseries},
  guide/.style={font=\small,text=figureInk},
  lab/.style={font=\normalsize,anchor=north,inner sep=0pt}
]
\foreach \x in {0,5.6,11.2}{
  \begin{scope}[xshift=\x cm]
    \draw[rest] (.05,-1.25) rectangle (4.95,.3);
    \node[font=\normalsize,text=figureInk] at (2.5,-.91) {rest of $G$};
    \node[font=\small] at (1.25,1.69) {twins};
    \node[font=\small] at (3.75,1.69) {twins};
  \end{scope}
}

\begin{scope}
  \node[heading] at (2.5,2.65) {No common vertex};
  \node[guide] at (2.5,2.20)
    {Section~\ref{sec:disjoint}};
  \draw[bandA] (.2,-.53) rectangle (2.3,.17);
  \draw[bandB] (2.7,-.53) rectangle (4.8,.17);
  \draw[figureTwinA,line width=.6pt,rounded corners=4pt] (.2,-.53) rectangle (2.3,.17);
  \draw[figureTwinB,line width=.6pt,rounded corners=4pt] (2.7,-.53) rectangle (4.8,.17);
  \foreach \i/\x in {1/.55,2/1.25,3/1.95}
    \coordinate (s\i) at (\x,0);
  \foreach \i/\x in {1/3.05,2/3.75,3/4.45}
    \coordinate (t\i) at (\x,0);
  \foreach \x in {.85,1.65}{
    \foreach \v in {s1,s2,s3} \draw[edgeA] (\x,1.22)--(\v);
    \node[twinA] at (\x,1.22) {};
  }
  \foreach \x in {3.35,4.15}{
    \foreach \v in {t1,t2,t3} \draw[edgeB] (\x,1.22)--(\v);
    \node[twinB] at (\x,1.22) {};
  }
  \foreach \i in {1,2,3}{
    \node[vertex] at (s\i) {};
    \node[lab] at ([yshift=-.12cm]s\i) {$s_\i$};
    \node[vertex] at (t\i) {};
    \node[lab] at ([yshift=-.12cm]t\i) {$t_\i$};
  }
\end{scope}

\begin{scope}[xshift=5.6cm]
  \node[heading] at (2.5,2.65) {Two common vertices};
  \node[guide] at (2.5,2.20) {Section~\ref{sec:two-common}};
  \draw[bandA] (.28,-.53) rectangle (3.49,.17);
  \draw[bandB] (1.51,-.56) rectangle (4.72,.2);
  \draw[figureTwinA,line width=.6pt,rounded corners=4pt] (.28,-.53) rectangle (3.49,.17);
  \draw[figureTwinB,line width=.6pt,rounded corners=4pt] (1.51,-.56) rectangle (4.72,.2);
  \coordinate (c) at (.65,0); \coordinate (a) at (1.85,0);
  \coordinate (b) at (3.15,0); \coordinate (d) at (4.35,0);
  \foreach \x in {.85,1.65}{
    \foreach \v in {c,a,b} \draw[edgeA] (\x,1.22)--(\v);
    \node[twinA] at (\x,1.22) {};
  }
  \foreach \x in {3.35,4.15}{
    \foreach \v in {a,b,d} \draw[edgeB] (\x,1.22)--(\v);
    \node[twinB] at (\x,1.22) {};
  }
  \foreach \v in {c,a,b,d}{
    \node[vertex] at (\v) {};
    \node[lab] at ([yshift=-.12cm]\v) {$\v$};
  }
\end{scope}

\begin{scope}[xshift=11.2cm]
  \node[heading] at (2.5,2.65) {One common vertex};
  \node[guide] at (2.5,2.20) {Section~\ref{sec:one-common}};
  \draw[bandA] (.28,-.53) rectangle (2.82,.17);
  \draw[bandB] (2.18,-.56) rectangle (4.72,.2);
  \draw[figureTwinA,line width=.6pt,rounded corners=4pt] (.28,-.53) rectangle (2.82,.17);
  \draw[figureTwinB,line width=.6pt,rounded corners=4pt] (2.18,-.56) rectangle (4.72,.2);
  \coordinate (b) at (.65,0); \coordinate (c) at (1.575,0);
  \coordinate (a) at (2.5,0);
  \coordinate (d) at (3.425,0); \coordinate (e) at (4.35,0);
  \foreach \x in {.85,1.65}{
    \foreach \v in {b,c,a} \draw[edgeA] (\x,1.22)--(\v);
    \node[twinA] at (\x,1.22) {};
  }
  \foreach \x in {3.35,4.15}{
    \foreach \v in {a,d,e} \draw[edgeB] (\x,1.22)--(\v);
    \node[twinB] at (\x,1.22) {};
  }
  \foreach \v in {b,c,a,d,e}{
    \node[vertex] at (\v) {};
    \node[lab] at ([yshift=-.12cm]\v) {$\v$};
  }
\end{scope}
\end{tikzpicture}}}
\endgroup
\caption{An illustration of the three cases.}
\label{fig:three-cases}
\end{figure}
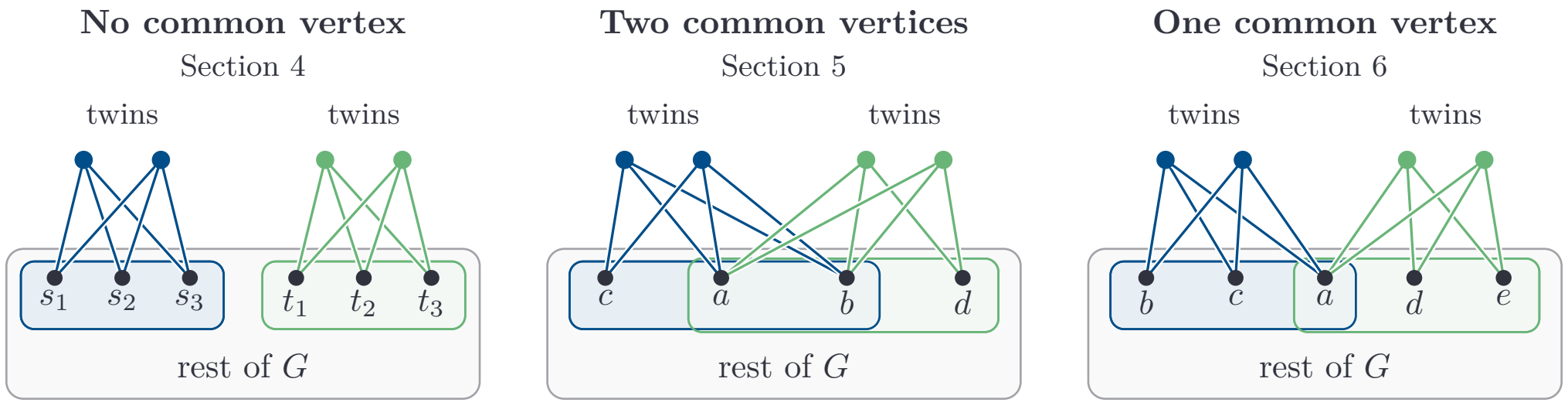

An independent triple that is the neighborhood of at least two twins is
a \emph{twin triple}. Two distinct twin triples share at most two
vertices, so the twin triples of Theorem~\ref{thm:two-pages} give three
cases (Figure~\ref{fig:three-cases}). If some two twin triples are
disjoint, Theorem~\ref{thm:disjoint} applies. If no two are disjoint but
some two share two vertices, Theorem~\ref{thm:two-common} applies.
Otherwise every two distinct twin triples share exactly one vertex, and
Theorem~\ref{thm:one-common} applies. In each case we fix two twin triples
and two twins with each neighborhood, and write $J_1,J_2$ for their
twin subgraphs. In Section~\ref{sec:disjoint}, three paths joining
$J_1$ to $J_2$ give the graph $\Psi_1$. In
Sections~\ref{sec:two-common} and~\ref{sec:one-common}, we use their
union as $\Psi_2$ and $\Psi_3$, respectively. We study the rest of
$G$ through bridges. Sections~\ref{core:three-root-section}--\ref{sec:common-construction}
develop the tools used for these bridges.

\subsection{Bridges with three attachments}\label{core:three-root-section}

Let $B$ be a graph with a set $A$ of \emph{attachment} vertices, no two
of them adjacent. We assume that its interior is nonempty and connected,
and every attachment has a neighbor in it. Every bridge of $G$ with a
nonempty interior has these properties. 

Adding edges between attachments to $B$ gives a \emph{completion}
of $B$, and when $|A|=3$, the completion $B\cup K_3[A]$ is its \emph{triangle completion}. An \emph{allowed operation} deletes an interior vertex or
an edge, or contracts an edge with an interior end and deletes any
edges this creates between attachments. The result must retain
the labeled attachments, a nonempty connected interior, and an
interior neighbor at every attachment. Graphs obtained by these
operations are \emph{allowed minors}. Planar completions
remain planar under allowed minor operations, and replacing a bridge of $G$
by an allowed minor gives a minor of $G$.

\begin{lemma}\label{core:relevant}
Let $B$ be a bridge of $G$ with nonempty interior, and let $\mathcal E$ be a set of sets of edges between its attachments. Suppose that, for every allowed minor $B'$ of $B$, whether the graph obtained from $G$ by replacing $B$ with $B'$ is toroidal depends only on which of the completions $B'+E$ with $E\in\mathcal E$ are planar. Then every proper allowed minor $B'$ of $B$ has $B'+E$ planar for some $E\in\mathcal E$ with $B+E$ nonplanar.
\end{lemma}
\begin{proof}
Otherwise $B'+E$ is planar exactly when $B+E$ is, for every $E\in\mathcal E$. Replacing $B$ with $B'$ then gives a nontoroidal graph, since $G$ is nontoroidal. This is a proper minor of $G$, contrary to minimality.
\end{proof}

When the hypothesis of Lemma~\ref{core:relevant} holds, we call completions $B+E$ with $E\in\mathcal E$ the \emph{relevant completions} of $B$. We use the following lemma to describe the bridges resulting from a three-cut when it has some planar relevant~completion.

\begin{lemma}\label{core:three-root-setup}
Let $B$ be a graph with attachment set $A=\{x,y,z\}$ as above, and suppose $B\cup K_3[A]$ is three-connected and $B+E$
is planar for some $E\subseteq E(K_3[A])$. Let $X$ be a three-cut of
$B\cup K_3[A]$ such that $(B\cup K_3[A])-X$ has at least three components. Then
$X\ne A$, and the following hold.
\begin{enumerate}[label=\textup{(\roman*)}]
\item $(B\cup K_3[A])-X$ has exactly three components, $D_0,D_1,D_2$, with $V(D_0)\supseteq A-X$
and $V(D_1),V(D_2)\subseteq \operatorname{Int}(B)$.
\item Every component of $B-X$ inside $D_0$ meets $A$.
\item For $i=1,2$, $N_B(D_i)=X$ and $B[V(D_i)\cup X]\cup K_3[X]$ is planar.
\item No bridge on $X$ in $B+E$ other than those of $D_1,D_2$ is
adjacent to all of $X$.
\item If $B$ is a bridge of $G$ on $A$, then $D_1$ and $D_2$ are single
vertices $p$ and $q$ with neighborhood $X$.
\item For every $E'\subseteq E(K_3[A])$, a bridge on $X$ in $B+E'$ containing
an edge of $E'$ has at least two attachments. Also, $B+E'$ is nonplanar
exactly when such a bridge is adjacent to all of $X$, or has two attachments
$\alpha,\beta$ and becomes nonplanar on adding $\alpha\beta$, and in the
latter case the completed bridge is a minor of $B\cup K_3[A]$. In particular, if $u$ and $v$ are distinct vertices of $A\cap X$, then $B+E'+uv$ is planar if and only if~$B+E'$~is~planar.
\end{enumerate}
If $B\cup K_3[A]$ is nonplanar and has no $K_5$ minor, such a three-cut $X$ exists.
\end{lemma}
\begin{proof}
Deleting $A$ from $B\cup K_3[A]$ leaves $B[\operatorname{Int}(B)]$,
which is connected, so $X\ne A$. Since $|A|=|X|=3$, $A-X$ is
nonempty and is a clique in $B\cup K_3[A]$, so it lies in one
component, say $D_0$, of $(B\cup K_3[A])-X$.

Every other component of $(B\cup K_3[A])-X$ avoids $A$. The added
edges of $B+E$ and $B\cup K_3[A]$ join attachments, so the graphs
$B-X$, $(B+E)-X$ and $(B\cup K_3[A])-X$ have the same components
avoiding $A$. A component of $B-X$ avoiding $A$ is therefore a component
of $(B\cup K_3[A])-X$ other than $D_0$. Hence every
component of $B-X$ inside $D_0$ meets $A$, proving (ii).

Let $D$ be a component of $(B\cup K_3[A])-X$ other than $D_0$.
Three-connectivity of $B\cup K_3[A]$ gives $N_B(D)=X$, so $D$ is a
component of $(B+E)-X$ whose bridge on $X$ is adjacent to all of $X$.
By Lemma~\ref{core:assembly-general}(i), the planar graph $B+E$ has at
most two such bridges. Since $(B\cup K_3[A])-X$ has at least three
components, there are exactly two, $D_1$ and $D_2$, which proves (i)
and (iv). By Lemma~\ref{core:assembly-general}(ii), the bridge of each
$D_i$ plus $K_3[X]$ is planar, and this graph is
$B[V(D_i)\cup X]\cup K_3[X]$, proving (iii). For (v), suppose $B$ is a bridge of $G$ on $A$. It is planar since $B+E$
is, and $D_1,D_2$ are components of $B-X$ disjoint from $A$, so
Lemma~\ref{lem:avoiding-components} makes them single vertices with
neighborhood $X$.
For (vi), suppose a bridge on $X$ in $B+E'$ containing an edge of $E'$ has
at most one attachment. Then its interior contains an end of that edge and
has at most two neighbors in $B\cup K_3[A]$, consisting of at most one
attachment and at most one further vertex of $A$. This contradicts
three-connectivity. Every bridge not containing an edge of $E'$ is unchanged
from $B$, which is planar since $B+E$ is. If such a bridge has two
attachments $\alpha,\beta$, an $\alpha$--$\beta$ path through $D_1$ makes it
plus $\alpha\beta$ a minor of $B$, so this completion is also planar. As $B$ is
planar, Lemma~\ref{core:assembly-general}(i) allows no such bridge other than
those of $D_1,D_2$ to be adjacent to all of $X$. The bridges of $D_1,D_2$ are
adjacent to all of $X$ by (iii), so Lemma~\ref{core:assembly-general}(i),(ii)
makes $B+E'$ nonplanar exactly when a bridge containing an edge of $E'$ is
adjacent to all of $X$, or has two attachments $\alpha,\beta$ and becomes
nonplanar on adding $\alpha\beta$. In the latter case, the completed bridge is
a minor of $B\cup K_3[A]$, using an $\alpha$--$\beta$ path through $D_1$.

Finally, if $B\cup K_3[A]$ is nonplanar and has no $K_5$ minor, Halin's
theorem gives such a three-cut, since this graph contains a triangle.
\end{proof}

\begin{figure}[htbp]
\centering
\tikzsetnextfilename{three-root-cut-bottom}%
\begin{tikzpicture}[
  font=\small,text=figureInk,
  vertex/.style={circle,fill=figureInk,inner sep=1.55pt},
  root/.style={circle,fill=figurePath,inner sep=1.9pt},
  edge/.style={draw=figureInk!75,line width=.75pt,line cap=round,line join=round},
  every label/.style={inner sep=1.2pt}
]
\begin{scope}
  \node at (-.45,1.34) {$T(z)$};
  \node[vertex,label=left:{$p$}] (d1) at (-.95,0) {};
  \node[vertex,label=right:{$q$}] (d2) at (.95,0) {};
  \node[root,label=above:{$x$}] (x) at (0,.68) {};
  \node[root,label=above:{$y$}] (y) at (0,0) {};
  \node[vertex,label=below:{$u$}] (u) at (0,-.68) {};
  \node[root,label=left:{$z$}] (z) at (-1.85,-.68) {};
  \foreach \w in {x,y,u} {\draw[edge] (d1)--(\w)--(d2);}
  \draw[edge] (u)--(z);
\end{scope}
\begin{scope}[xshift=5cm]
  \node at (-.65,1.34) {$S(x;y)$};
  \node[vertex,label=above:{$p$}] (d1) at (-.91,0) {};
  \node[vertex,label=right:{$q$}] (d2) at (.95,0) {};
  \node[root,label=above:{$x$}] (x) at (0,.68) {};
  \node[root,label=above:{$y$}] (y) at (0,0) {};
  \node[vertex,label=below:{$u$}] (u) at (0,-.68) {};
  \node[vertex,label=below:{$f$}] (f) at (-1.75,0) {};
  \node[root,label=left:{$z$}] (z) at (-2.35,0) {};
  \foreach \w in {x,y,u} {\draw[edge] (d1)--(\w)--(d2);}
  \draw[edge] (x) .. controls (-.76,.74) and (-1.55,.64) .. (f)
    .. controls (-1.55,-.64) and (-.76,-.74) .. (u);
  \draw[edge] (f)--(z);
\end{scope}
\begin{scope}[xshift=9.25cm]
  \node at (0,1.34) {$Q(x;y,z)$};
  \node[vertex,label=above:{$p$}] (d1) at (0,.68) {};
  \node[vertex,label=below:{$q$}] (d2) at (0,-.68) {};
  \node[vertex,label=above:{$u$}] (u) at (-.93,0) {};
  \node[vertex,label=above:{$v$}] (v) at (.93,0) {};
  \node[root,label=right:{$x$}] (x) at (0,0) {};
  \node[root,label=left:{$y$}] (y) at (-1.73,0) {};
  \node[root,label=right:{$z$}] (z) at (1.73,0) {};
  \draw[edge] (u)--(d1)--(v)--(d2)--(u);
  \draw[edge] (d1)--(x)--(d2) (y)--(u) (v)--(z);
\end{scope}
\end{tikzpicture}
\caption{The graphs $T(z)$, $S(x;y)$ and $Q(x;y,z)$.}
\label{fig:three-root-cut}
\end{figure}
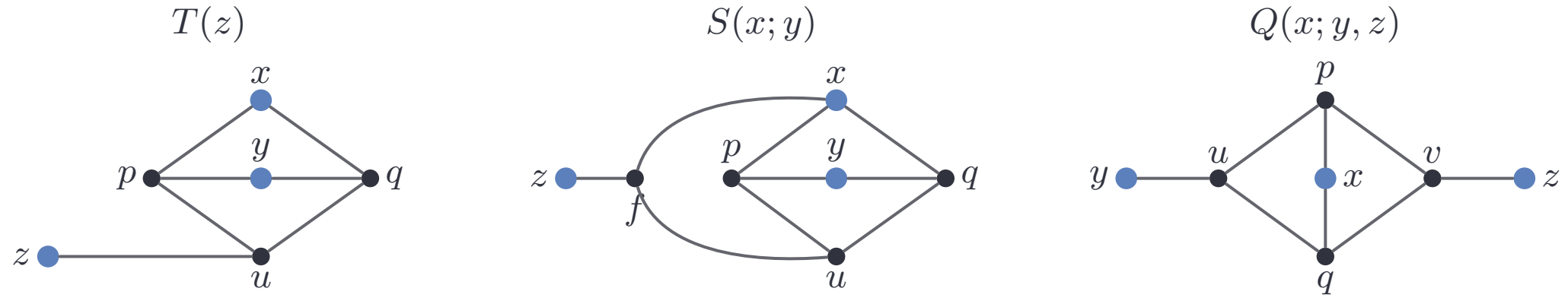

The next lemma treats the case where $B+xy+xz$ is planar and the triangle
completion of $B$ is nonplanar. It gives one of the three graphs in
Figure~\ref{fig:three-root-cut} as an allowed minor of $B$. We use the
same notation after relabeling~the~attachments.

\begin{lemma}\label{core:three-root-cut}
Let $B$ be a bridge of $G$ on $A=\{x,y,z\}$, and suppose
$B+xy+xz$ is planar. Suppose the triangle completion $B\cup K_3[A]$ is three-connected, nonplanar and without a $K_5$ minor.
By Lemma~\ref{core:three-root-setup}, this graph has a three-cut leaving
at least three components. Let $X$ be any such cut, and let $D_0$, $p$
and $q$ be as in that lemma with $E=\{xy,xz\}$.
Then one of the following holds, where
$u,v\in\operatorname{Int}(B)$.
\begin{enumerate}[label=\textup{(\roman*)}]
\item $X=\{x,u,v\}$. Up to interchanging $u,v$, the graph
$B[V(D_0)]$ has exactly two components, one containing $y$
and adjacent to $u$ but not $v$, and one containing $z$
and adjacent to $v$ but not $u$. Moreover, $B+xy+yz$ and
$B+xz+yz$ are nonplanar, and $B$ has $Q(x;y,z)$ as an allowed minor.

\item $X=\{x,y,u\}$. No vertex of $D_0-z$ is adjacent to $y$,
and every component of $B[V(D_0)\setminus\{z\}]$ has
neighborhood $\{z,x,u\}$ in $B$. If $V(D_0)=\{z\}$, then
$B=T(z)$. Otherwise,
$B+xy+yz$ and $B+xz+yz$ are nonplanar, and $B$ has
$S(x;y)$ as an allowed minor.

\item $X=\{x,z,u\}$. The claims of (ii) hold with $y,z$
interchanged.
\end{enumerate}
\end{lemma}
\begin{proof}
Applying Lemma~\ref{core:three-root-setup} with $E=\{xy,xz\}$,
the components of $(B\cup K_3[A])-X$ are $D_0,\{p\},\{q\}$,
with $A-X\subseteq V(D_0)$ and $X\ne A$.
The connections from $D_0$ to the vertices of $X$ necessarily exist
in $B\cup K_3[A]$, since $N_{B\cup K_3[A]}(D_0)\subseteq X$
and, by three-connectivity, $N_{B\cup K_3[A]}(D_0)=X$.
However, they need not be edges of $B$.

By Lemma~\ref{core:three-root-setup}(v), $p,q$ have neighborhood
$X$, so $X$ is independent by
Lemma~\ref{lem:independent-neighborhood}(a).
Suppose for contradiction that $x\notin X$. We first show that
$D_0$ is also a component of $(B+E)-X$.
Since $x\in A-X\subseteq V(D_0)$, it is enough to connect each
$w\in V(D_0)$ to $x$ within $(B+E)-X$.
By Lemma~\ref{core:three-root-setup}(ii), the component of $w$
in $B-X$ contains an attachment in $A$. If this attachment is $x$,
we are done. Otherwise it is $y$ or $z$, and the edge $xy$ or $xz$
connects $w$ to $x$.
Also, $N_{B+E}(D_0)=X$, since edges of $B$ are retained and every
attachment in $X$ is joined to $x\in V(D_0)$.
Thus the bridge with interior $V(D_0)$ is a third bridge adjacent
to all of $X$, contradicting planarity by
Lemma~\ref{core:assembly-general}(i).
Hence $x\in X$, and since $X\ne A$, its intersection with $A$
is $\{x\}$, $\{x,y\}$ or $\{x,z\}$, with
$X-A\subseteq\operatorname{Int}(B)$.

In case (i), $X=\{x,u,v\}$.
We have shown that $x\notin X$ gives a $K_{3,3}$ minor of
$B+xy+xz$. Interchanging $x$ with $y$ or $z$ shows that
$B+xy+yz$ and $B+xz+yz$ are nonplanar, since $y,z\notin X$.
Now, if $B[V(D_0)]$ were connected, its bridge on $X$ in
$B+xy+xz$ would be adjacent to $x$ through $xy$ and to $u,v$
through edges of $B$. This would form a third bridge adjacent
to all of $X$, contradicting planarity by
Lemma~\ref{core:assembly-general}(i).
Since $B[V(D_0)]+yz$ is connected, $B[V(D_0)]$ therefore has
exactly two components, $D_y$ containing $y$ and $D_z$
containing $z$.
Each of $D_y,D_z$ is adjacent to $x$ in $B+xy+xz$, so neither
can be adjacent to both $u,v$, as this would give a third bridge
adjacent to all of $X$.
If $D_y$ were adjacent to neither, its neighborhood in
$B\cup K_3[A]$ would be contained in $\{x,z\}$, contradicting
three-connectivity. The same argument applies to $D_z$ with
$\{x,y\}$. Thus each is adjacent to exactly one of $u,v$.
Since $D_0$ is adjacent to both, after interchanging $u,v$ if
necessary, $D_y$ is adjacent to $u$ but not $v$, and $D_z$
to $v$ but not $u$. There is a path from $u$ to $y$ through $D_y$, and a path from
$v$ to $z$ through $D_z$. Contract their internal vertices into
$u$ and $v$, respectively, retaining the edges $uy,vz$.
Delete any remaining vertices of $D_y-\{y\}$ and $D_z-\{z\}$,
starting with those furthest from $u$ or $v$, respectively,
with distances measured within the interior, to keep the interior
connected. Delete any edges $xu$ or $xv$ that may have been formed.
The graph that remains is $Q(x;y,z)$.
This minor is allowed, since the interior remains nonempty and
connected, the attachments remain distinct, and each keeps an
interior neighbor, as $x$ keeps $p,q$, and $y,z$ keep $u,v$.

Now consider case (ii), $X=\{x,y,u\}$.
Each added triangle edge has an end in $X$, so
\[
(B\cup K_3[A])-X=B-X,
\]
and $D_0$ is a component of $B-X$. We claim that no
vertex of $D_0-z$ is adjacent to $y$.
If some vertex of $D_0-z$ were adjacent to $y$, then $D_0$
already has a neighbor at $u$, since its neighborhood in the
triangle completion is $X=\{x,y,u\}$ and none of the added
triangle edges touches $u$. The added edge $xz$ makes $D_0$
adjacent to $x$. Contracting $D_0$ to a single vertex $m$ gives
a $K_{3,3}$ minor with parts $\{x,y,u\}$ and $\{p,q,m\}$,
contradicting the planarity of $B+xy+xz$.

If $D_0-z$ is empty, then $N_{B\cup K_3[A]}(D_0)=X$ forces
$uz\in E(B)$, and since $X$ is independent, $B=T(z)$.
Otherwise, choose a connected component $C$ of $D_0-z$.
Inside $D_0$, $z$ is the only possible neighbor of $C$ outside
itself. Since $D_0$ is a component of $B-X$, every neighbor
outside $D_0$ lies in $X$. No vertex of $D_0-z$ is adjacent
to $y$, so $N_B(C)\subseteq\{x,u,z\}$.
Since $C$ contains no attachments, it has the same neighborhood
in $B\cup K_3[A]$. Thus $N_B(C)=\{x,u,z\}$, otherwise
the triangle completion would not be three-connected. This argument for $x\notin X$ also shows that $B+xz+yz$ is
nonplanar, since $z\notin X$.

\looseness=-1
In $B+xy+yz$, the connected graph $D_0$ is adjacent to $x,u$
through edges from $C$ and to $y$ through $yz$.
Contracting $D_0$ again gives a $K_{3,3}$ minor, so this graph
is also nonplanar.
Now contract $C$ to a single vertex $f$ and delete the other
components of $D_0-z$, starting with vertices furthest from $u$
in the interior, and keeping the edges $fx,fu,fz$.
Delete $uz$ if it is present.
This retains all three attachments and keeps the
interior nonempty and connected. Also, $x,y$ retain their
neighbors $p,q$, and $z$ retains its neighbor $f$, so the
operations are allowed. The resulting graph is $S(x;y)$.

In case (iii), $X=\{x,z,u\}$, and the same argument applies
with $y,z$ swapped.
\end{proof}

We next determine which completions of $T$, $S$ and $Q$ are planar.

\begin{lemma}\label{core:three-graph-completions}
For $E\subseteq\{xy,xz,yz\}$, $T(z)+E$ is planar if and only if
$E$ does not contain both $xz$ and $yz$, $S(x;y)+E$ is planar if
and only if $yz\notin E$, and $Q(x;y,z)+E$ is planar if and only if
$E\subseteq\{xy,xz\}$ or $E=\{yz\}$.
\end{lemma}
\begin{proof}
PAGE computed the genera of all $24$ completions, finding genus zero
exactly as stated~\supp{E30}.
\end{proof}

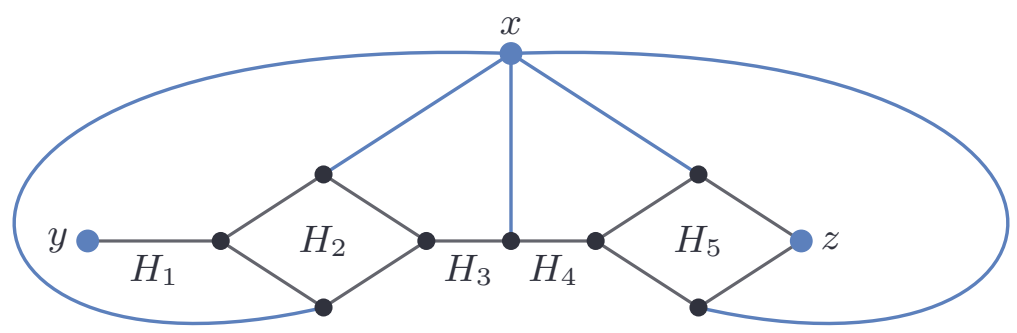
\begin{figure}[htbp]
\centering
\scalebox{.8}{%
\tikzsetnextfilename{three-root-chain}%
\begin{tikzpicture}[
  font=\small,text=figureInk,
  vertex/.style={circle,fill=figureInk,inner sep=1.5pt},
  root/.style={circle,fill=figurePath,inner sep=1.9pt},
  edge/.style={draw=figureInk!75,line width=.75pt,
    line cap=round,line join=round},
  xedge/.style={edge,draw=figurePath,line width=.8pt},
  every label/.style={inner sep=1.5pt}
]
  \coordinate (x) at (3.5,1.55);
  \coordinate (y) at (0,0);
  \coordinate (s1) at (1.1,0);
  \coordinate (p1) at (1.95,.55);
  \coordinate (q1) at (1.95,-.55);
  \coordinate (s2) at (2.8,0);
  \coordinate (w) at (3.5,0);
  \coordinate (s3) at (4.2,0);
  \coordinate (p2) at (5.05,.55);
  \coordinate (q2) at (5.05,-.55);
  \coordinate (z) at (5.9,0);
  \draw[edge] (y)--(s1) (s1)--(p1)--(s2)--(q1)--cycle (s2)--(w)--(s3)
    (s3)--(p2)--(z)--(q2)--cycle;
  \draw[xedge] (x)--(p1) (x)--(w) (x)--(p2);
  \draw[xedge] (x) .. controls (-1.9,1.75) and (-1.5,-1.35) .. (q1);
  \draw[xedge] (x) .. controls (8.9,1.75) and (8.5,-1.35) .. (q2);
  \foreach \v in {s1,p1,q1,s2,w,s3,p2,q2} \node[vertex] at (\v) {};
  \node[root,label=above:{$x$}] at (x) {};
  \node[root,label=left:{$y$}] at (y) {};
  \node[root,label=right:{$z$}] at (z) {};
  \node at (.55,-.24) {$H_1$};
  \node at (1.95,0) {$H_2$};
  \node at (3.15,-.24) {$H_3$};
  \node at (3.85,-.24) {$H_4$};
  \node at (5.05,0) {$H_5$};
\end{tikzpicture}%
}
\caption{A bridge $B$ as in Proposition~\ref{prop:three-root-blocks}, where $B-x$ has five blocks.}
\label{fig:three-root-chain}
\end{figure}

Now we will describe the bridges with three attachments and a
planar two-edge completion (as in Figure~\ref{fig:three-root-chain}).

\begin{proposition}\label{prop:three-root-blocks}
Let $B$ be a bridge of $G$ on $A=\{x,y,z\}$ whose triangle completion
is three-connected and a minor of $G$, and suppose $B+xy+xz$ is
planar. Then the blocks of $B-x$ have the following properties.
\begin{enumerate}[label=\textup{(\roman*)}]
\item The blocks are edges or four-cycles. They can be
ordered as $H_1,\dots,H_n$ so that $y$ lies only in $H_1$, $z$ lies
only in $H_n$, consecutive blocks share exactly one vertex, and
nonconsecutive blocks are vertex-disjoint.
\item Each four-cycle $H_j$ has the form $apbqa$, where $a=y$ for
$j=1$ and $a\in V(H_{j-1})$ for $j>1$, while $b=z$ for $j=n$ and
$b\in V(H_{j+1})$ for $j<n$. The vertices $p,q$ have degree three
in $B$ and are adjacent to $x$.
\item No three consecutive blocks are edges.
\item The neighbors of $x$ in $B$ are the vertices $p,q$ of the
four-cycles and the vertices shared by two edge blocks.
\item If every block is an edge, then $B$ has exactly one
interior vertex, adjacent to all three attachments $x,y,z$.
\end{enumerate}
We say that the blocks form a \emph{chain} if they can be ordered
so that consecutive blocks share exactly one vertex and
nonconsecutive blocks are vertex-disjoint.
\end{proposition}
\begin{proof}
We prove (i) and the first sentence of (ii) by induction on the
number of interior vertices of $B$. Let $H=B+xy+xz$, which is
planar by hypothesis. Suppose first that $B\cup K_3[A]$ is planar.
If $G-A$ were connected, every vertex outside $A$ would be interior
to $B$. Then every edge of $G$ outside $B$ would join two vertices
of $A$, so $G$ would be a subgraph of the planar graph $B\cup K_3[A]$,
a contradiction. Thus $A$ is a three-cut. By
Lemma~\ref{lem:flat-replacement}, the interior of $B$ is a single
vertex $w$. Each attachment is adjacent to $w$, so $B-x$ is the
path $y-w-z$.

Now suppose $B\cup K_3[A]$ is nonplanar. It is three-connected
by hypothesis and has no $K_5$ minor since it is a minor of $G$.
By Lemmas~\ref{core:three-root-setup} and~\ref{core:three-root-cut},
it has a three-cut $X$ of the form $\{x,u,v\}$, $\{x,y,u\}$, or
$\{x,z,u\}$. The components of $(B\cup K_3[A])-X$ are $D_0$ and
two single vertices $p,q$, both interior to $B$ and both with
neighborhood $X$. Write $X=\{x,a,b\}$. Since $X$ is the
neighborhood of $p$ in $G$, it is independent in $G$ by
Lemma~\ref{lem:independent-neighborhood}(a). Thus $a,p,b,q$
induce the four-cycle $C=apbqa$ in $B-x$. Deleting $x$ and the
vertices of $C$ from $B$ leaves $B[V(D_0)]$.

In case~(i) of Lemma~\ref{core:three-root-cut}, $B[V(D_0)]$ has
two components. After interchanging $u,v$ if necessary, the
component containing $y$ has $u$ as its only neighbor on $C$,
and the one containing $z$ has $v$ as its only neighbor on $C$.
In case~(ii), $X=\{x,y,u\}$. Deleting $X$ removes all the added
triangle edges, so $D_0$ is a component of $B-X$ and contains
$z$. Its only possible neighbors in $B-x$ are $y,u$. By
Lemma~\ref{core:three-root-cut}(ii), $y$ has no neighbor in
$D_0-z$, and $yz\notin E(B)$. There must be an edge from $D_0$
to $u$ in the triangle completion, since otherwise deleting $x,y$
would separate $D_0$ from $u$. This edge belongs to $B$ since
$u\notin A$. Thus $u$ is the only neighbor of $D_0$ on $C$.
For Case (iii) interchange $y,z$.

In all three cases, each component of $(B-x)-V(C)$ has just one
neighbor on $C$, so $C$ is a block of $B-x$. Interchange $a,b$
if necessary so that $a=y$ when $y$ lies on $C$, and otherwise
$a$ is the neighbor on $C$ of the component containing $y$.
If $a=y$, let $H_L$ be the single vertex $y$. Otherwise let
$H_L$ consist of the component containing $y$ together with $a$
and all its edges to that component. Define $H_R$ in the same
way with $z,b$ in place of $y,a$. Then $H_L$ and $H_R$ are
disjoint, meet $C$ only at $a$ and $b$, respectively, and
together with $C$ form $B-x$.

If $H_R$ has at most two vertices, it is either the single vertex
$z$ or the edge $bz$. Otherwise, let $D$ be a component of
$H_R-\{b,z\}$. Its vertices are interior to $B$, and its only
possible neighbors in $G$ are $x,b,z$. Since $p$ lies outside
$V(D)\cup\{x,b,z\}$, three-connectivity gives $N_G(D)=\{x,b,z\}$.
Thus $D$ is a component of $G-\{x,b,z\}$. There are at most two
components of $H_R-\{b,z\}$, since three would give a $K_{3,3}$
minor of $H$.

If there are two components, contract either one onto $b$ in $H$.
This gives the edges $xb,bz$, and $xz$ is already in $H$. Thus
the triangle completion of the bridge of the other component is
planar. By Lemma~\ref{lem:flat-replacement}, both components are
single vertices. Their common neighborhood $\{x,b,z\}$ is
independent by Lemma~\ref{lem:independent-neighborhood}(a), so
$bz\notin E(B)$. Hence $H_R$ is a four-cycle with $b,z$ opposite.

Suppose now that $D$ is the only component of $H_R-\{b,z\}$.
Let $B'$ be its bridge in $G$ on $\{x,b,z\}$. Contracting $xp$
in $H$ gives $B'+xb+xz$ as a planar minor. If $bz$ were an edge,
the same contraction would give the triangle completion of $B'$
as a planar minor. Then $D$ would be a single vertex by
Lemma~\ref{lem:flat-replacement}, and its neighbors $b,z$ would
be adjacent, contradicting
Lemma~\ref{lem:independent-neighborhood}(a). Thus $bz\notin E(B)$
and $H_R=B'-x$.

In $B\cup K_3[A]$, follow $bqa$, then an $a$--$y$ path in $H_L$,
then the edge $yz$ to obtain a $b$--$z$ path. This path and $xpb$
meet only at $b$, and their internal vertices lie outside $B'$.
Contracting each path to an edge and retaining $xz$ gives the
triangle completion of $B'$ as a minor of $B\cup K_3[A]$, hence
of $G$. It is three-connected by Lemma~\ref{lem:triangle-completion}.
We have already shown that $B'+xb+xz$ is planar. Every interior
vertex of $B'$ is interior to $B$, and $B'$ omits $p,q$. By
induction, the blocks of $H_R=B'-x$ satisfy (i) and the first
sentence of (ii) with $b$ in place of $y$.

Interchanging $y,z$ and $a,b$ gives the same argument for $H_L$.
Take the blocks of $H_L$ from $y$ to $a$, then $C$, then the
blocks of $H_R$ from $b$ to $z$, omitting either side if it is a
single vertex. The last block of $H_L$, when present, meets $C$
only at $a$, and the first block of $H_R$, when present, meets
$C$ only at $b$. Since $a,b$ are opposite on $C$, this proves
(i) and the first sentence of (ii).

The vertices $p,q$ of each four-cycle and the vertices shared by
two edge blocks have degree two in $B-x$. They are interior to
$B$, so they have degree at least three in $B$. Thus each is
adjacent to $x$ and has degree three in $B$. In each four-cycle
$apbqa$, the neighborhood $\{x,a,b\}$ of $p$ is independent in
$G$ by Lemma~\ref{lem:independent-neighborhood}(a), so neither
$a$ nor $b$ is adjacent to $x$. This accounts for every vertex
except possibly $y,z$, and $xy,xz\notin E(B)$ since $x,y,z$ are
attachments. This proves (ii) and (iv).

\enlargethispage{\baselineskip}
If three consecutive blocks were edges, let $st$ be the middle
edge. Both $s$ and $t$ are shared by two edge blocks, so both
are adjacent to $x$. Then $s$ has degree three and its neighbors
$x,t$ are adjacent, contradicting
Lemma~\ref{lem:independent-neighborhood}(a). This proves (iii).
If every block is an edge, there are at most two blocks by (iii).
A single block would be the edge $yz$, which is absent in $B$.
Thus $B-x$ is the path $y-w-z$, and (iv) gives $wx\in E(B)$.
This proves (v).
\end{proof}

\subsection{Nonplanar completions}\label{sec:nonplanar-completions}

By Lemma~\ref{core:relevant}, no proper allowed minor of a bridge of $G$
keeps all its nonplanar relevant completions nonplanar. So if a bridge of $G$ has an allowed minor that does, the bridge must itself be equal to it.

\begin{lemma}\label{lem:tripod-contraction}
Contracting each attachment path of a tripod into its end in the
attachment triple $A$, and each of the six paths of its $K_{2,3}$
subdivision to an edge, gives two vertices with neighborhood $A$.
\end{lemma}

\begin{lemma}\label{endpoint:primitive}
Let $B$ and $A=\{x,y,z\}$ be as in Section~\ref{core:three-root-section},
with $B\cup K_3[A]$ three-connected and nonplanar. Suppose $B$ lies in a
graph with no $K_5$ minor that also contains three paths joining the
pairs of vertices of $A$, internally disjoint from each other and from
$B$. Then $B$ has $T(r)$ as an allowed minor for some $r\in A$.
\end{lemma}
\begin{proof}
Since $B\cup K_3[A]$ is nonplanar, it has no embedding in a closed disk
with its triangle $C$ on $A$ as the boundary. Applying
Theorem~\ref{thm:mohar}, we see that outcome (a) would require $A$
to have four vertices. Also, the auxiliary graph
$\operatorname{Aux}(B\cup K_3[A],C)$ is three-connected, since
$B\cup K_3[A]$ is and the added vertex has three neighbors. Thus it has
no three-connected component other than the one containing $C$, so
outcome (c) cannot occur. Therefore $(B\cup K_3[A])-E(C)=B$ contains a
tripod.

Call the two degree-three vertices of the tripod's $K_{2,3}$ subdivision
$p,q$. Suppose all three of its attachment paths were trivial, so that
$x,y,z$ correspond to the three degree-two vertices of $K_{2,3}$. Contract each
of its six paths to a single edge, placing all contractions in
$\operatorname{Int}(B)$. Then $p,q$ are each adjacent to $x,y,z$.
Since $B[\operatorname{Int}(B)]$ is connected, it remains so after these
contractions, so it must still contain a $p$--$q$ path that avoids
$x,y,z$. Shorten this path to an edge, and do the same with the three
paths joining the pairs of vertices of $A$ that exist by hypothesis.
This gives a $K_5$ minor, a contradiction.

Thus at least one attachment path has positive length. After relabeling
$A$, let it end at $z$. Contract the rest of the interior of $B$ onto
the tripod and delete the edges outside it, so that every later step
keeps the interior connected and adjacent to all of $A$. Shorten the
chosen attachment path to an edge $uz$, contract the other two
attachment paths into their ends $x,y$, and shorten the six paths of
the $K_{2,3}$ subdivision to edges. The remaining graph has $p,q$ each
adjacent to $x,y,u$, together with the edge $uz$. This is $T(z)$.
\end{proof}

\begin{lemma}\label{core:two-cap-critical}
Let $B$ and $A=\{x,y,z\}$ be as in Section~\ref{core:three-root-section},
and suppose that
\begin{enumerate}[label=\textup{(\roman*)}]
\item $B\cup K_3[A]$ is three-connected and has no $K_5$ minor,
\item $B+xy$ is planar, but $B+xy+yz$ and $B+xy+xz$ are nonplanar, and
\item no proper allowed minor of $B$ keeps both $B+xy+yz$ and
$B+xy+xz$ nonplanar.
\end{enumerate}
Then $B=Q(z;x,y)$, or, after possibly interchanging $x,y$, $B-A$ is
the four-cycle $a_1b_1a_2b_2a_1$ and the remaining edges of $B$ are
$za_1,za_2,yb_1,yb_2,xa_2$. If $B+xz+yz$ is also nonplanar, then
$B\ne Q(z;x,y)$.
\end{lemma}
\begin{proof}
Contract $xy$ in $B\cup K_3[A]$ and call the merged vertex $r$.
By (ii) and Lemma~\ref{r2:planar-splitting} applied to $B+xy$,
the resulting graph is nonplanar. It has no $K_5$ minor by (i),
so it contains a $K_{3,3}$ subdivision. Fix such a subdivision.
It uses $rz$, since deleting $rz$ leaves the planar graph
$(B+xy)/xy$.

Deleting $r,z$ from the subdivision leaves a connected subgraph
of $B-A$ containing at least four corners. To see this, $r,z$
are consecutive on one branch path. Deleting consecutive vertices
on a single branch path leaves the remaining corners connected,
since $K_{3,3}$ stays connected after deleting an edge or at most
two vertices.

In $B-A$, contract each component outside the vertex set of the
subdivision into one of its neighbors on the subdivision. For each
edge $rv$ of the subdivision other than $rz$, choose $xv$ or $yv$,
using both $x$ and $y$ if possible. If one of $x,y$ has no chosen
edge, choose one more edge from it to the interior. Let $B'$ consist
of the chosen edges and the edges of the subdivision not incident
with $r$. Then $B'-A$ is the subdivision with $r,z$ deleted,
and every attachment has an interior neighbor, so $B'$ is an
allowed minor of $B$.

In $B'+xy+yz$, contracting $xy$ restores the chosen edges at $r$
and turns $yz$ into $rz$, giving back the subdivision. The same
holds in $B'+xy+xz$, since $xz$ also becomes $rz$. Thus both
completions of $B'$ are nonplanar. By (iii), $B'$ cannot be a
proper allowed minor, so the vertices and edges we kept are all
of $B$.

Suppose one of $r,z$ has degree two in the subdivision, and let
$s$ be its neighbor outside $\{r,z\}$. If the degree-two vertex
is $r$, contract the edge of $B$ chosen for $rs$. If it is $z$,
contract $zs$. The attachment on the edge we contract has $s$
as its only interior neighbor.

If another attachment also had $s$ as its only interior neighbor,
deleting $s$ and the third attachment would separate the first two
from the rest of $B\cup K_3[A]$, contradicting (i). Thus the other
two attachments keep an interior neighbor, and the merged attachment
does too since $s$ has a neighbor in $B-A$. The remaining interior
is connected, since it is obtained from the subdivision by deleting
$r,z,s$, three consecutive vertices on one branch path. Thus the
contraction is allowed. In either completion, contracting $xy$ now
leaves a $K_{3,3}$ subdivision with that branch path one edge shorter.
Both completions are therefore nonplanar, contradicting (iii). Hence
$r,z$ are both corners.

Any degree-two vertex of the subdivision is now interior to $B$.
It has an interior neighbor, since otherwise its two neighbors
would be $r,z$, forming a triangle with $rz$. Contracting the edge
to that neighbor keeps the interior connected and all three
attachments adjacent to it. This again shortens a branch path and
leaves both completions nonplanar. By (iii) there is no such vertex,
so the subdivision is $K_{3,3}$ itself.

Since $r,z$ are adjacent, write the two parts as $\{r,a_1,a_2\}$
and $\{z,b_1,b_2\}$. Then $B-A$ is the four-cycle $a_1b_1a_2b_2a_1$,
and $z$ is adjacent to $a_1,a_2$. If the edges chosen for $rb_1,rb_2$
use both $x$ and $y$, then $B=Q(z;x,y)$. Otherwise, after
interchanging $x,y$ if necessary, the chosen edges are $yb_1,yb_2$.
Since we chose edges from both attachments whenever possible, $x$
is adjacent to neither $b_1$ nor $b_2$. Its only edge therefore
joins it to $a_1$ or $a_2$, and we may call that vertex $a_2$.
Finally, $Q(z;x,y)+xz+yz$ is planar by
Lemma~\ref{core:three-graph-completions}, so $B\ne Q(z;x,y)$ if
$B+xz+yz$ is nonplanar.
\end{proof}

We now use minimality of $G$ to identify
the bridges themselves for the following completion tests.

\begin{proposition}\label{prop:three-root-disjoint}
Let $B$ be a bridge of $G$ on $A=\{x,y,z\}$ whose triangle completion
is three-connected, nonplanar and a minor of $G$. Suppose its relevant
completions are $B\cup K_3[A]$, $B+xy+xz$, $B+xy+yz$ and possibly $B+xy$,
at least one is planar, and every proper allowed minor of $B$ makes at least one nonplanar relevant completion planar. Then $B$ is $T(x)$, $T(y)$, $T(z)$, $S(x;y)$ or $S(y;x)$, or
$B+xy$ is relevant and $B$ is $Q(z;x,y)$ or a four-cycle form of
Lemma~\ref{core:two-cap-critical}.
\end{proposition}
\begin{proof}
Suppose first that at least one of $B+xy+xz$ and $B+xy+yz$ is
planar. Interchange $x,y$ if necessary so that $B+xy+xz$ is
planar. If $B$ is $T(y)$ or $T(z)$, we are done. Otherwise,
Lemma~\ref{core:three-root-cut} gives an allowed minor $Q(x;y,z)$,
$S(x;y)$ or $S(x;z)$. Contracting $vz$ into $z$ in $Q(x;y,z)$
gives $T(y)$. Contracting $fy$ into $y$ in $S(x;z)$ and deleting
the resulting edge $xy$ also gives $T(y)$. Thus $B$ has $T(y)$
or $S(x;y)$ as an allowed minor.

For both $T(y)$ and $S(x;y)$, adding $xy,yz$ gives a nonplanar
graph by Lemma~\ref{core:three-graph-completions}, so their triangle
completions are also nonplanar. Since $B+xy+xz$ and $B+xy$ are
planar, the allowed minor keeps every nonplanar relevant completion
of $B$ and hence equals $B$ by minimality. Hence,
the possibilities are $T(x)$, $T(y)$, $T(z)$, $S(x;y)$, $S(y;x)$.

Now suppose both $B+xy+xz$ and $B+xy+yz$ are nonplanar. Since at
least one relevant completion is planar, $B+xy$ must be relevant
and planar. A proper allowed minor keeping both two-edge completions
nonplanar would also keep the triangle completion nonplanar,
contradicting minimality. Thus Lemma~\ref{core:two-cap-critical}
gives $B=Q(z;x,y)$ or one of the two four-cycle forms.
\end{proof}

When two relevant completions are nonplanar, we look for a reduction of the bridge
that keeps both nonplanar. The next lemma shows we can find a $K_{3,3}$
subdivision for which replacing $xy$ by $xz,yz$ still gives
a nonplanar graph.

\begin{lemma}\label{r2:three-mark-entrance}
Let $B$ and $A=\{x,y,z\}$ be as in Section~\ref{core:three-root-section},
with $B$ planar. Suppose $B\cup K_3[A]$ is two-connected, has no
$K_5$ minor, and every two-cut contains $z$. If $B+xy$ and
$B+xz+yz$ are nonplanar, then $B+xy$ contains a $K_{3,3}$
subdivision $H$ using $xy$ for which $H-xy+xz+yz$ is nonplanar.
\end{lemma}
\begin{proof}
Since $B$ is planar, every $K_{3,3}$ subdivision in $B+xy$ uses
$xy$. For any such subdivision, let $C$ be the cycle formed by
the four branch paths avoiding the corner ends of the branch
through $xy$. By Lemma~\ref{r2:position-of-third}, it suffices to
find a subdivision with $z\notin V(C)$.

We use induction on $|V(B)|$. Suppose first that $B\cup K_3[A]$
has a two-cut $\{z,w\}$. Deleting two attachments leaves the
connected interior and the third attachment, so $w\notin A$.
The edge $xy$ keeps $x,y$ in the same component of
$(B\cup K_3[A])-\{z,w\}$. Choose another component $D$.
By two-connectivity, $D$ has neighbors at both $z,w$, and all
its vertices and incident edges belong to $B$. If $z$ has no
neighbor in $B$ outside $D$, the planar graph $B[V(D)\cup\{z,w\}]$
meets the rest of $B+xy$ only at $w$. Thus $(B+xy)-V(D)-z$ is
nonplanar. It has no $K_5$ minor, so it contains a $K_{3,3}$
subdivision avoiding $z$, as required.

Otherwise $B$ has a $z$--$w$ path avoiding $D$. Indeed, $z$ has an
interior neighbor outside $D$, and $B-A$ is connected with $w$ the
only neighbor of $D$ in the rest of $B-A$. Shortening this path
to $zw$ shows that $B[V(D)\cup\{z,w\}]+zw$ is a planar minor of
$B$. Contract $D$ onto $w$ to obtain $B_0$. This contraction is
within the interior, so the interior stays connected and each
attachment keeps an interior neighbor. Suppose $B_0+xy$ or
$B_0+xz+yz$ has a plane drawing. Take a small disk on one side of
$zw$ that meets the drawing exactly in the edge $zw$. Draw
$B[V(D)\cup\{z,w\}]+zw$ inside this disk using the same edge $zw$
on its boundary. We can do this since the graph is planar and we
may take a face incident with $zw$ as its outer face. This restores
$D$ and all its incident edges. Delete $zw$ if $zw\notin E(B)$.
We now have a plane drawing of $B+xy$ in the first case or of
$B+xz+yz$ in the second, contradicting the hypothesis that both
are nonplanar. Thus $B_0+xy$ and $B_0+xz+yz$ are both nonplanar.
Also, $B_0\cup K_3[A]$ is a minor of $B\cup K_3[A]$, and any cut
of at most two vertices in $B_0\cup K_3[A]$ also separates
$B\cup K_3[A]$, since the only neighbors of $D$ are $z,w$, now
joined by an edge. Thus $B_0$ satisfies the hypotheses. By induction
it has the required subdivision $H$ in $B_0+xy$. If $H$ uses a new
edge $zw$, replace it by a path through $D$. This subdivides the
same edge in $H$ and in $H-xy+xz+yz$, so the latter remains
nonplanar.

We may now assume $B\cup K_3[A]$ is three-connected. Since it is
nonplanar, Lemma~\ref{core:three-root-setup} with $E=\varnothing$
gives a three-cut $X\ne A$ and three components $D_0,D_1,D_2$ of
$(B\cup K_3[A])-X$. The components $D_1,D_2$ avoid $A$, and their
bridges in $B$ are the only bridges on $X$ meeting all of $X$.
By part~(vi) of that lemma, the bridge of $B+xy$ on $X$ containing
$xy$ either meets all of $X$ or has two attachments $\alpha,\beta$
and becomes a nonplanar minor of $B\cup K_3[A]$ when $\alpha\beta$
is added.

In the two-attachment case, the bridge containing $xy$ avoids $z$.
Otherwise it contains all of $A$, so even in $B\cup K_3[A]$ deleting
$\alpha,\beta$ separates its interior from $D_1$, contradicting
three-connectivity. Joining $\alpha,\beta$ gives a nonplanar graph with no $K_5$ minor,
so it contains a $K_{3,3}$ subdivision. If the subdivision uses
$\alpha\beta$, replace that edge by an $\alpha$--$\beta$ path
through $D_1$. The resulting subdivision lies in $B+xy$ and
avoids $z$.

Suppose the bridge containing $xy$ meets all of $X$, as do the
bridges of $D_1,D_2$. Each of these three bridges contains a tree
whose only vertices of degree one are the vertices of $X$, since
its interior is connected and has a neighbor at each vertex of $X$.
Together the three trees form a $K_{3,3}$ subdivision in $B+xy$.
Its branch through $xy$ lies in the first bridge, so the four paths
forming $C$ lie in the other two. Since $D_1,D_2$ avoid $A$, this proves the claim when
$z\notin X$.

It remains to consider $X=\{z,u,v\}$. The edges $xz,yz$ belong
to distinct bridges of $B+xz+yz$ on $X$. Indeed, an edge with
both ends in $X$ is a bridge by itself. If $x,y$ are both outside
$X$, they lie in different components of $B-X$, since otherwise
adding $xy$ would leave the attachments of their bridge unchanged
and it would already meet all of $X$ in $B$. Since $B+xz+yz$ is
nonplanar, Lemma~\ref{core:three-root-setup}(vi) gives a bridge
containing $xz$ or $yz$ that meets all of $X$ or becomes nonplanar
when its two attachments are joined. Interchange $x,y$ if necessary
so that this bridge contains $xz$. A single edge has neither
property, so $x\notin X$.

If this bridge meets all of $X$, the component of $B-X$ containing
$x$ has neighbors at $u,v$ and none at $z$, since its bridge in $B$
cannot meet all three. As the bridge containing $xy$ meets $z$, the
edge to $z$ must come from a different component containing $y$.
Take a $u$--$v$ path through the component containing $x$, and
join $x$ to that path within the component, stopping at the first
vertex on the path. From this vertex we have paths to $u,v$ along
the path and a third path back to $x$, across $xy$, and through the
component of $y$ to $z$. Using these with the six paths through
$D_1,D_2$ gives a subdivision whose branch through $xy$ ends at
$z$. Hence $z\notin V(C)$.

Otherwise write $\{z,u\}$ for the attachments of the bridge
containing $xz$. Joining them gives a nonplanar minor of
$B\cup K_3[A]$ by Lemma~\ref{core:three-root-setup}(vi), so it
contains a $K_{3,3}$ subdivision. This subdivision uses $xz$,
since deleting $xz$ from the completed bridge leaves a minor of
the planar graph $B$, using a $z$--$u$ path through $D_1$ for $zu$.
Since the bridge containing $xy$ meets $v$ and the bridge containing
$xz$ does not, either $y=v$ or the component of $B-X$ containing
$y$ has a neighbor at $v$. Replace $xz$ by $xy$, then a $y$--$v$
path through that component (trivial if $y=v$), then a $v$--$z$
path through $D_2$. If the subdivision uses $zu$, replace it by a
$z$--$u$ path through $D_1$. The replacement paths have disjoint
interiors outside the bridge containing $xz$, so they give a
subdivision in $B+xy$. Its branch through $xy$ contains $z$, and
hence $z\notin V(C)$.
\end{proof}

\subsection{Reducing a bridge while keeping a \texorpdfstring{$K_{3,3}$}{K3,3} subdivision}\label{sec:common-construction}

Once a completion is known to be nonplanar, the next lemma reduces any
bridge with at most five attachments while keeping a $K_{3,3}$
subdivision of that completion, and gives a bound on its number of interior vertices.

\begin{lemma}\label{core:retained-subdivision}
Let $B$ be a graph with an independent set $A$ of attachments, $1\le|A|\le5$, such that $B-A$ is nonempty, connected, and every vertex of $A$ has a neighbor in $B-A$. Let $H_0$ be a planar graph with $V(B)\cap V(H_0)=A$, such that every vertex of $H_0-A$ has all its neighbors in $A$. Let $H\subseteq B\cup H_0$ be a $K_{3,3}$ subdivision. Then $B$ has an allowed minor $B'$ such that $B'\cup H_0$ contains a $K_{3,3}$ subdivision $H'$ with the following properties.
\begin{enumerate}[label=\textup{(\roman*)}]
\item The same vertices of $A$ are corner vertices of $H$ and $H'$. Corresponding paths between corner vertices, with no corner vertices in their interiors, contain the same vertices of $A$ in the same order.
\item Every vertex of $B'-A$ lies on $H'$. On each such path in $H'$, the corner vertices and vertices of $A$ divide the path into segments. Each segment is the edge of $H_0$ joining its ends, if there is one. Otherwise it has exactly one vertex in its interior if both ends lie in $A$, and is a single edge if not.
\item The edges of $B'$ outside $H'$ consist of a minimal set of edges in $B'-A$ connecting the components of $H'-V(H_0)$, and one edge to $B'-A$ from each attachment with no such edge in $H'$. Deleting any edge of the first set would disconnect $B'-A$.
\end{enumerate}
In particular, $B'-A$ has at most $|A|+5$ vertices. If $B_0\cup H_0$ is planar for every graph $B_0$ obtained from $B$ by one allowed operation, then $B'=B$, so $B$ itself has this form and bound.
\end{lemma}
\begin{proof}
Split each branch path of $H$ at its internal vertices in $A$.
Whenever the ends of a segment are joined by an edge of $H_0$, use
that edge. It cannot occur elsewhere in $H$, since each branch path
is simple and two distinct branch paths share at most one vertex.
Every vertex of $H_0-A$ has all its neighbors in $A$, so a segment
meeting $H_0-A$ is already a single edge or a two-edge path between
attachments. For each remaining segment in $B$, contract its interior
into an end outside $A$ if it has one. If both ends lie in $A$, the
segment has an internal vertex since $A$ is independent in $B$.
In this case contract its interior to one vertex. This gives a
subdivision $H'$ with the same corners and the same attachments in
the same order on each branch path, proving (i).

All contractions so far are within $B-A$, so the interior stays
connected and each attachment keeps an interior neighbor. Since
$H_0$ is planar, $H'$ meets the interior. Thus each component of
the interior outside $H'$ has a neighbor in $H'-V(H_0)$. Contract
each such component into one of these neighbors. Keep the edges of
$H'$ in the contracted graph and a minimal set of interior edges
joining the components of $H'-V(H_0)$. For each attachment with no
edge to the interior in $H'$, keep one such edge. Delete every other
edge of the contracted graph to obtain $B'$. The edges kept in $B'$
keep its interior connected and give each attachment an interior
neighbor throughout the deletions. Thus all operations are allowed,
and (ii) and (iii) follow.

For the bound, each vertex of $B'-A$ other than a corner lies in a
segment between consecutive attachments, with at most one such vertex
per segment. Let $i$ be the number of corners in $A$. A branch path
with $k$ attachments in its interior has at most $k$ such segments
unless both ends lie in $A$, when it has $k+1$. At most $3i/2$ paths
have both ends in $A$, since each uses two of the $3i$ branch ends
at these corners. There are at most $|A|-i$ attachments in branch
interiors and at most $6-i$ corners in $B'-A$. Thus if $i>0$, then
\[
|V(B'-A)|\le (6-i)+(|A|-i)+3i/2=|A|+6-i/2<|A|+6.
\]
If $i=0$, each branch path meeting $A$ has one fewer such segment
than attachments. This gives at most $|A|-1$ segments in all, also
when $H'$ misses $A$ since $|A|\ge1$. Adding the at most six corners
gives $|V(B'-A)|\le |A|+5$ in this case as well.

Suppose $B_0\cup H_0$ is planar for every $B_0$ obtained from $B$
by one allowed operation. If we used any operations, the union with
$H_0$ would be planar after the first and would stay planar after
every further deletion or contraction. But $B'\cup H_0$ contains
the $K_{3,3}$ subdivision $H'$, so $B'=B$.
\end{proof}

\begin{figure}[htbp]
\centering
\tikzsetnextfilename{reduced-form-example}%
\begin{tikzpicture}[
  font=\small,text=figureInk,
  vertex/.style={circle,fill=figureInk,inner sep=1.5pt},
  root/.style={circle,fill=figurePath,inner sep=1.9pt},
  edge/.style={draw=figureInk!75,line width=.75pt,
    line cap=round,line join=round},
  added/.style={edge,draw=figureBridge,line width=.9pt},
  extra/.style={edge,dash pattern=on 2.2pt off 1.6pt},
  every label/.style={inner sep=1.5pt}
]
  \coordinate (a)  at (-2.75,0);
  \coordinate (v2) at (-1.38,.55);
  \coordinate (u2) at (1.38,.55);
  \coordinate (v1) at (2.75,0);
  \coordinate (u3) at (1.38,-.55);
  \coordinate (v3) at (-1.38,-.55);
  \coordinate (b) at (-2.30,1.30);
  \coordinate (w) at (-.15,1.65);
  \coordinate (c) at (2.30,1.30);
  \coordinate (d) at (2.72,-.95);
  \draw[edge] (a)--(v2)--(u2)--(v1)--(u3)--(v3)--cycle;
  \draw[edge] (v2)--(u3);
  \draw[edge,preaction={draw=white,line width=2.8pt}] (v3)--(u2);
  \draw[added] (a)--(b);
  \draw[edge] (b)--(w)--(c)--(v1);
  \draw[extra] (w)--(v1) (d)--(u3);
  \foreach \v in {u2,u3,v1,v2,v3} \node[vertex] at (\v) {};
  \node[vertex,label=above:{$w$}] at (w) {};
  \node[root,label=left:{$a$}] at (a) {};
  \node[root,label=above left:{$b$}] at (b) {};
  \node[root,label=above right:{$c$}] at (c) {};
  \node[root,label=right:{$d$}] at (d) {};
\end{tikzpicture}
\caption{A candidate from Construction~\ref{common:marked-construction} with attachments $a,b,c,d$. The edge $ab$ of $H_0$ is purple, and the edges added in step (iii) are dashed.}
\label{fig:reduced-form-example}
\end{figure}
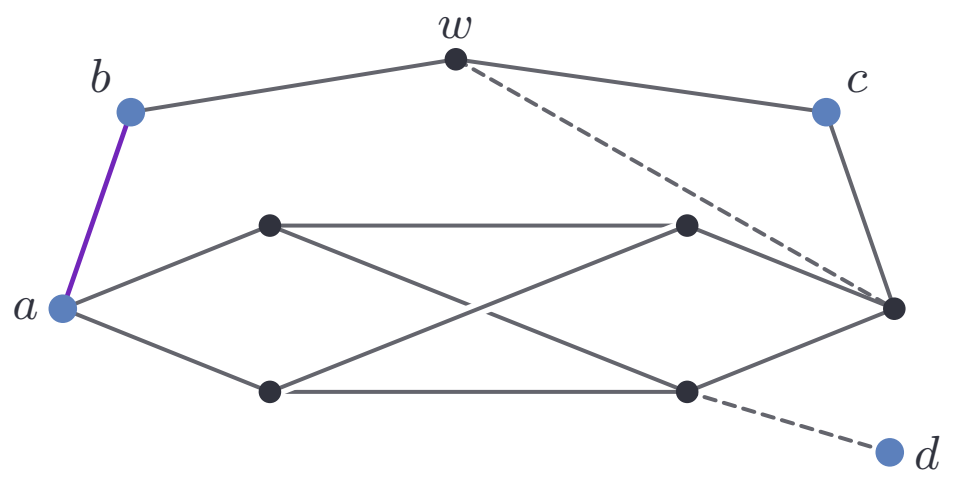

The following construction lists candidates for the reduced bridge
in Lemma~\ref{core:retained-subdivision}.
Figure~\ref{fig:reduced-form-example} shows one with four attachments.

\begin{construction}\label{common:marked-construction}
Let $B,A,H_0$ be as in Lemma~\ref{core:retained-subdivision}.
With $H_0$ fixed, we list candidates for the reduced bridge $B'$
by first choosing a $K_{3,3}$ subdivision $H'$ in $B'\cup H_0$.

\begin{enumerate}[label=\textup{(\roman*)}]
\item Choose which vertices of $H_0$ lie on $H'$. Place them at
distinct corners or internal vertices of its nine branch paths in
all possible ways, including every order along each path. Every
edge of $H'$ incident with a vertex of $H_0-A$ or joining two
attachments must belong to $H_0$.

Divide each branch path into segments at its internal vertices
in $A$. If the ends of a segment are adjacent in $H_0$, use that
edge. For each remaining segment, use a path of length two if
both ends lie in $A$ and a single edge otherwise. Keep only
subdivisions with at most $|A|+5$ vertices outside $H_0$.

\item To obtain the part of $H'$ that $B'$ must supply, take a copy
of $H'\cup H_0$ and delete the vertices of $H_0-A$ and all edges of
$H_0$. This keeps every vertex of $A$, including those outside $H'$.
At least one vertex of $H'$ remains outside $A$, since $H_0$ is
planar.

\item The interior of $B'$ must be connected. If the interior of
the graph from (ii) has $r$ components, we need $r-1$ edges to
connect them. We take every choice of $r-1$ edges that joins
these components in a tree. Each attachment of $B'$ must also
have an interior neighbor. For each attachment that still has
none, we add one edge to an interior vertex, considering all
choices of these neighbors.
\end{enumerate}

The union of each candidate with $H_0$ contains $H'$, and
Lemma~\ref{core:retained-subdivision} shows that every possible
reduced bridge $B'$ occurs. The bound in (i) makes the list finite.
We keep one representative up to isomorphism fixing each attachment
and any further vertices specified in the application.

Each application states which edges the subdivision must use and
which edges may belong to the original graph. If an original edge
has the same ends as an edge supplied by $H_0$, step (ii) keeps
the original edge. If a segment between two attachments may be an
original edge or a path of length two, we carry out (ii) and (iii)
first with the edge and then with a two-edge path in its place.

When $a,b$ have been identified to $x$, we restore them before
testing the candidates. In Lemma~\ref{r2:32}, each retained edge
$xv$ is replaced by $av$ or $bv$, with one further edge to the
interior if either attachment still has no interior neighbor.
In Lemma~\ref{r2:ab-cut-apart}(iii), we replace each $xv$ by
$av$ alone, $bv$ alone, or both edges, taking all combinations
of these choices over the original neighbors of $x$.

When the interior may have several components, we partition the
interior components from (ii) according to which will lie in the
same component after the edges are added. We take every partition
with the required vertices together or apart and the number of
components within the stated bounds. Within each class of the
partition we use a tree as in (iii), adding one original edge
for each pair of components adjacent in that tree.
If the case permits further components outside $H'$, we consider
every number of such components from zero to the stated bound.
For each choice, we add one vertex per component.
For each such vertex we try each subset of $A$ that it may have
as its neighborhood. Finally, we consider each set of further
original edges that supplies the required paths and attachment
neighbors, makes the required subgraphs connected, and leaves
vertices in different classes of the partition in different
components after deleting $A$. We keep the graph if no proper
subset of these further edges does the same.

\end{construction}

In Figure~\ref{fig:reduced-form-example}, $A=\{a,b,c,d\}$ and $H_0$ consists of $A$ and the edge $ab$. In step~(i), one branch path starts at the corner $a$, follows the
edge $ab$ of $H_0$ and the path $bwc$, and ends with the edge from
$c$ to the next corner. Every other branch path is a single edge,
and $d$ lies outside $H'$. Step~(ii) deletes $ab$. The vertex $d$
is isolated, and deleting $A$ leaves two components, one consisting
of $w$ and the other of the five remaining corners. Step~(iii) adds
the two dashed edges. One joins $w$ to the corner adjacent to $c$,
and the other gives $d$ an interior neighbor.

The next corollary gives conditions under which $B$ itself is
among the candidates.

\begin{corollary}\label{common:reduced-forms}
Let $B$ be a bridge of $G$ with nonempty interior and attachment set $A$, where $1\le|A|\le5$, and let $H_0$ be as in Lemma~\ref{core:retained-subdivision}, with $B\cup H_0$ nonplanar and without a $K_5$ minor. Suppose that $B_0\cup H_0$ is planar for every graph $B_0$ obtained from $B$ by one allowed operation, or that replacing $B$ in $G$ by an allowed minor $B'$ gives a nontoroidal graph whenever $B'\cup H_0$ is nonplanar. Then $B$ is one of the candidates of Construction~\ref{common:marked-construction} for $A$ and $H_0$, and has at most $|A|+5$ interior vertices.
\end{corollary}
\begin{proof}
The graph $B\cup H_0$ contains a $K_{3,3}$ subdivision, so Lemma~\ref{core:retained-subdivision} gives an allowed minor $B'$ of $B$ for which $B'\cup H_0$ also contains one. Under the first assumption, that lemma gives $B'=B$. Under the second, replacing $B$ by $B'$ gives a nontoroidal minor of $G$, so $B'=B$ by minimality. The construction lists every such $B'$.
\end{proof}

The constructions, counts and checks for the enumerations below
are given in~\supp{E1--E29}. Proposition~\ref{prop:membership} in Section~\ref{sec:completion} verifies that $\mathfrak{T}_{1},\ldots,\mathfrak{T}_{74}$ are obstructions and that each graph listed in Theorems~\ref{thm:disjoint}, \ref{thm:two-common} and~\ref{thm:one-common} has the properties stated there. We use this throughout without repeatedly referencing it each time.

\section{Two twin triples with no common vertex}\label{sec:disjoint}

\begin{figure}[htbp]
\centering
\begingroup
\scaledfigure{1}{%
\tikzsetnextfilename{disjoint-neighborhoods}%
\begin{tikzpicture}[
  font=\small,text=figureInk,
  vertex/.style={circle,fill=figureInk,inner sep=1.45pt},
  twinA/.style={circle,fill=figureTwinA,inner sep=1.7pt},
  twinB/.style={circle,fill=figureTwinB,inner sep=1.7pt},
  edgeA/.style={draw=figureTwinA,line width=.65pt,
    preaction={draw=white,line width=1.6pt}},
  edgeB/.style={draw=figureTwinB,line width=.65pt,
    preaction={draw=white,line width=1.6pt}},
  pathedge/.style={draw=figurePath,line width=.8pt},
  lab/.style={anchor=north,inner sep=0pt}
]
\begin{scope}[yshift=-.1cm]
\draw[draw=figureInk!45,fill=figureInk!3,line width=.6pt,rounded corners=4pt]
  (.05,-.92) rectangle (4.95,.17);
\draw[draw=figureTwinA,fill=figureTwinA!9,line width=.6pt,rounded corners=3pt]
  (.2,-.43) rectangle (2.3,.15);
\draw[draw=figureTwinB,fill=figureTwinB!9,line width=.6pt,rounded corners=3pt]
  (2.7,-.43) rectangle (4.8,.15);
\node at (2.5,-.69) {rest of $G$};
\foreach \i/\x in {1/.55,2/1.25,3/1.95}
  \coordinate (s\i) at (\x,0);
\foreach \i/\x in {1/3.05,2/3.75,3/4.45}
  \coordinate (t\i) at (\x,0);
\foreach \v/\x in {p/.85,q/1.65}{
  \foreach \u in {s1,s2,s3} \draw[edgeA] (\x,1.05)--(\u);
  \node[twinA,label=above:{$\v$}] at (\x,1.05) {};
}
\foreach \v/\x in {r/3.35,z/4.15}{
  \foreach \u in {t1,t2,t3} \draw[edgeB] (\x,1.05)--(\u);
  \node[twinB,label=above:{$\v$}] at (\x,1.05) {};
}
\foreach \i in {1,2,3}{
  \node[vertex] at (s\i) {};
  \node[lab] at ([yshift=-.1cm]s\i) {$s_\i$};
  \node[vertex] at (t\i) {};
  \node[lab] at ([yshift=-.1cm]t\i) {$t_\i$};
}
\end{scope}
\begin{scope}[xshift=6.15cm]
\draw[draw=figureTwinA,fill=figureTwinA!9,line width=.6pt,rounded corners=3pt]
  (1.08,-1.54) rectangle (1.92,1.34);
\draw[draw=figureTwinB,fill=figureTwinB!9,line width=.6pt,rounded corners=3pt]
  (3.33,-1.54) rectangle (4.17,1.34);
\foreach \i/\y in {1/1.15,2/0,3/-1.15}{
  \coordinate (s\i) at (1.35,\y);
  \coordinate (t\i) at (3.9,\y);
  \draw[pathedge] (s\i)--node[midway,above=3pt] {$P_\i$} (t\i);
}
\foreach \v/\y in {p/.55,q/-.55}{
  \foreach \u in {s1,s2,s3} \draw[edgeA] (.05,\y)--(\u);
  \node[twinA,label=left:{$\v$}] at (.05,\y) {};
}
\foreach \v/\y in {r/.55,z/-.55}{
  \foreach \u in {t1,t2,t3} \draw[edgeB] (5.2,\y)--(\u);
  \node[twinB,label=right:{$\v$}] at (5.2,\y) {};
}
\foreach \i in {1,2,3}{
  \node[vertex] at (s\i) {};
  \node[anchor=north west,inner sep=0pt] at ([xshift=.1cm,yshift=-.1cm]s\i) {$s_\i$};
  \node[vertex] at (t\i) {};
  \node[anchor=north east,inner sep=0pt] at ([xshift=-.1cm,yshift=-.1cm]t\i) {$t_\i$};
}
\node[vertex] at (2.625,1.15) {};
\node[vertex] at (2.2,0) {};
\node[vertex] at (3.05,0) {};
\end{scope}
\end{tikzpicture}%
}
\endgroup
\caption{Left, the disjoint twin triples. Right, an example of $\Psi_1$, with the paths $P_i$ in blue.}
\label{fig:disjoint-neighborhoods}
\end{figure}
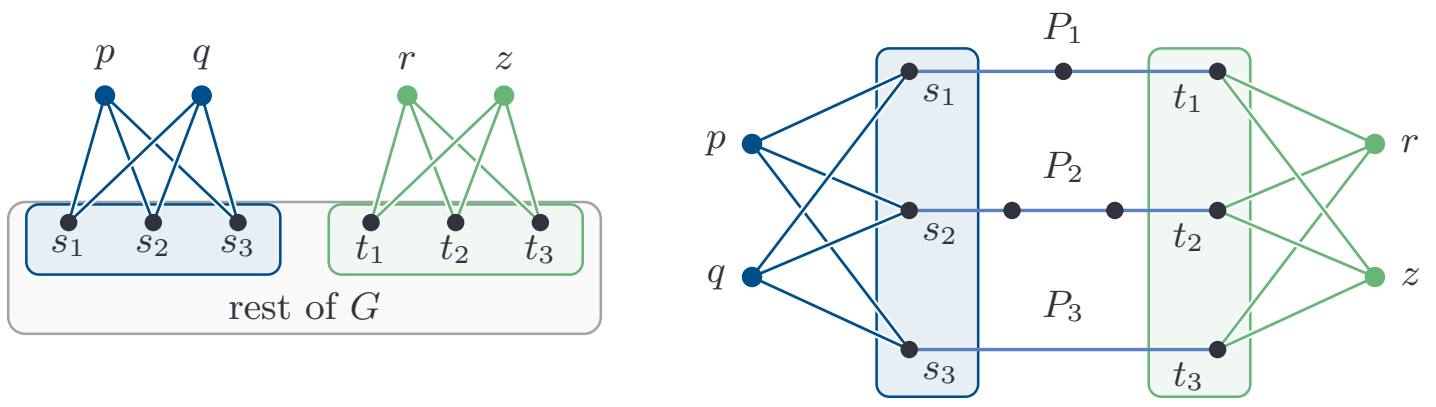

\looseness=-1
In this section we consider the case where some two twin triples in an obstruction are entirely disjoint. We show that $G$ must then
contain a subdivision $\Psi_1$ of a graph formed from them, and
determine $G$'s structure from the $\Psi_1$-bridges. Denote $G$'s twins as $p,q$ with $N(p)=N(q)=A_1=\{s_1,s_2,s_3\}$ and twins $r,z$ with $N(r)=N(z)=A_2=\{t_1,t_2,t_3\}$, where $A_1\cap A_2=\varnothing$. Write $J_1,J_2$ for the twin subgraphs
of $p,q$ and $r,z$, respectively. These ten vertices are distinct, as a vertex
adjacent to its twin would be its own neighbor, and twins with different
neighborhoods are not adjacent by
Lemma~\ref{lem:independent-neighborhood}(b). As $G$ is three-connected,
Menger's theorem gives, after possibly relabeling $A_2$, disjoint paths $P_i$ from
$s_i$ to $t_i$ meeting $A_1\cup A_2$ only at their ends. They avoid
$p,q,r,z$, as an internal vertex $p$, say, would have two neighbors on
the path, both in $A_1$. Choose the paths with minimum total length, and
let $\Psi_1=J_1\cup J_2\cup P_1\cup P_2\cup P_3$. Then
\begin{equation}\label{dis:shortest}
 P_i\ne s_it_i\quad\Longrightarrow\quad s_it_i\notin E(G),
\end{equation}
since otherwise replacing $P_i$ by the edge $s_it_i$ would give a
shorter choice. $\Psi_1$ is a subdivision of two copies of $K_{2,3}$ joined by
the matching $s_1t_1,s_2t_2,s_3t_3$, as illustrated in
Figure~\ref{fig:disjoint-neighborhoods}.

\begin{theorem}\label{thm:disjoint}
A graph is a three-connected obstruction in which some two twin triples are disjoint if and only if it is isomorphic to one of
\[
 \mathfrak{T}_{4},\mathfrak{T}_{5},\mathfrak{T}_{21},\mathfrak{T}_{32},\mathfrak{T}_{38},\mathfrak{T}_{43},\mathfrak{T}_{45},\mathfrak{T}_{48},
 \mathfrak{T}_{58},\mathfrak{T}_{61},\mathfrak{T}_{63},\mathfrak{T}_{66},\mathfrak{T}_{67},\mathfrak{T}_{68},\mathfrak{T}_{69}.
\]
\end{theorem}

Section~\ref{dis:tests} finds the faces of $\Psi_1$ in which the
$\Psi_1$-bridges lie and reduces toroidality to planarity tests.
Sections~\ref{sec:path-choice} and~\ref{dis:active-support} restrict the
$\Psi_1$-bridges, and Proposition~\ref{dis:types} collects the result. Any further vertex with neighborhood $A_1$ or $A_2$ gives
$G\cong\mathfrak{T}_5$ by Lemma~\ref{dis:noncofacial}.
Section~\ref{dis:finite-method} reduces the bridges and bounds their number, and
Sections~\ref{dis:endpoint-support}--\ref{dis:path-four} handle
the remaining cases.

\subsection{Faces and embedding tests of \texorpdfstring{$\Psi_1$}{Psi1}}\label{dis:tests}

Any torus embedding of $G$ would restrict to an embedding of $\Psi_1$,
with each $\Psi_1$-bridge drawn in one of its faces
(Section~\ref{sec:extension}). We therefore determine the possible
face boundaries of $\Psi_1$ and test whether all its bridges can be
drawn together in those faces. 

\begin{lemma}\label{dis:faces}
Every torus embedding of $\Psi_1$ is cellular and has five faces.
For some $k\in\{1,2,3\}$, write
$\{i,j\}=\{1,2,3\}\setminus\{k\}$.
Omitting the twins from the boundary walks, these faces are
\begin{enumerate}[label=\textup{(\roman*)}]
\item two faces, each following
$s_iP_it_i,t_jP_j^{-1}s_j$ and meeting no other path.
\item one face following
$s_kP_kt_k,t_y,t_kP_k^{-1}s_k,s_x$, where
$x,y\in\{i,j\}$ may be equal.
\item two other faces meeting the paths exactly at
$\{s_k,s_{x'}\}$ and $\{t_k,t_{y'}\}$, respectively,
where $\{x,x'\}=\{y,y'\}=\{i,j\}$.
\end{enumerate}
All twelve choices of $k$, $x$ and $y$ occur.
\end{lemma}

\begin{figure}[htbp]
\centering
\begingroup
\scaledfigure{0.9}{\tikzsetnextfilename{F-torus}%
\begin{tikzpicture}[
  scale=1.4,
  font=\small,
  vertex/.style={circle,fill=figureInk,inner sep=1.6pt},
  edge/.style={draw=figureInk!75,line width=.8pt},
  pathedge/.style={draw=figurePath,line width=1.1pt},
  side/.style={draw=figureInk!35,line width=.6pt}
]
\foreach \dx in {0,4.8}{
  \begin{scope}[xshift=\dx cm]

  \fill[faceRepeated] (0,0) rectangle (4,3);
  \fill[faceSimple] (0,1) rectangle (4,2);
  \fill[faceOther] (1.5,1)--(3.5,1)--(3,.45)--cycle;
  \fill[faceOther] (1.5,2)--(3.5,2)--(3,2.55)--cycle;
  \draw[side] (0,0) rectangle (4,3);
  \draw[side,-{Stealth[length=5pt]}] (0,0.25)--(0,0.65);
  \draw[side,-{Stealth[length=5pt]}] (4,0.25)--(4,0.65);
  \draw[side,-{Stealth[length=5pt] Stealth[length=5pt]}]
    (1.7,0)--(2.3,0);
  \draw[side,-{Stealth[length=5pt] Stealth[length=5pt]}]
    (1.7,3)--(2.3,3);
  \end{scope}}
\begin{scope}
  \coordinate (s1) at (0.5,1);
  \coordinate (p) at (1.5,1);
  \coordinate (s2) at (2.5,1);
  \coordinate (q) at (3.5,1);
  \coordinate (t1) at (0.5,2);
  \coordinate (r) at (1.5,2);
  \coordinate (t2) at (2.5,2);
  \coordinate (z) at (3.5,2);
  \coordinate (s3) at (3.0,0.45);
  \coordinate (t3) at (3.0,2.55);
  \draw[edge] (0,1)--(s1)--(p)--(s2)--(q)--(4,1);
  \draw[edge] (0,2)--(t1)--(r)--(t2)--(z)--(4,2);
  \draw[edge] (p)--(s3)--(q);
  \draw[edge] (r)--(t3)--(z);
  \draw[pathedge] (s1)--(t1);
  \draw[pathedge] (s2)--(t2);
  \draw[pathedge] (s3)--(3,0);
  \draw[pathedge] (3,3)--(t3);
  \node[vertex,label=below:{$s_1$}] at (s1) {};
  \node[vertex,label=above:{$p$}] at (p) {};
  \node[vertex,label=above left:{$s_2$}] at (s2) {};
  \node[vertex,label=above:{$q$}] at (q) {};
  \node[vertex,label=right:{$s_3$}] at (s3) {};
  \node[vertex,label=above:{$t_1$}] at (t1) {};
  \node[vertex,label=below:{$r$}] at (r) {};
  \node[vertex,label=below left:{$t_2$}] at (t2) {};
  \node[vertex,label=below:{$z$}] at (z) {};
  \node[vertex,label=right:{$t_3$}] at (t3) {};
  \node at (1.5,1.5) {(i)};
  \node at (3.5,1.5) {(i)};
  \node at (1.1,0.45) {(ii)};
  \node at (1.1,2.55) {(ii)};
  \node at (2.75,0.8) {(iii)};
  \node at (2.75,2.2) {(iii)};
\end{scope}
\begin{scope}[xshift=4.8cm]
  \coordinate (s1) at (0.5,1);
  \coordinate (p) at (1.5,1);
  \coordinate (s2) at (2.5,1);
  \coordinate (q) at (3.5,1);
  \coordinate (t2) at (0.5,2);
  \coordinate (r) at (1.5,2);
  \coordinate (t1) at (2.5,2);
  \coordinate (z) at (3.5,2);
  \coordinate (s3) at (3.0,0.45);
  \coordinate (t3) at (3.0,2.55);
  \draw[edge] (0,1)--(s1)--(p)--(s2)--(q)--(4,1);
  \draw[edge] (0,2)--(t2)--(r)--(t1)--(z)--(4,2);
  \draw[edge] (p)--(s3)--(q);
  \draw[edge] (r)--(t3)--(z);
  \draw[pathedge] (s1)--(t1);
  \draw[pathedge] (s2)--(4,1.75);
  \draw[pathedge] (0,1.75)--(t2);
  \draw[pathedge] (s3)--(3,0);
  \draw[pathedge] (3,3)--(t3);
  \node[vertex,label=below:{$s_1$}] at (s1) {};
  \node[vertex,label=below:{$p$}] at (p) {};
  \node[vertex,label=above left:{$s_2$}] at (s2) {};
  \node[vertex,label=below right:{$q$}] at (q) {};
  \node[vertex,label=right:{$s_3$}] at (s3) {};
  \node[vertex,label=above:{$t_2$}] at (t2) {};
  \node[vertex,label=above:{$r$}] at (r) {};
  \node[vertex,label=below right:{$t_1$}] at (t1) {};
  \node[vertex,label=above:{$z$}] at (z) {};
  \node[vertex,label=right:{$t_3$}] at (t3) {};
  \node at (2.05,1.45) {(i)};
  \node at (3.6,1.3) {(i)};
  \node at (1.1,0.45) {(ii)};
  \node at (1.1,2.55) {(ii)};
  \node at (2.6,0.75) {(iii)};
  \node at (2.6,2.25) {(iii)};
\end{scope}
\end{tikzpicture}}
\caption{The two torus embeddings of $\Psi_1$ up to automorphism, with additional endpoints $s_1,t_1$ (left) and $s_1,t_2$ (right). Blue marks simple faces (i), tan the repeated face (ii), and gray end faces (iii).}
\label{fig:F-torus}
\endgroup
\end{figure}
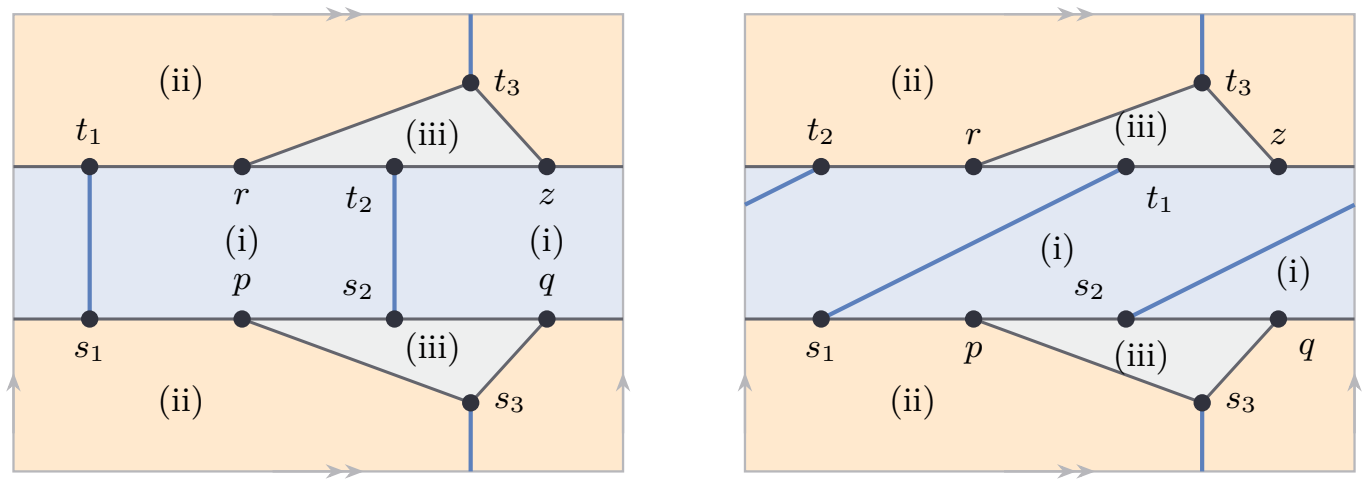

\begin{proof}
The graph $\Psi_1$ is connected and nonplanar, as $\Psi_1-z$ is a subdivision of
$K_{3,3}$ with parts $A_1$ and $\{p,q,r\}$, so its torus embeddings are
cellular by Youngs's theorem. Subdividing an edge does not change the faces,
so we may take each $P_i$ to be an edge. Of the $1{,}024$ rotation
systems of this graph, $48$ have genus one. Up to automorphisms, these
give the two embeddings in Figure~\ref{fig:F-torus}~\supp{R1}. Their boundary walks,
after permuting the indices, give (i)--(iii), with $x=y$ in one and
$x\ne y$ in the other. Permuting the indices of both triples
simultaneously is an automorphism of $\Psi_1$, so the first gives the six
choices with $x=y$ and the second the six with $x\ne y$.
\end{proof}

We call the faces in (i) the \emph{simple faces}, the face in (ii)
the \emph{repeated face}, and those in (iii) the \emph{end faces},
since they meet the paths only at their ends. In each embedding,
write $P_{\mathrm R}=P_k$ for the \emph{repeated path}, traversed twice
by the repeated face boundary walk. Lemma~\ref{dis:faces} shows that internal vertices of $P_{\mathrm R}$ only appear on the boundary walk of the repeated face, and internal vertices of the other paths only lie in the boundary walks of the simple faces. A $\Psi_1$-bridge attaching inside
$P_{\mathrm R}$ must therefore use the repeated face, and one attaching inside another path must use a simple face. We use the tuple $(P_{\mathrm R};s_x,t_y)$ to record the repeated path and
the two additional endpoints on its face. There are three possible
repeated paths. For $P_{\mathrm R}=P_k$, independently choose
$s_x\in\{s_i,s_j\}$ and $t_y\in\{t_i,t_j\}$, giving
$3\cdot2\cdot2=12$ combinations of repeated path and additional
endpoints. 

For each choice of repeated path and additional endpoints, we first determine which face boundaries contain each bridge's attachments. We then check whether the bridges assigned to each face can be drawn there together without crossings. For the simple faces and the repeated face, the next lemma gives planarity conditions.

\begin{lemma}\label{dis:joint-test}
Let $H$ be a graph containing $\Psi_1$ in which $p,q,r,z$ have degree three. For a choice $(P_{\mathrm R};s_x,t_y)$ from Lemma~\ref{dis:faces}, write $P_{\mathrm R}=P_k$ and $\{i,j\}=\{1,2,3\}\setminus\{k\}$. Consider an assignment of each $\Psi_1$-bridge to the pair of simple faces, the repeated face or one of the end faces, with all its attachments on the corresponding boundary~walks.

Let $H_{\mathrm{simp}}$ consist of the bridges assigned to the pair of simple faces together with $(P_i\cup P_j)+s_is_j+t_it_j$. Let $H_{\mathrm{rep}}$ consist of the bridges assigned to the repeated face together with $P_{\mathrm R}+s_ks_x+t_kt_y$. Keep each added edge separate from any bridge edge with the same ends. Then $H$ is toroidal if and only if one of the twelve choices has an assignment satisfying the following conditions.
\begin{enumerate}[label=\textup{(\alph*)}]
\item $H_{\mathrm{simp}}$ is planar.
\item $H_{\mathrm{rep}}$ is planar.
\item The bridges assigned to each end face can be drawn together without crossings in a disk, with their attachments on its boundary in the order given by that face's boundary walk.
\end{enumerate}
If $A_1$ and $A_2$ are independent in $H$ and every $\Psi_1$-bridge with nonempty interior has at least three attachments, no bridge can lie in an end face. Thus $H$ is toroidal exactly when, for one of the twelve choices, the bridges can be assigned to the pair of simple faces and the repeated face with their attachments on the corresponding boundary walks and both $H_{\mathrm{simp}}$ and $H_{\mathrm{rep}}$ planar.
\end{lemma}

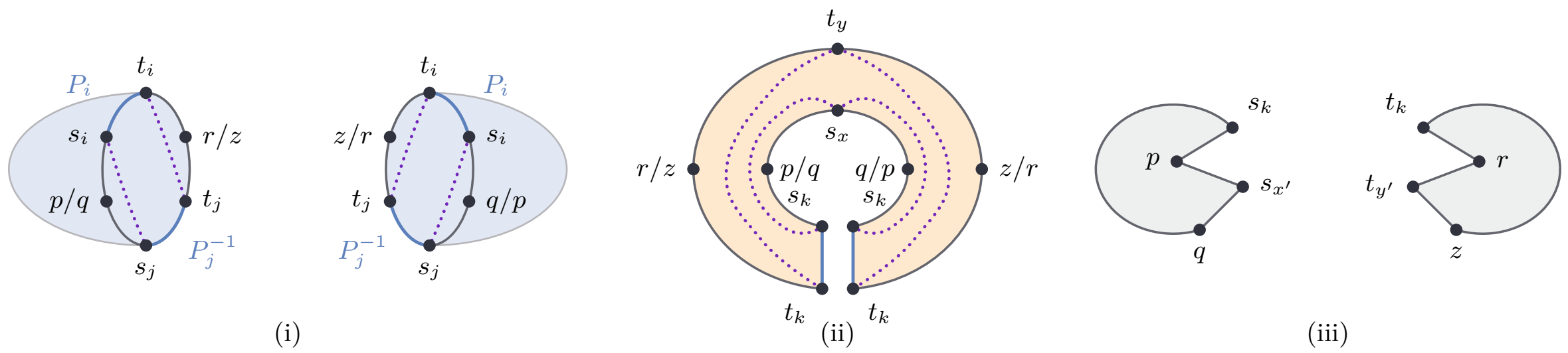
\begin{figure}[htbp]
\centering
\tikzset{external/export next=false}%
\scaledfigure{1.03}{%
\begin{tikzpicture}[
  font=\small,
  xscale=.85,yscale=.71,
  vertex/.style={circle,fill=figureInk,inner sep=1.5pt},
  twin/.style={circle,fill=figureInk,inner sep=1.5pt},
  edge/.style={draw=figureInk!75,line width=.8pt},
  pathedge/.style={draw=figurePath,line width=1.1pt},
  aux/.style={draw=figureBridge,line width=1.1pt,line cap=round,
    dash pattern=on 0pt off 2.6pt},
  bd/.style={draw=figureInk!35,line width=.6pt},
  face/.style={fill=figureInk!8}
]

\begin{scope}[xshift=1.95cm]
  \fill[faceSimple] (0,1.3) arc[start angle=90,end angle=270,x radius=1.95,y radius=1.3]
    arc[start angle=-90,end angle=90,x radius=0.62,y radius=1.3];
  \draw[bd] (0,1.3) arc[start angle=90,end angle=270,x radius=1.95,y radius=1.3];
  \draw[pathedge] (155:0.62 and 1.3) arc[start angle=155,end angle=90,x radius=0.62,y radius=1.3];
  \draw[edge] (90:0.62 and 1.3) arc[start angle=90,end angle=-25,x radius=0.62,y radius=1.3];
  \draw[pathedge] (-25:0.62 and 1.3) arc[start angle=-25,end angle=-90,x radius=0.62,y radius=1.3];
  \draw[edge] (270:0.62 and 1.3) arc[start angle=270,end angle=155,x radius=0.62,y radius=1.3];
  \draw[aux] (90:0.62 and 1.3) .. controls (0.22,0.75) and (0.42,0.0) .. (-25:0.62 and 1.3);
  \draw[aux] (155:0.62 and 1.3) .. controls (-0.42,0.0) and (-0.22,-0.75) .. (270:0.62 and 1.3);
  \node[vertex,label=above:{$t_i$}] at (90:0.62 and 1.3) {};
  \node[twin,label=right:{$r/z$}] at (25:0.62 and 1.3) {};
  \node[vertex,label=right:{$t_j$}] at (-25:0.62 and 1.3) {};
  \node[vertex,label=below:{$s_j$}] at (270:0.62 and 1.3) {};
  \node[twin,label=left:{$p/q$}] at (205:0.62 and 1.3) {};
  \node[vertex,label=left:{$s_i$}] at (155:0.62 and 1.3) {};
  \node[text=figurePath] at (-.95,1.43) {$P_i$};
  \node[text=figurePath] at (.95,-1.42) {$P_j^{-1}$};
\end{scope}
\begin{scope}[xshift=5.98cm]
  \fill[faceSimple] (0,1.3) arc[start angle=90,end angle=-90,x radius=1.95,y radius=1.3]
    arc[start angle=270,end angle=90,x radius=0.62,y radius=1.3];
  \draw[bd] (0,1.3) arc[start angle=90,end angle=-90,x radius=1.95,y radius=1.3];
  \draw[pathedge] (25:0.62 and 1.3) arc[start angle=25,end angle=90,x radius=0.62,y radius=1.3];
  \draw[edge] (90:0.62 and 1.3) arc[start angle=90,end angle=205,x radius=0.62,y radius=1.3];
  \draw[pathedge] (205:0.62 and 1.3) arc[start angle=205,end angle=270,x radius=0.62,y radius=1.3];
  \draw[edge] (270:0.62 and 1.3) arc[start angle=270,end angle=385,x radius=0.62,y radius=1.3];
  \draw[aux] (90:0.62 and 1.3) .. controls (-0.22,0.75) and (-0.42,0.0) .. (205:0.62 and 1.3);
  \draw[aux] (25:0.62 and 1.3) .. controls (0.42,0.0) and (0.22,-0.75) .. (270:0.62 and 1.3);
  \node[vertex,label=above:{$t_i$}] at (90:0.62 and 1.3) {};
  \node[twin,label=left:{$z/r$}] at (155:0.62 and 1.3) {};
  \node[vertex,label=left:{$t_j$}] at (205:0.62 and 1.3) {};
  \node[vertex,label=below:{$s_j$}] at (270:0.62 and 1.3) {};
  \node[twin,label=right:{$q/p$}] at (-25:0.62 and 1.3) {};
  \node[vertex,label=right:{$s_i$}] at (25:0.62 and 1.3) {};
  \node[text=figurePath] at (.95,1.43) {$P_i$};
  \node[text=figurePath] at (-.95,-1.42) {$P_j^{-1}$};
\end{scope}
\node at (3.965,-2.82) {(i)};

\begin{scope}[xshift=11.78cm]
  \fill[faceRepeated] (263.84:2.05) arc[start angle=263.84,end angle=-83.84,radius=2.05]
    -- (-77.29:1) arc[start angle=-77.29,end angle=257.29,radius=1]
    -- cycle;
  \draw[edge] (263.84:2.05) arc[start angle=263.84,end angle=-83.84,radius=2.05];
  \draw[edge] (257.29:1) arc[start angle=257.29,end angle=-77.29,radius=1];
  \draw[pathedge] (257.29:1)--(263.84:2.05) (-77.29:1)--(-83.84:2.05);
  \draw[aux] plot[smooth,tension=.75] coordinates
    {(90:2.05) (122:1.66) (180:1.58) (228:1.62) (263.84:2.05)};
  \draw[aux] plot[smooth,tension=.75] coordinates
    {(90:2.05) (58:1.66) (0:1.58) (-48:1.62) (-83.84:2.05)};
  \draw[aux] (90:1)
    .. controls (-.30,1.34) and (-.60,1.20) .. (130:1.25)
    arc[start angle=130,end angle=230,radius=1.25]
    .. controls (-.62,-1.14) and (-.39,-1.13) .. (257.29:1);
  \draw[aux] (90:1)
    .. controls (.30,1.34) and (.60,1.20) .. (50:1.25)
    arc[start angle=50,end angle=-50,radius=1.25]
    .. controls (.62,-1.14) and (.39,-1.13) .. (-77.29:1);
  \node[vertex,label=above:{$t_y$}] at (90:2.05) {};
  \node[twin,label=left:{$r/z$}] at (180:2.05) {};
  \node[twin,label=right:{$z/r$}] at (0:2.05) {};
  \node[vertex,label=below left:{$t_k$}] at (263.84:2.05) {};
  \node[vertex,label=below right:{$t_k$}] at (-83.84:2.05) {};
  \node[vertex,label=below:{$s_x$}] at (90:1) {};
  \node[twin,label={[inner sep=1pt,xshift=1pt]right:{$p/q$}}] at (180:1) {};
  \node[twin,label={[inner sep=1pt,xshift=-1pt]left:{$q/p$}}] at (0:1) {};
  \node[vertex,label={[xshift=2pt,yshift=2pt]above left:{$s_k$}}] at (257.29:1) {};
  \node[vertex,label={[xshift=-2pt,yshift=2pt]above right:{$s_k$}}] at (-77.29:1) {};
  \node at (0,-2.82) {(ii)};
\end{scope}

\begin{scope}[xshift=16.55cm]
  \coordinate (sk) at (40:1.1);   \coordinate (p) at (0.05,0.13);
  \coordinate (sx) at (0.99,-0.3); \coordinate (q) at (290:1.1);
  \fill[faceOther] (sk) arc[start angle=40,end angle=290,radius=1.1]
    -- (sx) -- (p) -- cycle;
  \draw[edge] (sk) arc[start angle=40,end angle=290,radius=1.1];
  \draw[edge] (q)--(sx)--(p)--(sk);
  \node[vertex,label=above right:{$s_k$}] at (sk) {};
  \node[twin,label=left:{$p$}] at (p) {};
  \node[vertex,label=right:{$s_{x'}$}] at (sx) {};
  \node[twin,label=below:{$q$}] at (q) {};
\end{scope}
\begin{scope}[xshift=20.95cm]
  \coordinate (tk) at (140:1.1);  \coordinate (r) at (-0.05,0.13);
  \coordinate (ty) at (-0.99,-0.3); \coordinate (z) at (250:1.1);
  \fill[faceOther] (tk) arc[start angle=140,end angle=-110,radius=1.1]
    -- (ty) -- (r) -- cycle;
  \draw[edge] (tk) arc[start angle=140,end angle=-110,radius=1.1];
  \draw[edge] (z)--(ty)--(r)--(tk);
  \node[vertex,label=above left:{$t_k$}] at (tk) {};
  \node[twin,label=right:{$r$}] at (r) {};
  \node[vertex,label=left:{$t_{y'}$}] at (ty) {};
  \node[twin,label=below:{$z$}] at (z) {};
\end{scope}
\node at (18.75,-2.82) {(iii)};
\end{tikzpicture}}
\caption{Left, the simple faces. Middle, the repeated face before gluing its two copies of $P_{\mathrm R}$. Right, the end faces. Dotted edges replace paths through twins.}
\label{fig:joint-test}
\end{figure}

\begin{proof}
Suppose first that $H$ has a torus embedding. Its restriction
to $\Psi_1$ has the faces described in Lemma~\ref{dis:faces}
with repeated path $P_{\mathrm R}=P_k$ and additional endpoints
$s_x,t_y$. Assign each $\Psi_1$-bridge to
the face containing its drawing away from its attachments
(Section~\ref{sec:extension}). Since $p,q,r,z$ have degree
three in $H$, no bridge attaches there. Represent each
face by a separate closed disk whose boundary traces its
boundary walk.

The boundary of each simple facial disk follows both $P_i$
and $P_j$. Replace its path through $p$ or $q$ by the
auxiliary edge $s_is_j$, and its path through $r$ or $z$
by $t_it_j$. Both disks now have boundary
$C=s_iP_it_it_jP_j^{-1}s_js_i$. Gluing them along their
entire boundary $C$ gives a sphere, with $C$ as an equator.
The bridge drawings in the two disks give a sphere drawing
of $H_{\mathrm{simp}}$, proving (a).

The repeated facial disk has two copies of $P_{\mathrm R}$ on its
boundary, traversed in opposite directions. Identifying
these copies gives an annulus, with $s_k,s_x$ on one
boundary component and $t_k,t_y$ on the other
(Figure~\ref{fig:joint-test}). The drawing of $P_{\mathrm R}$ joins
the two components, and the assigned bridges remain drawn
without crossings. Lemma~\ref{lem:annulus-embedding}, applied
with these two pairs of boundary vertices, gives the
planarity of $H_{\mathrm{rep}}$, proving (b). The drawings in the two end faces give (c).

Conversely, suppose (a)--(c) hold for one of the twelve combinations
$(P_{\mathrm R};s_x,t_y)$. Lemma~\ref{dis:faces} gives a torus embedding
of $\Psi_1$ with this repeated path and these additional endpoints. Let $C$ be the cycle in $H_{\mathrm{simp}}$ formed by
$P_i$, $P_j$ and the auxiliary edges $s_is_j,t_it_j$. By (a), take a sphere embedding of $H_{\mathrm{simp}}$. Cutting along $C$
gives two disks. The drawing of each bridge assigned to $H_{\mathrm{simp}}$
lies in one of these disks, since its drawing away from
its attachments is connected. On each disk boundary,
replace $s_is_j$ and $t_it_j$ by the two-edge paths through
the twins on the corresponding simple face boundary.
This gives disk drawings for the bridges in the two
simple faces. By Lemma~\ref{lem:annulus-embedding}, the planarity of $H_{\mathrm{rep}}$
gives an annulus drawing of $P_{\mathrm R}$ and the bridges assigned
to the repeated face, with $s_k,s_x$ on one boundary
component and $t_k,t_y$ on the other. The drawing of $P_{\mathrm R}$
meets the boundary only at $s_k$ and $t_k$. Cutting along
it gives a disk with two copies of $P_{\mathrm R}$ on its boundary.
Restoring the boundary paths through the twins gives a
disk drawing for the repeated face. Together with the
disk drawings in (c), these drawings extend the embedding
of $\Psi_1$ to $H$ (Lemma~\ref{lem:embedding-extension}).

For the last claim, a bridge assigned to an end face
has all its attachments in $\{s_k,s_{x'}\}$ or
$\{t_k,t_{y'}\}$. It therefore has at most two attachments,
both in $A_1$ or both in $A_2$. It is not an edge, as $A_1$ and
$A_2$ are independent, so it has nonempty interior, contradicting that every such bridge has at least three attachments.
\end{proof}

\subsection{Restrictions on the attachments of \texorpdfstring{$\Psi_1$-bridges}{Psi1}}\label{sec:path-choice}

We now use these tests to restrict the $\Psi_1$-bridges. For $u,v\in V(P_i)$, write $P_i[u,v)$ for
$P_i[u,v]$ without $v$. For a path $P$, write $\operatorname{int}P$ for its set of internal vertices.

\begin{lemma}\label{dis:same-side}
Suppose that paths with internal vertices outside $\Psi_1$ join internal
vertices of $P_1$ to $s_2$ and to $s_3$. Then $G\cong\mathfrak{T}_{4}$.
The same is true for $t_2,t_3$ in place of $s_2,s_3$ by exchanging $A_1$ and $A_2$.   
\end{lemma}
\begin{proof}
If the two paths lie in different bridges, shorten both to edges and
contract the segment of $P_1$ between their ends. If they lie in one
bridge, contract its interior into one of its attachments in the
interior of $P_1$. Deleting everything else outside $\Psi_1$ and
shortening the remaining segments gives $P_1=s_1ut_1$, the other
paths edges, and edges $us_2,us_3$.
This minor is isomorphic to $\mathfrak{T}_{4}$~\supp{W2}.
Thus, $G\cong\mathfrak{T}_{4}$.
\end{proof}
\begin{lemma}\label{dis:noncofacial}
If the attachments of a $\Psi_1$-bridge with nonempty interior lie on no
common face in any torus embedding of $\Psi_1$, then
$G\cong\mathfrak{T}_5$. If $G\not\cong\mathfrak{T}_5$, then
$p,q$ are the only vertices with neighborhood $A_1$, and
$r,z$ are the only vertices with neighborhood $A_2$.
\end{lemma}
\begin{proof}
Let $B$ be the bridge. By Lemma~\ref{dis:faces}, every set
$V(P_\ell)\cup V(P_m)$ or $V(P_\ell)\cup\{s_m,t_n\}$ with
$m,n\ne\ell$ lies on a face boundary in some torus embedding
of $\Psi_1$. So the attachments of $B$ lie in no such set.
Let $u_i$ denote an internal vertex of $P_i$. After permuting the path indices and, if necessary, interchanging $A_1,A_2$, the attachments contain one of
\[
\{u_1,u_2,u_3\},\qquad
\{u_1,u_2,s_3\},\qquad
\{u_1,s_2,s_3\},\qquad A_1.
\]
Indeed, if they include internal vertices of all three paths,
the first case holds. If of exactly two, say $P_1,P_2$, they
also include an end of $P_3$, giving the second case. If of
exactly one, say $P_1$, they must include both $s_2,s_3$ or
both $t_2,t_3$, giving the third case. If there are no internal
attachments, they contain $A_1$ or $A_2$. Otherwise, each of
$A_1,A_2$ contains at most two attachments, so all attachments
lie in $\{s_\ell,s_m,t_\ell,t_n\}$ for some
$\ell,m,n\in\{1,2,3\}$ with $m,n\ne\ell$.

Contract the connected interior of $B$ to a vertex $w$,
retain only its edges to the chosen triple, and delete
the other $\Psi_1$-bridges. We show the first three attachment
sets cannot occur. In each, contract $wu_1$. For each $u_j$
in the triple with $j=2,3$, contract $P_j[s_j,u_j]$ into
$s_j$. Shortening the remaining segments gives
$P_1=s_1u_1t_1$, the other paths edges, and edges
$u_1s_2,u_1s_3$. This is the $\mathfrak{T}_4$ minor
in the proof of Lemma~\ref{dis:same-side}. It is proper
because $wu_1$ was contracted, and it is nontoroidal,
contrary to minimality. Hence the attachments contain $A_1$.
Shortening each $P_i$ to an edge gives $\Psi_1$ with a vertex
$w$ adjacent to $A_1$,
which is $\mathfrak{T}_5$~\supp{W3}. This minor is
nontoroidal, so $G\cong\mathfrak{T}_5$.

Any vertex $w\notin\{p,q\}$ with neighborhood $A_1$ lies
outside $\Psi_1$, since every vertex of $\Psi_1$ other
than $p,q$ has a neighbor outside $A_1$. The vertex $w$
and its three incident edges therefore form a bridge
attached at $A_1$, which lies on no common face boundary
by Lemma~\ref{dis:faces}, so the first part gives
$G\cong\mathfrak{T}_5$. The same argument applies to $A_2$.
\end{proof}

We next use the absence of a $K_5$ minor to exclude two further
attachment possibilities for a $\Psi_1$-bridge.
\begin{lemma}\label{dis:skew}
A connected set outside $\Psi_1$ cannot attach to $s_i,s_j,t_i,t_k$ for
three distinct indices $i,j,k$.
\end{lemma}
\begin{proof}
Suppose a connected set $D$ outside $\Psi_1$ is adjacent
to $s_i,s_j,t_i,t_k$, where $i,j,k$ are distinct.
Relabel so that $(i,j,k)=(1,2,3)$.
Retain $\Psi_1$, the graph induced by $D$, and one edge
from $D$ to each of $s_1,s_2,t_1,t_3$.
Contract $D$ to a vertex $x$ and shorten each $P_i$
to the edge $s_it_i$. Thus $x$ is adjacent to
$s_1,s_2,t_1,t_3$. The five sets
\[
\{p,q,s_1\},\qquad
\{r,s_2,t_2\},\qquad
\{s_3,t_3\},\qquad
\{t_1,z\},\qquad
\{x\}
\]
are disjoint and connected, with an edge between every
two sets. Contracting each to a vertex gives a $K_5$
minor.
\end{proof}

\begin{figure}[htbp]
\centering
\begingroup
\scaledfigure{1.1}{\tikzsetnextfilename{opposite-order-minor-r3}%
\begin{tikzpicture}[
  font=\small,
  vertex/.style={circle,draw=figureInk,line width=.3pt,inner sep=2.6pt},
  dot/.style={circle,fill=figureInk,inner sep=1.4pt},
  edge/.style={draw=figureInk!75,line width=.7pt,
    preaction={draw=white,line width=1.7pt}},
  pathedge/.style={draw=figurePath,line width=1pt},
  link/.style={draw=figureInk!75,line width=1pt},
  new/.style={draw=figureBridge,line width=1pt,
    preaction={draw=white,line width=2pt}}
]
\foreach \dx in {-10.2,-6.4}{
  \begin{scope}[xshift=\dx cm]
  \draw[pathedge] (0,.6)--(2.4,.6) (0,-.6)--(2.4,-.6);
  \end{scope}
}
\begin{scope}[xshift=-10.2cm]
  \draw[link] (.75,.6)--(.75,-.6) (1.65,.6)--(1.65,-.6);
  \node at (1.2,-1.35) {same order};
\end{scope}
\begin{scope}[xshift=-6.4cm]
  \draw[new] (.75,.6)--node[pos=.25,left=2pt,inner sep=1pt] {$P$} (1.65,-.6);
  \draw[new] (1.65,.6)--node[pos=.25,right=2pt,inner sep=1pt] {$P'$} (.75,-.6);
  \node at (1.2,-1.35) {opposite order};
\end{scope}
\foreach \dx in {-10.2,-6.4}{
  \begin{scope}[xshift=\dx cm]
  \foreach \xx/\yy in {0/.6,2.4/.6,0/-.6,2.4/-.6,
      .75/.6,1.65/.6,.75/-.6,1.65/-.6}
    \node[dot] at (\xx,\yy) {};
  \node[left=2pt,inner sep=0pt] at (0,.6) {$s_i$};
  \node[right=2pt,inner sep=0pt] at (2.4,.6) {$t_i$};
  \node[left=2pt,inner sep=0pt] at (0,-.6) {$s_j$};
  \node[right=2pt,inner sep=0pt] at (2.4,-.6) {$t_j$};
  \node[above=3pt,inner sep=0pt] at (.75,.6) {$u$};
  \node[above=3pt,inner sep=0pt] at (1.65,.6) {$v$};
  \node[below=3pt,inner sep=0pt] at (.75,-.6) {$w$};
  \node[below=3pt,inner sep=0pt] at (1.65,-.6) {$h$};
  \end{scope}
}
\draw[-{Stealth[length=4pt]},figureInk!65,line width=.6pt]
  (-3.15,0)--node[above=2pt,font=\footnotesize] {contract} (-1.95,0);

\foreach \v/\xx/\yy in {
  s1/0/1,s2/0/0,s3/0/-1,t1/2.6/1,t2/2.6/0,t3/2.6/-1,
  p/-1.2/.55,q/-1.2/-.55,r/3.8/.55,z/3.8/-.55}
  \coordinate (\v) at (\xx,\yy);
\foreach \a in {p,q}
  \foreach \b in {s1,s2,s3} \draw[edge] (\a)--(\b);
\foreach \a in {r,z}
  \foreach \b in {t1,t2,t3} \draw[edge] (\a)--(\b);
\draw[pathedge] (s1)--(t1) (s2)--(t2) (s3)--(t3);
\draw[new] (s1)--(t2);
\draw[new] (s2)--(t1);
\node[vertex,fill=branchA,label=left:{$p$}] at (p) {};
\node[vertex,fill=branchB,label=left:{$q$}] at (q) {};
\node[vertex,fill=branchC,label=right:{$r$}] at (r) {};
\node[vertex,fill=branchA,label=right:{$z$}] at (z) {};
\node[vertex,fill=branchD,label=above:{$s_1$}] at (s1) {};
\node[vertex,fill=branchB,label=below right:{$s_2$}] at (s2) {};
\node[vertex,fill=branchA,label=below:{$s_3$}] at (s3) {};
\node[vertex,fill=branchC,label=above:{$t_1$}] at (t1) {};
\node[vertex,fill=branchE,label=below left:{$t_2$}] at (t2) {};
\node[vertex,fill=branchA,label=below:{$t_3$}] at (t3) {};
\node[text=branchA,anchor=west] at (4.75,.9) {$W_1=\{p,s_3,t_3,z\}$};
\node[text=branchB!55!figureInk,anchor=west] at (4.75,.45) {$W_2=\{q,s_2\}$};
\node[text=branchC,anchor=west] at (4.75,0) {$W_3=\{r,t_1\}$};
\node[text=branchD,anchor=west] at (4.75,-.45) {$W_4=\{s_1\}$};
\node[text=branchE,anchor=west] at (4.75,-.9) {$W_5=\{t_2\}$};
\end{tikzpicture}}
\endgroup
\caption{Disjoint paths joining $P_i$ to $P_j$ in the same order and in
opposite order.}
\label{fig:opposite-order}
\end{figure}
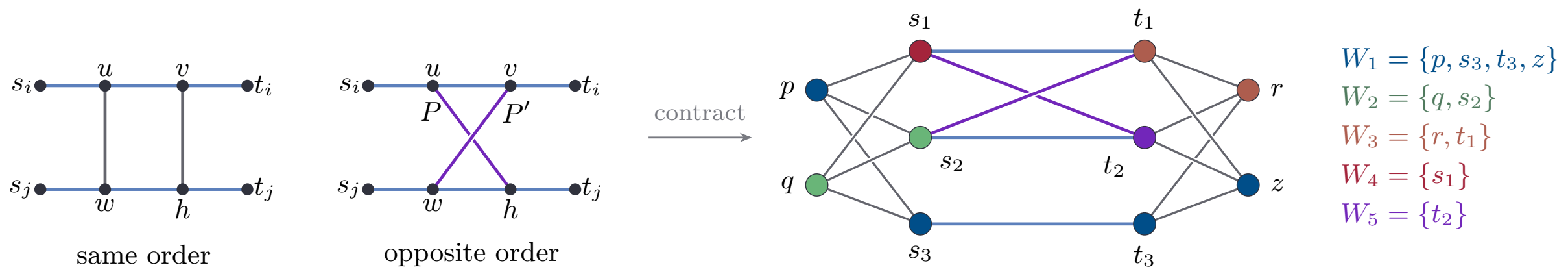

\begin{lemma}\label{dis:opposite-order}
Let $i\ne j$. At most one of $s_it_j$ and $s_jt_i$ is an edge of $G$,
and neither is an edge if a path with internal vertices outside
$\Psi_1$ joins internal vertices of $P_i$ and $P_j$. More generally,
order $P_i$ and $P_j$ from $A_1$ to $A_2$, and let $u$ precede $v$ on
$P_i$ and $w$ precede $h$ on $P_j$. Then there are no disjoint paths
from $u$ to $h$ and from $v$ to $w$ with internal vertices outside
$\Psi_1$.
\end{lemma}
\begin{proof}
Suppose there are disjoint paths $P$ from $u$ to $h$ and $P'$ from $v$
to $w$, with internal vertices outside $\Psi_1$. Delete everything
outside $\Psi_1\cup P\cup P'$. Contract $P_i[s_i,u]$ into $s_i$ and
$P_i[v,t_i]$ into $t_i$, and contract $P_j[s_j,w]$ into $s_j$ and
$P_j[h,t_j]$ into $t_j$. Now $P$ joins $s_i$ to $t_j$, and $P'$ joins
$t_i$ to $s_j$. Shortening $P$, $P'$ and the three paths of $\Psi_1$ to
edges gives $\Psi_1$ with every path an edge, together with $s_it_j$
and $s_jt_i$. Let $k$ be the third index. The five sets
\[
\{p,s_k,t_k,z\},\qquad \{q,s_j\},\qquad
\{r,t_i\},\qquad \{s_i\},\qquad \{t_j\}
\]
are disjoint and connected, with an edge between every two.
Contracting each to a vertex gives a $K_5$ minor (Figure~\ref{fig:opposite-order}). The edges $s_it_j$ and $s_jt_i$ would be such paths, with $u=s_i$,
$v=t_i$, $w=s_j$ and $h=t_j$, so at most one of them is an edge of $G$.
Now let $P_0$ join an internal vertex $x$ of $P_i$ to an internal vertex
$y$ of $P_j$, with internal vertices outside $\Psi_1$. Then $P_0$
contains none of $s_i,t_i,s_j,t_j$. So the edge $s_it_j$ and $P_0$
would be such paths, with $u=s_i$, $v=x$, $w=y$ and $h=t_j$.
Interchanging $i$ and $j$ excludes $s_jt_i$ in the same way.
\end{proof}

We next exclude three configurations in which $P_2=s_2t_2$,
$P_3=s_3t_3$, and all $\Psi_1$-bridge attachments lie in
$V(P_1)\cup\{s_2,t_2\}$. In the first two, $P_1$ is also
an edge and there is one bridge, attached at
$\{s_1,t_1,s_2,t_2\}$ or $\{s_1,t_1,s_2\}$, respectively.
In the third, every bridge attaches in the interior of $P_1$.
The next lemma shows that, in each configuration, a torus
embedding of $H/s_3t_3$ gives one of $H$.

\begin{lemma}\label{dis:matching-contraction}\label{dis:aligned-contraction}
Let $H$ contain $\Psi_1$ and have no $K_5$ minor, with $p,q,r,z$ of
degree three, $P_2=s_2t_2$ and $P_3=s_3t_3$. Suppose every
$\Psi_1$-bridge has all its attachments in $V(P_1)\cup\{s_2,t_2\}$, and
one of the following holds.
\begin{enumerate}
\item[\textup{(a)}]
$P_1=s_1t_1$, and there is exactly one $\Psi_1$-bridge, with attachment
set $\{s_1,t_1,s_2,t_2\}$.

\item[\textup{(b)}]
$P_1=s_1t_1$, and there is exactly one $\Psi_1$-bridge, with attachment
set $\{s_1,s_2,t_1\}$.

\item[\textup{(c)}]
Every $\Psi_1$-bridge has an attachment in the interior of $P_1$.
\end{enumerate}
If $H/s_3t_3$ is toroidal, then so is $H$.
\end{lemma}

\begin{proof}
Suppose $H/s_3t_3$ is toroidal. Take a torus embedding of it and consider the induced embedding on $\Psi_1/s_3t_3$, which is cellular by Youngs's theorem since this graph is connected and nonplanar. Write $a=s_1$, $b=t_1$, $c=s_2$, and $d=t_2$.

In (a), the unique bridge $B$ lies in a facial disk $F$ whose boundary
contains $a,b,c,d$. With the exception of $abcd$, every attachment
sequence of such a face of $\Psi_1/s_3t_3$ also occurs for $\Psi_1$,
up to cyclic shift and reversal~\supp{R2}. When the sequences agree,
a homeomorphism of the two disks matching the boundary appearances
carries the drawing of $B$ to the corresponding face of $\Psi_1$.

\looseness=-1
It remains to consider the order $a,b,c,d$. A simple face of
$\Psi_1$ has order $a,b,d,c$, so it suffices to draw $B$ in this order. Put $D=B-\{a,b,c,d\}$. The graph $D$ is nonempty and connected. Every edge of $B$ has an end in $D$, and every attachment has a neighbor there. Thus, if $B$ has no disk drawing with boundary order $a,b,d,c$, Lemma~\ref{dis:four-disk-order} gives vertex-disjoint $a$--$d$ and $b$--$c$ paths with nonempty interiors in $D$. Take a shortest path in $D$ joining their interiors.

Contract the interiors of the first two paths into $a,b$, respectively, and shorten the joining path to $ab$. Retaining these edges, $J_1\cup J_2$ and $P_2,P_3$ gives $\Psi_1$ with all three paths edges, together with $ad=s_1t_2$ and $bc=s_2t_1$. The contractions in the proof of Lemma~\ref{dis:opposite-order} give a $K_5$ minor, a contradiction. Thus the required drawing exists and gives a torus embedding of $H$.

In (b), the bridge uses only $a,b,c$. Omitting $d$ from the
facial walks gives the same sequences containing $a,b,c$ for
$\Psi_1/s_3t_3$ and $\Psi_1$, up to cyclic shift and reversal~\supp{R2}.
Thus the bridge can be drawn in the corresponding face of $\Psi_1$.

For (c), if there are no bridges, then $H=\Psi_1$ is toroidal. Otherwise every bridge has an attachment to an internal vertex of $P_1$. Let $H_1$ consist of $P_1$ and all such bridges. We combine the faces incident with $P_1$ in the induced embedding into one region containing the drawing of $H_1$.

To construct this region, represent each such face by a disk whose boundary follows its facial walk, with repeated appearances of a vertex kept separate. If the two sides of $P_1$ lie in different faces, identify the copies of $P_1$ on the boundaries of their two disks. This gives a disk. If both sides lie in the same face, identify the two copies of $P_1$ on that disk's boundary. This gives an annulus. In either case, keep all other boundary appearances separate. The resulting region contains the drawing of $H_1$, with $P_1$ inside it apart from its ends.

The regions along $P_1$ in the two graphs are compared in~\supp{R2}.
If the region is an annulus, its boundary sequences are $ac\mid b$,
$ac\mid bd$, $a\mid b$ or $a\mid bd$. All four occur for $\Psi_1$, so
Lemma~\ref{lem:annulus-embedding} gives a torus embedding of $H$.

If the region is a disk, its sequence is $ab$, $abc$, $abd$, $abcd$
or $acdbdc$. The sequences $ab$ and $acdbdc$ also occur for $\Psi_1$.
In the disk with sequence $acdbdc$, the appearances after $b$ give
either $abc$ or $abd$ by using just $c$ or just $d$. Thus only $abcd$
with all four attachments remains.

Again we get a $K_5$ minor if the required drawing in order $a,b,d,c$ does not exist. Put $D=H_1-\{a,b,c,d\}$. The internal vertices of $P_1$ connect all the bridge interiors, so $D$ is connected. Every edge of $H_1$ has an end in $D$, and every attachment has a neighbor there. As in (a), Lemma~\ref{dis:four-disk-order} gives a drawing in order $a,b,d,c$, or disjoint $a$--$d$ and $b$--$c$ paths through $D$. Joining these through $D$ gives the same $K_5$ minor as in (a), with the joining path supplying $ab$. Hence the drawing in order $a,b,d,c$ exists. Using the appearances of $d,c$ after $b$ in the disk above gives a torus embedding of $H$.
\end{proof}

\subsection{\texorpdfstring{$\Psi_1$}{Psi1}-bridges meeting the interior of two paths}\label{dis:active-support}

We now show that no $\Psi_1$-bridge can meet the interiors of two
of the paths.

\begin{lemma}\label{dis:active-criterion}
Let $H$ consist of $\Psi_1$ together with some of its bridges
in $G$. Suppose every such bridge has all its attachments in
$V(P_1)\cup V(P_2)$, and at least one meets the interiors of
both paths. Let $H_{\mathrm{simp}}$ consist of $P_1,P_2$, these bridges, and
the edges $s_1s_2,t_1t_2$. Then $H_{\mathrm{simp}}$ has no $K_5$
minor, and $H$ is toroidal if and only if $H_{\mathrm{simp}}$ is planar.
\end{lemma}

\begin{proof}
Contracting $ps_1,rt_1$ and deleting $q,z$ and all vertices
of $P_3$ gives $H_{\mathrm{simp}}$ from $H$. Thus $H_{\mathrm{simp}}$ is a minor of $G$
and has no $K_5$ minor. Suppose first that $H$ is toroidal. By Lemma~\ref{dis:faces},
a bridge meeting the interiors of $P_1,P_2$ lies in a simple
face for these paths, so $P_3$ is the repeated path.
The repeated face meets $P_1\cup P_2$ only in $s_x,t_y$,
and each end face meets it in one vertex. Since $G$
is three-connected, every retained bridge with nonempty
interior has at least three attachments. Thus only an
edge $s_xt_y$ could lie outside the two simple faces. Both $P_1,P_2$ have internal vertices, so
\eqref{dis:shortest} excludes $s_1t_1,s_2t_2$.
The bridge meeting both interiors contains a path between
them with internal vertices outside $\Psi_1$, so
Lemma~\ref{dis:opposite-order} excludes $s_1t_2,s_2t_1$.
Hence every bridge lies in a simple face, and
Lemma~\ref{dis:joint-test} gives that $H_{\mathrm{simp}}$ is planar. Conversely, suppose $H_{\mathrm{simp}}$ is planar. Apply
Lemma~\ref{dis:joint-test} with $(P_{\mathrm R};s_x,t_y)=(P_3;s_1,t_1)$,
assigning every bridge to the pair of simple faces.
The graph $H_{\mathrm{rep}}=P_3+s_3s_1+t_3t_1$ is a path, and no
bridges are assigned to the end faces. Thus
$H$ is toroidal.
\end{proof}

\begin{lemma}\label{dis:sole-two-rails}
If every $\Psi_1$-bridge has all its attachments in
$V(P_1)\cup V(P_2)$, then no $\Psi_1$-bridge meets the
interiors of both paths.
\end{lemma}
\begin{proof}
Suppose a bridge contains a path $P$ from the interior of $P_1$ to that of $P_2$, with internal vertices outside $\Psi_1$. No bridge attaches to $P_3$, so minimum degree gives $P_3=s_3t_3$.

Let $H_1$ consist of $P_1$, $P_2$, and every $\Psi_1$-bridge, and let $H_2=(J_1\cup J_2)+s_3t_3$. Write $A=\{s_1,t_1,s_2,t_2\}$ and $E=\{s_1s_2,t_1t_2\}$. Then $G=H_1\cup H_2$ and $V(H_1)\cap V(H_2)=A$.

If a component of $G-V(\Psi_1)$ attaches only at path ends, it has at least three attachments in $A$ by three-connectivity, and therefore has either an $s_1$--$t_2$ or an $s_2$--$t_1$ path disjoint from $P$. This path along with $P$ contradicts Lemma~\ref{dis:opposite-order}. Therefore every such component attaches inside $P_1$ or $P_2$, and $P$ makes $H_1-A$ nonempty and connected. Both $P_1$ and $P_2$ have internal vertices, so each vertex of $A$
has an interior neighbor on its path. Neither $s_1t_1$ nor $s_2t_2$
is an edge, since either would shorten its path. Applying
Lemma~\ref{dis:opposite-order} to $P$ rules out $s_1t_2$ and $s_2t_1$.
The remaining pairs $s_1s_2$ and $t_1t_2$ lie within the independent
triples $A_1$ and $A_2$, respectively. Thus $H_1$ is an $H_2$-bridge
with independent attachment set $A$.

By Lemma~\ref{dis:active-criterion}, $H_1+E$ is nonplanar with no $K_5$ minor. Let $H_1'$ be an allowed minor of $H_1$ for which $H_1'+E$ is nonplanar, and suppose $H_2\cup H_1'$ has a torus embedding. Since $H_1'-A$ is connected, it lies in one face of $H_2$. Each of
the four attachments has a neighbor in $H_1'-A$, so all four lie on
the boundary of that face.

Let $\mathcal N$ be a closed regular neighborhood of $H_2$ in the
torus. Its genus is zero or one. In either case, no boundary walk
contains all four attachments~\supp{R12}. The face containing
$H_1'-A$ must therefore have at least two boundary walks. By
Lemma~\ref{lem:torus-complement}, this face is an annulus and
$\mathcal N$ has genus zero.

In each genus-zero rotation, exactly one pair of boundary walks
contains all four attachments between them. One meets $s_1,s_2$
and the other meets $t_1,t_2$, each attachment appearing once~\supp{R12}.
Cap the two boundary circles with disks and draw $s_1s_2$ in one
disk and $t_1t_2$ in the other. This gives a sphere embedding of
$H_1'+E$, contradicting its nonplanarity.

Thus replacing $H_1$ by $H_1'$ gives a nontoroidal graph.
By Corollary~\ref{common:reduced-forms}, $H_1$ is among the
candidates with attachment set $A$ and added edges $E$.
Of the $26{,}633$ candidates, $14$ are excluded because they have no
vertex-disjoint $s_1$--$t_1$ and $s_2$--$t_2$ paths. For every other
candidate, adjoining $H_2$ gives a graph with a $K_5$ minor, a cut
of at most two vertices, or a proper nontoroidal minor~\supp{E1}.
Hence no bridge contains such a path.
\end{proof}

We now show further that even allowing bridges to attach outside of only $P_1\cup P_2$, no $\Psi_1$-bridge meets the interior of two of the paths $P_\ell$. 

\begin{lemma}\label{dis:normalize}
No $\Psi_1$-bridge meets the interiors of two of the paths $P_1,P_2,P_3$.
\end{lemma}
\begin{proof}
Suppose a bridge contains a path $P$ from
$u\in\operatorname{int}P_1$ to $v\in\operatorname{int}P_2$,
with internal vertices outside $\Psi_1$.
By Lemma~\ref{dis:faces}, drawing $P$ requires
$P_{\mathrm R}=P_3$.

No bridge meets both the interior of $P_1$ and $P_3$.
Otherwise, choose a path $P'$ from
$u'\in\operatorname{int}P_1$ to $w\in V(P_3)$,
with internal vertices outside $\Psi_1$.
If it lies in a different bridge from $P$, contract
$P_1[u,u']$ into $u$. Otherwise contract into $u$ the
vertices outside $\Psi_1$ of a tree joining $u,v,w$.
Shortening the resulting paths gives $uv,uw$.
Interchanging $A_1,A_2$ if necessary, take $w\ne t_3$.
Contract $P_2[s_2,v]$ into $s_2$ and $P_3[s_3,w]$ into
$s_3$. The edges $us_2,us_3$ give $\mathfrak T_4$ by
the contractions in the proof of
Lemma~\ref{dis:same-side}. This nontoroidal minor is
proper since $v$ is internal to $P_2$, a contradiction.
Interchanging $P_1,P_2$ also excludes a bridge meeting
the interior of $P_2$ and $P_3$.

Since a bridge meets the interiors of both $P_1$ and $P_2$,
Lemma~\ref{dis:sole-two-rails} gives a bridge with an attachment
outside $P_1\cup P_2$. This attachment lies on $P_3$, since all
edges at $p,q,r,z$ already belong to $\Psi_1$.
Deleting all such bridges, retaining their attachments,
gives a toroidal proper minor.
Lemma~\ref{dis:active-criterion} therefore makes the graph
$H_{\mathrm{simp}}$ formed by $P_1,P_2$, the remaining
bridges and $s_1s_2,t_1t_2$ planar.
By Lemma~\ref{dis:joint-test}, the remaining bridges
can therefore be drawn in the simple faces for each
of the four pairs of additional endpoints when
$P_{\mathrm R}=P_3$. The bridges meeting $P_3$ must use
the repeated face, since the simple faces miss $P_3$
and three-connectivity and independence of $A_1,A_2$
exclude the end faces by the same lemma.
Since $G$ is nontoroidal, these bridges cannot be drawn
together there for any of the four pairs.
Retain only $\Psi_1$, $P$ and the bridges meeting $P_3$.
The path $P$ still requires $P_{\mathrm R}=P_3$, and
these bridges and the repeated-face boundaries are
unchanged. The retained graph is therefore nontoroidal,
so minimality makes it equal to $G$.
Minimum degree now gives $P=uv$, $P_1=s_1ut_1$ and
$P_2=s_2vt_2$.

Take a torus embedding of the proper minor $G/us_1$.
The bridges meeting $P_3$ avoid $u,v$, so their
attachments are unchanged. The only new edges are
$s_1t_1,s_1v$, so $A_1,A_2$ remain independent, and
every bridge with nonempty interior retains at least
three attachments. Lemma~\ref{dis:joint-test} therefore
excludes the end faces. With $P_{\mathrm R}=P_3$,
the repeated face still has the same boundary and
bridges, so they still cannot be drawn together there.
If $P_{\mathrm R}=P_1$, the only faces containing $v$
are the simple faces, which do not contain $s_1$, so the edge
$s_1v$ cannot be drawn. Thus $P_2$ is repeated,
and its repeated face contains $s_1$.

A bridge meeting $P_3$ and either $s_2$ or $t_2$ must
use this repeated face. That face meets $P_3$ at most
at $t_3$, so the bridge meets $P_3$ only at $t_3$.
Applying the same argument to $G/ut_1$ forces it
to meet $P_3$ only at $s_3$, a contradiction.
Interchanging indices $1,2$ likewise excludes
attachments at $s_1,t_1$. Thus every bridge meeting
$P_3$ attaches only to that path. There is at least
one such bridge, so $\{s_3,t_3\}$ separates $p$
from the internal vertices of $P_3$ and these bridges,
contrary to three-connectivity.
\end{proof}

We may apply Lemma~\ref{dis:normalize} to any three disjoint paths
$P_i$ from $s_i$ to $t_i$ that meet $A_1\cup A_2$ only at their
ends, provided we take $P_i=s_it_i$ whenever $s_it_i$ is an edge
of $G$, as in~\eqref{dis:shortest}.

To apply the bridge results of Section~\ref{sec:core}, we now
consider the $(J_1\cup J_2)$-bridges. The internal vertices of the
paths then lie in the bridge interiors. Since the twins have degree
three, these are exactly the bridges on $A_1\cup A_2$ in
$G-\{p,q,r,z\}$.

By Lemma~\ref{dis:normalize}, no $\Psi_1$-bridge meets the interiors
of two paths. Thus each $P_m$ with an internal vertex, together with
all $\Psi_1$-bridges meeting its interior, forms one
$(J_1\cup J_2)$-bridge. Let $D_m$ be the graph induced by its
interior and put $W_m=N(D_m)-\{s_m,t_m\}$, the set of attachments
other than the ends of $P_m$. The only other bridges are the
$\Psi_1$-bridges attached only at path ends and the paths $P_i$ that
are single edges. For each such single-edge path, put
$W_i=\varnothing$. The next proposition bounds the attachment sets
of bridges with nonempty interior and lists the cases with four
attachments.

\begin{proposition}\label{dis:types}
Suppose $G\not\cong\mathfrak{T}_4$ and
$G\not\cong\mathfrak{T}_5$. Let $B$ be a $(J_1\cup J_2)$-bridge
with nonempty interior. Then $B$ has either
three or four attachments, and its interior meets at most one
of the paths $P_i$. If $B$ has exactly four attachments,
one of the following holds.
\begin{enumerate}[label=\textup{(\alph*)}]
\item It lies outside the paths, with attachment set
$\{s_x,t_x,s_y,t_y\}$ for some $x\ne y$.
\item It contains $P_m$, with attachment set
$\{s_m,t_m,s_x,t_x\}$ for some $x\ne m$.
\item It contains $P_m$, with attachment set
$\{s_m,t_m,s_x,t_y\}$, where
$\{m,x,y\}=\{1,2,3\}$.
\end{enumerate}
\end{proposition}

\begin{proof}
By three-connectivity, $B$ has at least three attachments. If its interior met two of the
paths $P_i$, take a shortest path within it between
interiors of distinct $P_i$. Its internal vertices
lie outside $\Psi_1$, so it lies in a $\Psi_1$-bridge
meeting both path interiors, contrary to
Lemma~\ref{dis:normalize}.

If the interior of $B$ avoids all the paths, then
$B$ is a $\Psi_1$-bridge. Since
$G\not\cong\mathfrak{T}_5$,
Lemma~\ref{dis:noncofacial} places its attachments
on one face boundary in some torus embedding of
$\Psi_1$. By Lemma~\ref{dis:faces}, they therefore
lie in some $\{s_\ell,t_\ell,s_n,t_k\}$ with
$n,k\ne\ell$. Hence there are at most four.
With four, Lemma~\ref{dis:skew} forces $n=k$,
giving type (a).

Otherwise, $B$ contains $P_m$. For each $d\in W_m$,
connectivity of $D_m$ gives a shortest path from
$d$ to the interior of $P_m$ with internal
vertices in $D_m$. These vertices lie outside
$\Psi_1$, since the path first meets $P_m$ at its
end and $D_m$ meets no other path. As
$G\not\cong\mathfrak{T}_4$,
Lemma~\ref{dis:same-side} allows at most one
vertex of $A_1$ and one of $A_2$ in $W_m$.
Together with $s_m,t_m$, this gives at most
four attachments. With four,
$W_m=\{s_x,t_y\}$ for $x,y\ne m$. If $x=y$,
we have type (b). Otherwise,
$\{m,x,y\}=\{1,2,3\}$, giving type (c).
\end{proof}

\subsection{Bridge reductions and bounds}\label{dis:finite-method}

Apart from $\mathfrak{T}_4,\mathfrak{T}_5$, Proposition~\ref{dis:types} leaves three cases. If every $P_i$ is an edge, Section~\ref{dis:endpoint-support} applies. Otherwise, if no bridge has four attachments and at least two interior vertices, every bridge with at least two interior vertices has three attachments, and Section~\ref{dis:three-construction} applies. If some bridge has four attachments and at least two interior vertices, it has type~(a), (b) or~(c), and Section~\ref{dis:path-four} applies. From here on, unless stated otherwise, bridges are
$(J_1\cup J_2)$-bridges other than the single-edge paths $P_i$,
and $G\not\cong\mathfrak{T}_4$ and $G\not\cong\mathfrak{T}_5$. By Lemma~\ref{dis:noncofacial}, only $p,q$ have neighborhood
$A_1$, and only $r,z$ have neighborhood $A_2$. 

To use the bridge results of Section~\ref{sec:core} to restrict the possible
forms of the bridges, we first need their triangle completions
to be three-connected and have no $K_5$ minor. 

\begin{lemma}\label{dis:exterior-triangles}
Let $B$ be a $\Psi_1$-bridge with attachment set
$A=\{w,s_j,t_k\}$, where $w\in V(P_i)$, and either $i,j,k$ are
distinct, or $j\ne i$, $w=s_i$ and $k=i$. Then $B\cup K_3[A]$ is a
three-connected minor of $G$, and if it is nonplanar, $B$ has an
allowed minor $T(v)$ for some $v\in A$. This also holds with $A_1,A_2$
interchanged.
\end{lemma}
\begin{proof}
First suppose $i,j,k$ are distinct. Join $w$ to $s_j$
along $wP_is_ips_j$, $w$ to $t_k$ along $wP_it_irt_k$,
and $s_j$ to $t_k$ along $s_jP_jt_jzt_k$. In the second
case, $A=\{s_i,s_j,t_i\}$. The paths $s_ips_j$, $P_i$
and $s_jP_jt_jrt_i$ join $s_i$ to $s_j$, $s_i$ to
$t_i$, and $s_j$ to $t_i$, respectively. In either case, the three paths have disjoint interiors
in $\Psi_1-A$, so they meet $B$ only at their ends.
Keeping $B$ and contracting these paths to edges,
without identifying their ends, gives $B\cup K_3[A]$
as a minor of $G$.
Since $B[\operatorname{Int}(B)]$ is a component of
$G-A$, Lemma~\ref{lem:triangle-completion} also gives
that $B\cup K_3[A]$ is three-connected. If $B\cup K_3[A]$ is nonplanar, apply
Lemma~\ref{endpoint:primitive} in $G$ using these
three paths. As $G$ has no $K_5$ minor, this gives
an allowed minor $T(v)$ of $B$ for some $v\in A$.
The same argument applies with $A_1,A_2$ interchanged.
\end{proof}

Lemmas~\ref{dis:faces} and~\ref{dis:joint-test} give the following
test for when an embedding of $\Psi_1$ extends to one of $G$.

\begin{corollary}\label{dis:choice-fails}
Fix $(P_{\mathrm R};s_x,t_y)$ as in
Lemma~\ref{dis:joint-test}, where $P_{\mathrm R}=P_k$
and $\{i,j,k\}=\{1,2,3\}$.
A bridge containing $P_k$ must be assigned to
$H_{\mathrm{rep}}$, and one containing $P_i$ or $P_j$
to $H_{\mathrm{simp}}$.
Every other bridge must have all its attachments in
$\{s_i,t_i,s_j,t_j\}$ to be assigned to $H_{\mathrm{simp}}$,
or in $\{s_k,t_k,s_x,t_y\}$ to be assigned to
$H_{\mathrm{rep}}$.
Only the edge $s_xt_y$, with $x\ne y$, can have both assignments. Then the embedding extends to $G$ exactly when every bridge
has an assignment and both test graphs are planar,
for at least one assignment of $s_xt_y$ if it is present.
\end{corollary}
\begin{proof}
Since $G$ is three-connected, every $\Psi_1$-bridge
with nonempty interior has at least three attachments.
Also, $A_1,A_2$ are independent by
Lemma~\ref{lem:independent-neighborhood}(a).
Thus Lemma~\ref{dis:joint-test} excludes the end
faces, leaving only the simple faces and the
repeated face. Only the repeated face contains
internal vertices of $P_k$, and only the simple
faces contain internal vertices of the other two
paths. Hence a bridge containing $P_m$ must be
assigned to $H_{\mathrm{rep}}$ if $m=k$, and to
$H_{\mathrm{simp}}$ otherwise. For a bridge outside the paths, its attachments
must lie in $\{s_i,t_i,s_j,t_j\}$ to be assigned
to $H_{\mathrm{simp}}$, or in
$\{s_k,t_k,s_x,t_y\}$ to be assigned to
$H_{\mathrm{rep}}$. These sets have only $s_x,t_y$
in common, so a bridge with at least three
attachments cannot have both assignments.
The only edge that could have both is $s_xt_y$.
If $x=y$, this edge, when present, is $P_x$
itself by \eqref{dis:shortest} and is not a
bridge. Thus only $s_xt_y$ with $x\ne y$
can have both assignments. The end faces are unavailable, so a bridge with no assignment
cannot be added. Otherwise, Lemma~\ref{dis:joint-test} requires both test graphs
to be planar, with $s_xt_y$, if present, added to one of them.
\end{proof}

The next
lemma lists the relevant completions for bridges with three attachments, allowing
us to apply Proposition~\ref{prop:three-root-disjoint} to
restrict their possible forms.

\begin{lemma}\label{dis:relevant}
Let $B$ be a bridge with three attachments. Then its triangle completion $B\cup K_3[A]$ is a three-connected minor of $G$. If this completion is planar, the interior of $B$ is a single vertex.
\begin{enumerate}[label=\textup{(\roman*)}]
\item If $B$ lies outside the paths, its attachments are $c_1,c_2$ in one of $A_1,A_2$ and $h$ in the other. Its relevant completions are
\[
B\cup K_3[A],\qquad B+c_1c_2+c_1h,\qquad B+c_1c_2+c_2h.
\]
\item If $B$ contains $P_m$, let $e$ be its third attachment. Label the ends of $P_m$ so that $g$ lies in the same twin triple as $e$ and $h$ lies in the other. Its relevant completions are
\[
B\cup K_3[A],\qquad B+ge+gh,\qquad B+ge+eh,\qquad B+ge.
\]
\end{enumerate}
If the triangle completion is nonplanar and some relevant completion is planar, Proposition~\ref{prop:three-root-disjoint} applies to $B$, with $c_1,c_2,h$ in case (i), or $g,e,h$ in case (ii), in the roles of $x,y,z$.
\end{lemma}

\begin{proof}
If $B$ lies outside the paths, its attachment set cannot be $A_1$ or $A_2$, since Lemmas~\ref{dis:faces} and~\ref{dis:noncofacial} would give $G\cong\mathfrak T_5$. Thus two attachments $c_1,c_2$ lie in one triple and the third $h$ lies in the other. Lemma~\ref{dis:exterior-triangles} shows the triangle completion is a three-connected minor of $G$. If $B$ contains $P_m$, we may interchange $A_1,A_2$ and assume its attachments are $s_m,t_m,s_x$. The paths
\[
s_mps_x,\qquad s_xP_xt_xrt_m,\qquad
s_mqs_wP_wt_wzt_m,
\]
where $w$ is the third index, join the three attachment pairs and have disjoint interiors outside $B$. Contracting each to an edge gives the triangle completion as a minor of $G$, and Lemma~\ref{lem:triangle-completion} makes it three-connected. In either case, if the triangle completion is planar, Lemma~\ref{lem:flat-replacement} makes the interior of $B$ a single vertex.

Now consider a test graph $H_{\mathrm{simp}}$ or $H_{\mathrm{rep}}$ containing $B$. Delete the interior and edges of $B$, retaining $A$, and call the remaining graph $H_1$. Consider the bridges on $A$ in $H_1$. If at least two attach at all three vertices, they and $B$ give a $K_{3,3}$ minor, so the test graph is nonplanar. If exactly one does, Lemma~\ref{core:assembly-general} requires the triangle completion of $B$ to be planar. If none does, the condition on $B$ is that $B+E$ is planar, where $E$ contains the edge $\alpha\beta$ whenever a bridge on $A$ in $H_1$ has attachments $\{\alpha,\beta\}$. The rest of the test checks the bridges on $A$ in $H_1$, with the
added edges prescribed by Lemma~\ref{core:assembly-general}.

Suppose first that $B$ lies outside the paths, and consider the case when no bridge on $A$ in $H_1$ attaches at all three vertices. In $H_{\mathrm{simp}}$, $H_1$ contains a cycle through $c_1,c_2,h$. These three attachments divide the cycle into a $c_1$--$c_2$ path, a $c_2$--$h$ path and an $h$--$c_1$ path, each with no attachment in its interior. Each path lies in a bridge whose attachments are its two ends, so $c_1c_2,c_2h,c_1h$ all belong to $E$. In $H_{\mathrm{rep}}$, $H_1$ contains a path through $c_1,c_2,h$ or $c_2,c_1,h$, in that order. In the first case, $c_1c_2,c_2h\in E$, and another bridge can only add $c_1h$. In the second, $c_1c_2,c_1h\in E$, and another bridge can only add $c_2h$. Thus the completion is one of the two listed two-edge completions or the triangle completion, giving (i).

If $B$ contains $P_m$, $H_1$ contains a path through $g,e,h$ in $H_{\mathrm{simp}}$. 
If no bridge on $A$ in $H_1$ attaches at all three vertices, the
$g$--$e$ and $e$--$h$ portions of the path lie in bridges with
attachment sets $\{g,e\}$ and $\{e,h\}$, respectively. Thus $E$
contains $ge$ and $eh$. The completion is $B+ge+eh$, or the triangle
completion if $gh$ also belongs to $E$. In $H_{\mathrm{rep}}$, $H_1$ contains the edge $ge$, so $ge\in E$. Either, both or neither of $gh,eh$ may also belong to $E$. This gives precisely the four completions in (ii).

Finally, let $B'$ be an allowed minor of $B$. Its interior is nonempty and connected, and every attachment still has an interior neighbor. If $B$ contains $P_m$, choose a $g$--$h$ path through the interior of $B'$ to replace $P_m$ in $\Psi_1$. The attachments and the graph outside $B$ are unchanged. For each assignment of the bridges, whether $H_{\mathrm{simp}}$ and $H_{\mathrm{rep}}$ are planar after the replacement therefore depends only on which of the listed completions of $B'$ are planar. By Lemma~\ref{dis:joint-test}, the same holds for whether the resulting graph is toroidal. Lemma~\ref{core:relevant} now shows that every proper allowed minor of $B$ makes a nonplanar relevant completion planar. Hence Proposition~\ref{prop:three-root-disjoint} applies whenever the triangle completion is nonplanar and some relevant completion is planar.
\end{proof}

Forming the candidate graphs in Sections~\ref{dis:endpoint-support} and~\ref{dis:three-construction} also requires a bound on how many
bridges can have the same attachments.

\begin{lemma}\label{dis:counts}
At most three bridges outside the paths have the same attachment set. If a bridge has one interior vertex and three attachments,
its attachments are independent. Deleting any bridge outside the paths, with its
attachments retained, leaves a graph to which an embedding of
$\Psi_1$ extends.
\end{lemma}
\begin{proof}
Deleting any bridge outside the paths and retaining its attachments
gives a proper minor of $G$ containing $\Psi_1$. A torus embedding
of this minor restricts to the required embedding of $\Psi_1$.

Now suppose four bridges outside the paths have the same attachment
set $A$. By three-connectivity, $|A|\ge3$. Delete one
bridge and take a torus embedding as above. The other three have the
same attachments, so Corollary~\ref{dis:choice-fails} puts them in
the same planar test graph. Contracting their interiors gives a
$K_{3,3}$ minor using three vertices of $A$, a contradiction.

If a bridge has one interior vertex and three attachments, those
attachments are the neighbors of a vertex of degree three. They
are therefore independent by Lemma~\ref{lem:independent-neighborhood}(a).
\end{proof}
Lemma~\ref{dis:relevant} concerns three-attachment bridges
with a planar relevant completion. If a bridge has none,
the next corollary and minor-minimality shows that it is the \emph{only} bridge and
every path outside it is an edge.

\begin{corollary}\label{dis:lone-bridge}
If $B$ has no planar relevant completion, then $G=\Psi_1\cup B$
and every path outside $B$ is an edge.
\end{corollary}
\begin{proof}
By Lemma~\ref{dis:joint-test}, $\Psi_1\cup B$ is nontoroidal, so
$G=\Psi_1\cup B$ by minimality. Any internal vertex of a path
outside $B$ would have degree two, contrary to three-connectivity.
\end{proof}

Suppose every bridge with at least two interior vertices has
three attachments and a planar relevant completion. By
Lemma~\ref{dis:relevant}, these bridges have the forms in
Proposition~\ref{prop:three-root-disjoint}. We combine them with
bridges with one interior vertex and the permitted edges, keeping
the interiors disjoint and using the attachment restrictions and
Lemma~\ref{dis:counts}. Four bridges with the same four attachments
would give a $K_{4,4}$ and hence a $K_5$ minor, so we allow at most
three. We use~\supp{E3} when all three paths are edges
and~\supp{E4} when some path has an internal vertex.

We use the next lemma when contracting the paths gives all four
twins the attachment triple of a bridge as their common neighborhood.

\begin{lemma}\label{lem:four-twins}
Let $B$ have three independent attachments $x,y,t$, with $B-\{x,y,t\}$ nonempty, connected and adjacent to all three. If adding four vertices with neighborhood $\{x,y,t\}$ to $B$ gives a toroidal graph, then at least one of $B+xy+xt$, $B+xy+yt$ and $B+xt+yt$ is planar.
\end{lemma}
\begin{proof}
Restrict the torus embedding to the $K_{3,4}$ on $x,y,t$ and the
four added vertices. Since $K_{3,4}$ is nonplanar, its regular
neighborhood has genus one, so its faces are disks by
Lemma~\ref{lem:torus-complement}. Since $B-\{x,y,t\}$ is connected, it lies in one face, with
$x,y,t$ on its boundary. Euler's formula gives five faces with
total boundary length $24$. Each of the other four boundary walks
has length at least four, so this one has length at most eight.

The walk alternates between attachments and added vertices, so
there are at most four attachment appearances. Consecutive
attachments are distinct, since otherwise the walk would reverse
an edge at an added vertex. As all three attachments occur, the
sequence is $xyt$ or $xyxt$ after relabeling. In both sequences
$x$ is consecutive with both $y$ and $t$, so $B+xy+xt$ is planar
by Lemma~\ref{cor:boundary-completions}.
\end{proof}

\subsection{When every path \texorpdfstring{$P_i$}{Pi} is a single edge}\label{dis:endpoint-support}

Here $P_i=s_it_i$ for each $i$, so every attachment of a
$\Psi_1$-bridge lies in $A_1\cup A_2$. We first consider graphs
with exactly one bridge. By Corollary~\ref{dis:lone-bridge}, a
bridge with no planar relevant completion must be the only bridge. 

\begin{lemma}\label{dis:sole}
If every path is an edge and $G$ has exactly one bridge, then
$G$ is isomorphic to one of $\mathfrak{T}_{5},\mathfrak{T}_{68},\mathfrak{T}_{69}$.
\end{lemma}
\begin{proof}
Let $B$ be the unique bridge, with attachment set $A$. If $B$ is
an edge, its ends are path ends, since all edges at the twins
belong to $\Psi_1$. By Lemma~\ref{dis:faces}, we can choose a
repeated path containing neither end and draw $B$ in a simple
face. This would embed $G$ in the torus. Thus $B$ has nonempty
interior and, by three-connectivity, at least three attachments.
We may assume some face in some torus embedding of $\Psi_1$
contains all of $A$, since otherwise
Lemma~\ref{dis:noncofacial} gives $G\cong\mathfrak T_5$.

By Lemma~\ref{dis:faces}, $A$ has at most four vertices and meets
each of $A_1,A_2$ in at most two. If $|A|=4$, the set
$\{s_i,t_i,s_j,t_j\}$ is excluded by
Lemma~\ref{dis:aligned-contraction}(a), and the set
$\{s_i,t_i,s_j,t_k\}$ with distinct indices is excluded by
Lemma~\ref{dis:skew}. Lemma~\ref{dis:aligned-contraction}(b)
also excludes a triple containing both $s_i,t_i$. After
interchanging $A_1,A_2$ and relabeling the paths if necessary,
we have $A=\{s_1,s_2,t_3\}$.

By Lemma~\ref{dis:faces}, only the repeated faces for
$(P_1;s_2,t_3)$ and $(P_2;s_1,t_3)$ contain $A$. Their attachment
sequences are $s_1s_2s_1t_3$ and $s_1s_2t_3s_2$, respectively.
By Lemma~\ref{cor:boundary-completions}, planarity of
$B+s_1s_2+s_1t_3$ or $B+s_1s_2+s_2t_3$ would let us draw $B$
in the corresponding face. Thus both completions are nonplanar,
as is the triangle completion. The latter is a three-connected
minor of $G$ by Lemma~\ref{dis:exterior-triangles}.

Contract $P_1,P_2,P_3$ into $s_1,s_2,t_3$, respectively. In this
toroidal proper minor of $G$, the four twins have common
neighborhood $A$. Lemma~\ref{lem:four-twins} gives a planar
two-edge completion of $B$, which must be $B+s_1t_3+s_2t_3$.

Apply Lemma~\ref{core:three-root-cut} with $(x,y,z)=(t_3,s_1,s_2)$.
It gives $B=T(s_1)$ or $T(s_2)$, or an allowed minor
$S(t_3;s_1)$, $S(t_3;s_2)$ or $Q(t_3;s_1,s_2)$. By
Lemma~\ref{core:three-graph-completions}, each $T$ form has a
planar completion containing $s_1s_2$, ruling out both forms.
The same lemma shows that each $S$ or $Q$ minor keeps both
completions containing $s_1s_2$ nonplanar. Replacing $B$ by
that minor therefore gives a nontoroidal graph, so minimality
forces $B$ to equal the minor. With $S(t_3;s_1)$ or
$Q(t_3;s_1,s_2)$, the resulting graph is isomorphic to
$\mathfrak T_{68}$ or $\mathfrak T_{69}$, respectively~\supp{W4}.
Interchanging the first two paths gives the other $S$ form.
\end{proof}

Write $a=s_1$, $b=t_1$, $c=s_2$, $d=t_2$, $e=s_3$,
$f=t_3$ and $A=\{a,b,c,d\}$. The boundary of each simple
face for $P_1,P_2$ contains $a,b,d,c$ in this cyclic order.

\begin{lemma}\label{dis:four-root-isolation}
Let $D$ be a component of $G-V(\Psi_1)$ with $N(D)=A$.
Then no other component of $G-V(\Psi_1)$ has all its
attachments in $A$, and neither $ad$ nor $bc$ is an edge.
Let $B$ be the $\Psi_1$-bridge formed by $D$ and its attaching
edges. If $B$ has a disk drawing with $a,b,d,c$ on the boundary in this
cyclic order, then $D$ is a single vertex.
\end{lemma}
\begin{proof}
Suppose another component $D'$ of $G-V(\Psi_1)$
has all its attachments in $A$. By three-connectivity,
it has at least three, so they include $a,d$ or
$b,c$. There is therefore a path joining that pair
with internal vertices in $D'$. Since $D$ is connected
and $N(D)=A$, there is a path joining the other pair
with internal vertices in $D$. These paths are
disjoint, contrary to Lemma~\ref{dis:opposite-order}.
An edge $ad$ or $bc$ would give the same contradiction
with a path through $D$ joining the other pair.

Now suppose $B$ has the stated disk drawing and
$D$ has more than one vertex. Contract $D$ to a
vertex $v$. The resulting proper minor is toroidal,
and $v$ lies in a face $F$ of the induced embedding
of $\Psi_1$ whose boundary contains $A$. By
Lemma~\ref{dis:faces}, $F$ is a simple face for
$P_1,P_2$ or the repeated face for
$(P_1;s_2,t_2)$ or $(P_2;s_1,t_1)$.
Its boundary meets $A_1\cup A_2$ only in $A$,
in one of the cyclic orders
\[
a,b,d,c,\qquad a,c,a,b,d,b,\qquad a,c,d,b,d,c,
\]
up to reversal. No other bridge with nonempty interior can lie
in $F$, since all its attachments would belong
to $A$, which we have just excluded. Also no edge bridge can lie there. The edges $ac,bd$
are excluded by the independence of $A_1,A_2$,
the edges $ad,bc$ were excluded above, and
$ab,cd$ already belong to $\Psi_1$. Each listed order contains $a,b,d,c$ as a cyclic
subsequence. Remove $v$ and its incident edges, and represent
$F$ by a disk with distinct boundary appearances
kept separate. Choose appearances of $a,b,d,c$
in that order. The Jordan--Schoenflies theorem
lets us insert the assumed drawing of $B$
through them, giving a torus embedding of $G$,
a contradiction.
\end{proof}

We now consider a component $D$ of $G-V(\Psi_1)$ with at
least four attachments. If $G\not\cong\mathfrak{T}_5$, we
show that, after relabeling, $N(D)=A$ and its bridge can
be drawn in a disk with $a,b,d,c$ in this cyclic order
on the boundary. Then Lemma~\ref{dis:four-root-isolation}
implies that $D$ is a single vertex.

\begin{proposition}\label{dis:four-attachments}
If a component of $G-V(\Psi_1)$ has at least four attachments, then
$G\cong\mathfrak{T}_{5}$, or the component is a single vertex with
neighborhood $\{s_i,t_i,s_j,t_j\}$ for some $i\ne j$.
\end{proposition}
\begin{proof}
Suppose $G\not\cong\mathfrak T_5$. By Lemmas~\ref{dis:noncofacial}
and~\ref{dis:skew}, the component has neighborhood
$\{s_i,t_i,s_j,t_j\}$ for distinct $i,j$. Relabel these indices as
$1,2$, so its neighborhood is $A=\{a,b,c,d\}$. Call the component
$D$, and let $B$ be the $\Psi_1$-bridge formed by $D$ and its
attaching edges.

Every path $P_i$ is the edge $s_it_i$, and the attachments of the $\Psi_1$-bridges all lie in $A_1\cup A_2$. By Lemma~\ref{dis:four-root-isolation}, it suffices to show that $B$ has a disk drawing in the order $a,b,d,c$. For a contradiction suppose it does not. Then $|V(D)|>1$. By Lemma~\ref{dis:faces}, in any torus embedding of a graph containing $\Psi_1\cup B$, either $P_1$ is the repeated path and its repeated face also contains $s_2,t_2$, or $P_2$ is the repeated path and its repeated face also contains $s_1,t_1$. In either case $B$ lies in that repeated face.

Suppose there is another component $D'$ of $G-V(\Psi_1)$. If three of its attachments have distinct indices, Lemma~\ref{dis:opposite-order} rules out any triple containing both $a$ and $d$ or both $b$ and $c$. Lemma~\ref{dis:noncofacial} rules out $A_1$ and $A_2$, leaving $\{e,b,d\}$~and~$\{f,a,c\}$.

For $\{e,b,d\}$, contract $D,D'$ to single vertices. A face containing $e,b,d$ must be the repeated face by Lemma~\ref{dis:faces}. $P_3=ef$ cannot be the repeated path, since it could contain only one additional $t$-vertex, not both $b=t_1$ and $d=t_2$. Hence the embedding choice is either $(P_1;s_3,t_2)$ or $(P_2;s_3,t_1)$. Neither choice has a face containing all of $A$, so the contracted vertex representing $D$ cannot be placed in a toroidal embedding, contradicting minor minimality, since $|V(D)|>1$. Swapping $A_1$ and $A_2$ gives the same contradiction for $\{f,a,c\}$.

Therefore the attachments of $D'$ must lie in one of $\{a,b,c,d\}$, $\{a,b,e,f\}$, or $\{c,d,e,f\}$. Lemma~\ref{dis:four-root-isolation} excludes the first set. Any extra edges outside $\Psi_1$ between its vertices must join two of $a,b,c,d,e,f$. Since $A_1,A_2$ are independent and by Lemma~\ref{dis:four-root-isolation}, only $af,eb,cf,ed$ are possible.

Every remaining attachment triple contains one of $ab,cd,ef$,
so it has two vertices adjacent in $G$. If its bridge could be drawn in a simple face, its triangle completion would be planar by Lemma~\ref{cor:boundary-completions}, and Lemma~\ref{dis:relevant} would make its interior a single vertex. That vertex has two adjacent neighbors, contradicting Lemma~\ref{lem:independent-neighborhood}(a). Thus none of these three-attachment bridges can be drawn in a simple face. $D'$ therefore has either three attachments and a bridge that cannot be drawn in a simple face, or four attachments whose bridge may or may not be drawable in a simple face.

Suppose the bridge $B'$ associated with $D'$ cannot be drawn in a simple face. Since $D'$ has at least three attachments, it cannot be placed in an end face, which contains only two possible attachments. In either choice available to $B$ discussed earlier, the repeated face contains only $a,b,c,d$ among $A_1\cup A_2$. But $D'$ cannot have all its attachments in $\{a,b,c,d\}$, by Lemma~\ref{dis:four-root-isolation}. Thus $\Psi_1\cup B\cup B'$ is nontoroidal, and anything else in $G$ can be deleted without making it toroidal.

Suppose instead at least one edge from each of $\{af,eb\}$ and $\{cf,ed\}$ is present in $G$. No end face can contain any of these edges, as its vertices are in only one of $A_1,A_2$. In $(P_1;s_2,t_2)$ neither the repeated face nor either simple face can contain an edge from $\{af,eb\}$. In $(P_2;s_1,t_1)$ neither the repeated face nor either simple face can contain an edge from $\{cf,ed\}$.

In either of these two cases, retain $\Psi_1\cup B$ together with
the other bridge, or one edge from each pair, and delete everything
else. The resulting graph is still nontoroidal, so minimality makes
it equal to $G$.

Apply Theorem~\ref{thm:mohar} to $B\cup C_0$, where $C_0=abdca$.
The auxiliary graph is three-connected by
\cite[Lemma~4.2]{Mohar1994}, applied to $B$ in $G+E(C_0)$.
Lemma~\ref{dis:opposite-order} rules out crossing paths, so $B$
contains a tripod. By interchanging the path indices $1,2$ or
the triples $A_1,A_2$ as needed, take its attachments to be
$a,b,d$. If all three attachment paths have length zero, the two
trees left after deleting $a,b,d$ each have a neighbor at all
three attachments. Join the trees by a path in $D$, contract
each tree to a vertex, and shorten the joining path to an edge.
Shortening $ab$, $brd$ and $dzfeqa$ to edges then gives a $K_5$
minor. Thus some attachment path has positive length.

Take a shortest path from $c$ to the tripod through $D$ and
shorten the tripod as in Lemma~\ref{endpoint:primitive}. If the
path from $c$ ends on an attachment path, shorten that attachment
path to $ug$, contracting the meeting vertex to $u$. Otherwise
contract the meeting vertex to the degree-three end of its branch
path and shorten any positive-length attachment path to $ug$.
Shorten the path from $c$ to an edge and contract the rest of $D$
onto the retained interior. This gives $T(g)$ on $\{a,b,d\}$,
with twins $x,y$, together with an edge from $c$ to $u$, $x$ or $y$.

A disk drawing of this bridge in the order $a,b,d,c$ would let
us add the triangle on $a,b,d$ outside the disk. This contradicts
the nonplanarity of the triangle completion of $T(g)$ by
Lemma~\ref{core:three-graph-completions}. The other bridge or
edges still rule out the two repeated faces that contain
$a,b,c,d$. The reduced graph is therefore nontoroidal, so
minimality~gives~$B$~itself.

If $g=d$, one of $x,y$ has degree three with adjacent neighbors
$a,b$, contradicting Lemma~\ref{lem:independent-neighborhood}(a).
If $cx$ is present, use the disjoint paths $auyd,bxc$ for $g=a$
and $ayd,buxc$ for $g=b$. Lemma~\ref{dis:opposite-order} excludes
both pairs, and interchanging $x,y$ also excludes $cy$. If $cu$
is present and $g=b$, the same lemma excludes $axd,buc$. Hence
$B=T(a)+cu$.

In $G-ua$ the remaining bridge is $T(c)$ on $\{b,c,d\}$, and in
$G-uc$ it is $T(a)$ on $\{a,b,d\}$. Their triangle completions are nonplanar, so neither bridge fits
in a simple face. In both $G-ua$ and $G-uc$, the second bridge
or the two extra edges still rule out $(P_1;s_2,t_2)$ and
$(P_2;s_1,t_1)$.
By Lemma~\ref{dis:faces}, the torus embeddings must therefore
use $(P_2;s_3,t_1)$ for $G-ua$ and $(P_1;s_3,t_2)$ for $G-uc$.

If there is a second bridge, it must use the repeated face in
both embeddings. Its attachments therefore lie in both
$\{b,c,d,e\}$ and $\{a,b,d,e\}$. By three-connectivity its
attachment set is $\{b,d,e\}$, which was already excluded.
In the edge case, the first embedding rules out $cf$ and the
second rules out $af$, leaving $eb,ed$. Delete $y$, contract
$bx$, and delete $ab,cd$ and the new edge $bd$. This gives
twins $p,q$ with neighborhood $\{a,c,e\}$, twins $d,f$ with
neighborhood $\{r,z,e\}$, and the edges $ua,uc,ub,br,bz,be$.
The resulting graph is the proper nontoroidal minor
$\mathfrak T_6$, a contradiction.

It now remains only to consider $D'$ with exactly four attachments and $B'$ drawable in a simple face. By Lemma~\ref{dis:four-root-isolation}, $D'$ must be a single vertex $h$. Up to interchanging the indices $1,2$ of the paths $P_i$, assume $N(h)=\{a,b,e,f\}$. Contract $ha$. The edge $hb$ becomes the existing edge $ab$, $hf$ becomes $af$, and $he$ becomes $ae$. Delete $ae$. The resulting graph is a proper minor and hence toroidal. As $B$ is unchanged, its embedding uses one of $(P_1;s_2,t_2),(P_2;s_1,t_1)$. In the former no face contains both $a,f$, a contradiction. Thus it uses $(P_2;s_1,t_1)$, and $af$ lies in a simple face for $P_1,P_3$, whose boundary contains $a,b,e,f$.

Applying Lemma~\ref{dis:four-root-isolation} to $h$ excludes another component having all its attachments in $\{a,b,e,f\}$ and excludes $af,eb$ in the original graph. The edges $ae,bf$ are absent by independence, while $ab,ef$ already belong to $\Psi_1$. Thus $af$ is the only addition inside this face. Erase $af$, place $h$ in the now empty disk and draw its four edges to $a,b,e,f$, restoring $G$ without crossings, a contradiction. Therefore there are no extra components $D'$.

Only the original bridge $B$ and possible extra edges remain outside $\Psi_1$. There is not an edge from both $\{af,eb\}$ and $\{cf,ed\}$, and Lemma~\ref{dis:opposite-order} says both edges cannot occur in either pair, so there is at most one extra edge. If there are no extra edges, $B$ is the unique bridge. The proper minor
$G/ef$ is toroidal, so Lemma~\ref{dis:aligned-contraction}(a), with
$H=G$, would imply that $G$ is toroidal. Thus exactly one extra edge remains, and after interchanging indices or exchanging $A_1,A_2$ or both, assume it is $af$. Hence
\[
G=(\Psi_1\cup B)+af.
\]
For every allowed minor $B'$ of $B$, we have
\[
(\Psi_1\cup B')+af\text{ is toroidal}
\quad\Longleftrightarrow\quad B'+ac+cd+bd\text{ is planar}.
\]
In a torus embedding, $B'$ lies in a face containing $a,b,c,d$.
If it lies in a simple face, Lemma~\ref{dis:joint-test} requires
$B'+ab+cd+ac+bd$ to be planar, so $B'+ac+cd+bd$ is planar.
By Lemma~\ref{dis:faces}, the only repeated faces containing
$a,b,c,d$ are those for $(P_1;s_2,t_2)$ and $(P_2;s_1,t_1)$.
The first choice has no face containing both $a$ and $f$, and
the second requires $B'+ac+cd+bd$ to be planar.

Conversely, if $B'+ac+cd+bd$ is planar, use $(P_2;s_1,t_1)$.
Assign $B'$ to the repeated face and $af$ to a simple face for
$P_1,P_3$. The simple-face test is the cycle $abfea$ with diagonal
$af$, so both tests are planar. Lemma~\ref{dis:joint-test} gives
a torus embedding.

Taking $B'=B$ shows that $B+ac+cd+bd$ is nonplanar. It is a minor
of $G$, obtained by retaining $B$ and $cd$ and shortening $apc$
and $brd$ to $ac$ and $bd$, so it has no $K_5$ minor. By the
equivalence, Corollary~\ref{common:reduced-forms} puts $B$ among
the candidates with attachments $a,b,c,d$ and added edges
$ac,cd,bd$. For each of the $15{,}626$ candidates, the graph
$(\Psi_1\cup B)+af$ has an edge whose deletion or contraction is
nontoroidal~\supp{E2}, a contradiction. Thus $B$ has a disk drawing
in the order $a,b,d,c$, and Lemma~\ref{dis:four-root-isolation}
gives that $D$ is a single vertex.
\end{proof}

We have determined $\mathfrak T_5$ and the graphs with only one
bridge. By Corollary~\ref{dis:lone-bridge} and
Proposition~\ref{dis:four-attachments}, every remaining bridge
with nonempty interior either has three attachments and a planar
relevant completion or has one interior vertex and four
attachments. We now bound the number of three-attachment bridges,
then form the candidates by adding these two kinds of bridges
and the edges $s_it_j$ with $i\ne j$ to $\Psi_1$.

\begin{theorem}\label{dis:endpoint-only}
A three-connected obstruction containing $p,q,r,z$ with
neighborhoods $A_1,A_2$, where $A_1\cap A_2=\varnothing$, and
the three edges $s_1t_1,s_2t_2,s_3t_3$ is isomorphic to one of
\[
\mathfrak T_5,\quad \mathfrak T_{48},\quad \mathfrak T_{58},
\quad \mathfrak T_{68},\quad \mathfrak T_{69}.
\]
\end{theorem}
\begin{proof}
Assume $G\not\cong\mathfrak T_5$. By three-connectivity and
Lemma~\ref{dis:noncofacial}, a bridge with nonempty interior has
at least three attachments on one face boundary. Hence its
attachment set meets both $A_1$ and $A_2$. By
Proposition~\ref{dis:four-attachments}, a bridge with four or
more attachments consists of one interior vertex with neighborhood
$\{s_i,t_i,s_j,t_j\}$ for some $i\ne j$. If a three-attachment
bridge has no planar relevant completion, it is the only bridge
by Corollary~\ref{dis:lone-bridge}, and Lemma~\ref{dis:sole} gives
$G\cong\mathfrak T_{68}$ or $\mathfrak T_{69}$. We may thus assume
every three-attachment bridge has a planar relevant completion.

Take a torus embedding of a proper minor obtained by deleting
bridges and retaining their attachments and $\Psi_1$. If a
three-attachment bridge $B$ lies in a simple face, its attachments
are three vertices of $C=s_it_it_js_js_i$, so they include both
ends of $P_i$ or $P_j$. These ends are adjacent, so
Lemma~\ref{lem:independent-neighborhood}(a) rules out a single
interior vertex. By Lemma~\ref{dis:relevant}, the triangle completion
of $B$ is nonplanar. Yet contracting an edge of $C$ incident with
the vertex outside the attachment set gives this completion as
a minor of the planar simple-face test, a contradiction.

Thus all $m$ three-attachment bridges lie in the repeated face.
In its planar test graph, contract each bridge interior to one
vertex and keep just these vertices, the four possible attachments,
and the three edges from each contracted vertex to its attachments.
The resulting simple bipartite graph is planar, so
$3m\le 2(m+4)-4$ and $m\le4$. Deleting one three-attachment bridge
of $G$ therefore shows that $G$ has at most five. If it has five,
deleting any other bridge would leave all five in a toroidal proper
minor. Thus in that case there are no four-attachment bridges or
edges $s_it_j$ with $i\ne j$.

Up to interchanging $A_1,A_2$, a three-attachment bridge has
attachments $c_1,c_2\in A_1$ and $h\in A_2$. By
Lemma~\ref{dis:relevant} and Proposition~\ref{prop:three-root-disjoint},
it has one interior vertex or is one of
\[
T(c_1),\quad T(c_2),\quad T(h),\quad S(c_1;c_2),\quad S(c_2;c_1).
\]
If the triple contains a path edge, label it $c_1h$. The
single-vertex form, $T(c_2)$ and $S(c_1;c_2)$ each have a vertex
of degree three adjacent to both $c_1$ and $h$, so
Lemma~\ref{lem:independent-neighborhood}(a) leaves only
$T(c_1),T(h),S(c_2;c_1)$.

Of the eighteen triples meeting both $A_1,A_2$, twelve contain
a path edge. This gives $12\cdot3+6\cdot6=72$ labeled forms.
Add these bridges to $\Psi_1$ with disjoint interiors, using at
most three on each triple by Lemma~\ref{dis:counts} and at most
five in total. By Lemma~\ref{dis:four-root-isolation}, allow at
most one four-attachment bridge on each $\{s_i,t_i,s_j,t_j\}$.
We also allow the six edges $s_it_j$ with $i\ne j$, subject to
the restriction above when there are five three-attachment bridges.

We keep the graphs satisfying Lemmas~\ref{lem:independent-neighborhood}(a)
and~\ref{dis:four-root-isolation}. By Lemma~\ref{dis:opposite-order},
we exclude graphs with disjoint $s_i$--$t_j$ and $s_j$--$t_i$ paths
through distinct bridges. By~\supp{E3}, each remaining graph is a
subgraph of a recorded toroidal graph or contains a recorded
nontoroidal graph. Since $G$ is a minor obstruction, it equals one
of the latter. All but $\mathfrak T_{48}$ and $\mathfrak T_{58}$ have a
proper nontoroidal minor, and neither has a four-attachment bridge
in this case. Together with $\mathfrak T_5$,
$\mathfrak T_{68}$ and $\mathfrak T_{69}$, these give the stated list.
\end{proof}

\subsection{Three-attachment bridges and four-attachment bridges with one interior vertex}
\label{dis:three-construction}

For the remainder of Section~\ref{sec:disjoint}, we suppose now that some path $P_i$ is not an edge. For this subsection, we consider the case where every bridge with four attachments has one interior vertex. By
Proposition~\ref{dis:types}, every bridge with at least two
interior vertices then has three attachments. The next corollary of
Corollary~\ref{dis:lone-bridge} and Lemma~\ref{dis:aligned-contraction}
shows that these bridges have planar relevant completions. 

\begin{corollary}\label{dis:remaining-completions}
Every bridge with three attachments has a planar relevant completion.
\end{corollary}
\begin{proof}
A bridge with one interior vertex has a planar triangle completion.
Suppose that a bridge $B$ has no planar relevant completion.
By Corollary~\ref{dis:lone-bridge}, $G=\Psi_1\cup B$ and the paths
outside $B$ are edges. Since some path is not an edge, $B$
contains a path $P_m$. Its third attachment belongs to
$\{s_x,t_x\}$ for some $x\ne m$. Every $\Psi_1$-bridge in $B$ has an attachment in the interior
of $P_m$ and all its attachments in
$V(P_m)\cup\{s_x,t_x\}$. Let $w$ be the third index.
The graph $G/s_wt_w$ is a toroidal proper minor, so
Lemma~\ref{dis:aligned-contraction}(c), after relabeling the
indices, makes $G$ toroidal.
\end{proof}
\begin{proposition}\label{dis:remaining-classification}
If some $P_i$ is not an edge and every bridge with four
attachments has one interior vertex, then $G$ is isomorphic
to one of
\[
\mathfrak T_{32},\quad \mathfrak T_{38},\quad
\mathfrak T_{43},\quad \mathfrak T_{45},\quad
\mathfrak T_{48},\quad \mathfrak T_{61},\quad
\mathfrak T_{66},\quad \mathfrak T_{67}.
\]
\end{proposition}

\begin{proof}
By Corollary~\ref{dis:remaining-completions}, every three-attachment
bridge has a planar relevant completion. On each of the eighteen
triples $\{c_1,c_2,h\}$, with $c_1,c_2$ in one of $A_1,A_2$ and $h$
in the other, Lemma~\ref{dis:relevant} and
Proposition~\ref{prop:three-root-disjoint} allow a single interior
vertex or one of the five forms
\[
T(c_1),\quad T(c_2),\quad T(h),\quad S(c_1;c_2),\quad S(c_2;c_1).
\]
A bridge containing $P_m$ may also have the form $Q(h;g,e)$ or either
four-cycle form of Lemma~\ref{core:two-cap-critical}, with $g,e,h$ as
in Lemma~\ref{dis:relevant}(ii). These three further forms are possible
on the twelve triples containing both ends of some $P_m$. This gives
$18\cdot6+12\cdot3=144$ labeled three-attachment forms.

Every four-attachment bridge has one interior vertex. By
Proposition~\ref{dis:types}, its neighborhood is one of the three
sets $\{s_i,t_i,s_j,t_j\}$ or one of the six sets
$\{s_m,t_m,s_x,t_y\}$ of type (c). Together with the nine possible
edges between $A_1$ and $A_2$, this gives $144+9+9=162$ choices of
an edge or bridge to add to $J_1\cup J_2$.

We combine these with disjoint interiors. By Lemma~\ref{dis:counts},
at most three bridges outside the paths have the same attachment
triple. If the triple contains $s_m,t_m$ and $s_mt_m$ is absent, we
allow one more bridge for $P_m$. At most three vertices have the same
four neighbors, since four would give a $K_{4,4}$ and hence a $K_5$
minor. For each $m$, at most one of the $Q$, four-cycle or type (c)
forms containing $P_m$ can occur, and it requires $s_mt_m$ to be absent.

By Lemmas~\ref{lem:independent-neighborhood}(a) and~\ref{dis:opposite-order},
we exclude a degree-three interior vertex with adjacent neighbors and
any choice containing both $s_it_j$ and $s_jt_i$ for $i\ne j$. If $s_mt_m$ is absent, then $P_m$ has an internal vertex. The bridge containing $P_m$ therefore has both $s_m$ and $t_m$ as attachments. At least one of the three edges $s_mt_m$ is absent
by~\eqref{dis:shortest}.

These choices include $G$ without having to specify the paths inside
each candidate. By the coverage check in~\supp{E4}, each graph in
this family is a subgraph of a recorded toroidal graph or contains
a recorded nontoroidal graph. Since $G$ is nontoroidal, by minimality
it equals any such recorded nontoroidal graph. With at least
one $s_mt_m$ absent, the checks leave $\mathfrak T_{32},\mathfrak
T_{43},\mathfrak T_{61}$ when every bridge has at most one interior
vertex, and $\mathfrak T_{38},\mathfrak T_{45},\mathfrak T_{48},
\mathfrak T_{66},\mathfrak T_{67}$ when some bridge has form $T,S$
or $Q$. No graph with a four-cycle form remains.
\end{proof}

\subsection{A bridge with four attachments and at least two interior vertices}
\label{dis:path-four}

We now consider a bridge with four attachments and at least two interior vertices. We treat types (c), (a) and (b) of Proposition~\ref{dis:types} in that order. 

\begin{proposition}\label{dis:type-c}
If $G$ has a bridge of type (c) with at least two interior vertices, then $G\cong\mathfrak{T}_{21}$ or $G\cong\mathfrak{T}_{63}$.
\end{proposition}
\begin{proof}
Let $B$ be such a bridge. Number the paths so that $B$ contains
$P_1$ and its other attachments are $s_2,t_3$. Exchange $d$ and $f$
in the notation of Section~\ref{dis:endpoint-support}. The ends of
$P_1,P_2,P_3$ are now $(a,b)$, $(c,f)$ and $(e,d)$, respectively, and
the attachment set of $B$ is $A=\{a,b,c,d\}$.
The bridge $B$ consists of $P_1$ and all $\Psi_1$-bridges meeting its
interior. We first determine the embedding choice that could extend
to $G$. If $P_2$ were repeated, the interior of $P_1$ would occur
only on the boundaries of the simple faces following $P_1,P_3$.
Neither contains $c$, so the bridge joining $c$ to the interior of
$P_1$ could not be drawn. Similarly, if $P_3$ were repeated, the
simple faces following $P_1,P_2$ would miss $d$. Therefore $P_1$
must be repeated, and the bridges joining its interior to $c,d$
require both vertices on the repeated-face boundary. By
Lemma~\ref{dis:faces}, the choice must be $(P_1;s_2,t_3)$. Figure~\ref{fig:type-c-embedding} shows this embedding.

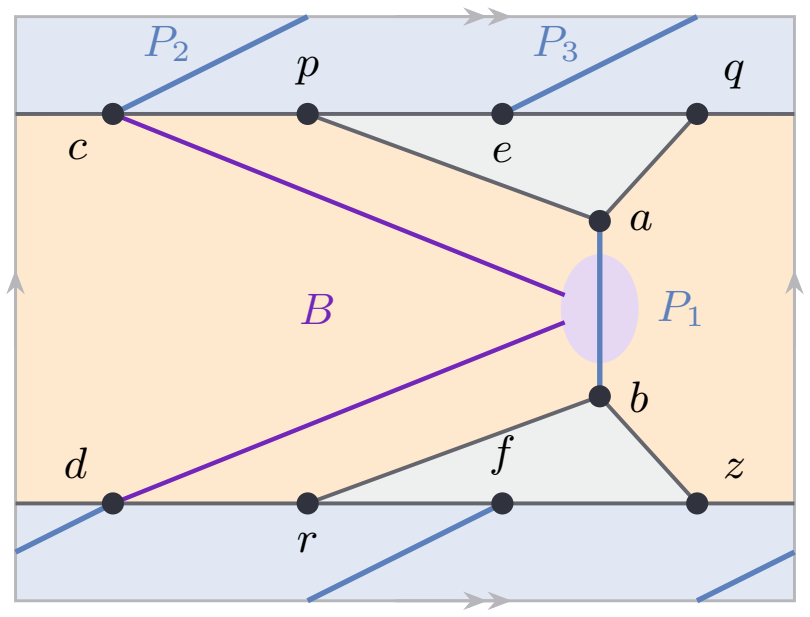
\begin{figure}[htbp]
\centering
\begingroup
\scaledfigure{1.15}{\tikzsetnextfilename{type-c-embedding}%
\begin{tikzpicture}[
  scale=1.4,
  font=\small,
  vertex/.style={circle,fill=figureInk,inner sep=1.6pt},
  edge/.style={draw=figureInk!75,line width=.8pt},
  pathedge/.style={draw=figurePath,line width=1.1pt},
  bridge/.style={draw=figureBridge,line width=.9pt},
  side/.style={draw=figureInk!35,line width=.6pt}
]

  \fill[faceSimple] (0,0) rectangle (4,3);
  \fill[faceRepeated] (0,.5) rectangle (4,2.5);
  \fill[faceOther] (1.5,2.5)--(3.5,2.5)--(3,1.95)--cycle;
  \fill[faceOther] (1.5,.5)--(3.5,.5)--(3,1.05)--cycle;
  \draw[side] (0,0) rectangle (4,3);
  \draw[side,-{Stealth[length=5pt]}] (0,1.3)--(0,1.7);
  \draw[side,-{Stealth[length=5pt]}] (4,1.3)--(4,1.7);
  \draw[side,-{Stealth[length=5pt] Stealth[length=5pt]}] (1.95,0)--(2.55,0);
  \draw[side,-{Stealth[length=5pt] Stealth[length=5pt]}] (1.95,3)--(2.55,3);
  \coordinate (c) at (0.5,2.5);
  \coordinate (p) at (1.5,2.5);
  \coordinate (e) at (2.5,2.5);
  \coordinate (q) at (3.5,2.5);
  \coordinate (d) at (0.5,.5);
  \coordinate (r) at (1.5,.5);
  \coordinate (f) at (2.5,.5);
  \coordinate (z) at (3.5,.5);
  \coordinate (a) at (3.0,1.95);
  \coordinate (b) at (3.0,1.05);
  \fill[figureBridge!18] (3,1.5) ellipse (.2 and .28);
  \draw[bridge] (2.82,1.57)--(c);
  \draw[bridge] (2.82,1.43)--(d);
  \draw[edge] (0,2.5)--(c)--(p)--(e)--(q)--(4,2.5);
  \draw[edge] (0,.5)--(d)--(r)--(f)--(z)--(4,.5);
  \draw[edge] (p)--(a)--(q);
  \draw[edge] (r)--(b)--(z);
  \draw[pathedge] (a)--(b);
  \draw[pathedge] (c)--(1.5,3);
  \draw[pathedge] (1.5,0)--(f);
  \draw[pathedge] (e)--(3.5,3);
  \draw[pathedge] (3.5,0)--(4,.25);
  \draw[pathedge] (0,.25)--(d);
  \node[vertex,label=below left:{$c$}] at (c) {};
  \node[vertex,label=above:{$p$}] at (p) {};
  \node[vertex,label=below:{$e$}] at (e) {};
  \node[vertex,label=above right:{$q$}] at (q) {};
  \node[vertex,label=above left:{$d$}] at (d) {};
  \node[vertex,label=below:{$r$}] at (r) {};
  \node[vertex,label=above:{$f$}] at (f) {};
  \node[vertex,label=above right:{$z$}] at (z) {};
  \node[vertex,label=right:{$a$}] at (a) {};
  \node[vertex,label=right:{$b$}] at (b) {};
  \node[text=figureBridge] at (1.55,1.5) {$B$};
  \node[text=figurePath] at (3.42,1.5) {$P_1$};
  \node[text=figurePath] at (.78,2.86) {$P_2$};
  \node[text=figurePath] at (2.78,2.86) {$P_3$};
\end{tikzpicture}}
\caption{The torus embedding of $\Psi_1$ for the choice $(P_1;s_2,t_3)$, with $B$ in the repeated face.}
\label{fig:type-c-embedding}
\endgroup
\end{figure}

For this choice, $H_{\mathrm{rep}}$ contains $B+ac+bd$, while
$H_{\mathrm{simp}}$ contains $P_2,P_3$ and the added edges $ce,fd$.
The boundaries of the repeated face and the simple faces share
only $c,d$ among the vertices where a bridge can attach.
Since a bridge with nonempty interior has at least
three attachments, only the possible edge $cd$ can be assigned
to either test. Leave it out of both for now.

Contract the interior of $B$ to a single vertex $w$. This gives
a proper minor, since $B$ has at least two interior vertices,
and hence a toroidal graph. The path $a-w-b$ replaces $P_1$,
and the edges $wc,wd$ still force the same embedding choice by
the preceding argument. Since $B$ is assigned to the repeated face, the contraction
leaves $H_{\mathrm{simp}}$ unchanged. Applying
Lemma~\ref{dis:joint-test}(a) to the toroidal minor shows that
$H_{\mathrm{simp}}$ is planar. Also, another
bridge of type (c) containing $P_2$ or $P_3$ would attach at $a$
or $b$. It would have to lie in a simple face, but neither simple
face contains $a$ or $b$. Thus no such bridge occurs.

\looseness=-1
Suppose for a contradiction that $cd\in E(G)$. The proper minor
$G-cd$ is toroidal, so both $H_{\mathrm{simp}}$ and
$H_{\mathrm{rep}}$ are planar. If $H_{\mathrm{simp}}+cd$ were
planar, Lemma~\ref{dis:joint-test} would give a torus embedding
of $G$. Hence $H_{\mathrm{simp}}+cd$ must be nonplanar.

If a bridge joins $x\in V(P_2)-\{c\}$ to
$y\in V(P_3)-\{d\}$, it contains an $x$--$y$ path with internal
vertices outside $\Psi_1$. This path is disjoint from $cd$. Its end $x$ follows $c$ on
$P_2$, whereas $y$ precedes $d$ on $P_3$, contrary to
Lemma~\ref{dis:opposite-order}. Therefore no bridge joins
$P_2-c$ to $P_3-d$. Each bridge in $H_{\mathrm{simp}}$ has all
its attachments in one of $V(P_2)\cup\{d\}$ and
$V(P_3)\cup\{c\}$.

Let $L$ consist of $P_2+fd$ and the bridges with attachments in
$V(P_2)\cup\{d\}$, and let $R$ consist of $P_3+ce$ and the bridges
with attachments in $V(P_3)\cup\{c\}$. Then $H_{\mathrm{simp}}=L\cup R$,
and $L,R$ meet only at $c,d$. Each contains a $c$--$d$ path. Keeping
one part and shortening such a path in the other to $cd$ gives
$L+cd$ and $R+cd$ as minors of the planar graph $H_{\mathrm{simp}}$.
Thus $H_{\mathrm{simp}}+cd$ is planar as the edge amalgamation
of $L+cd$ and $R+cd$ along $cd$. Since $H_{\mathrm{rep}}$ is also planar, assigning $cd$ to
the simple faces gives a torus embedding of $G$ by
Lemma~\ref{dis:joint-test}. Thus $cd\notin E(G)$.

Since $H_{\mathrm{simp}}$ is planar, $H_{\mathrm{rep}}$ is nonplanar.
Deleting any bridge assigned to the simple faces, or shortening
$P_2$ or $P_3$ without identifying its ends, would give a proper
minor with the same forced choice and the same nonplanar
repeated-face test. By minimality there is no such bridge and
neither path has an internal vertex. Thus
\[
P_2=cf,\qquad P_3=ed,\qquad H_{\mathrm{simp}}=cfdec.
\]

If another bridge with nonempty interior attaches at both $a,b$,
replace $P_1$ by an $a$--$b$ path through that bridge. The interior
of $B$ then lies outside the new subdivision and attaches at
$s_1,t_1,s_2,t_3$, which Lemma~\ref{dis:skew} forbids. Every other
bridge with nonempty interior has at least three attachments among
$a,b,c,d$, so its attachment set is $\{a,c,d\}$ or $\{b,c,d\}$.

Contract the interior of $B$ to $w$ again. In the planar
repeated-face test of this proper minor, keep one of the other
bridges and the three edges from $w$ to its attachments. Replacing
these three edges by a triangle gives a planar triangle completion
of the bridge. Its interior is therefore a single vertex by
Lemma~\ref{lem:flat-replacement}. Two such vertices on the same
triple would give a $K_{3,3}$ with $w$ in the planar test. Thus there
is at most one vertex $h_a$ with neighborhood $\{a,c,d\}$ and at
most one vertex $h_b$ with neighborhood $\{b,c,d\}$.

Of the edges on $A$, independence excludes $ac,bd$,
\eqref{dis:shortest} excludes $ab$, and we have excluded $cd$.
If $h_a$ is present, Lemma~\ref{lem:independent-neighborhood}(a)
excludes $ad$. Likewise, $h_b$ excludes $bc$. Interchanging $a,b$ and $c,d$ together exchanges $ad$ with $bc$
and $h_a$ with $h_b$. If neither edge is present, there may be
no further bridge, just the bridge at $h_a$, or bridges at both
$h_a$ and $h_b$. If exactly one edge is present, take it to be
$ad$. Then the only possible further bridge is the one at $h_b$.
If both edges are present, they exclude both $h_a$ and $h_b$.

For each of the six choices, let $H_0$ consist of $A$, the chosen
edges and vertices with their incident edges, and $ac,bd$. Then
$H_{\mathrm{rep}}=B\cup H_0$. Replacing $ac,bd$ by the paths
$apc,brd$ gives a subgraph of $G$, so $H_{\mathrm{rep}}$ is a
minor of $G$. It is nonplanar and has no $K_5$ minor, so it
contains a $K_{3,3}$ subdivision.

Let $B_0$ be obtained from $B$ by one allowed operation. Choose
an $a$--$b$ path through its interior to replace $P_1$. Since
$B_0-A$ is connected and both $c,d$ have neighbors in it, each
has a path to the interior of $P_1$ through $B_0-A$. This again
forces $(P_1;s_2,t_3)$. The proper minor with $B_0$ in place of
$B$ is toroidal, so $B_0\cup H_0$ is planar. In all six cases,
$H_0$ is planar and every vertex of $H_0-A$ has all its neighbors
in $A$. Thus $B$ itself has the form given by
Lemma~\ref{core:retained-subdivision}, with at most $|A|+5=9$
interior vertices.

Delete any one of $ad,bc,h_ac,h_ad,h_bc,h_bd$ that is present.
The resulting proper minor is toroidal, and $B$ still forces
$(P_1;s_2,t_3)$. Deleting an edge at $h_a$ or $h_b$ leaves a
two-edge path. The only such paths that can lie in an end face
are $a h_a c$ and $b h_b d$, which we can draw beside $ac$ and
$bd$, respectively, in the repeated-face test. Thus deleting any
of the listed edges makes $H_{\mathrm{rep}}$ planar, so every
$K_{3,3}$ subdivision in $H_{\mathrm{rep}}$ uses each listed edge
that is present. If $h_a$ is present, contracting $ah_a$ in the
planar graph $H_{\mathrm{rep}}-h_ac$ gives
$H_{\mathrm{rep}}-h_a+ad$. Likewise, $H_{\mathrm{rep}}-h_b+bc$
is planar when $h_b$ is present.

Suppose neither $h_b$ nor $bc$ is present. Take a torus embedding
of $G/cf$ and identify the two boundary copies of $P_1$ in its
incident facial disks, keeping all other boundary appearances
distinct. The resulting region contains $P_1$ and every bridge
meeting its interior, hence all of $B$. Among the regions along
$P_1$ in $\Psi_1/cf$, those whose boundary contains $c,d$ are an
annulus with sequences $ac\mid bd$ and disks with sequence $acbdc$
or $acbcdc$. In each case no other face contains both $a,d$~\supp{R3},
so the region also contains $h_a$ or $ad$ if present.

In the annulus, capping the two boundary circles and drawing $ac$
in one and $bd$ in the other would give a plane drawing of
$H_{\mathrm{rep}}$. In either disk, join consecutive attachment
appearances along the boundary. Join the copies of $c$ by a tree
in the complementary disk and contract the tree. For $acbdc$ this
gives a plane drawing of $H_{\mathrm{rep}}+bc+cd$, again
impossible. Thus the order is $acbcdc$, which gives a plane
drawing of $H_{\mathrm{rep}}-bd+bc+cd$. In particular,
$H_{\mathrm{rep}}-bd$ is planar, so every $K_{3,3}$ subdivision
uses $bd$. Similarly, if neither $h_a$ nor $ad$ is present,
contracting $ed$ shows that $H_{\mathrm{rep}}-ac$ is planar
and every such subdivision uses $ac$.

Apply Construction~\ref{common:marked-construction} to each of
the six choices of $H_0$, using the required subdivision edges
and planarity conditions above. Adding $J_1\cup J_2$, the edges
$cf,ed$, and the chosen additions to each candidate for $B$ gives
$1{,}586$ candidates for $G$. Each has a $K_5$ minor, a cut of at
most two vertices, or a proper nontoroidal minor, except for
$\mathfrak T_{21}$ and $\mathfrak T_{63}$~\supp{E5}.
\end{proof}

We now consider bridges of types~(a) and~(b) with at least two
interior vertices. For type~(a), let $B_a$ have attachments
$\{s_x,t_x,s_y,t_y\}$ and put
\[
C=s_xt_xt_ys_ys_x,\qquad C'=s_xt_yt_xs_ys_x.
\]
For type~(b), let $B_b$ contain $P_m$ and write $a=s_m$, $b=t_m$,
$c=s_x$ and $d=t_x$ for its four attachments. In either case, let
$w$ be the index of the path with neither end in the attachment set.

Let $B$ be either bridge and let $A$ be its attachment set.
After deleting at most two vertices from $B\cup K_4[A]$, the
remaining attachments still form a clique. Any component disjoint
from them would also be separated in $G$ by the same deletion.
Thus $B\cup K_4[A]$ is three-connected.

The next lemma uses the cyclic order of the attachments around
a tree obtained from $B$. If adding the four-cycle in that order
makes a planar graph, we can draw $B$ in place of the tree.
\begin{lemma}\label{dis:four-completions}
Let $B$ be $B_a$ or $B_b$, with attachment set $A$, and let $T$ be
an allowed minor of $B$ that is a tree whose vertices of degree one
are exactly the four attachments. Consider a torus embedding of the
graph obtained from $G$ by replacing $B$ by $T$. A thin closed disk
around $T$ meets the rest of the embedding precisely in $A$. We call
the cyclic order of the attachments on its boundary the
\emph{order around} $T$.
\begin{enumerate}[label=\textup{(\roman*)}]
\item If $T$ consists of four edges at a vertex $h$, the order
around $T$ is that of $C$ or $C'$ for $B_a$, and of $acdba$ or~$acbda$~for~$B_b$.
\item If the order around $T$ is that of a four-cycle $C_0$ on $A$
for which $B\cup C_0$ is planar, then $G$ is toroidal.
\end{enumerate}
\end{lemma}
\begin{proof}
For (i), use the face walks in Lemma~\ref{dis:faces}. For $B_a$,
the four edges at $h$ lie in one face. In a simple face, the
attachments occur in the order $s_x,t_x,t_y,s_y$. In the repeated
face, they occur in the order $s_x,s_y,s_x,t_x,t_y,t_x$, after
interchanging $x,y$ if necessary. Whichever appearances the edges
use, $s_x,s_y$ are consecutive around $h$, as are $t_x,t_y$.
This gives the order of $C$ or $C'$.

For $B_b$, take $ahb$ as $P_m$. If $hc$ and $hd$ leave $h$ on
opposite sides of $ahb$, the order is $a,c,b,d$. Suppose they
leave on the same side. The boundary walk passes $a,h,b,d,c$ in
a simple face or $a,h,b,d,b,h,a,c$ in the repeated face. In either
case, following the boundary from the chosen appearance of $h$
through $a$ reaches $c$ before $d$. Thus $hc$ comes before $hd$
on that side from $ha$ to $hb$, giving the order $a,c,d,b$.

For (ii), since $B\cup C_0$ is planar,
Lemma~\ref{cor:boundary-completions} gives a disk drawing of $B$
with its attachments in the order of $C_0$. This is their order
around $T$, so the Jordan--Schoenflies
theorem lets us replace $T$
by this drawing, matching the four attachments as in the proof of
Lemma~\ref{lem:flat-replacement}. This gives a torus embedding of $G$.
\end{proof}

\begin{lemma}\label{dis:four-cycles}
The graphs $B_a\cup C$, $B_a\cup C'$, $B_b\cup acdba$ and $B_b\cup acbda$ are nonplanar.
\end{lemma}
\begin{proof}
Let $B$ be either bridge, with attachment set $A$. Contract its
interior to a single vertex. Since the interior has at least two
vertices, this gives a proper minor of $G$, so it has a torus
embedding. By Lemma~\ref{dis:four-completions}(i), the four edges
at the contracted vertex occur in the order of one of the two
listed cycles, say $C_0$. By part (ii), $B\cup C_0$ is nonplanar,
since $G$ is not toroidal. Suppose $B\cup C_1$ is planar for the
other listed cycle $C_1$.

In both bridge types, the two attachments in $A_1$ are consecutive
on each cycle. Read both cycles starting with these two attachments
in the same order. The two attachments in $A_2$ then occur in
opposite orders. By Lemma~\ref{cor:boundary-completions}, $B$ has
a disk drawing in the order of $C_1$ and none in the order of $C_0$. Thus by
Lemma~\ref{dis:four-disk-order} there are two disjoint paths
through $B$ joining the two pairs of opposite vertices of $C_0$.
Each path joins $A_1$ to $A_2$ and has an internal vertex, since
$B$ has no edges between its attachments.

Since $B-A$ is connected, extend the interiors of the two paths
to two disjoint connected sets containing all of $B-A$. Contract these sets
to $u,v$, keeping $uv$ and the two edges from each vertex to the
ends of its path. This gives an allowed minor $T$ of $B$. Suppose $T\ne B$. Replacing $B$ by $T$ would give a proper minor of
$G$. Hence taking a torus embedding of $G$ with $B$ replaced by $T$, in a thin disk around $T$,
the paths through $u$ and $v$ are disjoint, so their ends cannot
alternate on its boundary. This rules out the order of $C_0$.
Contracting $uv$ preserves the order around $T$, so
Lemma~\ref{dis:four-completions}(i) leaves only $C_1$. Part (ii)
then gives a torus embedding of $G$, a contradiction. Thus~$B=T$.

Therefore $u$ and $v$ have degree three in $G$. By
Lemma~\ref{lem:independent-neighborhood}(a), the two attachments
at each are nonadjacent. Use the two paths of length two through
$u,v$ and the unchanged $P_w$ as the three connecting paths,
relabeling the two $t$-ends if necessary. They satisfy
\eqref{dis:shortest}, so the remark after Lemma~\ref{dis:normalize}
allows us to apply that lemma to the new subdivision. But $uv$ is
a bridge joining the interiors of two of its paths, a contradiction.
\end{proof}

\begin{proposition}\label{dis:no-type-a}
If some $P_i$ is not an edge, then no bridge of type (a) has at least two interior vertices.
\end{proposition}
\begin{proof}
Suppose $B$ is such a bridge, with attachments
$A=\{s_x,t_x,s_y,t_y\}$, and put $H=G-(V(B)\setminus A)$.
A path through another $\Psi_1$-bridge from $s_x$ to $t_y$ or to
the interior of $P_y$ would be disjoint from a $t_x$--$s_y$ path
through $B$, contrary to Lemma~\ref{dis:opposite-order}. By symmetry,
\[
s_xt_y,s_yt_x\notin E(G),\qquad W_x,W_y\subseteq\{s_w,t_w\}.
\]
Every three vertices of $A$ contain one of the two pairs
$s_x,t_y$ and $s_y,t_x$. Thus no other bridge outside the paths
with nonempty interior attaches only in $A$, since it has at least
three attachments. There is no edge bridge on $A$ either, since
$A_1,A_2$ are independent and $s_it_i$ belongs to $P_i$ whenever
it is present.

Choose some $P_k$ with an internal vertex. By three-connectivity,
$W_k\ne\varnothing$, so some $\Psi_1$-bridge contains a path $Q$
from the interior of $P_k$ to an end of another path $P_\ell$,
with its internal vertices outside $\Psi_1$. By the restrictions
on $W_x,W_y$, we may interchange $x,y$ so that
$\{k,\ell\}=\{y,w\}$. In the choice $(P_y;s_x,t_x)$, no face
contains both ends of $Q$, by Lemma~\ref{dis:faces}. Thus this
choice is already impossible for $H$.

Retain $\Psi_1\cup B\cup Q$ and shorten $Q$ and the three
connecting paths to edges. On $P_k$, contract the end of $Q$ into
$t_k$ if its other end is $s_\ell$, and into $s_k$ if its other
end is $t_\ell$. This gives a proper minor with $B$ unchanged
and with the extra edge $s_yt_w$ or $s_wt_y$.
Take a torus embedding of this minor. If $P_w$ is repeated, $B$
lies in a simple face and its test contains $B\cup C$, which is
nonplanar by Lemma~\ref{dis:four-cycles}. If $P_y$ is repeated,
the attachments of $B$ force $(P_y;s_x,t_x)$, and no face contains
both ends of the new edge. Thus the choice must be $(P_x;s_y,t_y)$,
whose repeated-face test $B\cup(C-s_yt_y)$ is planar by
Lemma~\ref{dis:joint-test}.

We show that this is also not possible. If $P_x$ had an internal vertex, its bridge would attach to
$s_w$ or $t_w$, and neither lies in the repeated face. Thus $P_x$
would be an edge. The repeated face would then contain no bridge
of $H$, since its only possible attachments are in $A$. Adding
$B$ there gives precisely the planar repeated-face test above,
so Lemma~\ref{dis:joint-test} would embed $G$ in the torus.

Now replace $B$ by any $B_0$ obtained by one allowed operation.
This proper minor is toroidal and still contains $H$. By Lemma~\ref{dis:faces}, the four attachments of $B_0$ force
$(P_x;s_y,t_y)$ if $P_x$ is repeated and $(P_y;s_x,t_x)$ if
$P_y$ is repeated. Both choices are impossible for $H$, so $P_w$
is repeated. The interiors of $P_x,P_y$ then lie only on simple
faces. If either path had an internal vertex, a bridge of $H$
would join its interior to $s_w$ or $t_w$, but neither end of
$P_w$ lies on a simple face. Thus
$P_x,P_y$ are edges, $k=w$, and $Q$ joins the interior of $P_w$
to $s_y$ or $t_y$. The simple-face test gives $B_0\cup C$ planar
for every such $B_0$.

The paths $P_x,P_y,s_xps_y,t_xrt_y$ give $B\cup C$ as a minor
of $G$, so it has no $K_5$ minor. Every $K_{3,3}$ subdivision
in $B\cup C$ uses $s_yt_y$, since deleting that edge leaves the
planar completion above. Corollary~\ref{common:reduced-forms}
with $H_0=C$ and $s_yt_y$ required in the subdivision gives
$1{,}378$ candidate bridges. Here $B\cup K_4[A]$ is
three-connected, Lemma~\ref{dis:opposite-order} excludes disjoint
$s_x$--$t_y$ and $s_y$--$t_x$ paths through $B$, and every
degree-three vertex has independent neighbors by
Lemma~\ref{lem:independent-neighborhood}(a). The checks
in~\supp{E6} leave two labeled bridges. Up to interchanging $A_1,A_2$,
the interior is $\{p',q',v\}$, with
\[
N(p')=N(q')=\{s_x,s_y,v\},
\]
and its remaining edges are $vt_x,vt_y$.

Delete $vt_x$. The remaining bridge attaches at $s_x,s_y,t_y$.
If $P_x$ is repeated, these attachments force $(P_x;s_y,t_y)$,
already excluded for $H$. If $P_y$ is repeated, the only faces
meeting the interior of $P_w$ are the simple faces following
$P_x,P_w$, and neither contains the other end of $Q$. Thus $P_w$
must be repeated. Its simple-face test contains $(B-vt_x)\cup C$.
Contracting $t_xt_y$ to $t$ and deleting $s_xs_y$ in this graph
gives $K_{3,3}$ with parts $\{p',q',t\}$ and $\{s_x,s_y,v\}$.
Hence $G-vt_x$ is a proper nontoroidal minor of $G$, a contradiction.
\end{proof}

It remains to consider bridges of type (b). Let $B=B_b$, and put $e=s_w$ and $f=t_w$.

\begin{lemma}\label{dis:type-b-forms}
The path $P_x$ is the edge $cd$. Up to interchanging $A_1,A_2$, the
bridge $B$ is one of the five forms below, and one of the following holds.
\begin{enumerate}[label=\textup{(\roman*)}]
\item $P_w=ef$, and $G=(\Psi_1\cup B)+eb$ or $G=(\Psi_1\cup B)+af$.
\item Every other $(J_1\cup J_2)$-bridge attaches within
$\{c,d,e,f\}$, and $H_{\mathrm{simp}}$ for $(P_m;c,d)$ is nonplanar.
\end{enumerate}
In each row, $p',q'$ have the displayed neighborhood. The other interior
vertices are $u,v$ in $B_1,\ldots,B_4$ and $v$ in $B_5$.
\[
\scaledtable{1}{$\displaystyle
\begin{array}{c|c|l}
 & N(p')=N(q') & \text{remaining edges}\\ \hline
B_1 & \{a,c,v\} & au,\ bu,\ du,\ uv\\
B_2 & \{a,c,v\} & au,\ bu,\ dv,\ uv\\
B_3 & \{a,u,v\} & cu,\ bv,\ dv\\
B_4 & \{d,u,v\} & cu,\ av,\ bv\\
B_5 & \{a,b,v\} & cv,\ dv
\end{array}
$}
\]
\end{lemma}
\begin{proof}
Put $A=\{a,b,c,d\}$. The bridge $B$ contains a path from $d$ to the
interior of $P_m$ with internal vertices outside $\Psi_1$. A path
through another $\Psi_1$-bridge from $b$ to $c$ or to the interior
of $P_x$ would be disjoint from it. These paths would join $P_m,P_x$
in opposite order, contrary to Lemma~\ref{dis:opposite-order}.
Interchanging $A_1,A_2$ gives the same exclusion from $a$ to $d$ or
to the interior of $P_x$. Thus $ad,bc\notin E(G)$ and
$W_x\subseteq\{e,f\}$. No other bridge outside the paths attaches
only in $A$, since any three vertices of $A$ contain $a,d$ or $b,c$,
and the only possible edges on $A$ belong to the connecting~paths.

Let $D$ be $B$ or an allowed minor of $B$, and take an $a$--$b$ path
$P$ through its interior in place of $P_m$. Since $D-A$ is connected,
every $\Psi_1$-bridge in $D$ for this choice meets the interior of
$P$, and both $c,d$ have paths to its interior. If $P_x$ were
repeated, neither simple face following $P,P_w$ would contain $c$
or $d$. By Lemma~\ref{dis:faces}, either $P$ is repeated with choice
$(P;c,d)$, or $P_w$ is repeated. In the first case the repeated-face
test contains $D+ac+bd$. In the second, the simple-face test contains
$(D\cup P_x)+ac+bd$.

Let $L$ be obtained from $\Psi_1\cup B$ by shortening $P_x,P_w$ to
$cd,ef$. The proper minor $L/ef$ is toroidal. Every $\Psi_1$-bridge
in $L$ meets the interior of $P_m$ and attaches within
$V(P_m)\cup\{c,d\}$, so Lemma~\ref{dis:aligned-contraction}(c) gives
a torus embedding of $L$. Both choices above require $B+ac+bd$ to be
planar.

Now let $B_0$ be obtained from $B$ by one allowed operation. Replacing
$B$ by $B_0$ gives a proper minor of $G$, so take a torus embedding
of it. For $(P;c,d)$, no bridge outside $B_0$ can use the repeated
face, since its only possible attachments there are in $A$. Restoring
$B$ would change only the repeated-face test, to the planar graph
$B+ac+bd$. By Lemma~\ref{dis:joint-test} this would embed $G$ in the
torus. Thus $P_w$ is repeated. If $P_x$ had an internal vertex,
three-connectivity and $W_x\subseteq\{e,f\}$ would give a bridge
joining its interior to $e$ or $f$. Neither $e$ nor $f$ lies on a simple face, so $P_x=cd$.

The simple faces now follow $P$ and $cd$. No bridge outside $B_0$
can use them, so their test is $B_0\cup H_0$, where $H_0$ is the
path $acdb$. Thus $B_0\cup H_0$ is planar after every single allowed
operation. Restoring $B$ changes only this test, so $B\cup H_0$ is
nonplanar.

Suppose another $(J_1\cup J_2)$-bridge $K$ attaches at $a$ or $b$,
say at $a$ after interchanging $A_1,A_2$. If $K$ contains $P_w$,
contracting its interior into $f$
gives $af$ and $ef$. Otherwise $K$ lies outside the paths. If it
attaches at $f$, it gives an $a$--$f$ path. If it misses $f$, it
must attach at $e$, since its attachments cannot all lie in $A$.
It misses $d$ by the first paragraph, and it cannot also attach at
$c$, since no face contains $\{a,c,e\}=A_1$
(Lemma~\ref{dis:noncofacial}). Since $ae$ is absent,
three-connectivity gives a third attachment at $b$, so $K$ gives
an $e$--$b$ path.

Thus $L+af$ or $L+eb$ is a minor of $G$. For $(P_m;c,d)$ no face
contains both ends of the extra edge, and for $P_w$ repeated the
simple-face test contains $B\cup H_0$. Hence the minor is
nontoroidal, and minimality gives (i). If no other bridge attaches
at $a$ or $b$, all the others attach within $\{c,d,e,f\}$ and can
be assigned to the simple faces for $(P_m;c,d)$. The repeated-face
test $B+ac+bd$ is planar, so $H_{\mathrm{simp}}$ is nonplanar,
giving (ii).

The paths $apc,brd$ and the edge $cd$ give $B\cup H_0$ as a minor
of $G$, so it has no $K_5$ minor. Every $K_{3,3}$ subdivision in it
uses $cd$, since $B+ac+bd$ is planar. Corollary~\ref{common:reduced-forms}
with $cd$ required gives $2{,}104$ candidates. Here $B\cup K_4[A]$
is three-connected, $L$ has no $K_5$ minor, and no degree-three
interior vertex is adjacent to both ends of $cd$ by
Lemma~\ref{lem:independent-neighborhood}(a). With the planarity
conditions above, the checks in~\supp{E7} leave nine labeled
bridges, giving the five displayed forms up to interchanging $A_1,A_2$.
\end{proof}

\begin{proposition}\label{dis:type-b-reduction}
If $G$ has a bridge of type (b) with at least two interior vertices,
then $G\cong\mathfrak T_{66}$, or the twins and paths can be chosen
as at the start of this section so that one of the following holds.
\begin{enumerate}[label=\textup{(\roman*)}]
\item All three paths are edges.
\item There is a bridge of type (c) with at least two interior vertices.
\item Every bridge with nonempty interior has three attachments.
\end{enumerate}
\end{proposition}
\begin{proof}
In case (i) of Lemma~\ref{dis:type-b-forms}, the nine labeled bridges
and the two choices of added edge give $18$ graphs. Sixteen have a
proper nontoroidal minor and the other two are isomorphic to
$\mathfrak T_{66}$~\supp{E7}. Suppose (ii) holds. Replacing $ce,df$
in the nonplanar simple-face test by $cqe,dzf$ gives a connected
nonplanar subgraph $S$ of $G-aq$. In a torus embedding of this proper
minor, every face of $S$ is a disk by Lemma~\ref{lem:torus-complement}.

If $B=B_1$ or $B_2$, then $a,v$ and their common neighbors $p',q',u$
give a $K_{2,3}$ disjoint from $S$. It lies in one face of $S$, and
$p'c,q'c,ubrd$ join its three degree-two vertices to the boundary
without meeting the $K_{2,3}$ again. Thus all three lie on the outer
face of the $K_{2,3}$ in this disk. Adding a vertex there adjacent to
all three gives a plane drawing of $K_{3,3}$, a contradiction.
For $B_3$ and $B_5$, use the $K_{2,3}$ on $p',q'$ and their three
neighbors. The paths to $S$ are $apc,uc,vbrd$ for $B_3$ and
$apc,brd,vc$ for $B_5$, giving the same contradiction.

Thus $B=B_4$. Choose $p,q$ and $p',q'$ as the twins, with triples
$\{a,c,e\}$ and $\{v,u,d\}$, and choose connecting paths of minimum
total length. Outside the four twins, $a$ has only the neighbor $v$
and $u$ only the neighbor $c$. Thus two paths are $av,cu$ and the
third joins $e$ to $d$. If $ed$ is an edge, this gives (i).

Suppose $ed$ is absent. Outside the two new twin subgraphs, let $D$
be the component containing $b,r,z,f$. It attaches at $v$ through
$b$, at $d$ through $r$, and at $e$ along the old path $P_w$.
Neither $a$ nor $u$ has a neighbor in any component, and $b$ is
the only such neighbor of $v$. Every other component therefore
attaches within $\{c,d,e\}$ and has exactly three attachments by
three-connectivity. If $D$ misses $c$, this gives (iii). Otherwise
its attachments are $e,d,c,v$, both ends of the $e$--$d$ path and
one end of each of $cu,av$. By Proposition~\ref{dis:types} its
bridge is of type (c), with $b,r,z,f$ in its interior, giving (ii).
\end{proof}

\begin{proof}[Proof of Theorem~\ref{thm:disjoint}]
Both $\mathfrak T_4$ and $\mathfrak T_5$ are in the list, so suppose
$G$ is neither. If all three paths are edges,
Theorem~\ref{dis:endpoint-only} applies. Otherwise, if every bridge
with four attachments has one interior vertex,
Proposition~\ref{dis:remaining-classification}~applies.

It remains to consider a bridge with four attachments and at least
two interior vertices. By Propositions~\ref{dis:types}
and~\ref{dis:no-type-a}, it has type (b) or (c). Type (c) gives
$\mathfrak T_{21}$ or $\mathfrak T_{63}$ by Proposition~\ref{dis:type-c}.
For type (b), Proposition~\ref{dis:type-b-reduction} gives
$\mathfrak T_{66}$ or a new choice of twins and paths. For that
choice, either all three paths are edges, there is a bridge of type
(c) with at least two interior vertices, or every bridge with
nonempty interior has three attachments. These cases were just
considered. The converse follows from Proposition~\ref{prop:membership}.
\end{proof}

\section{Two twin triples with two common vertices}\label{sec:two-common}

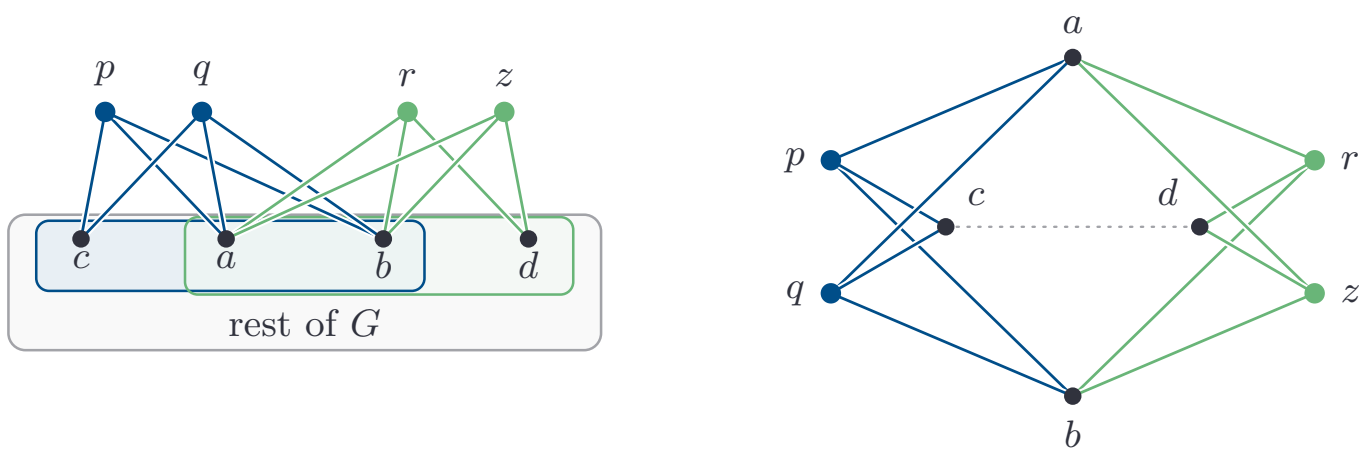
\begin{figure}[htbp]
\centering
\begingroup
\scaledfigure{1}{%
\tikzsetnextfilename{two-common-neighborhoods}%
\begin{tikzpicture}[
  font=\small,text=figureInk,
  vertex/.style={circle,fill=figureInk,inner sep=1.45pt},
  twinA/.style={circle,fill=figureTwinA,inner sep=1.7pt},
  twinB/.style={circle,fill=figureTwinB,inner sep=1.7pt},
  edgeA/.style={draw=figureTwinA,line width=.65pt,
    preaction={draw=white,line width=1.6pt}},
  edgeB/.style={draw=figureTwinB,line width=.65pt,
    preaction={draw=white,line width=1.6pt}},
  lab/.style={anchor=north,inner sep=0pt}
]
\begin{scope}[yshift=-.1cm]
  \draw[draw=figureInk!45,fill=figureInk!3,line width=.6pt,
    rounded corners=4pt]
    (.05,-.92) rectangle (4.95,.2);
  \fill[figureTwinA!9,rounded corners=3pt]
    (.28,-.43) rectangle (3.49,.15);
  \fill[figureTwinB!9,fill opacity=.6,rounded corners=3pt]
    (1.51,-.46) rectangle (4.72,.18);
  \draw[figureTwinA,line width=.6pt,rounded corners=3pt]
    (.28,-.43) rectangle (3.49,.15);
  \draw[figureTwinB,line width=.6pt,rounded corners=3pt]
    (1.51,-.46) rectangle (4.72,.18);
  \node at (2.5,-.69) {rest of $G$};

  \coordinate (c) at (.65,0);
  \coordinate (a) at (1.85,0);
  \coordinate (b) at (3.15,0);
  \coordinate (d) at (4.35,0);
  \coordinate (p) at (.85,1.05);
  \coordinate (q) at (1.65,1.05);
  \coordinate (r) at (3.35,1.05);
  \coordinate (z) at (4.15,1.05);

  \foreach \v in {p,q}
    \foreach \u in {a,b,c} \draw[edgeA] (\v)--(\u);
  \foreach \v in {r,z}
    \foreach \u in {a,b,d} \draw[edgeB] (\v)--(\u);
  \foreach \v in {p,q}
    \node[twinA,label=above:{$\v$}] at (\v) {};
  \foreach \v in {r,z}
    \node[twinB,label=above:{$\v$}] at (\v) {};
  \foreach \v in {c,a,b,d}{
    \node[vertex] at (\v) {};
    \node[lab] at ([yshift=-.1cm]\v) {$\v$};
  }
\end{scope}

\begin{scope}[xshift=8.85cm]
  \coordinate (a) at (0,1.4);
  \coordinate (b) at (0,-1.4);
  \coordinate (c) at (-1.05,0);
  \coordinate (d) at (1.05,0);
  \coordinate (p) at (-2,.55);
  \coordinate (q) at (-2,-.55);
  \coordinate (r) at (2,.55);
  \coordinate (z) at (2,-.55);

  \draw[draw=figureInk!45,line width=.55pt,dotted,line cap=round]
    (c)--(d);
  \foreach \v in {p,q}
    \foreach \u in {a,b,c} \draw[edgeA] (\v)--(\u);
  \foreach \v in {r,z}
    \foreach \u in {a,b,d} \draw[edgeB] (\v)--(\u);

  \node[vertex,label=above:{$a$}] at (a) {};
  \node[vertex,label=below:{$b$}] at (b) {};
  \node[vertex,label=above right:{$c$}] at (c) {};
  \node[vertex,label=above left:{$d$}] at (d) {};
  \node[twinA,label=left:{$p$}] at (p) {};
  \node[twinA,label=left:{$q$}] at (q) {};
  \node[twinB,label=right:{$r$}] at (r) {};
  \node[twinB,label=right:{$z$}] at (z) {};
\end{scope}
\end{tikzpicture}%
}
\endgroup
\caption{Left, two twin triples meeting at $a,b$.
Right, the graph $\Psi_2$, with the possible edge $cd$ dotted.}
\label{fig:two-common-neighborhoods}
\end{figure}

In this section, we consider the case where an obstruction has no two disjoint twin triples but has some that share two vertices. Let the two triples be $A_1=\{a,b,c\}$ and $A_2=\{a,b,d\}$,
where $a,b,c,d$ are distinct, and call them the $c$ triple and the $d$ triple.
Choose twins $p,q$ with neighborhood $A_1$ and $r,z$ with neighborhood
$A_2$, and write $J_1,J_2$ for their twin subgraphs. Put
$\Psi_2=J_1\cup J_2$. Figure~\ref{fig:two-common-neighborhoods} shows the twin triples
and the graph $\Psi_2$. We determine $G$ from the bridges of $\Psi_2$.

Only $cd$ can join $a,b,c,d$, by
Lemma~\ref{lem:independent-neighborhood}(a). Moreover, every twin has degree
three, so each component of $G-V(\Psi_2)$ has its neighbors among
$a,b,c,d$. It has at least three of them, as $G$ is three-connected.
The $\Psi_2$-bridges are these components with their attaching edges,
and the edge $cd$ when present. Since $\Psi_2=J_1\cup J_2$, these
are also exactly the bridges on $A_1\cup A_2$ in $G-\{p,q,r,z\}$.
Here the two bridge descriptions used in Section~\ref{sec:disjoint}
coincide.

\begin{theorem}\label{thm:two-common}
A graph is a three-connected obstruction in which some two twin triples share two vertices, and no two twin triples are disjoint, if and only if it is isomorphic to one of
\[
\begin{gathered}
\mathfrak{T}_{2},\mathfrak{T}_{9},\mathfrak{T}_{10},\mathfrak{T}_{11},\mathfrak{T}_{15},\mathfrak{T}_{16},\mathfrak{T}_{17},\mathfrak{T}_{18},\mathfrak{T}_{19},
\mathfrak{T}_{22},\mathfrak{T}_{24},\mathfrak{T}_{26},\\
\mathfrak{T}_{30},\mathfrak{T}_{34},\mathfrak{T}_{39},\mathfrak{T}_{40},\mathfrak{T}_{51},\mathfrak{T}_{52},\mathfrak{T}_{53},\mathfrak{T}_{56},\mathfrak{T}_{57},
\mathfrak{T}_{72},\mathfrak{T}_{73},\mathfrak{T}_{74}.
\end{gathered}
\]
\end{theorem}

We split the proof into cases depending on whether some choice of the two twin
pairs gives a $\Psi_2$-bridge attached to all four vertices $a,b,c,d$,
or no choice does. We first give the embedding tests and bridge reductions used in
both cases. We then treat the case of a bridge with four attachments
in Sections~\ref{r2:four-neighbors}--\ref{r2:planar-four} and the
remaining case in Sections~\ref{r2:no-four}--\ref{r2:planar-bridges}.

\begin{lemma}\label{r2:unique-four}
Suppose $B$ is a $\Psi_2$-bridge with attachments $a,b,c,d$.
Then $B$ is the only such bridge, $cd$ is absent, and every other
$\Psi_2$-bridge has a single interior vertex with neighborhood
$\{a,b,c\}$ or $\{a,b,d\}$. If at least four vertices have the same neighborhood
$\{a,b,c\}$ or $\{a,b,d\}$, then $G\cong\mathfrak{T}_{2}$.
\end{lemma}

\begin{proof}
Suppose another $\Psi_2$-bridge meets both $c,d$, or $cd$
is an edge. Then there is a $c$--$d$ path $P$ whose
internal vertices avoid $\Psi_2$ and the interior of $B$.
Contract the interior of $B$ to a vertex and shorten $P$ to $cd$.
The resulting graph has a $K_5$ minor~\supp{W5}, a contradiction.
Thus $B$ is the only bridge with all four attachments,
$cd$ is absent, and every other component of
$G-V(\Psi_2)$ has neighborhood $\{a,b,c\}$ or $\{a,b,d\}$
by three-connectivity.

Suppose such a component $D$ with neighborhood $\{a,b,c\}$
has nonplanar triangle completion. Its completion is three-connected by
Lemma~\ref{lem:triangle-completion}. The paths $arb,apc,bqc$
join the three attachment pairs with disjoint interiors
outside the bridge containing $D$, so Lemma~\ref{endpoint:primitive} gives
$T(t)$ as an allowed minor of this bridge for some
$t\in\{a,b,c\}$. Contracting the unique edge of $T(t)$
incident with $t$ gives two further vertices with
neighborhood $\{a,b,c\}$. Contracting the interior of $B$
to a vertex then gives a proper $\mathfrak T_2$
minor~\supp{W6}, contrary to minimality.

Thus every such triangle completion is planar. The same argument applies to $\{a,b,d\}$ after interchanging $c,d$ and the two twin pairs.
Lemma~\ref{lem:flat-replacement} makes each component a
single vertex. If at least four vertices have either
neighborhood, contracting the interior of $B$ and retaining
four vertices with that neighborhood and two with the
other gives a $\mathfrak T_2$ minor~\supp{W7}.
Minimality therefore gives $G\cong\mathfrak T_2$.
\end{proof}

If $B$ exists and $G\not\cong\mathfrak T_2$, then
Lemma~\ref{r2:unique-four} shows that $G$ consists of $B$,
$n_c$ vertices whose neighborhood is the $c$ triple and $n_d$ whose
neighborhood is the $d$ triple. Each triple is the neighborhood of its two chosen twins and,
by the lemma, of at most one other vertex. Thus, after
interchanging $c,d$ if needed, $(n_c,n_d)$ is $(2,2)$, $(3,2)$
or $(3,3)$.

\subsection{Faces and embedding tests of \texorpdfstring{$\Psi_2$}{Psi2}}\label{r2:tests}

The next lemma gives the embedding tests for these three cases.

\begin{lemma}\label{r2:four-tests}
Let $B$ have four independent attachments $a,b,c,d$. Suppose
$B-\{a,b,c,d\}$ is nonempty and connected, and each attachment
has a neighbor in it. For $(n_c,n_d)\in\{(2,2),(3,2),(3,3)\}$,
form $G$ from $B$ by adding $n_c$ vertices with neighborhood
$\{a,b,c\}$ and $n_d$ vertices with neighborhood $\{a,b,d\}$.
Then $G$ is toroidal if and only if at least one of the graphs
in the corresponding row is planar.
\[
\scaledtable{.9}{$\displaystyle
\begin{array}{c|l}
(n_c,n_d)&\text{planarity tests}\\\hline
(3,3)&B+ac+bc+bd+ad\\
(3,2)&B+ab+ac+bc+ad,\quad B+ab+ac+bc+bd\\
(2,2)&B+ac+bd,\quad B+ad+bc,\\
&B+ab+ac+ad,\quad B+ab+bc+bd.
\end{array}
$}
\]
\end{lemma}

\begin{proof}
Let $K$ consist of the added vertices, their incident edges and
the four attachments. Since the interior of $B$ is connected,
it must lie in one face of $K$ in any torus embedding of $G$.
We therefore need only consider faces whose boundary meets all
four attachments. For each rotation of $K$, \supp{R4} gives the
attachments in the order they occur on each boundary walk. We give
these sequences up to cyclic shift and reversal.
For the disk faces that occur, Lemma~\ref{cor:boundary-completions}
gives the planarity test by joining consecutive distinct
attachments on the boundary walk.

When $n_c=3$, the three vertices with neighborhood $\{a,b,c\}$
give a $K_{3,3}$ in $K$, so every face in a torus embedding of
$K$ is a disk by Youngs's theorem. For $(3,3)$ the only boundary
sequence meeting all four attachments is $acbd$, giving the
first row. For $(3,2)$ the sequences are $abcad$, $abdbc$ and
$acbd$. The first two sequences give the graphs in the $(3,2)$ row. If
$B\cup acbda$ is planar, its connected subgraph $B-\{a,b,c,d\}$
lies on one side of the cycle in a plane drawing, so we can draw
$ab$ on the other side. Deleting $bd$ gives a plane drawing of
the first graph in the row, and deleting $ad$ gives one of the second.

For $(2,2)$ we have $K=\Psi_2$. In a genus-zero rotation the
boundary walks, restricted to the attachments, are
$ab,ab,ac,bc,ad,bd$. No single walk meets all four attachments,
so by Lemma~\ref{lem:torus-complement} the face containing $B$
must be an annulus. The only pairs meeting all four attachments are $ac\mid bd$ and
$ad\mid bc$. For either pair, adding a handle between the two
corresponding faces of the sphere embedding of $K$ makes them
a single annular face. Lemma~\ref{lem:annulus-embedding} gives the
tests $B+ac+bd$ and $B+ad+bc$.

In genus one all faces are disks. Of the seven boundary sequences
meeting all four attachments, five have either $ac,bd$ or $ad,bc$
as consecutive pairs. Thus if $B$ can be drawn in any of these five faces, then
$B+ac+bd$ or $B+ad+bc$ is planar. The remaining sequences
$abacad$ and $abcbdb$ give $B+ab+ac+ad$ and $B+ab+bc+bd$,
respectively.

Conversely, if one of the graphs in the $(n_c,n_d)$ row is
planar, the disk or annulus test lets us draw $B$ in the
corresponding face of a torus embedding of $K$. This gives
a torus embedding of $G$.
\end{proof}

We now show that the
completions used in these tests are minors of $G$,
so they also have no $K_5$ minor.

\begin{lemma}\label{r2:test-minors}
Let $B$ be a $\Psi_2$-bridge of $G$ with attachments $a,b,c,d$,
and let $H=B\cup acbda$ be the completion in the $(3,3)$ row of
Lemma~\ref{r2:four-tests}. Both $H$ and $H+ab$ are three-connected
$K_5$-minor-free minors of $G$.
\end{lemma}
\begin{proof}
In $B\cup\Psi_2$, contract one twin from each pair into $a$ and
the other into $b$. This gives $H+ab$, which contains every
completion in the table. Deleting $ab$ gives $H$, so both graphs
are $K_5$-minor-free minors of $G$.

Delete at most two vertices from $H$. Every remaining component
contains an attachment, since otherwise the same deletion would
disconnect $G$. If at most one attachment was deleted, the cycle
$acbda$ still joins the others. If two were deleted, the interior
of $B$ is still connected and joins the remaining two. Thus $H$
is three-connected, and so is $H+ab$.
\end{proof}

\begin{figure}[htbp]
\centering
\scaledfigure{1.4}{%
\tikzsetnextfilename{psi2-annulus-and-restoration-r2}%
\begin{tikzpicture}[
  scale=.5, font=\scriptsize,
  every label/.append style={inner sep=.9pt},
  label distance=.5pt,
  vertex/.style={circle,fill=figureInk,inner sep=1.6pt},
  edge/.style={draw=figureInk!80,line width=.8pt},
  boundary/.style={draw=figurePath,line width=1pt},
  seam/.style={draw=figurePath,line width=.7pt,dashed},
  newedge/.style={draw=figureBridge,line width=1.2pt},
  contractedge/.style={draw=figureBridge,line width=1.2pt,dashed},
  arrow/.style={draw=figureInk!70,line width=.75pt,
    -{Stealth[length=5pt]}},
  panel/.style={anchor=north,align=center,inner sep=0pt}
]

\begin{scope}
  \coordinate (a) at (0,2);
  \coordinate (b) at (0,-1.8);
  \coordinate (p) at (-2.8,0);
  \coordinate (c) at (-1.8,0);
  \coordinate (q) at (-.8,0);
  \coordinate (r) at (.8,0);
  \coordinate (d) at (1.8,0);
  \coordinate (z) at (2.8,0);

  \fill[faceRepeated] (a)--(p)--(c)--(q)--cycle;
  \fill[faceRepeated] (a)--(r)--(d)--(z)--cycle;
  \foreach \v in {p,q,r,z}
    \draw[edge] (a)--(\v)--(b);
  \draw[edge] (p)--(c)--(q) (r)--(d)--(z);

  \foreach \x/\mark in {-1.02/I,1.02/II}{
    \filldraw[fill=white,seam] (\x,.75) circle (.36);
    \node at (\x,.75) {\mark};
  }
  \node[vertex,label=above:{$a$}] at (a) {};
  \node[vertex,label=below:{$b$}] at (b) {};
  \node[vertex,label=left:{$p$}] at (p) {};
  \node[vertex,label=below:{$c$}] at (c) {};
  \node[vertex,label=above right:{$q$}] at (q) {};
  \node[vertex,label=above left:{$r$}] at (r) {};
  \node[vertex,label=below:{$d$}] at (d) {};
  \node[vertex,label=right:{$z$}] at (z) {};
  \node[panel] at (0,-2.72)
    {(a) Remove\\the two disks};
\end{scope}

\draw[arrow] (3.2,.35)--(3.85,.35);

\begin{scope}[xshift=5cm,yshift=.35cm]
  \fill[faceSimple] (-.72,-.72) rectangle (.72,.72);
  \fill[faceSimple] (.72,0) ellipse (.29 and .72);
  \draw[figureInk!65,line width=.7pt]
    (-.72,.72)--(.72,.72)
    (-.72,-.72)--(.72,-.72);
  \filldraw[fill=white,seam] (-.72,0) ellipse (.29 and .72);
  \draw[seam] (.72,0) ellipse (.29 and .72);
  \node[above=1pt,inner sep=.6pt] at (-.72,.72) {I};
  \node[above=1pt,inner sep=.6pt] at (.72,.72) {II};
  \node[panel] at (0,-3.07)
    {(b) Join them\\by a cylinder};
\end{scope}

\draw[arrow] (6.1,.35)--(7,.35);

\begin{scope}[xshift=9cm,yshift=.2cm]
  \fill[faceSimple] (0,0) circle (1.6);
  \fill[white] (0,0) circle (.65);
  \draw[seam] (0,.07) circle (.87) (0,.07) circle (1.07);
  \draw[boundary] (0,0) circle (.65) (0,0) circle (1.6);

  \node[vertex,label=above left:{$a_1$}] at (140:1.6) {};
  \node[vertex,label=below left:{$d$}] at (-140:1.6) {};
  \node[vertex] at (-.65,0) {};
  \node[vertex] at (.65,0) {};
  \node[inner sep=0pt] at (-.15,.23) {$a_2$};
  \node[inner sep=0pt] at (.27,-.21) {$c$};
  \node[inner sep=0pt] at (0,-1.30) {$B_1$};
  \node[panel] at (0,-2.92)
    {(c) Insert $B_1$\\$a_1,a_2\mapsto a$};
\end{scope}

\draw[figureInk!25,line width=.4pt]
  (11.45,2.25)--(11.45,-2.4);

\begin{scope}[xshift=14cm]
\begin{scope}
  \fill[faceRepeated] (0,0) circle (1.35);
  \draw[boundary] (0,0) circle (1.35);
  \coordinate (aone) at (140:1.35);
  \coordinate (cone) at (60:1.35);
  \coordinate (atwo) at (0:1.35);
  \coordinate (ctwo) at (-60:1.35);
  \coordinate (dd) at (-140:1.35);

  \draw[contractedge] (cone)
    .. controls (2.65,1.65) and (2.65,-1.65) ..
    node[pos=.30,right=1pt,inner sep=1pt,text=figureBridge] {$e$} (ctwo);
  \draw[newedge] (cone)
    .. controls (1.65,1.4) and (1.8,.55) .. (atwo);
  \draw[newedge] (aone)
    .. controls (-2.05,1.3) and (-2.05,-1.3) .. (dd);

  \node at (-.05,0) {$B_1$};
  \node[vertex,label=above left:{$a_1$}] at (aone) {};
  \node[vertex,label=above left:{$c_1$}] at (cone) {};
  \node[vertex,label=below left:{$a_2$}] at (atwo) {};
  \node[vertex,label=below:{$c_2$}] at (ctwo) {};
  \node[vertex,label=below left:{$d$}] at (dd) {};
  \node[panel] at (.2,-2.72)
    {(d) Add\\$e,c_1a_2,a_1d$};
\end{scope}

\draw[arrow] (2.85,0)--(3.8,0);

\begin{scope}[xshift=6cm]
  \fill[faceSimple,even odd rule]
    (0,0) circle (1.35)
    (.75,0) circle (.6);
  \draw[boundary]
    (0,0) circle (1.35)
    (.75,0) circle (.6);
  \coordinate (aone) at (140:1.35);
  \coordinate (dd) at (-140:1.35);
  \coordinate (atwo) at (.15,0);
  \coordinate (cc) at (1.35,0);

  \draw[newedge] (atwo)--(cc);
  \draw[newedge] (aone)
    .. controls (-2.05,1.3) and (-2.05,-1.3) .. (dd);
  \node at (-.62,-.12) {$B_1$};
  \node[vertex,label=above left:{$a_1$}] at (aone) {};
  \node[vertex,label=below left:{$d$}] at (dd) {};
  \node[vertex] at (atwo) {};
  \node[inner sep=0pt] at (.72,.26) {$a_2$};
  \node[vertex,label=above right:{$c$}] at (cc) {};
  \node[panel] at (-.15,-2.72)
    {(e) Contract $e$\\and recover $c$};
\end{scope}

\draw[arrow] (7.95,0)--(9.45,0)
  node[midway,below=2pt,align=center]
    {Lemma\\\ref{lem:annulus-embedding}};

\begin{scope}[xshift=11.2cm]
  \fill[faceSimple,even odd rule]
    (0,0) circle (1.35)
    (0,0) circle (.6);
  \draw[boundary] (0,0) circle (1.35) (0,0) circle (.6);
  \node[vertex,label=above left:{$a_1$}] at (140:1.35) {};
  \node[vertex,label=below left:{$d$}] at (-140:1.35) {};
  \node[vertex] at (-.6,0) {};
  \node[vertex] at (.6,0) {};
  \node[inner sep=0pt] at (-.13,.23) {$a_2$};
  \node[inner sep=0pt] at (.25,-.20) {$c$};
  \node at (0,-.98) {$B_1$};
  \node[panel] at (0,-2.72)
    {(f) $a_2c\mid a_1d$};
\end{scope}
\end{scope}
\end{tikzpicture}%
}
\caption{The proof of Lemma~\ref{r2:degree-four}. A torus embedding of $G$ in (a)--(c) and an annular drawing of $B_1$ in (d)--(f).}
\label{fig:psi2-annulus-construction}
\label{fig:psi2-twin-restoration}
\end{figure}
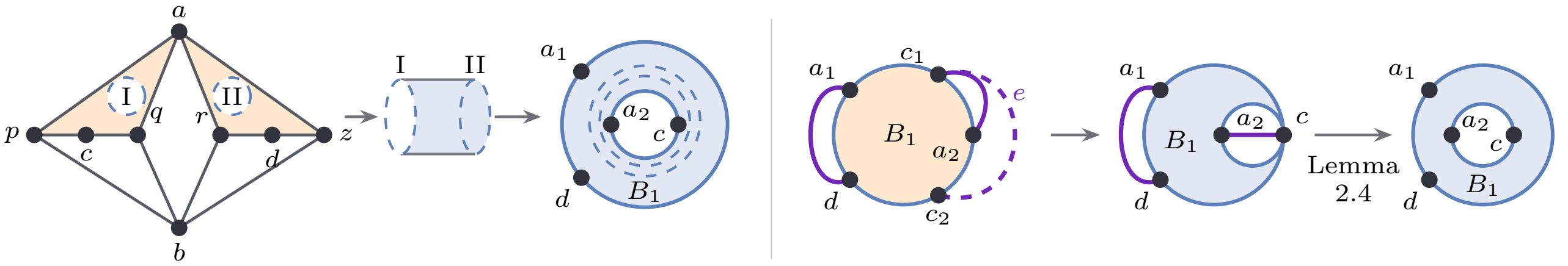

\begin{lemma}\label{r2:degree-four}
No three-connected obstruction has twin pairs with neighborhoods $\{a,b,c\},\{a,b,d\}$ and $\deg(b)=4$.
\end{lemma}
\begin{proof}
Suppose $\deg(b)=4$. Every component outside $\Psi_2$
then has neighborhood $\{a,c,d\}$ by three-connectivity.
Let $B_1$ consist of their bridges and $cd$, if present,
so $G=\Psi_2\cup B_1$. By minimality, $G-z$ has a torus
embedding. All bridges in $B_1$ lie in the same face $F$
of $H=\Psi_2-z$, since each meets $d$ and $\deg_H(d)=1$. Record only $a,c,d$ on the boundary walks of $F$.
If consecutive appearances are equal, move their incident
edges to one appearance along the intervening boundary arc,
preserving their order. Since $B_1$ meets both $c,d$,
the rotation table~\supp{R5} and
Lemma~\ref{lem:torus-complement} leave only the disk sequences
$acd,acad,acacd$ and the annular pairs $ad\mid ac,ad\mid c$. We show that $B_1$ can be drawn in an annulus with $a_1,d$
on one boundary and $a_2,c$ on the other, allowing the edges
at $a$ to use $a_1$ or $a_2$. The annular cases already give
this, with $a_2$ unused for $ad\mid c$. In a disk case,
separate the appearances and use the boundary order
$a_1,c_1,a_2,c_2,d$, leaving positions unused where needed.
Draw $c_1c_2,c_1a_2,a_1d$ outside the disk without crossings,
then contract $c_1c_2$ to recover $c$. The resulting graph
is planar, with $a_1,a_2$ distinct and $a_1d,a_2c$ added.
Lemma~\ref{lem:annulus-embedding} gives the required drawing (Figure~\ref{fig:psi2-annulus-construction}(d)--(f)). A planar embedding of $\Psi_2$ has faces $apcqa$ and $ardza$.
Remove an open disk from each and join the new boundary
circles by a cylinder. This gives a torus embedding with
an annular face $ac\mid ad$. Insert the drawing of $B_1$,
with the edges at $a_2$ reaching $a$ between $ap,aq$ and
those at $a_1$ reaching $a$ between $ar,az$.
This embeds $G$ on the torus (Figure~\ref{fig:psi2-annulus-construction}(a)--(c)), a contradiction.
\end{proof}

Moreover, we can exclude the case of $(n_c,n_d)=(3,3)$ before reducing the bridges. 

\begin{lemma}\label{r2:33-empty}
If $G$ has a $\Psi_2$-bridge $B$ with attachments $a,b,c,d$, then
$(n_c,n_d)\ne(3,3)$.
\end{lemma}
\begin{proof}
Suppose $(n_c,n_d)=(3,3)$,
and put $C=acbda$ and $H=B\cup C$.
By Lemmas~\ref{r2:test-minors} and~\ref{r2:four-tests},
$H$ is three-connected and nonplanar.
The graph $\operatorname{Aux}(H,C)$ is also
three-connected, since the added vertex has four
neighbors. Thus outcome (c) of
Theorem~\ref{thm:mohar} is impossible, so
$H-E(C)=B$ contains crossing paths or a tripod.

Suppose first that $B$ contains crossing paths
$P$ from $a$ to $b$ and $Q$ from $c$ to $d$.
Since neither $ab$ nor $cd$ is an edge of $B$, both $P$ and $Q$ have internal vertices. Let $R$ be a shortest path in $B-\{a,b,c,d\}$ from an internal vertex of $P$ to an internal vertex of $Q$. Let $W_1=V(Q)$ and let
$W_2$ consist of the internal vertices of $P$
and the vertices of $R$ outside $W_1$.
Then $W_1,W_2$ are disjoint and connected.
Each of $\{a\},\{b\},W_1$ is adjacent to $W_2$
and to each of the six vertices with neighborhoods
$\{a,b,c\},\{a,b,d\}$, giving a $K_{3,7}$ minor.
Since both path interiors are nonempty,
$|V(G)|\ge12>10$, so this minor is proper,
contradicting minimality.

Now suppose $B$ contains a tripod on
$A\subseteq\{a,b,c,d\}$.
Lemma~\ref{lem:tripod-contraction} gives a minor
$H_A$ consisting of $a,b,c,d$, the six vertices
with neighborhoods $\{a,b,c\},\{a,b,d\}$, and two
vertices with neighborhood $A$. This is
$\mathfrak{T}_{13}$ if $a,b\in A$ and
$\mathfrak{T}_{9}$ otherwise~\supp{W8}.
If $|\operatorname{Int}(B)|>2$, then $H_A$ has fewer
vertices than $G$. Otherwise the connected interior
consists of two adjacent vertices, which are nonadjacent
in $H_A$. Thus $H_A$ is a proper minor in both cases,
contradicting minimality.
\end{proof}

\subsection{Reducing bridges while preserving the tests}\label{sec:three-root-lists}

Proposition~\ref{prop:three-root-eight} lists the bridges that
cannot be reduced without changing which of their eight
completions are planar, provided one completion by two edges
is planar.

\begin{proposition}\label{prop:three-root-eight}
Let $B$ be a bridge of $G$ on $A=\{x,y,t\}$ whose triangle completion
is three-connected and a minor of $G$, and suppose $B+xy+xt$ is planar.
Assume that every proper allowed minor $B'$ has $B'+E$ planar and
$B+E$ nonplanar for some $E\subseteq E(K_3[A])$.
Then $B$ consists of one interior vertex adjacent to $x,y,t$, or is one of
$T(y),T(t),Q(x;y,t),S(x;y),S(x;t)$.
\end{proposition}
\begin{proof}
By Proposition~\ref{prop:three-root-blocks}, $B-x$
is a chain of edges and four-cycles. If there is
no four-cycle, $B$ has one interior vertex.
With exactly one four-cycle, the paths from $y,t$
to the cycle have lengths at most two. Lengths
$(0,0)$ disconnect the interior, while
$(0,1),(1,0),(1,1)$ give $T(t),T(y),Q(x;y,t)$. Otherwise, there is one four-cycle with an
end path of length two, or at least two four-cycles.
After interchanging $y,t$ if necessary, contract
the end path at $y$ into $y$, leaving the first
cycle as $ypuqy$. If this is the only cycle,
let $f$ be the middle vertex of the length-two
end path at $t$. Otherwise, let $f$ be one of
the two vertices of the next cycle other than
its connecting vertices. Proposition~\ref{prop:three-root-blocks}(ii),(iv)
gives $xf\in E(B)$. Keep a $u$--$t$ path through
$f$, the cycle $ypuqy$, and the edges $xp,xq,xf$.
Delete the other interior vertices, starting with those furthest from the interior vertices of this path and cycle. Keep only this path, the cycle and the edges $xp,xq,xf$. Contract the internal
vertices between $u,f$ into $u$ and those between
$f,t$ into $t$. This keeps the interior connected
and a neighbor at each attachment, giving the
allowed minor $S(x;y)$. By Lemma~\ref{core:three-graph-completions}, every
completion of this minor containing $yt$ is
nonplanar, so the corresponding completion of
$B$ is also nonplanar. Every completion of $B$
not containing $yt$ is a subgraph of the planar
graph $B+xy+xt$. Thus the minor has exactly the
same planar completions as $B$ and cannot be
proper by the hypothesis. Hence $B=S(x;y)$,
or $B=S(x;t)$.
\end{proof}

\subsection{A bridge with four attachments: nonplanar \texorpdfstring{$B$}{B} or \texorpdfstring{$(n_c,n_d)=(3,2)$}{(nc,nd)=(3,2)}}\label{r2:four-neighbors}

We now classify the graphs for which some choice of the twin pairs
gives a $\Psi_2$-bridge $B$ attached to all four vertices $a,b,c,d$.
Apart from $\mathfrak T_2$,
Lemmas~\ref{r2:unique-four} and~\ref{r2:33-empty} show that
$(n_c,n_d)=(2,2)$ or $(3,2)$.
\begin{lemma}\label{r2:32}
If $G$ has a $\Psi_2$-bridge $B$ with attachments $a,b,c,d$ and
$(n_c,n_d)=(3,2)$, then adding the triangle $abc$ to $B$ gives
a nonplanar graph.
\end{lemma}
\begin{proof}
Suppose $H=B+ab+ac+bc$ is planar. The two tests in the $(3,2)$
row of Lemma~\ref{r2:four-tests} are $H+ad$ and $H+bd$, so both
are nonplanar. Contract $ab$ to a vertex $x$ and put
$Q=(B+ab)/ab$. Now $ad$ and $bd$ both become $xd$, so the two
tests give the same graph $Q+xc+xd$. By
Lemma~\ref{r2:planar-splitting}, this graph is nonplanar. It is
a $K_5$-minor-free minor of $G$ by Lemma~\ref{r2:test-minors},
so it contains a $K_{3,3}$ subdivision.

Apply Lemma~\ref{core:retained-subdivision} to $Q$ with
attachments $x,c,d$ and fixed edges $xc,xd$. Since $xc$ and $xd$ are fixed edges, every segment of the
shortened $K_{3,3}$ subdivision is a single edge except for a
possible $c$--$d$ segment of length two. At most one branch path
contains both $c$ and $d$. Thus at most seven vertices of the
reduced bridge lie in its interior, at most six corners and one
other vertex on the $c$--$d$ segment.

We make the same contractions in the interior of $B$, leaving
$a,b$ distinct. We keep the edges of the reduced bridge not
incident with $x$. For each edge $xv$ of the reduced bridge, we
keep one of $av,bv$ present after the contractions. Both $a$ and
$b$ still have interior neighbors after the contractions, so if
the chosen edges miss one of them, we keep one more edge from
it to the interior. This gives an allowed minor $B'$ of $B$.
Contracting $ab$ in either test for $B'$ gives a graph containing
the shortened subdivision, so both tests are nonplanar. By Lemma~\ref{r2:four-tests},
replacing $B$ by $B'$ gives a nontoroidal minor of $G$. Hence
$B'=B$ by minimality.

Construction~\ref{common:marked-construction} lists the reduced
bridges on $x,c,d$. Splitting $x$ as above gives $2{,}953$ candidates
for $B$ up to isomorphisms fixing the four attachments. We add three
vertices with neighborhood $\{a,b,c\}$ and two with neighborhood
$\{a,b,d\}$ to each. Of the resulting graphs, $2{,}159$ have a $K_5$
minor, $92$ have a cut of at most two vertices, and $702$ have a proper
nontoroidal minor~\supp{E9}. None can be $G$.
\end{proof}

\begin{proposition}[Nonplanar completions with four attachments]\label{r2:nonplanar-list}
If $B$ is nonplanar and $G=B\cup\Psi_2$, then $G$ is isomorphic to one of
$\mathfrak{T}_{17},\mathfrak{T}_{18},\mathfrak{T}_{19},\mathfrak{T}_{22},\mathfrak{T}_{24},\mathfrak{T}_{26},\mathfrak{T}_{30},\mathfrak{T}_{34}, \mathfrak{T}_{39},\mathfrak{T}_{48},\mathfrak{T}_{52},\mathfrak{T}_{57},\mathfrak{T}_{66},\mathfrak{T}_{74}.$
If $(n_c,n_d)=(3,2)$, $G\cong\mathfrak{T}_{11}$.
\end{proposition}

Of the graphs in Theorem~\ref{thm:classification}, only
$\mathfrak T_{74}$ is absent from the collection of Myrvold and
Woodcock~\cite{MyrvoldWoodcock2018}, by the comparison in~\supp{O4}.

\begin{proof}
Take $E=\varnothing$ in the first case and $E=\{ab,ac,bc\}$
in the second. In either case, $B+E$ is nonplanar by
the hypothesis or Lemma~\ref{r2:32}, and is contained
in every applicable test of Lemma~\ref{r2:four-tests}.
By Lemma~\ref{r2:test-minors}, it is a $K_5$-minor-free minor
of $G$. Every applicable test for an allowed minor $B'$ of $B$
contains $B'+E$. By Lemma~\ref{r2:four-tests}, replacing $B$ by $B'$
therefore gives a nontoroidal graph whenever $B'+E$ is nonplanar.
Corollary~\ref{common:reduced-forms}, with attachments $a,b,c,d$ and
added edges $E$, lists every possible $B$.

For $E=\varnothing$, the construction gives $57{,}479$ bridges up to
isomorphisms fixing the attachments. The enumeration in~\supp{E8}
checks each with $(n_c,n_d)=(2,2),(3,2),(3,3)$. It leaves the
fourteen graphs in the first list, all with $(n_c,n_d)=(2,2)$.
For $E=\{ab,ac,bc\}$, the same enumeration checks $14{,}638$ bridges
with $(n_c,n_d)=(3,2)$ and leaves only $\mathfrak T_{11}$.
Of the $3\cdot57{,}479+14{,}638=187{,}075$ candidates, every other
graph has a $K_5$ minor, a cut of at most two vertices, or a proper
nontoroidal minor.
\end{proof}

\subsection{A bridge with four attachments: planar \texorpdfstring{$B$}{B} and \texorpdfstring{$(n_c,n_d)=(2,2)$}{(nc,nd)=(2,2)}}\label{r2:planar-four}

It remains to consider $G=B\cup\Psi_2$ with $B$ planar and
$(n_c,n_d)=(2,2)$. By Lemma~\ref{r2:four-tests}, all four graphs
\[
B+ac+bd,\qquad B+ad+bc,\qquad
B+ab+ac+ad,\qquad B+ab+bc+bd
\]
are nonplanar. Let $H_0$ be the graph on $a,b,c,d$ with edges
$ab,ac,bc,ad,bd$. The graph $B\cup H_0$ is a three-connected,
$K_5$-minor-free minor of $G$ by Lemma~\ref{r2:test-minors}.
It is nonplanar since it contains $B+ac+bd$, and it contains
the triangle $abc$.
Halin's theorem therefore gives a three-cut $X$ such that
$(B\cup H_0)-X$ has at least three components.

\begin{lemma}\label{r2:cut-locations}
The cut $X$ has, up to interchanging $a,b$ or $c,d$, one of
the forms $\{a,b,w\}$, $\{a,u,v\}$, $\{a,c,w\}$, where
$u,v,w\in\operatorname{Int}(B)$.
In the latter two cases, exactly two components of $(B\cup H_0)-X$
contain none of $a,b,c,d$. Both are single vertices $p',q'$
with neighborhood $X$ in $G$. In each test, $p',q'$ then give two bridges on $X$ adjacent to all of $X$, each with planar triangle completion, and the test is nonplanar exactly when another bridge is adjacent to all of $X$, has two attachments $\alpha,\beta$ and is nonplanar after adding $\alpha\beta$, or has at most one attachment and is nonplanar.
\end{lemma}
\begin{proof}
Put $A=\{a,b,c,d\}$. If $X$ contains three attachments, then $(B\cup H_0)-X$ consists of the connected interior of $B$ and the remaining attachment, which has a neighbor there. Thus $(B\cup H_0)-X$ is connected, a contradiction. If $X\cap A=\{a,b\}$, we have the first form.

Otherwise, if a component contains an attachment then it is connected to all of $A-X$ via the added edges. Thus there are at least two components that avoid $A$. These are also components of $B-X$, and then
Lemma~\ref{lem:avoiding-components} implies they are exactly two single vertices $p',q'$ with neighborhood $X$ in $G$. Thus $X$ is a twin triple and, by the assumption of this
section, meets both $\{a,b,c\}$ and $\{a,b,d\}$.

Suppose $X$ contains two attachments. If one is in $\{a,b\}$
and the other in $\{c,d\}$, then after possibly interchanging
$a,b$ or $c,d$, we have $X=\{a,c,w\}$. The only remaining
possibility is $X=\{c,d,w\}$. Then $p',q'$ are leaves of
$B[\operatorname{Int}(B)]$, both adjacent to $w$.
Deleting them leaves a connected graph still adjacent
to $a,b$, since neither attachment is adjacent to
$p'$ or $q'$. Contract this graph to $w$.
The five sets
\[
\{a,q,r\},\quad \{b,p\},\quad \{c,q'\},\quad
\{d,p',z\},\quad \{w\}
\]
are connected and pairwise adjacent, giving a $K_5$ minor,
a contradiction.

If $X$ contains one attachment, it must belong to both
$\{a,b,c\}$ and $\{a,b,d\}$, so it is $a$ or $b$.
After possibly interchanging $a,b$, we have $X=\{a,u,v\}$. If $X$ contains no attachments, it meets neither triple,
a contradiction. The last claim follows from
Lemma~\ref{core:assembly-general}.
\end{proof}

By Lemma~\ref{r2:four-tests}, the four test graphs are the relevant completions of $B$. By Lemma~\ref{core:relevant}, a reduction of $B$ to an allowed minor therefore leaves $B$ unchanged if it keeps all four tests nonplanar. Several of the reductions we will use below retain a $K_{3,3}$ subdivision in a component of $B-X$, by this next lemma.

\begin{lemma}\label{r2:retained-components}
Let $D$ be a component of $B-X$, and let $M$ be the set of
attachments of $B$ in $D$. Suppose $H_1$ is a $K_{3,3}$ subdivision
in $B[V(D)\cup X]$ after adding edges among $M\cup X$ and possibly
identifying $a$ with $b$. Fix any $k$ components of $D-M$.
By deletions and contractions within $D-M$, we can keep these
components and shorten $H_1$ as in the proof of
Lemma~\ref{core:retained-subdivision}, leaving at most $k$ vertices
of $D-M$ outside the shortened subdivision. Each remaining
component keeps its neighbors in $M\cup X$.
\end{lemma}
\begin{proof}
Split each branch path at its internal vertices in $M\cup X$.
Leave any segment that is an edge of $B$ unchanged and shorten
the others as in Lemma~\ref{core:retained-subdivision}. Only edges within $D-M$ are
contracted, so we can do this in $B$ with $a,b$ distinct even
when they were identified in the subdivision.

In each component meeting the shortened subdivision, repeatedly
contract a vertex outside the subdivision into a neighbor on it.
Since the component is connected, we can continue until none of
its vertices remain outside the subdivision. Of the components
disjoint from the subdivision, we keep only those among the $k$
chosen components. Contracting each of these to a single vertex
leaves at most $k$ vertices of $D-M$ outside the subdivision.

All contractions are within components of $D-M$, so each
component we keep remains connected and has the same neighbors
in $M\cup X$. To retain a path between two
vertices of $M$ through $D-M$, we keep a component adjacent to both.
\end{proof}

After applying Lemma~\ref{r2:retained-components}, delete the edges of $D$ outside the subdivision, except for a minimal set preserving the required neighborhoods and connectivity. Each proof below shows that these neighborhoods, this connectivity and the nonplanarity of the retained completion keep all four tests nonplanar. Minimality therefore leaves $B$ unchanged, and Construction~\ref{common:marked-construction} lists every possible form.

\begin{samepage}
\begin{lemma}\label{r2:one-attachment}
Let $D$ be a component of $B-X$ containing exactly one of
$a,b,c,d$, say $g$.
\begin{enumerate}[label=\textup{(\roman*)}]
\item If $N_B(D)=\{s\}$, then $D=\{g\}$ and $gs\in E(B)$.
\item If $N_B(D)=\{s,t\}$, then each component of $D-g$ is a
component of $G-\{g,s,t\}$ with neighborhood $\{g,s,t\}$. There are
at most two, and if there are two, both are single vertices.
\end{enumerate}
\end{lemma}
\end{samepage}
\begin{proof}
The vertices of $D-g$ are interior vertices of $B$, so each component
of $D-g$ has all its neighbors outside it in $\{g\}\cup N_B(D)$. In
(i), three-connectivity leaves no such component, so $D=\{g\}$ and
$gs\in E(B)$. In (ii), each component of $D-g$ is a component of
$B-\{g,s,t\}$ disjoint from $\{a,b,c,d\}$, and
Lemma~\ref{lem:avoiding-components} applies.
\end{proof}

\subsubsection{A cut of the form \texorpdfstring{$\{a,b,w\}$}{\{a,b,w\}}}\label{r2:ab-cut}

This subsection shows that a cut of this form gives $\mathfrak T_{56},\mathfrak T_{72},\mathfrak T_{73}$. We first determine the possible forms of the components of $B-X$ and then combine them in Proposition~\ref{r2:ab-cut-outcome}.

Let $X=\{a,b,w\}$. Every edge of $H_0$ meets $a$ or $b$, so $(B\cup H_0)-X=B-X$ has at least three components. At most two contain $c$ or $d$, so at least one contains neither. By Lemma~\ref{lem:avoiding-components}, there are at most two such components. Each is a component of $G-X$ with neighborhood $X$, and if there are two, both are single vertices. Since the interior of $B$ is connected, every component of $B[\operatorname{Int}(B)]-w$ has a neighbor at $w$. Every component of $B-X$ also has a neighbor at $w$, since $c,d$ each have an interior neighbor.

\begin{lemma}\label{r2:ab-cut-single}
Suppose exactly one component $D_0$ of $B-X$ contains neither $c$ nor $d$, and let $B_0$ be its bridge on $X$. Then $c,d$ lie in separate components $D_c,D_d$ of $B-X$, and $B_0$ has one interior vertex or is one of
\[
T(a),\ T(b),\ T(w),\ S(a;b),\ S(a;w),\ S(b;a),\ S(b;w),\
Q(a;b,w),\ Q(b;a,w),\ Q(w;a,b).
\]
\end{lemma}
\begin{proof}
Since $B-X$ has at least three components and only $D_0$ contains neither $c$ nor $d$, these vertices lie in separate components $D_c,D_d$. The edge $cd$ is absent by Lemma~\ref{r2:unique-four}, so the components of $G-X$ are $D_0$, $D_c\cup\{p,q\}$ and $D_d\cup\{r,z\}$. Lemma~\ref{lem:triangle-completion} therefore makes the triangle completion of $B_0$ a three-connected minor of $G$.

In each of the four tests, $B_0$ is a bridge on $X$ adjacent to all of $X$. By Lemma~\ref{core:assembly-general}, the conditions on $B_0$ depend only on which graphs $B_0+E$, with $E\subseteq\{ab,aw,bw\}$, are planar. Replacing $B_0$ by an allowed minor gives an allowed minor of $B$, since every component of $D_c-c$ and $D_d-d$ remains adjacent to $w$. These graphs are therefore the relevant completions of $B_0$, and by Lemma~\ref{core:relevant} every proper allowed minor of $B_0$ makes at least one previously nonplanar completion planar.

\looseness=-1
At least one of $B_0+ab+aw$, $B_0+ab+bw$ and $B_0+aw+bw$ is planar. Otherwise, contract a $c$--$w$ path through $D_c$ and a $d$--$w$ path through $D_d$ into $w$. Delete the remaining vertices of $D_c,D_d$ and any edges $aw,bw$. What remains is $B_0$ together with $p,q,r,z$, each now having neighborhood $\{a,b,w\}$. This is a proper minor of $G$, since $c,d$ were contracted, so it is toroidal. Lemma~\ref{lem:four-twins} then makes at least one of the three completions planar, a contradiction.

Now suppose $B_0+ab$ is nonplanar. Then $B_0+ab+aw$ and $B_0+ab+bw$ are nonplanar, so $B_0+aw+bw$ is planar. Apply Proposition~\ref{prop:three-root-eight} with $x=w$. The single-vertex form, $T(a),T(b)$ and $Q(w;a,b)$ remain planar after adding $ab$, so $B_0=S(w;a)$ or $S(w;b)$. Suppose $B_0=S(w;a)$, with interior vertices $p',q',u,f$, where $p',q'$ have neighborhood $\{w,a,u\}$ and $N(f)=\{w,u,b\}$. Contract $fb$ and retain the resulting edge $bw$. This replaces $B_0$ by $T(b)$ together with $bw$. By Lemma~\ref{core:three-graph-completions}, $T(b)+bw+E$ is planar exactly when $ab\notin E$, which is also the condition for $B_0+E$ to be planar. The new edge $bw$ is itself a bridge on $X$. Thus, by
Lemma~\ref{core:assembly-general}, all four tests remain
nonplanar after the contraction. The path $wp'u$ joins the only two neighbors of $f$ in the
interior of $B$, so the interior stays connected when we contract
$fb$. The contraction gives $bu$, while $a$ is still adjacent to
$p'$ and $c,d$ keep their interior neighbors. Since $w$ is an interior vertex of $B$, retaining $bw$ creates no edge between its attachments $a,b,c,d$. Thus this is a proper allowed minor of $B$ keeping all four tests nonplanar, contrary to minimality. Interchanging $a,b$ excludes $S(w;b)$.

Hence $B_0+ab$ is planar, and at least one completion obtained by adding two edges is planar. Apply Proposition~\ref{prop:three-root-eight} with $x$ the attachment incident with both added edges. The forms $S(w;a),S(w;b)$ are excluded because adding $ab$ makes them nonplanar. Taking the three possible choices of $x$ leaves the eleven forms in the statement.
\end{proof}

The next two lemmas give the possible forms of the components of $B-X$ containing $c$ or $d$.

\begin{lemma}\label{r2:ab-cut-apart}
Suppose $c,d$ lie in different components $D_c,D_d$ of $B-X$, and put
\[
H_c=B[V(D_c)\cup X]\cup K_3[X].
\]
If $H_c+ac+bc$ is planar, then $D_c=\{c\}$ with the edge $cw$, or
$D_c-c$ has at most two vertices, each adjacent to $c,w$.
Otherwise Construction~\ref{common:marked-construction} gives every
possible form of $D_c$ with one of the following choices.
\begin{enumerate}[label=\textup{(\roman*)}]
\item Use attachments $a,b,c,w$ and add the triangle on $X$ and both
$ac,bc$. The subdivision uses both $ac,bc$.
\item Use the same attachments and triangle with just one of $ac,bc$.
The subdivision uses that edge.
\item Identify $a,b$ to $x$ and use attachments $x,c,w$ with added
edges $xw,xc$. The subdivision uses $xc$. Then split $x$ back into
$a,b$, assigning each original neighbor of $x$ to $a$, to $b$, or
to both.
\end{enumerate}
In each case a $c$--$w$ segment may be the edge $cw$ or a path of
length two. At most two vertices of $D_c-c$ lie outside the shortened
subdivision. The same holds for $D_d$ with $d$ in place of $c$.
\end{lemma}
\begin{proof}
Each of the four tests in Lemma~\ref{r2:four-tests} adds exactly one
of $ac,bc$. In a test adding $ac$, the bridge containing $D_c$
attaches at $a,w$. If it also attaches at $b$, its triangle
completion is $H_c+ac$. Otherwise its completion adds just $aw$,
and the path $abw$ can be drawn beside $aw$ to give $H_c+ac$.
The same argument gives $H_c+bc$ for the tests adding $bc$.
Since $B-X$ has at least three components, choose a component $D_0$
containing neither $c$ nor $d$. Its bridge on $X$ attaches at
$a,b,w$ in every test. By Lemma~\ref{core:assembly-general}, it
therefore suffices to keep $D_c$ connected, preserve $N_B(D_c)$,
and preserve which of $H_c+ac,H_c+bc$ are planar.

In a plane drawing of $B$, contract $D_0$ to one vertex and replace
the three edges from this vertex to $a,b,w$ by the triangle on $X$.
Restricting the drawing to $H_c$ shows that $H_c$ is planar.
In $B\cup H_0$, contract $D_0$ into $w$ to obtain $aw,bw$.
Since $H_0$ contains $ab,ac,bc$, the graph $H_c+ac+bc$ is a
$K_5$-minor-free minor of $B\cup H_0$.

If neither $a$ nor $b$ has a neighbor in $D_c$, then
Lemma~\ref{r2:one-attachment}(i) gives $D_c=\{c\}$ with the edge
$cw$. Otherwise, for each of $a,b$ with a neighbor in $D_c$, choose
one component of $D_c-c$ containing such a neighbor. Each component
of $D_c-c$ is adjacent to $c$ since $D_c$ is connected, and to $w$
since it is a component of $B[\operatorname{Int}(B)]-w$.
If $H_c+ac+bc$ is planar, contract each chosen component to one
vertex and delete the others. This leaves at most two vertices
besides $c$, each adjacent to $c,w$, and both $H_c+ac$ and
$H_c+bc$ remain planar.

Suppose $H_c+ac+bc$ is nonplanar. By Kuratowski's theorem, each
nonplanar minor of this graph contains a $K_{3,3}$ subdivision.
If $H_c+ac$ and $H_c+bc$ are both planar, a subdivision in
$H_c+ac+bc$ must use both added edges, giving (i). If just one of
$H_c+ac$ and $H_c+bc$ is nonplanar, take a subdivision in that
graph. It uses the added edge since $H_c$ is planar, giving (ii).

If both are nonplanar, contracting $ab$ in either gives the same
graph $(H_c/ab)+xc$, where $x$ is the vertex obtained by identifying
$a,b$. This graph is nonplanar by Lemma~\ref{r2:planar-splitting}.
A $K_{3,3}$ subdivision in it uses $xc$ since $H_c/ab$ is planar.
This gives (iii), and retaining the subdivision in the common minor
keeps both single-edge completions nonplanar.

In (i)--(iii), apply Lemma~\ref{r2:retained-components} to the
subdivision, keeping the components chosen for $a,b$ and making the
contractions with $a,b$ distinct. In both the planar and nonplanar
cases, $D_c$ stays connected with the same neighbors in $X$, and
its remaining interior components stay connected through $w$. Thus the
reduction is allowed and all four tests remain nonplanar. By
minimality $B$ is unchanged, so $D_c$ already has the stated form.
\end{proof}

\begin{lemma}\label{r2:ab-cut-same}
Suppose $c,d$ lie in the same component $D_1$ of $B-X$. Then, after interchanging $a,b$ if necessary, $N_B(D_1)=\{w\}$ or $\{a,w\}$, and Construction~\ref{common:marked-construction} gives every possible form of $D_1$ with attachments $a,c,d,w$ and added edges $aw,ac,ad$, requiring the subdivision to use at least one of $ac,ad$. For a segment with ends $c,w$ or $d,w$, either the original edge or a path with one interior vertex is allowed.
\end{lemma}
\begin{proof}
Since $B-X$ has at least three components, Lemma~\ref{lem:avoiding-components} gives exactly two other components, each a single vertex adjacent to all of $X$. If $N_B(D_1)=X$, contracting $D_1$ to a vertex would give a $K_{3,3}$ minor in the planar graph $B$. Since $D_1$ has a neighbor at $w$, we may interchange $a,b$ and assume $N_B(D_1)=\{w\}$ or $\{a,w\}$.

In $B+ac+bd$, the two added edges join both $a$ and $b$ to $D_1$. Contracting $D_1$ to a vertex then gives, with the two single vertices, a $K_{3,3}$ minor. The same holds in $B+ad+bc$. For the other two tests, let $B_1$ be the bridge of $D_1$ on $X$ in $B$, and put
\[
H_1=B_1+\{a,w\}+aw.
\]
In $B+ab+ac+ad$, both $ac,ad$ join $a$ to $D_1$, and $ab$ is a separate edge bridge on $X$. The bridge containing $D_1$ has exactly two attachments $a,w$, and adding $aw$ to it gives $H_1+ac+ad$. The two single-vertex bridges each have triangle completion $K_4$, so Lemma~\ref{core:assembly-general}(ii) makes this test nonplanar exactly when $H_1+ac+ad$ is nonplanar. In $B+ab+bc+bd$, if $N_B(D_1)=\{w\}$, the same argument gives the same completed graph with $b$ in place of $a$. If $N_B(D_1)=\{a,w\}$, adding $bc$ makes the bridge adjacent to all of $X$, again giving a $K_{3,3}$ minor. Thus it suffices to keep $D_1$ connected, preserve $N_B(D_1)$ and keep $H_1+ac+ad$ nonplanar.

An $a$--$w$ path through one of the single vertices supplies $aw$ outside $B_1$. Contracting one edge of this path gives $H_1$ as a minor of the planar graph $B$. The same contraction in $B\cup H_0$, retaining $ac,ad$, gives $H_1+ac+ad$, so it has no $K_5$ minor. Choose a $K_{3,3}$ subdivision in $H_1+ac+ad$. It uses at least one of $ac,ad$, since deleting~both~leaves~$H_1$.

Each component of $D_1-\{c,d\}$ has a neighbor at $w$ and at least one of $c,d$. Since $D_1$ is connected and $cd$ is absent, the interior of a $c$--$d$ path in $D_1$ lies in one such component. By Lemma~\ref{r2:retained-components}, we can shorten the
subdivision while keeping this component. If $a$ has a neighbor
in $D_1$ but none in this component, we also keep a component
containing a neighbor of $a$.
The path still joins $c$ to $d$, and each remaining component
has a neighbor at $c$ or $d$, so $D_1$ stays connected. Each
also has a neighbor at $w$, so the interior of $B$ stays
connected. Since $D_1$ keeps its neighbors in $X$ and the
shortened subdivision keeps $H_1+ac+ad$ nonplanar, all four
tests remain nonplanar. By minimality the reduction leaves $B$
unchanged, so $D_1$ has one of the forms in
Construction~\ref{common:marked-construction}.
\end{proof}

\begin{proposition}\label{r2:ab-cut-outcome}
If $X=\{a,b,w\}$, then $G$ is isomorphic to one of $\mathfrak T_{56},\mathfrak T_{72},\mathfrak T_{73}$.
\end{proposition}

\begin{proof}
We distinguish whether $c,d$ lie in different components of $B-X$ or in the same component. Lemmas~\ref{r2:ab-cut-apart} and~\ref{r2:ab-cut-same} give the
possible forms of $D_c,D_d$ and $D_1$. We discard a form if a
single allowed deletion or contraction leaves the component
connected with the same neighbors in $X$, keeps each interior
component adjacent to $w$, and preserves which of the completions
in the corresponding lemma are planar. This leaves $33$ forms
for $D_c$ and $54$ for $D_1$. Replacing $c$ by $d$ in the
forms for $D_c$ gives the forms for $D_d$.

When $c,d$ lie in different components, we pair each form for
$D_c$ with every form for $D_d$. If exactly one component avoids $c,d$, include one of the eleven bridges of Lemma~\ref{r2:ab-cut-single} and a permitted subset of $\{aw,bw\}$. If two components avoid $c,d$, include the two single vertices adjacent to all of $X$. In this case $aw,bw$ are absent by Lemma~\ref{lem:independent-neighborhood}(a). If $c,d$ lie together, take one retained form for $D_1$ and include those same two single vertices. These choices give $12{,}645$ graphs up to isomorphism, of which $1{,}471$ have $B$ planar. Among these the checks in~\supp{E10} leave only $\mathfrak T_{56},\mathfrak T_{72},\mathfrak T_{73}$. In each, exactly one component of $B-X$ avoids $c,d$.
\end{proof}

\subsubsection{A cut of the form \texorpdfstring{$\{a,u,v\}$}{\{a,u,v\}}}

Let $X=\{a,u,v\}$, and let $p',q'$ be as in Lemma~\ref{r2:cut-locations}, with $N(p')=N(q')=X$. This set is independent by Lemma~\ref{lem:independent-neighborhood}(a). Every other component of $B-X$ contains at least one of $b,c,d$ and has a neighbor at $u$ or $v$, since these attachments have interior neighbors and the interior of $B$ is connected. None has neighborhood $X$, since together with $p',q'$ it would give three bridges adjacent to all of $X$, contrary to the planarity of $B$.

We distinguish the cases according to whether $b,c,d$ lie in one, two or three components of $B-X$. The only edges of $H_0$ that survive deletion of $X$ are $bc$ and $bd$, so
\[
(B\cup H_0)-X=(B-X)+bc+bd.
\]
In $(B\cup H_0)-X$ the three vertices therefore lie in one component.

\begin{lemma}\label{r2:a-cut-one}
If one component $D$ of $B-X$ contains $b,c,d$, then $G$ is isomorphic to one of $\mathfrak T_{40},\mathfrak T_{53},\mathfrak T_{56},\mathfrak T_{73}$.
\end{lemma}

\begin{proof}
If $u$ had no neighbor in $D$, its only neighbors in $G$ would be
$p',q'$, contrary to three-connectivity. The same holds for $v$.
Since $D$ cannot have neighborhood $X$ in $B$, we have
$N_B(D)=\{u,v\}$.

Consider now $D-\{b,c,d\}$ and its connected components (possibly
only one) $C_1,\ldots,C_m$. Each is a component of the interior of
$B$ with $u,v$ deleted. Suppose that neither $u$ nor $v$ has an
edge to $C_i$. Since $C_i$ also has no edge to any other component
left after deleting $u,v$, it has no edge to the rest of the
interior of $B$. This contradicts the connectivity of the interior
of $B$. Therefore each $C_i$ has a neighbor at $u$ or $v$.

Contract each $C_i$ to one vertex $x_i$ and
keep a spanning tree of $D$. Retain an inclusion-minimal set of edges
from $D$ to $\{u,v\}$ such that each $x_i$ has a neighbor at $u$ or
$v$ and both $u$ and $v$ have a neighbor in $D$. The tree gives each
of $b,c,d$ an interior neighbor, and the interior stays connected
since each $x_i$ meets $u$ or $v$ and $p'$ joins $u$ to $v$.
This gives an allowed minor of $B$, and all four tests remain nonplanar
by Lemma~\ref{r2:cut-locations}. By minimality, each $C_i$ is a single vertex, $D$ is this tree, and its edges to $u$ and $v$ are precisely those retained in the reduction.

The tree $D$ is bipartite with parts $\{b,c,d\}$ and
$\{x_1,\ldots,x_m\}$. If $m=1$, then $D$ is the star with center
$x_1$ and leaves $b,c,d$. The center is adjacent to both $u,v$, or
to one of them while one of $b,c,d$ is adjacent to the other.
This gives $1+2\cdot3=7$ forms.

Now suppose $m\ge2$. If some $x_i$ were adjacent to both $u,v$,
take another $x_j$. Since $x_j$ has a neighbor at $u$ or $v$, we
may interchange $u,v$ so that $x_j$ is adjacent to $u$. Then $u$
already has $x_j$ as a neighbor in $D$, and $x_i$ still has $v$
as a neighbor if we delete $ux_i$. This contradicts the minimal
choice of edges.

Thus each $x_i$ has exactly one neighbor outside $D$ and hence at
least two in $D$, since $G$ has minimum degree three. Each edge
of $D$ has exactly one end among the $x_i$, so $D$ has at least
$2m$ edges. Since $D$ is a tree on $m+3$ vertices, it has $m+2$
edges. Thus $2m\le m+2$, and $m=2$.

Since $D$ has only four edges, both $x_1$ and $x_2$ must have
degree two in $D$. Their
neighborhoods in $\{b,c,d\}$ together contain all three vertices,
so they have exactly one common neighbor $g_2$. Thus $D$ is the
path $g_1x_1g_2x_2g_3$ for some ordering $g_1,g_2,g_3$ of $b,c,d$.

For each of the three choices of $g_2$, there are two ways to join
$x_1,x_2$ to different vertices of $\{u,v\}$. If both are adjacent
to the same vertex, we have two choices for that vertex and three
choices for the attachment in $\{b,c,d\}$ adjacent to the other.
This gives $3(2+2\cdot3)=24$ path forms. Together with the seven
stars, these are the $31$ forms checked in~\supp{E11}, leaving
exactly $\mathfrak T_{40},\mathfrak T_{53},\mathfrak T_{56},\mathfrak T_{73}$.
\end{proof}

\begin{lemma}\label{r2:a-cut-single}
Let $D_g$ be a component of $B-X$ containing exactly one of $b,c,d$,
say $g$, and suppose $D_g\ne\{g\}$. Then $N_B(D_g)=\{s,t\}$ for some $s,t\in X$. Either $D_g-g$ consists of two single vertices with neighborhood $\{g,s,t\}$, or it is connected and its bridge $B_g$ on $\{g,s,t\}$ is one of nine forms, namely one interior vertex, $T(g),T(s),T(t)$, $S(s;t),S(t;s)$, $Q(g;s,t)$ and the two four-cycle forms of Lemma~\ref{core:two-cap-critical}.
\end{lemma}
\begin{proof}
Since $D_g\ne\{g\}$, Lemma~\ref{r2:one-attachment}(i) excludes a
single neighbor in $X$. We already know that $D_g$ cannot have
neighborhood $X$, so $N_B(D_g)=\{s,t\}$ for some $s,t\in X$.
By Lemma~\ref{r2:one-attachment}(ii), each component of $D_g-g$
has neighborhood $\{g,s,t\}$ in $G$. There are at most two, and
if there are two, both are single vertices. This gives the first
outcome.

Suppose now that $D_g-g$ is connected, and let $B_g$ be its bridge
on $\{g,s,t\}$. Any original edges $gs,gt$ lie outside $B_g$.
Since $D_g-g$ is a component of $G-\{g,s,t\}$,
Lemma~\ref{lem:triangle-completion} makes the triangle completion
of $B_g$ three-connected. If it is planar,
Lemma~\ref{lem:flat-replacement} makes $D_g-g$ a single vertex.

\looseness=-1
Suppose the triangle completion is nonplanar. We first show that
it is a minor of $G$. One of the edges $bc,bd$ of $H_0$ joins $g$
to a vertex in another component of $B-X$. That component has a
neighbor at $u$ or $v$, say $u$, so $B\cup H_0$ contains a path
$P$ from $u$ to $g$ with all its internal vertices in that component.
Since $ag\in E(H_0)$, this gives the cycle
\[
g-a-p'-v-q'-u\xrightarrow{P}g.
\]
This cycle avoids the interior of $B_g$ and contains $g,s,t$
since it passes through all of $X$. Retain $B_g$ and this cycle
and shorten the three arcs between $g,s,t$ to single edges.
This gives the triangle completion as a minor of $B\cup H_0$,
hence of $G$.

In each test of Lemma~\ref{r2:four-tests}, the path $sp't$ joins
$s,t$ outside $B_g$. If $B_g$ is the only bridge adjacent to all of
$g,s,t$, the bridge containing $sp't$ has attachments $s,t$, so
Lemma~\ref{core:assembly-general}(iii) includes $st$ in the completion
of $B_g$. With exactly two bridges adjacent to all three vertices,
part (ii) requires the full triangle. Three such bridges already
make the test nonplanar by part (i). Thus it suffices to consider
$B_g+st+E$ with $E\subseteq\{gs,gt\}$. Retaining $B_g$ and the
path $sp't$ in $B$ and contracting $sp'$ gives $B_g+st$, so this
completion is planar.

An allowed reduction of $B_g$ keeps its interior connected and
adjacent to $g,s,t$. At least one of $s,t$ lies in $\{u,v\}$, so
the reduced interior stays joined to the rest of the interior of $B$.
Thus the reduction is allowed in $B$, and the four completions above
are relevant for $B_g$. By Lemma~\ref{core:relevant}, every proper
allowed minor makes at least one nonplanar relevant completion planar.
Proposition~\ref{prop:three-root-disjoint} now gives the other eight
forms with $s,t,g$ in the roles of $x,y,z$.
\end{proof}

The preceding lemma applies whether $b,c,d$ lie in two or three
components of $B-X$. We now consider the two-component case.

\begin{lemma}\label{r2:a-cut-table}
Suppose two of $b,c,d$ lie in a component $D_1$ of $B-X$ and the
remaining vertex $g$ lies in another component $D_2$. Let $B_1,B_2$
be their respective bridges on $X$. Up to interchanging $c,d$ or
$u,v$, the following rows give the possible neighborhoods and, in
each case, the condition for all four tests of
Lemma~\ref{r2:four-tests} to be nonplanar.
\end{lemma}

\noindent
\begin{minipage}{\linewidth}
\small
\begin{tabular*}{\linewidth}{@{\extracolsep{\fill}}rllll@{}}
\toprule
& $D_1$ contains & $N_B(D_1)$ & $N_B(D_2)$ & Condition \\
\midrule
(1) & $b,c$ & $\{u\}$ or $\{a,u\}$ & $\{u,v\}$ & always \\
(2) & & $\{u\}$ or $\{a,u\}$ & $\{v\}$ or $\{a,v\}$ &
  \begin{tabular}[t]{@{}l@{}}
  $B_2+av+ad$ is nonplanar, or both \\
  $H_1+bc$ and $H_1+ab+ac$ are nonplanar
  \end{tabular} \\
(3) & & $\{u,v\}$ & $\{u,v\}$ & always \\
(4) & & $\{u,v\}$ & $\{v\}$ or $\{a,v\}$ &
  $B_1+uv+bc$ or $B_2+av+ad$ is nonplanar \\
\addlinespace
(5) & $c,d$ & \multicolumn{2}{l}{at least one is $\{u,v\}$} & always \\
(6) & & $\{u\}$ or $\{a,u\}$ & $\{v\}$ or $\{a,v\}$ &
  $B_1+au+ac+ad$ or $B_2+av+ab$ is nonplanar \\
\bottomrule
\end{tabular*}
\par\smallskip
\noindent In row~(2), $H_1=B_1+\{a,u\}+au$.
\end{minipage}
\par
\begin{proof}
Each of $u,v$ has a neighbor in $D_1\cup D_2$, since its only
neighbors outside $D_1\cup D_2$ are $p',q'$ and $G$ has minimum
degree three. Also, $B-a$ is connected, since the interior of $B$
is connected and each of $b,c,d$ has an interior neighbor. Thus
each component of $B-X$ has a neighbor at $u$ or $v$.

The bridges of $p'$ and $q'$ meet all of $X=\{a,u,v\}$. By
Lemma~\ref{core:assembly-general}(i), neither $D_1$ nor $D_2$
can have neighborhood $X$ in the planar graph $B$. Each test of
Lemma~\ref{r2:four-tests} adds an edge from $a$ to $D_1\cup D_2$.
If it also joins $D_1$ and $D_2$, their bridge meets all of $X$,
so the test is nonplanar. If they remain separate, we apply
Lemma~\ref{r2:cut-locations}. The test is then nonplanar exactly
when one of their bridges meets all of $X$, has two attachments
and is nonplanar after adding the edge between them, or has one
attachment and is nonplanar.

Suppose first that $D_1$ contains $b,c$, so $D_2$ contains $d$.
The two tests containing $bd$ join $D_1$ and $D_2$. In the
remaining tests, $B+ad+bc$ and $B+ab+ac+ad$, the edge $ad$ joins
$a$ to $D_2$. If $D_2$ meets both $u,v$, its bridge therefore
meets all of $X$ in both tests. If $D_1$ meets just one of $u,v$,
call it $u$. Then $N_B(D_1)$ is $\{u\}$ or $\{a,u\}$, giving
row~(1). If $D_1$ meets both, then $N_B(D_1)=\{u,v\}$, giving
row~(3).

Suppose now that $D_2$ meets $v$ and misses $u$, interchanging
$u,v$ if necessary. Then $u$ must have a neighbor in $D_1$.
In both remaining tests, retain $D_2$, all edges from $a,v$ into
$D_2$ including $ad$, and the path $ap'v$. Contracting $ap'$
gives the bridge containing $D_2$ with $av$ added. If this
common minor is nonplanar, both tests are nonplanar. Suppose it
is planar. It remains to check the bridge containing $D_1$ in
each test.

\looseness=-1
Suppose $D_1$ misses $v$. In $B+ad+bc$, the edge $bc$ lies
within $D_1$. If $D_1$ meets $a$, its bridge has attachments
$a,u$, and adding $au$ gives $H_1+bc$. Otherwise $u$ is its
only attachment, and $H_1+bc$ is this bridge with the leaf $a$
added at $u$. In $B+ab+ac+ad$, the edges $ab,ac$ join $a$ to
$D_1$, so its bridge has attachments $a,u$, and adding $au$
gives $H_1+ab+ac$. Thus both tests are nonplanar exactly when
both $H_1+bc$ and $H_1+ab+ac$ are nonplanar, giving row~(2).

Suppose instead that $D_1$ meets both $u,v$. In $B+ab+ac+ad$,
the edge $ab$ makes its bridge adjacent to all of $X$, so this
test is nonplanar. In $B+ad+bc$, the edge $bc$ lies within
$D_1$, so its bridge still has just the attachments $u,v$.
Since the common minor is planar, Lemma~\ref{r2:cut-locations} shows that $B+ad+bc$ is nonplanar exactly when the bridge containing $D_1$ is nonplanar after adding $uv$. This is row~(4).

Now suppose $D_1$ contains $c,d$, so $D_2$ contains $b$. The
three tests containing $bc$ or $bd$ join $D_1$ and $D_2$. In the
remaining test $B+ab+ac+ad$, the edges $ac,ad$ join $a$ to $D_1$
and $ab$ joins $a$ to $D_2$. If either component meets both $u,v$,
its bridge now meets all of $X$, giving row~(5). Otherwise the
two components must meet different vertices of $\{u,v\}$, say
$D_1$ meets $u$ and $D_2$ meets $v$. Their bridges then have
attachments $a,u$ and $a,v$, respectively. The test is nonplanar
exactly when adding $au$ to the bridge containing $D_1$ or $av$
to the bridge containing $D_2$ gives a nonplanar graph. This
gives row~(6).
\end{proof}

We next rule out the possibility in Lemma~\ref{r2:a-cut-single}
that $D_2-g$ consists of two single vertices.

\begin{lemma}\label{r2:a-cut-second}
Let $D_2$ be the component containing $g$ in
Lemma~\ref{r2:a-cut-table}. If $D_2=\{g\}$, then $g$ is adjacent
in $B$ to $u$, to $v$, or to both. Otherwise $D_2-g$ is connected,
and the bridge with interior $D_2-g$ is as in Lemma~\ref{r2:a-cut-single}.
\end{lemma}
\begin{proof}
If $D_2=\{g\}$, then $g$ has an interior neighbor, and the only
interior vertices it can meet are $u,v$. Suppose $D_2\ne\{g\}$.
By Lemma~\ref{r2:a-cut-single}, it suffices to rule out the case
where $N_B(D_2)=\{s,t\}$ and $D_2-g$ consists of two single
vertices $x,y$, each with neighborhood $\{g,s,t\}$.

Delete $x$ from $B$. Since $y$ has the same neighbors as $x$,
the interior of $B$ stays connected and each attachment still has
an interior neighbor. Thus the deletion is allowed. Also, $D_2-x$
still has neighbors $s,t$ in $X$, so the row of
Lemma~\ref{r2:a-cut-table} stays the same.

In rows~(1), (3) and~(5), the neighborhoods alone make all four
tests nonplanar. In rows~(2), (4) and~(6), the only possibility
with two neighbors for $D_2$ is $\{s,t\}=\{a,v\}$. The graph
listed for $B_2$ is then $B_2+av+ag$. Before deleting $x$, this
graph is a subgraph of the planar graph $K_5-xy$ on $a,v,g,x,y$.
Thus the remaining condition in that row must already hold, and
it involves only the bridge containing $D_1$, which is unchanged.
Hence all four tests remain nonplanar after deleting $x$. By
Lemma~\ref{r2:four-tests}, this gives a proper nontoroidal minor
of $G$, contrary~to~minimality.
\end{proof}

In rows~(2), (4) and~(6) of Lemma~\ref{r2:a-cut-table}, we
distinguish two cases according to whether $B_2+av+ag$ is planar.
If it is nonplanar, the condition in the last column already holds.
We return to this case in Proposition~\ref{r2:a-cut-outcome}.
For now suppose instead that it is planar. We first handle rows~(4) and~(6), where a
single completion of $B_1$ must be nonplanar, and then row~(2),
where both $H_1+bc$ and $H_1+ab+ac$ must be nonplanar.

\begin{lemma}\label{r2:a-cut-row-four}
In row~(4) of Lemma~\ref{r2:a-cut-table}, suppose $B_2+av+ad$ is
planar. Then Construction~\ref{common:marked-construction} lists
every possible $B_1$ using only $K_{3,3}$ subdivisions in
$B_1+uv+bc$ that contain $bc$. At most two components of
$D_1-\{b,c\}$ miss the shortened subdivision, and each such
component is a single vertex.
\end{lemma}
\begin{proof}
Since $B_2+av+ad$ is planar, row~(4) requires $B_1+uv+bc$ to be
nonplanar. Contracting $up'$ on the path $up'v$ gives $B_1+uv$ as
a minor of the planar graph $B$. Keeping $bc$ as well gives
$B_1+uv+bc$ as a minor of the $K_5$-minor-free graph $B\cup H_0$.
Thus $B_1+uv+bc$ contains a $K_{3,3}$ subdivision, which must use
$bc$ since~$B_1+uv$~is~planar.

Choose a $b$--$c$ path in $D_1$. Since $bc\notin E(B)$, the interior of this path is nonempty and
lies in one component $C$ of $D_1-\{b,c\}$. By
Lemma~\ref{r2:retained-components}, we can shorten the subdivision
while keeping $C$ connected and adjacent to both $b$ and $c$, so
it still joins them in $D_1$. To stay in row~(4), $D_1$ must also
keep both $u$ and $v$ as neighbors. Every component of $D_1-\{b,c\}$ has a neighbor at
$u$ or $v$, since the interior of $B$ is connected. If $C$ misses
one of $u,v$, keep an edge from that vertex to $D_1$. Such an edge
exists since $N_B(D_1)=\{u,v\}$. If the edge ends at $b$ or $c$, we need no further component.
Otherwise we keep the component of $D_1-\{b,c\}$ containing its
end in $D_1$. Thus we keep at most two components disjoint from
the subdivision, each of which contracts to a single vertex.

Since $D_1$ was connected, each remaining component has a
neighbor at $b$ or $c$. Both $b$ and $c$ keep an interior
neighbor in $C$, so $D_1$ stays connected. Each
component also keeps a neighbor at $u$ or $v$, so
$(B_1-\{b,c\})+uv$ is connected. In $B$ the path $up'v$ joins
$u$ to $v$. Hence the reduction is allowed. It keeps
$N_B(D_1)=\{u,v\}$, and the shortened subdivision keeps
$B_1+uv+bc$ nonplanar. By row~(4), all four tests remain nonplanar. By minimality the
reduction cannot be proper, so $B_1$ already has one of the forms
in Construction~\ref{common:marked-construction}.
\end{proof}

\begin{lemma}\label{r2:a-cut-row-six}
In row~(6) of Lemma~\ref{r2:a-cut-table}, suppose $B_2+av+ab$ is
planar. Then Construction~\ref{common:marked-construction} lists
every possible $B_1$ using only $K_{3,3}$ subdivisions in
$B_1+au+ac+ad$ that contain at least one of $ac,ad$. Every
component of $D_1-\{c,d\}$ meets the shortened subdivision,
with the possible exception of a single vertex.
\end{lemma}
\begin{proof}
Since $B_2+av+ab$ is planar, row~(6) requires $B_1+au+ac+ad$ to
be nonplanar. Contracting $ap'$ on the path $ap'u$ gives $B_1+au$ as a minor
of the planar graph $B$. The same contraction in $B\cup H_0$,
keeping $ac,ad$, gives $B_1+au+ac+ad$ as a minor, so it has no
$K_5$ minor. It therefore contains a $K_{3,3}$ subdivision, which
uses at least one of $ac,ad$ since $B_1+au$ is planar.

As in row~(4), we shorten the subdivision by
Lemma~\ref{r2:retained-components}, keeping the component $C$ of
$D_1-\{c,d\}$ containing the interior of a $c$--$d$ path in
$D_1$. Every component of $D_1-\{c,d\}$ has a neighbor at $u$,
since the interior of $B$ is connected and $D_1$ has no neighbor
at $v$. Row~(6) allows both $\{u\}$ and $\{a,u\}$ as the
neighborhood of $D_1$, so we need no second component to keep a
neighbor at $a$. If any remaining component misses the shortened
subdivision, it is $C$ contracted to a single vertex.

Each remaining component has a neighbor at $c$ or $d$, and $C$
has neighbors at both, so $D_1$ stays connected. Their neighbors
at $u$ keep $B_1-\{a,c,d\}$ connected, and $a$ keeps its
interior neighbor $p'$. Thus the reduction is allowed. The
shortened subdivision keeps $B_1+au+ac+ad$ nonplanar, so all
four tests remain nonplanar by row~(6). By minimality $B_1$
already has one of the forms in
Construction~\ref{common:marked-construction}.
\end{proof}

\begin{lemma}\label{r2:a-cut-row-two}
In row~(2) of Lemma~\ref{r2:a-cut-table}, suppose $B_2+av+ad$ is
planar. Then $D_2=\{d\}$ and $dv$ is an edge of $B$.
Construction~\ref{common:marked-construction} lists every possible
$B_1$ using only $K_{3,3}$ subdivisions $H$ in $H_1+bc$ that
contain $bc$ and for which $H-bc+ab+ac$ is nonplanar.
\end{lemma}
\begin{proof}
Put $H_1=B_1+\{a,u\}+au$ as in Lemma~\ref{r2:a-cut-table}.
Since $B_2+av+ad$ is planar, row~(2) requires both $H_1+bc$ and
$H_1+ab+ac$ to be nonplanar. Suppose $D_2\ne\{d\}$. By
Lemma~\ref{r2:a-cut-second}, $D_2-d$ is connected and has
neighbors at $d$ and $v$. Contract $D_2-d$ into $v$ and delete
any edge $av$ created by the contraction. The interior of $B$
stays connected and $d$ keeps its interior neighbor $v$. Since
$\{v\}$ is an allowed neighborhood for $\{d\}$ in row~(2), and
both $H_1+bc$ and $H_1+ab+ac$ are unchanged, all four tests
remain nonplanar. This contradicts minimality, so $D_2=\{d\}$
with the edge $dv$.

To apply Lemma~\ref{r2:three-mark-entrance}, put
$H_2=H_1+ab+ac+bc$. The path $ap'u$ supplies $au$, so $H_1$ is
a planar minor of $B$. Keeping $ab,ac,bc$ as well gives $H_2$
as a $K_5$-minor-free minor of $B\cup H_0$. Every component of
$D_1-\{b,c\}$ has a neighbor at $u$, since the interior of $B$
is connected and $D_1$ misses $v$. Thus $H_1-\{a,b,c\}$ is
connected. The edge $au$ joins it to $a$, and $b,c$ each have
an interior neighbor in it. In particular, $H_2-a$ is connected.

After deleting at most two vertices of $H_2$ other than $a$,
each remaining vertex of $\{b,c,u\}$ is adjacent to $a$. Any component not containing $a$ would therefore miss $b,c,u$ as
well. All edges of $G$ incident with its vertices lie in $H_2$,
so the same deletion would disconnect $G$. Since $G$ is
three-connected and $H_2-a$ is connected, $H_2$ is two-connected
and every two-cut contains $a$. Lemma~\ref{r2:three-mark-entrance},
with $(x,y,z)=(b,c,a)$, now gives a $K_{3,3}$ subdivision $H$ in
$H_1+bc$ using $bc$ such that $H-bc+ab+ac$ is nonplanar.

We shorten the branch paths of $H$ as in
Lemma~\ref{core:retained-subdivision}, keeping the corners and
the vertices $a,b,c,u$ wherever they occur. For a segment from
$u$ to $b$ or $c$, we keep the edge if the segment has length
one, and otherwise keep one interior vertex. Write $H'$ for the
shortened subdivision. Both $b$ and $c$ keep an interior neighbor
along an edge of $H'$ other than $bc$.

Each component of $D_1-\{b,c\}$ has a neighbor at $u$, so we can
contract its vertices outside $H'$ into $u$ or into vertices of
$H'$ in the same component. We keep a minimal forest joining the
retained interior vertices to $u$ and delete any edge $au$ created
by contraction, since $ap'u$ still supplies it.

A branch path of $H'$ contains $a$ exactly when the corresponding
branch path of $H$ does. By
Lemma~\ref{r2:position-of-third}, $H'-bc+ab+ac$ is therefore
nonplanar. In all four tests the path $ap'u$ supplies $au$, so
the two tests containing $bc$ have $H'$ as a minor. The other
two both contain $ac$. In $B+ab+ac+ad$ the edge $ab$ is present,
and in $B+ac+bd$ the path $aq'vdb$ supplies $ab$. Thus both have
$H'-bc+ab+ac$ as a minor. All four tests remain nonplanar, so by minimality the reduction
cannot be proper.

Suppose two segments of length two join $a$ to $b$ or $c$. If
they have the same ends, they form a cycle. If one ends at $b$
and the other at $c$, they form a cycle with $bc$. In either
case the only possible corners on the cycle are $a,b,c$. Since
every cycle in $H'$ contains at least four corners, there is at
most one such segment.

Construction~\ref{common:marked-construction} considers every
position of $u$ on $H'$ and the case $u\notin V(H')$. For each
choice it includes every minimal forest joining the retained
interior vertices to $u$, so it lists every possible $B_1$.
\end{proof}

\begin{proposition}\label{r2:a-cut-outcome}
If $X=\{a,u,v\}$, then $G$ is isomorphic to one of $\mathfrak T_{40},\mathfrak T_{51},\mathfrak T_{53},\mathfrak T_{56},\mathfrak T_{72},\mathfrak T_{73}$.
\end{proposition}

\begin{proof}
If $b,c,d$ lie in one component of $B-X$, Lemma~\ref{r2:a-cut-one}
gives $\mathfrak T_{40},\mathfrak T_{53},\mathfrak T_{56}$ or
$\mathfrak T_{73}$.

Suppose $b,c,d$ lie in two components $D_1,D_2$, with $h,k$ in
$D_1$ and the remaining vertex $g$ in $D_2$. Let $B_1,B_2$ be
their respective bridges on $X$. Lemma~\ref{r2:a-cut-second}
gives the forms of $D_2$ in every row of Lemma~\ref{r2:a-cut-table}.
In rows~(2), (4) and~(6), when $B_2+av+ag$ is planar,
Lemmas~\ref{r2:a-cut-row-two}, \ref{r2:a-cut-row-four}
and~\ref{r2:a-cut-row-six} give the forms of $B_1$, respectively.

In the remaining cases, $D_1$ can be reduced to a path $hmk$
with the same neighbors in $\{u,v\}$. Since $D_1$ is connected
and $hk$ is absent, some component of $D_1-\{h,k\}$ meets both
$h,k$. Contract this component to $m$. It has a neighbor at $u$
or $v$ by connectedness of the interior of $B$. Each other component
meets $h$ or $k$ and can be contracted into one that it meets.
Deleting $hk$ if it was created, and all edges from $a$ to $D_1$,
leaves the path $hmk$. To preserve the neighbors in $\{u,v\}$,
it suffices to keep one edge from $m$ to $u$ or $v$ and, if $D_1$
meets both, one edge from the other vertex to $h,m$ or $k$.

In rows~(1), (3) and~(5), we may also contract $D_2$ to $g$,
keeping its neighbors in $\{u,v\}$ and deleting $ag$ if it was
created. In the other rows we keep $B_2$ unchanged. The interior
of $B$ remains connected through $u,v,p',q'$, and each attachment
still has an interior neighbor. The reductions preserve the row
of the table. In rows~(1), (3) and~(5) the neighborhoods alone
make all four tests nonplanar, and in the other rows under
consideration $B_2+av+ag$ is nonplanar. The resulting minor of $G$
is therefore nontoroidal by Lemma~\ref{r2:four-tests}, so by
minimality $B$ already has one of these forms.

As in Proposition~\ref{r2:ab-cut-outcome}, we discard a form
if a single allowed deletion or contraction preserves its neighborhood
in $X$, the required connectivity and the planarity conditions in its row.
Combining the remaining forms of $B_1,B_2$ with the permitted
original edges, $p',q'$ with neighborhood $X$, and $\Psi_2$
gives $1{,}229$ candidates up to isomorphism. Each is toroidal, has a
$K_5$ minor, has a cut of at most two vertices, or has a proper
nontoroidal minor~\supp{E13}.

Finally, suppose $b,c,d$ lie in distinct components $D_b,D_c,D_d$.
For each $g\in\{b,c,d\}$, if $D_g=\{g\}$, its neighborhood in $B$
is $\{u\}$, $\{v\}$ or $\{u,v\}$. Otherwise
Lemma~\ref{r2:a-cut-single} gives $N_B(D_g)=\{s,t\}$ for some
$s,t\in X$. Either $D_g-g$ consists of two single vertices with
neighborhood $\{g,s,t\}$, or its bridge on $\{g,s,t\}$ has one
of the nine forms of that lemma. In either case the original edges
$gs,gt$ may be present, except that $ag$ is absent.

Together with $p',q'$ and $\Psi_2$, these choices give the family
in~\supp{E12}. The coverage check shows that each graph in this
family is a subgraph of a recorded toroidal graph or contains a
recorded nontoroidal graph. Since $G$ is nontoroidal, it contains
one of the latter, and minimality gives equality. Among the recorded
graphs in this family, only $\mathfrak T_{51}$ and $\mathfrak T_{72}$
remain. In both, each
$D_g-g$ is empty or connected. Together with the one-component
case, these give the six graphs in the statement.
\end{proof}

\subsubsection{A cut of the form \texorpdfstring{$\{a,c,w\}$}{\{a,c,w\}}}

This subsection shows that a cut of the form $X=\{a,c,w\}$ cannot occur. By Lemma~\ref{r2:cut-locations}, the only
components of $B-X$ containing neither $b$ nor $d$ are the single
vertices $p',q'$, with $N(p')=N(q')=X$. Since $X$ is the
neighborhood of $p'$, the edges $aw,cw$ are absent by
Lemma~\ref{lem:independent-neighborhood}(a). The interior of $B$
is connected and each of $b,d$ has an interior neighbor, so
$B-\{a,c\}$ is connected. Deleting $w$ gives $B-X$, and hence
each component of $B-X$ has a neighbor~at~$w$.

Also, $c$ has a neighbor in a component containing $b$ or $d$.
Otherwise $N(c)=\{p,q,p',q'\}$, contrary to
Lemma~\ref{r2:degree-four} applied to the twin pairs with
neighborhoods $\{a,b,c\}$ and $\{a,c,w\}$. If $b,d$ lie in separate components of $B-X$,
Lemma~\ref{r2:ac-cut-separate} gives the possible neighborhoods
of those components. Lemma~\ref{r2:ac-cut-together} rules out
the case where $b,d$ lie in the same component, and
Proposition~\ref{r2:ac-cut-outcome} excludes the remaining forms.

\begin{lemma}\label{r2:ac-cut-separate}
If $b,d$ lie in separate components $D_b,D_d$ of $B-X$, then
$N_B(D_b)=\{c,w\}$ and $N_B(D_d)$ is $\{w\}$~or~$\{a,w\}$.
\end{lemma}
\begin{proof}
Suppose $D_d$ has a neighbor at $c$. Take $p,q$ and $p',q'$ as
the twin pairs defining $\Psi_2$, with neighborhoods $\{a,b,c\}$
and $\{a,c,w\}$. The bridge containing $D_d$ also contains $r,z$,
since both are adjacent to $d$. It meets $a,b$ through $r,z$ and
meets $c,w$ through $D_d$. By Lemma~\ref{r2:unique-four}, $bw$ is
absent and every other bridge has a single interior vertex with
neighborhood $\{a,b,c\}$ or $\{a,c,w\}$.

Each component of $D_b-b$ is the interior of one of these other
bridges. Since $D_b$ is connected, its single vertex is adjacent
to $b$ and hence has neighborhood $\{a,b,c\}$. Thus no vertex of
$D_b-b$ is adjacent to $w$, and $bw$ is absent. This contradicts
that $D_b$ has a neighbor at $w$.

Hence $D_d$ misses $c$, so the neighbor of $c$ lies in $D_b$.
If $D_b$ also met $a$, then $D_b,p',q'$ would give three bridges
adjacent to all of $X$ in the planar graph $B$, contrary to
Lemma~\ref{core:assembly-general}(i). Since both components meet
$w$, the stated neighborhoods follow.
\end{proof}

\begin{lemma}\label{r2:ac-cut-together}
The vertices $b,d$ lie in distinct components of $B-X$.
\end{lemma}
\begin{proof}
Suppose $b,d$ lie in one component $D_1$, and let $B_1$ be its
bridge on $X$. Since the only other components are $p',q'$, $D_1$
has neighbors at $c,w$. It has none at $a$, since otherwise
$D_1,p',q'$ give three bridges meeting all of $X$ in the planar
graph $B$. Thus deleting $a$ and the three pairs of twins from $G$
leaves $B_1$.

Now $B_1-\{c,w\}=D_1$, and $B_1-\{b,c,d\}$ is the connected
interior of $B$ with the two leaves $p',q'$ removed. Since $bc,cd$
are absent, the neighbor of $c$ in $D_1$ is neither $b$ nor $d$.
Hence
\begin{equation}\label{r2:cross-connectivity}
B_1-\{c,w\}\text{ is connected},\qquad
B_1-\{b,c,d\}\text{ is connected},\qquad
N_{B_1}(c)\setminus\{b,c,d,w\}\ne\varnothing.
\end{equation}

Three of the four tests add $ab$ or $ad$, giving a third bridge
meeting all of $X$. In the remaining test, $B+ac+bd$, the bridge
containing $D_1$ still has attachments $c,w$. Thus
Lemma~\ref{r2:cut-locations} requires $B_1+bd+cw$ to be nonplanar.
The paths $brd$ and $cp'w$ give this graph as a minor of $G$, while
$B_1+cw$ is a planar minor of $B$. Hence $B_1+bd+cw$ contains a
$K_{3,3}$ subdivision using $bd$.

Put $I=V(B_1)\setminus\{b,c,d,w\}$. By~\eqref{r2:cross-connectivity},
each component of $B_1[I]$ has a neighbor at $w$ and at least one
of $b,d$. Some component has neighbors at both $b,d$, since $D_1$
is connected and $bd$ is absent. Some component also has a neighbor
at $c$. Keeping these last two components, which may coincide,
Lemma~\ref{r2:retained-components} shortens the subdivision with
at most two vertices of $I$ outside it. The retained components
keep their neighbors, so~\eqref{r2:cross-connectivity} still holds.
Thus the reduction is allowed. The subdivision keeps $B+ac+bd$
nonplanar, and the other three tests still have a third bridge
meeting all of $X$. By minimality, $B$ already has the reduced form.

Construction~\ref{common:marked-construction} therefore includes
every possible $B_1$, with added edges $bd,cw$, required edge $bd$,
and possible original edges $bw,dw$. Restoring $a$ and the three
pairs of twins gives the $5{,}819$ candidates in~\supp{E16}.
Each has a $K_5$ minor, a cut of at most two vertices, or a proper
nontoroidal minor, a contradiction.
\end{proof}

\begin{proposition}\label{r2:ac-cut-outcome}
The cut $X$ cannot have the form $\{a,c,w\}$.
\end{proposition}
\begin{proof}
By Lemmas~\ref{r2:ac-cut-separate} and~\ref{r2:ac-cut-together},
$b,d$ lie in separate components $D_b,D_d$, with
$N_B(D_b)=\{c,w\}$ and $N_B(D_d)$ equal to $\{w\}$ or $\{a,w\}$.
For $g\in\{b,d\}$, write $B_g=B[V(D_g)\cup N_B(D_g)]$.
When $N_B(D_g)=\{s,t\}$, the interior $D_g-g$ is nonempty,
since $bc,ad$ are absent. By Lemma~\ref{r2:one-attachment}(ii),
each of its components has neighborhood $\{g,s,t\}$. If there are
two, both are single vertices and the triangle completion of $B_g$
is planar.

Suppose first that $N_B(D_d)=\{w\}$. Then $D_d=\{d\}$ with the
edge $dw$ by Lemma~\ref{r2:one-attachment}(i). In $B+ac+bd$ the
path $bdw$ joins $b$ to $w$. In $B+ad+bc$ the bridge containing
$d$ is the path $adw$, which is planar after adding $aw$.
Lemma~\ref{r2:cut-locations} therefore requires both
\[
B_b+cw+bw \qquad\hbox{and}\qquad B_b+cw+bc
\]
to be nonplanar. The path $cp'w$ gives $B_b+cw$ as a planar minor
of $B$, so $bw$ is absent. Thus $b,c,w$ are independent, and
$D_b-b$ is connected since the triangle completion is nonplanar.

The triangle completion of $B_b$ is three-connected by
Lemma~\ref{lem:triangle-completion} and is a minor of $G$ by the
paths $bpc$, $brdw$ and $cp'w$. Lemma~\ref{core:three-root-setup}
now gives a three-cut $Y$ of this completion that is a twin triple.
Since no two twin triples are disjoint, $Y$ meets both $\{a,b,d\}$
and $\{a,c,w\}$. Neither $a$ nor $d$ belongs to $B_b$, so $Y$
contains $b$ and some $s\in\{c,w\}$. By part~\textup{(vi)} of the
same lemma, adding $bs$ to the planar graph $B_b+cw$ preserves
planarity, contrary to the two tests above.

It remains to consider $N_B(D_d)=\{a,w\}$. In $B+ac+bd$, the
edge $bd$ joins $D_b,D_d$, giving a third bridge meeting all of
$X$. The edge $ab$ does the same for the bridge containing $D_b$
in the two tests that contain it. In the remaining test $B+ad+bc$,
Lemma~\ref{r2:cut-locations} requires at least one of
\[
B_b+cw+bc \qquad\hbox{and}\qquad B_d+aw+ad
\]
to be nonplanar. Both are minors of $G$, using $cp'w,bpc$ for the
first and $ap'w,ard$ for the second. Without $bc$ and $ad$,
respectively, they are planar minors of $B$.

Choose a nonplanar one of these two completions. Its bridge has
connected interior, and the completion contains a $K_{3,3}$
subdivision using $bc$ or $ad$, respectively.
By Lemma~\ref{r2:retained-components}, we can shorten the
subdivision while keeping the interior of its bridge connected
and adjacent to all three attachments. In the other bridge, each
interior component has neighbors at all three attachments. We
contract one such component to a vertex and keep only its three
edges to the attachments. The two interiors are then connected
through $w$, and $a,b,c,d$ each keep an interior neighbor. The subdivision keeps $B+ad+bc$ nonplanar, and the other
three tests still have a third bridge meeting all of $X$. Thus the
reduction is allowed and minimality leaves $B$ unchanged.

Construction~\ref{common:marked-construction} gives the possible
forms on either side. A segment with ends $b,w$ in $B_b$, or $d,w$
in $B_d$, may be an original edge or a path with one interior
vertex. The two choices of nonplanar completion give $156$
candidates~\supp{E15}. Each has a $K_5$ minor, a cut of at most two
vertices, or a proper nontoroidal minor. This excludes the remaining
case.
\end{proof}

\subsection{Bridges with three attachments: number and planar completions}\label{r2:no-four}

For the rest of this section, we suppose no choice of twin pairs gives a bridge with four attachments. Each component of $G-V(\Psi_2)$ then gives a bridge with three attachments. If all the triangle completions are planar, Lemma~\ref{r2:component-count} gives the three obstructions $\mathfrak T_9,\mathfrak T_{10},\mathfrak T_{16}$. If at least one triangle completion is nonplanar, its component has more than one vertex, and the same lemma leaves at most four components outside $\Psi_2$. We show that each of their bridges has a planar two-edge completion. Section~\ref{r2:planar-bridges} uses these completions to determine the possible bridges, and Proposition~\ref{r2:chain-outcome} gives $G\cong\mathfrak T_{15}$.

The triangle completion of each such bridge $B$ is three-connected by Lemma~\ref{lem:triangle-completion}. The cycle $apcqbzdra$ in $\Psi_2$ meets $B$ precisely at its three attachments. Shortening the three arcs between them to single edges gives the triangle completion as a minor of $G$, so it has no $K_5$ minor.

\begin{lemma}\label{r2:component-count}
If the bridge of every component outside $\Psi_2$ has planar triangle completion, then $G$ is isomorphic to one of $\mathfrak T_9,\mathfrak T_{10},\mathfrak T_{16}$. If some bridge has nonplanar triangle completion, then its component has more than one vertex and there are at most four components outside $\Psi_2$.
\end{lemma}
\begin{proof}
Put $A=\{a,b,c,d\}$. Each component of $G-A$ has exactly three neighbors in $A$ by three-connectivity and our assumption on the bridges. Four components are the single vertices $p,q,r,z$. The others are the components outside $\Psi_2$.

If $G-A$ has at least nine components, retain nine, including $p,q,r,z$, and contract each to a vertex. Delete $cd$ if present. The resulting bipartite minor has $13$ vertices and $27>2\cdot13$ edges, so it is nontoroidal by Euler's formula. By minimality it is $G$. Thus $G-A$ has at most nine components, with equality only if each is a single~vertex.

Suppose first that every triangle completion is planar. By Lemma~\ref{lem:flat-replacement}, each component outside $\Psi_2$ is a single vertex. There are at most five such vertices, each with neighborhood $\{a,b,c\}$, $\{a,b,d\}$, $\{a,c,d\}$ or $\{b,c,d\}$. There are $\binom{9}{4}=126$ choices for how many have each neighborhood. With $cd$ present or absent, this gives the $252$ choices checked in~\supp{E17}. Only $\mathfrak T_9,\mathfrak T_{10},\mathfrak T_{16}$ remain, all with $cd$ absent.

Now suppose some bridge has nonplanar triangle completion. Its component has more than one vertex, since a single vertex gives $K_4$. Thus $G-A$ has at most eight components, of which at most four lie outside $\Psi_2$.
\end{proof}

It remains to consider the case where at least one bridge has nonplanar triangle completion. For a bridge $B$ with attachments $x,y,t$, we show that the attachments can always be labeled so that $B+xy+xt$ is planar, even when adding the third edge $yt$ makes the completion nonplanar. 

\begin{lemma}\label{r2:cellular-center}
Suppose a graph consists of $\Psi_2$, bridges with three attachments among $a,b,c,d$, and possibly $cd$, with $cd$ or a bridge joining $c,d$. If it has a torus embedding in which a closed regular neighborhood of $\Psi_2$ has genus one, then each added bridge $B$ has a planar completion $B+xy+xt$ for some labeling of its attachments.
\end{lemma}
\begin{proof}
By Lemma~\ref{lem:torus-complement}, every complementary region of the closed regular neighborhood of $\Psi_2$ is a disk. The edge $cd$ or the bridge joining $c,d$ lies in one such disk, so some boundary walk contains both $c$ and $d$.

Let $B$ be an added bridge. Its interior is connected, so it lies in one disk face. Among the genus-one rotations with a walk containing $c,d$, no two labels alternate on a face meeting three or more attachments~\supp{R4}. Restrict the boundary sequence of the face containing $B$ to its three attachments. Combining consecutive equal appearances leaves $xyt$ or $xyxt$ up to relabeling. In both cases $x$ is consecutive with both $y$ and $t$, so Lemma~\ref{cor:boundary-completions} gives a plane drawing of $B+xy+xt$.
\end{proof}

\begin{lemma}\label{r2:two-edge-exists}
The attachments of every $\Psi_2$-bridge $B$ with nonempty interior can be labeled $x,y,t$ so that $B+xy+xt$ is planar.
\end{lemma}
\begin{proof}
Suppose some bridge $B$ has no planar two-edge completion. If $B$ attaches at $a,b$ and one of $c,d$, we may assume its third attachment is $c$. By three-connectivity there is a $c$--$d$ path in $G-\{a,b\}$. The path meets $B$ only at $c$, its only remaining attachment.
It also avoids the four twins, since each has degree one in
$G-\{a,b\}$. We contract the path into $c$ and keep only $B$
and the edges at the four twins. Since the contraction identifies
$d$ with $c$, all four now have neighborhood $\{a,b,c\}$. This proper minor is toroidal, so Lemma~\ref{lem:four-twins} gives a planar two-edge completion of $B$, a contradiction.

Thus $B$ attaches at $\{a,c,d\}$ or $\{b,c,d\}$. Some other component outside $\Psi_2$ has a different attachment triple, since otherwise the vertex of $\{a,b,c,d\}$ that $B$ does not attach to has degree four, contrary to Lemma~\ref{r2:degree-four}. Contract this component to a vertex $v$ and retain its three incident edges.

Write $x,y,t$ for the attachments of $B$. We first show that if $B$ is planar, some one-edge completion is nonplanar. Suppose all three are planar. The triangle completion is three-connected, nonplanar and has no $K_5$ minor, so Lemma~\ref{core:three-root-setup} gives a three-cut $X$ and two interior vertices with neighborhood $X$. By hypothesis the twin triple $X$ meets $\{a,b,c\}$. Every
vertex they share is an attachment of $B$, so we may label
the attachments so that $x\in X$.

In $B+xy+xt$, the edges $xy,xt$ must lie in the same bridge on $X$. Otherwise each lies in the same bridge as in its planar one-edge completion, so Lemma~\ref{core:three-root-setup}(vi) makes $B+xy+xt$ planar. Thus $y,t$ lie in one component $D$ of $B-X$. The only edge of the triangle that remains after deleting $X$ is
$yt$, and both its ends lie in $D$. Thus $D$ is still a component,
and by three-connectivity its neighborhood in the triangle
completion is $X$. Deleting $yt$ and $xt$ leaves this neighborhood
unchanged, since $xy$ still joins $D$ to $x$. Hence the bridge of
$D$ in the planar graph $B+xy$ meets all of $X$, contrary to
part~(iv) of the same lemma.

If $B$ is nonplanar, take a $K_{3,3}$ subdivision $H\subseteq B$ and put $E=\varnothing$. If $B$ is planar, label the attachments so that $B+xy$ is nonplanar and put $E=\{xy\}$. Since $B+xt+yt$ is also nonplanar, Lemma~\ref{r2:three-mark-entrance} gives a $K_{3,3}$ subdivision $H\subseteq B+xy$ using $xy$ for which $H-xy+xt+yt$ is nonplanar.

By Lemma~\ref{core:retained-subdivision}, $B$ has an allowed minor $B'$ with at most eight interior vertices and a shortened subdivision $H'\subseteq B'+E$. If $E=\varnothing$, all three two-edge completions of $B'$ contain $H'$. If $E=\{xy\}$, the two containing $xy$ contain $H'$, and the third contains $H'-xy+xt+yt$. This last graph is homeomorphic to $H-xy+xt+yt$, since the shortening preserves $xy$, the corners and the positions of the attachments on the branch paths. Thus all three completions of $B'$ are nonplanar.

Retain $\Psi_2$, $B'$ and the star at $v$. Suppose this minor is toroidal. Since $B'$ joins $c,d$, Lemma~\ref{r2:cellular-center} rules out a genus-one regular neighborhood of $\Psi_2$. In genus zero the boundary sequences are $ab,ab,ac,bc,ad,bd$~\supp{R4}. Since $B'$ and the star each have three attachments, Lemma~\ref{lem:torus-complement} puts them in the same annulus. Together they meet all four vertices, so the boundary sequences are $ac\mid bd$ or $ad\mid bc$. By Lemma~\ref{lem:annulus-embedding}, adding $ac,bd$ or $ad,bc$, respectively, gives a plane drawing of $B'$ and the star. The star and the two added edges form a connected graph meeting $B'$ in its three attachments. Contract a spanning tree of this graph to those three vertices, keeping them distinct. This gives a planar two-edge completion of $B'$, a contradiction.

The retained minor is therefore nontoroidal, so by minimality it is $G$. For each of the two possible attachment triples of $B'$, Construction~\ref{common:marked-construction} gives $869$ nonplanar and $120$ planar bridge forms with all three two-edge completions nonplanar. There are three remaining choices for $N(v)$. Each of the resulting $2\cdot3(869+120)=5{,}934$ graphs has a cut of at most two vertices, a $K_5$ minor, or a proper nontoroidal minor~\supp{E18}.
\end{proof}

\subsection{Possible forms and combinations of bridges with three attachments}\label{r2:planar-bridges}

By Lemma~\ref{r2:two-edge-exists}, one may label the attachments
of each bridge $B$ by $x,y,t$ so that $B+xy+xt$ is
planar. Proposition~\ref{prop:three-root-blocks} then makes
$B-x$ a chain of edges and four-cycles. We show that
minimality leaves eight possible bridge forms and
then determine their possible combinations. 

The following lemma lets us extend a drawing of a reduced bridge
to a drawing of the original bridge.

\begin{lemma}\label{r2:cycle-restoration}
Let $B$ satisfy the hypotheses of Proposition~\ref{prop:three-root-blocks}, with attachments $x,y,t$, and let $C$ be a four-cycle of $B-x$ containing neither $y$ nor $t$. Then $B'$ is obtained from $B$ by one of the following operations.
\begin{enumerate}[label=\textup{(\roman*)}]
\item Contract the path of one or two edges joining $C$ to the next four-cycle. If it has two edges, first delete the edge from its middle vertex to $x$.
\item For an end path $yws$ or $tws$ with $s\in V(C)$, delete $wx$ and contract $ws$.
\end{enumerate}
Suppose $B'$ is drawn in a plane graph $K$, possibly with its attachments split into copies. Assume that a connected subgraph $H$ of $K$ meets $B'$ precisely in the attachment copies and that no edge of $K$ outside $B'$ is incident with an interior vertex of $B'$. Then $B'$ can be replaced by $B$ in this drawing, using the same attachment copies and keeping the rest of $K$ fixed.
\end{lemma}
\begin{proof}
Let $s$ be the vertex created by the contraction, and write $C=spuqs$. Since $H$ is connected and disjoint from $C$, it lies on one side of $C$. Every edge at $s$ outside $C$ starts a path to $H$ avoiding the rest of $C$, through the next four-cycle in (i) or directly to $y$ or $t$ in (ii). Thus all these edges leave $C$ on the side containing $H$, so $sp,sq$ are consecutive around $s$. Undo the contraction in a small disk about $s$, keeping $s$ on $C$ and moving the edges outside $C$ to the other end of the restored path.

If an edge to $x$ was deleted, its other end $w$ is the neighbor of $s$ on the restored path. Let $x_p$ be the copy of $x$ adjacent to $p$. Both $sw$ and $px_p$ lie on the side of $C$ containing $H$. Since $s$ and $p$ have degree three, the face on this side of $sp$ has $wspx_p$ on its boundary. Draw $wx_p$ in this face close to $wspx_p$. This restores $B$ and leaves the rest of the drawing fixed.
\end{proof}

We now apply Lemma~\ref{r2:cycle-restoration} in the disk
and annulus tests to bound the number of four-cycles and
the lengths of the paths joining them.

\begin{figure}[htbp]
\centering
\begingroup
\scaledfigure{1.295}{%
\tikzsetnextfilename{psi2-eight-bridge-forms-row}%
\begin{tikzpicture}[
  xscale=.5,yscale=.8,
  font=\small,text=figureInk,
  vertex/.style={circle,fill=figureInk,inner sep=1.5pt},
  root/.style={circle,fill=figurePath,inner sep=1.9pt},
  edge/.style={draw=figureInk!75,line width=.75pt,
    line cap=round,line join=round},
  xedge/.style={edge,draw=figurePath,line width=.8pt},
  panel/.style={align=center,inner sep=0pt}
]

\begin{scope}
  \coordinate (x) at (1.5,1.65);
  \coordinate (y) at (0,0);
  \coordinate (z) at (3,0);
  \coordinate (w) at (1.5,0);
  \draw[edge] (y)--(w)--(z);
  \draw[xedge] (x)--(w);
  \node[vertex] at (w) {};
  \foreach \v/\lab in {x/x,y/y,z/t}
    \node[root,label=above:{$\lab$}] at (\v) {};
  \node[panel] at (1.5,-1.12) {star};
\end{scope}

\foreach \leftlen/\rightlen/\dx/\dy in {
  0/1/4.2/0,
  0/2/8.4/0,
  1/0/12.6/0,
  2/0/16.8/0,
  1/1/21/0
}{
  \begin{scope}[xshift=\dx cm,yshift=\dy cm]
    \coordinate (x) at (1.5,1.65);
    \coordinate (y) at (0,0);
    \coordinate (z) at (3,0);
    \pgfmathsetmacro{\sx}{.75*\leftlen}
    \pgfmathsetmacro{\tx}{3-.75*\rightlen}
    \pgfmathsetmacro{\mx}{(\sx+\tx)/2}
    \coordinate (s) at (\sx,0);
    \coordinate (t) at (\tx,0);
    \coordinate (p) at (\mx,.58);
    \coordinate (q) at (\mx,-.58);

    \draw[edge] (s)--(p)--(t)--(q)--cycle;
    \draw[xedge] (x)--(p);
    \ifnum\rightlen=0\relax
      \draw[xedge] (x) .. controls (3.95,1.65) and (3.95,-1.15) .. (q);
    \else
      \draw[xedge] (x) .. controls (-.95,1.65) and (-.95,-1.15) .. (q);
    \fi

    \ifnum\leftlen>0\relax
      \ifnum\leftlen=1\relax
        \draw[edge] (y)--(s);
      \else
        \coordinate (u) at (.75,0);
        \draw[edge] (y)--(u)--(s);
        \draw[xedge] (x)--(u);
        \node[vertex] at (u) {};
      \fi
      \node[vertex] at (s) {};
    \fi
    \ifnum\rightlen>0\relax
      \ifnum\rightlen=1\relax
        \draw[edge] (t)--(z);
      \else
        \coordinate (v) at (2.25,0);
        \draw[edge] (t)--(v)--(z);
        \draw[xedge] (x)--(v);
        \node[vertex] at (v) {};
      \fi
      \node[vertex] at (t) {};
    \fi

    \foreach \v in {p,q} \node[vertex] at (\v) {};
    \foreach \v/\lab in {x/x,y/y,z/t}
      \node[root,label=above:{$\lab$}] at (\v) {};
    \node[panel] at (1.5,-1.12) {$(\leftlen,\rightlen)$};
  \end{scope}
}

\foreach \len/\dx in {1/25.2,2/29.4}{
  \begin{scope}[xshift=\dx cm]
    \coordinate (x) at (1.5,1.65);
    \coordinate (y) at (0,0);
    \coordinate (z) at (3,0);
    \coordinate (s) at (1.1,0);
    \coordinate (t) at (1.9,0);
    \coordinate (p1) at (.55,.58);
    \coordinate (q1) at (.55,-.58);
    \coordinate (p2) at (2.45,.58);
    \coordinate (q2) at (2.45,-.58);

    \draw[edge] (y)--(p1)--(s)--(q1)--cycle;
    \draw[edge] (t)--(p2)--(z)--(q2)--cycle;
    \draw[xedge] (x)--(p1) (x)--(p2);
    \draw[xedge] (x) .. controls (-.95,1.65) and (-.95,-1.15) .. (q1);
    \draw[xedge] (x) .. controls (3.95,1.65) and (3.95,-1.15) .. (q2);
    \ifnum\len=1\relax
      \draw[edge] (s)--(t);
    \else
      \coordinate (w) at (1.5,0);
      \draw[edge] (s)--(w)--(t);
      \draw[xedge] (x)--(w);
      \node[vertex] at (w) {};
    \fi

    \foreach \v in {s,t,p1,q1,p2,q2}
      \node[vertex] at (\v) {};
    \foreach \v/\lab in {x/x,y/y,z/t}
      \node[root,label=above:{$\lab$}] at (\v) {};
    \node[panel] at (1.5,-1.12) {$(0,0),\;\ell=\len$};
  \end{scope}
}
\end{tikzpicture}%
}
\endgroup
\caption{The eight bridge forms in Lemma~\ref{r2:eight-bridge-forms}. The attachments and the edges
incident with $x$ are blue. In the last two drawings, $\ell$ is the length
of the path joining the two four-cycles.}
\label{fig:psi2-eight-bridge-forms}
\end{figure}
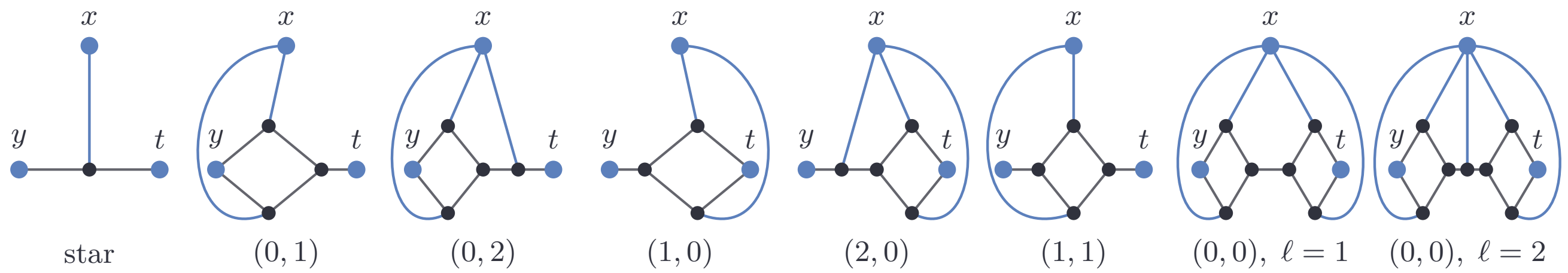

\begin{lemma}\label{r2:eight-bridge-forms}
Every $\Psi_2$-bridge $B$ with nonempty interior has one of the eight forms in Figure~\ref{fig:psi2-eight-bridge-forms} and hence has at most seven interior vertices. For the forms with a four-cycle, the end lengths are the distances in $B-x$ from $y$ and $t$, respectively, to the nearest four-cycle.
\end{lemma}
\begin{proof}
Suppose one of the shortenings in Lemma~\ref{r2:cycle-restoration} gives $B'$. It preserves the three attachments and gives a proper minor of $G$, so take a torus embedding of this minor. There is another bridge besides $B'$, possibly the edge $cd$. Indeed, $\deg(a),\deg(b)>4$ by Lemma~\ref{r2:degree-four}, and $\deg(c),\deg(d)\ge3$, so each of $a,b,c,d$ is incident with an edge outside $\Psi_2$.

If $B'$ lies in a disk face, let $D$ be the disk representing that face. Otherwise a regular neighborhood of $\Psi_2$ has genus zero by Lemma~\ref{lem:torus-complement}. No boundary walk then contains three attachments or both $c,d$~\supp{R4}, so every bridge meets both boundary components of the unique annulus. Cut the annulus along $cd$ if present, or along a path through another bridge joining the two boundary components. This gives a disk $D$ without cutting the interior of $B'$.

In either case, take the boundary cycle of $D$ as $H$ in Lemma~\ref{r2:cycle-restoration}. It contains all attachment copies of $B'$ and avoids its interior. Thus we can restore $B$ in $D$ and identify the boundary as before, giving a torus embedding of $G$. Hence neither shortening is possible.

By Proposition~\ref{prop:three-root-blocks}, the blocks of $B-x$ form a chain of edges and four-cycles with no three consecutive edge blocks. If two four-cycles shared a vertex $v$, the vertices $u,w$ opposite $v$ on the two cycles would give twin triples $\{x,u,v\}$ and $\{x,v,w\}$. But $v$ has degree four in $G$, contrary to Lemma~\ref{r2:degree-four}. Thus consecutive four-cycles are joined by a path of length one or two, and each end length is at most two.

If there is more than one four-cycle, any cycle avoiding both $y,t$ admits shortening~(i) toward a neighboring cycle. Thus there are exactly two, containing $y,t$ respectively, with joining path of length one or two. With one four-cycle, end lengths $(0,0)$ would disconnect the interior of $B$. If both end lengths are positive, neither can be two by shortening~(ii). This leaves $(0,1),(0,2),(1,0),(2,0),(1,1)$. With no four-cycle, $B$ is the three-edge star by Proposition~\ref{prop:three-root-blocks}(v). These are the eight forms. The largest has two four-cycles joined by a path of length two and has seven interior vertices.
\end{proof}

\begin{proposition}\label{r2:chain-outcome}
If some component outside $\Psi_2$ is not a singleton,
then $G\cong\mathfrak T_{15}$.
\end{proposition}
\begin{proof}
Suppose $\{s,t,x_i\}$ is a twin triple for $i=1,2,3$, with the $x_i$ distinct. In $G-\{s,t\}$, two chosen twins on each triple are leaves at $x_i$. Delete these six leaves from the connected graph $G-\{s,t\}$ and take a minimal tree joining $x_1,x_2,x_3$. Relabel so that $x_1,x_2$ are leaves of the tree. Choose the twin pairs on the first two triples for $\Psi_2$. The component outside $\Psi_2$ containing $x_3$ has neighbors at $x_1,x_2$ along the tree and at $s,t$ through the twins on the third triple. This gives a bridge with four attachments, a contradiction.

By Lemmas~\ref{r2:component-count} and~\ref{r2:eight-bridge-forms}, we add to $\Psi_2$ one to four of the eight bridge forms, with disjoint interiors and at least one form other than a star. On each attachment triple the eight forms give $19$ labeled bridges up to isomorphisms fixing the attachments. Thus there are $76$ labeled bridges, four of them stars. Choosing at most four with repetition, excluding choices consisting entirely of stars and allowing $cd$ to be present or absent, gives $2\bigl[\binom{80}{4}-\binom{8}{4}\bigr]=3{,}163{,}020$~combinations.

We discard a combination if two twin triples are disjoint, three distinct twin triples share two vertices, or some choice of twin pairs gives a bridge with four attachments. These exclusions also rule out any larger combination containing it. Indeed, only $a,b,c,d$ gain neighbors. The vertices $a,b$ already have degree at least four, and only $a,b$ could be twins of $c$ or $d$, since a twin of $c$ must be adjacent to $p,q$ and a twin of $d$ to $r,z$. Thus the twin pairs keep their neighborhoods, which must be independent in $G$ by Lemma~\ref{lem:independent-neighborhood}(a). A component adjacent to all four attachments remains in a component adjacent to all four.

Any combination containing a smaller nontoroidal one has a proper nontoroidal minor. Of the $3{,}163{,}020$ combinations, \supp{E19} records $2{,}174$ individually. Every other combination contains a smaller one that is nontoroidal or violates one of the conditions above. Of the $2{,}174$ recorded graphs, all but four are toroidal, have a cut of at most two vertices, or have a proper nontoroidal minor. The four remaining graphs are isomorphic to $\mathfrak T_{15}$.
\end{proof}

\begin{proof}[Proof of Theorem~\ref{thm:two-common}]
If some choice of twin pairs gives a bridge with four attachments, then either $G\cong\mathfrak T_2$ or $G$ is one of the graphs found in Sections~\ref{r2:four-neighbors} and~\ref{r2:planar-four}. Of these, $\mathfrak T_{48}$ and $\mathfrak T_{66}$ are excluded because they have disjoint twin triples.

If no choice gives such a bridge, Lemma~\ref{r2:component-count} and Proposition~\ref{r2:chain-outcome} give $\mathfrak T_9,\mathfrak T_{10},\mathfrak T_{15},\mathfrak T_{16}$. The converse follows from Proposition~\ref{prop:membership}.
\end{proof}

\section{Two twin triples with one common vertex}\label{sec:one-common}

\begin{figure}[htbp]
\centering
\begingroup
\scaledfigure{1}{%
\tikzsetnextfilename{one-common-neighborhoods}%
\begin{tikzpicture}[
  font=\small,text=figureInk,
  vertex/.style={circle,fill=figureInk,inner sep=1.45pt},
  twinA/.style={circle,fill=figureTwinA,inner sep=1.7pt},
  twinB/.style={circle,fill=figureTwinB,inner sep=1.7pt},
  edgeA/.style={draw=figureTwinA,line width=.65pt,
    preaction={draw=white,line width=1.6pt}},
  edgeB/.style={draw=figureTwinB,line width=.65pt,
    preaction={draw=white,line width=1.6pt}},
  lab/.style={anchor=north,inner sep=0pt}
]

\begin{scope}[yshift=-.1cm]
  \draw[draw=figureInk!45,fill=figureInk!3,line width=.6pt,
    rounded corners=4pt]
    (.05,-.92) rectangle (4.95,.2);
  \fill[figureTwinA!9,rounded corners=3pt]
    (.28,-.43) rectangle (2.83,.15);
  \fill[figureTwinB!9,fill opacity=.6,rounded corners=3pt]
    (2.17,-.46) rectangle (4.72,.18);
  \draw[figureTwinA,line width=.6pt,rounded corners=3pt]
    (.28,-.43) rectangle (2.83,.15);
  \draw[figureTwinB,line width=.6pt,rounded corners=3pt]
    (2.17,-.46) rectangle (4.72,.18);
  \node at (2.5,-.69) {rest of $G$};

  \coordinate (b) at (.65,0);
  \coordinate (c) at (1.575,0);
  \coordinate (a) at (2.5,0);
  \coordinate (d) at (3.425,0);
  \coordinate (e) at (4.35,0);
  \coordinate (p1) at (.85,1.05);
  \coordinate (p2) at (1.65,1.05);
  \coordinate (q1) at (3.35,1.05);
  \coordinate (q2) at (4.15,1.05);

  \foreach \i in {1,2}
    \foreach \u in {a,b,c} \draw[edgeA] (p\i)--(\u);
  \foreach \i in {1,2}
    \foreach \u in {a,d,e} \draw[edgeB] (q\i)--(\u);
  \foreach \i/\firsttwin/\secondtwin in {1/p/r,2/q/z}{
    \node[twinA,label=above:{$\firsttwin$}] at (p\i) {};
    \node[twinB,label=above:{$\secondtwin$}] at (q\i) {};
  }
  \foreach \v in {b,c,a,d,e}{
    \node[vertex] at (\v) {};
    \node[lab] at ([yshift=-.1cm]\v) {$\v$};
  }
\end{scope}

\begin{scope}[xshift=8.85cm]
  \coordinate (a) at (0,1.45);
  \coordinate (b) at (-1.05,.05);
  \coordinate (c) at (-1.55,-1.35);
  \coordinate (d) at (1.55,-1.35);
  \coordinate (e) at (1.05,.05);
  \coordinate (p1) at (-2,.5);
  \coordinate (p2) at (-2,-.65);
  \coordinate (q1) at (2,.5);
  \coordinate (q2) at (2,-.65);

  \foreach \u in {b,c}
    \foreach \v in {d,e}
      \draw[draw=figureInk!45,line width=.55pt,dotted,line cap=round]
        (\u)--(\v);
  \foreach \i in {1,2}
    \foreach \u in {a,b,c} \draw[edgeA] (p\i)--(\u);
  \foreach \i in {1,2}
    \foreach \u in {a,d,e} \draw[edgeB] (q\i)--(\u);

  \node[vertex,label=above:{$a$}] at (a) {};
  \node[vertex,label=above right:{$b$}] at (b) {};
  \node[vertex,label=below:{$c$}] at (c) {};
  \node[vertex,label=below:{$d$}] at (d) {};
  \node[vertex,label=above left:{$e$}] at (e) {};
  \node[twinA,label=left:{$p$}] at (p1) {};
  \node[twinA,label=left:{$q$}] at (p2) {};
  \node[twinB,label=right:{$r$}] at (q1) {};
  \node[twinB,label=right:{$z$}] at (q2) {};
\end{scope}
\end{tikzpicture}%
}
\endgroup
\caption{Left, two twin triples meeting only at $a$.
Right, the graph $\Psi_3$, with the possible cross edges dotted.}
\label{fig:one-common-neighborhoods}
\end{figure}
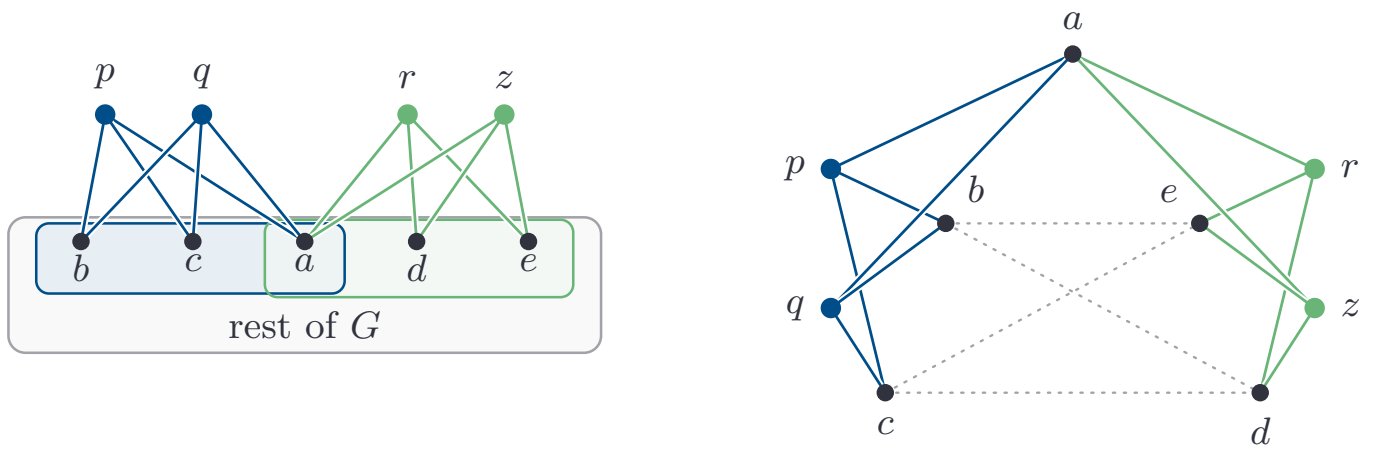

 Finally, in this section, we consider the case where $G$ is a three-connected obstruction with at least two distinct twin triples, every two of which have exactly one vertex in common. Choose twins $p,q$ with neighborhood $A_1=\{a,b,c\}$ and twins $r,z$ with neighborhood $A_2=\{a,d,e\}$. Write $J_1,J_2$ for their twin subgraphs and put $\Psi_3=J_1\cup J_2$, as in Figure~\ref{fig:one-common-neighborhoods}. Let $n_c,n_e$ count the vertices with neighborhoods $A_1,A_2$, respectively.

\begin{theorem}\label{thm:one-common}
A graph is a three-connected obstruction with at least two distinct twin triples, every two of which have exactly one vertex in common, if and only if it is isomorphic to one of
\[
\mathfrak{T}_{7},\mathfrak{T}_{14},\mathfrak{T}_{20},\mathfrak{T}_{23},
\mathfrak{T}_{29},\mathfrak{T}_{33},\mathfrak{T}_{35},\mathfrak{T}_{37},
\mathfrak{T}_{41},\mathfrak{T}_{49},\mathfrak{T}_{54},\mathfrak{T}_{55},
\mathfrak{T}_{59},\mathfrak{T}_{70}.
\]
\end{theorem}

Every component of $G-V(\Psi_3)$ has at least three neighbors, all in $A_1\cup A_2$, since $G$ is three-connected and the twins have degree three. The sets $A_1,A_2$ are independent by Lemma~\ref{lem:independent-neighborhood}(a), so the $\Psi_3$-bridges with empty interior are edges joining $\{b,c\}$ to $\{d,e\}$, which we call \emph{cross edges}. We first derive the embedding tests, then consider bridges with five attachments, four attachments excluding $a$, four attachments including $a$, and finally three attachments.

\begin{lemma}\label{r1:triangle-minor}
Let $U$ be a union of $\Psi_3$-bridges. The graph obtained from $U$ by adding the vertices $a,b,c,d,e$ and the edges of the triangles $abc$ and $ade$ is a minor of $G$, and so has no $K_5$ minor.
\end{lemma}
\begin{proof}
Contract $pb,qc,rd,ze$ and delete the bridges outside $U$.
\end{proof}

\subsection{Faces and embedding tests of \texorpdfstring{$\Psi_3$}{Psi3}}\label{r1:two-tests}

\begin{lemma}\label{r1:exterior-faces}
In a torus embedding of $\Psi_3$, let $\mathcal N$ be a closed regular
neighborhood of $\Psi_3$. If $\mathcal N$ has genus one and a face meets
$b,c,d,e$, that face is a disk with attachment sequence
$a,b,a,d,a,c,a,e$, up to cyclic shift and reversal, and the other
two faces have attachment sequences $bc,de$. If $\mathcal N$ has genus zero, any face
meeting $b,c,d,e$ is an annulus with attachment sequences $bc\mid de$, and
the other three faces are disks with attachment sequences $axay$, $ax'$ and
$ay'$, where $\{x,x'\}=\{b,c\}$ and $\{y,y'\}=\{d,e\}$.
\end{lemma}

\begin{proof}
For each of the $96$ rotation systems of $\Psi_3$, the table in~\supp{R6} gives the attachments in order along each facial walk. No rotation system of genus zero has a face containing all four of $b,c,d,e$ in its boundary walk. Among the genus-zero rotations, the only pair of attachment sequences for two facial walks that together meet $b,c,d,e$ is $bc,de$. Thus, by Lemma~\ref{lem:torus-complement}, a face meeting all four in a torus embedding with this rotation must be an annulus with attachment sequences $bc\mid de$. The other three facial walks of such a rotation have attachment sequences $axay$, $ax'$ and $ay'$, and bound disks by the same lemma. Of the genus-one rotations, the only facial attachment sequence containing all four is $abadacae$, up to cyclic shift and reversal. The other two faces have attachment~sequences~$bc,de$.
\end{proof}

\begin{figure}[htbp]
\centering
\begingroup
\scaledfigure{.85}{\tikzsetnextfilename{psi3-torus-minor-r1}%
\begin{tikzpicture}[
  scale=1.4,
  font=\small,
  vertex/.style={circle,fill=figureInk,inner sep=1.6pt},
  edge/.style={draw=figureInk!75,line width=.8pt},
  side/.style={draw=figureInk!35,line width=.6pt}
]

\begin{scope}
  \fill[faceRepeated]
    (4,0)--(2.8,0)--(2,.75)--(2,1.5)--(3,1.5)
    --(3.5,1.15)--(4,1.15)--cycle;
  \fill[faceRepeated]
    (4,3)--(4,1.85)--(3.5,1.85)--(3,1.5)--(2,1.5)
    --(2,2.25)--(2.8,2.6)--(2.8,3)--cycle;
  \fill[faceRepeated]
    (0,3)--(1.2,3)--(1.2,2.6)--(2,2.25)--(2,1.5)
    --(1,1.5)--(0,1.85)--cycle;
  \fill[faceRepeated]
    (0,0)--(0,1.15)--(1,1.5)--(2,1.5)--(2,.75)
    --(1.2,0)--cycle;

  \draw[side] (0,0) rectangle (4,3);
  \draw[side,-{Stealth[length=5pt]}] (0,.25)--(0,.65);
  \draw[side,-{Stealth[length=5pt]}] (4,.25)--(4,.65);
  \draw[side,-{Stealth[length=5pt] Stealth[length=5pt]}]
    (1.7,0)--(2.3,0);
  \draw[side,-{Stealth[length=5pt] Stealth[length=5pt]}]
    (1.7,3)--(2.3,3);

  \coordinate (a) at (2,1.5);
  \coordinate (p1) at (1,1.5);
  \coordinate (p2) at (3,1.5);
  \coordinate (q1) at (2,2.25);
  \coordinate (q2) at (2,.75);
  \coordinate (b) at (3.5,1.85);
  \coordinate (c) at (3.5,1.15);
  \coordinate (d) at (1.2,2.6);
  \coordinate (e) at (2.8,2.6);

  \draw[edge] (a)--(p1) (a)--(p2) (a)--(q1) (a)--(q2);
  \draw[edge] (p2)--(b) (p2)--(c) (q1)--(d) (q1)--(e);

  \draw[edge] (p1)--(0,1.85) (4,1.85)--(b);
  \draw[edge] (p1)--(0,1.15) (4,1.15)--(c);

  \draw[edge] (q2)--(1.2,0) (1.2,3)--(d);
  \draw[edge] (q2)--(2.8,0) (2.8,3)--(e);

  \node[vertex,label=below left:{$a$}] at (a) {};
  \node[vertex,label=below:{$p$}] at (p1) {};
  \node[vertex,label=below:{$q$}] at (p2) {};
  \node[vertex,label=-35:{$r$}] at (q1) {};
  \node[vertex,label=right:{$z$}] at (q2) {};
  \node[vertex,label=above:{$b$}] at (b) {};
  \node[vertex,label=below:{$c$}] at (c) {};
  \node[vertex,label=left:{$d$}] at (d) {};
  \node[vertex,label=right:{$e$}] at (e) {};

  \node[align=center,anchor=north,inner sep=1pt] at (2,-.25)
    {$a,b,a,d,a,c,a,e$};
\end{scope}

\begin{scope}[xshift=4.8cm]
  \fill[faceSimple] (0,0) rectangle (4,.8);
  \fill[faceSimple] (0,2.2) rectangle (4,3);

  \draw[side] (0,0) rectangle (4,3);
  \draw[side,-{Stealth[length=5pt]}] (0,.25)--(0,.65);
  \draw[side,-{Stealth[length=5pt]}] (4,.25)--(4,.65);
  \draw[side,-{Stealth[length=5pt] Stealth[length=5pt]}]
    (1.7,0)--(2.3,0);
  \draw[side,-{Stealth[length=5pt] Stealth[length=5pt]}]
    (1.7,3)--(2.3,3);

  \coordinate (a) at (2,1.5);
  \coordinate (p1) at (1.25,.8);
  \coordinate (c) at (2,.8);
  \coordinate (p2) at (2.75,.8);
  \coordinate (b) at (3.5,.8);
  \coordinate (q1) at (1.25,2.2);
  \coordinate (d) at (2,2.2);
  \coordinate (q2) at (2.75,2.2);
  \coordinate (e) at (3.5,2.2);

  \draw[edge] (p1)--(c)--(p2)--(b);
  \draw[edge] (q1)--(d)--(q2)--(e);
  \draw[edge] (a)--(p1) (a)--(p2) (a)--(q1) (a)--(q2);

  \draw[edge] (p1)--(0,.8) (4,.8)--(b);
  \draw[edge] (q1)--(0,2.2) (4,2.2)--(e);

  \node[vertex,label=left:{$a$}] at (a) {};
  \node[vertex,label=above left:{$p$}] at (p1) {};
  \node[vertex,label=below:{$c$}] at (c) {};
  \node[vertex,label=above right:{$q$}] at (p2) {};
  \node[vertex,label=below:{$b$}] at (b) {};
  \node[vertex,label=below left:{$r$}] at (q1) {};
  \node[vertex,label=above:{$d$}] at (d) {};
  \node[vertex,label=below right:{$z$}] at (q2) {};
  \node[vertex,label=above:{$e$}] at (e) {};
  \node[anchor=north,inner sep=1pt] at (2,-.25) {annulus $bc\mid de$};
\end{scope}
\end{tikzpicture}}
\endgroup
\caption{Embeddings of $\Psi_3$ in Lemma~\ref{r1:exterior-faces}, with a disk face shaded on the left and an annular face on the right.}
\label{fig:psi3-torus}
\end{figure}
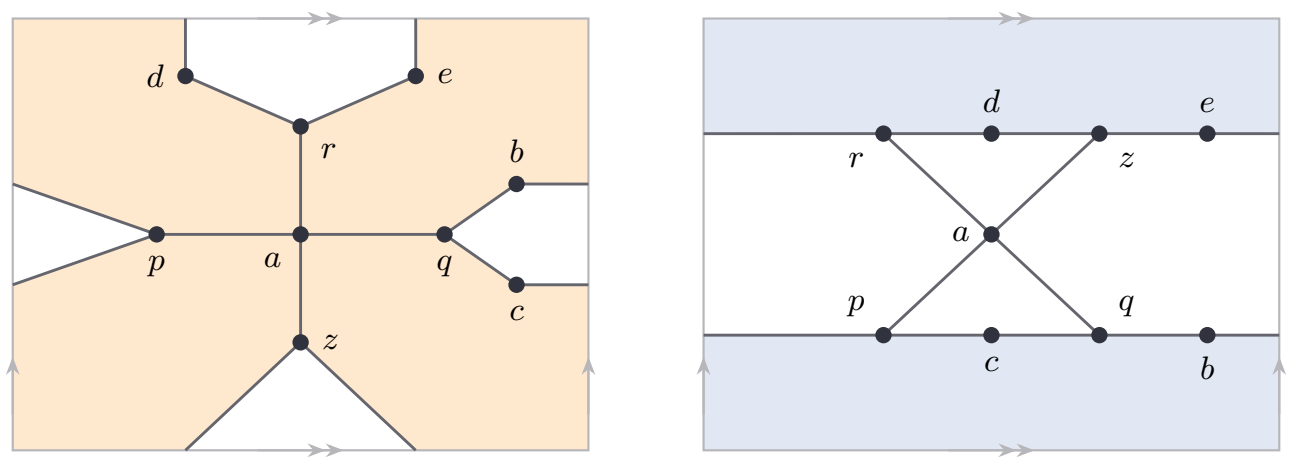

The disk boundary in Figure~\ref{fig:psi3-torus} follows the star $ab+ac+ad+ae$, with $ab,ac$ alternating with $ad,ae$ at $a$. The next lemma determines when a planar graph containing this star has an embedding with that order.

\begin{lemma}\label{r1:alternation}
Let $b,c,d,e$ be distinct neighbors of $a$ in a planar graph $H$, and suppose $H-a$ is connected. Then $H$ has a plane embedding in which $ab,ac$ alternate with $ad,ae$ around $a$ if and only if $H-a$ has no vertex-disjoint $b$--$c$ and $d$--$e$ paths.
\end{lemma}
\begin{proof}
Suppose $P,Q$ are vertex-disjoint $b$--$c$ and $d$--$e$ paths in $H-a$. In an alternating embedding, the cycle formed by $P$ and the edges $ab,ac$ separates $d$ from $e$, although $Q$ is disjoint from this cycle.

Conversely, suppose no such paths exist. Choose a plane embedding of $H$. If the four edges do not already alternate, we may assume their cyclic order is $ab,ad,ae,ac$ by relabeling within the two pairs. Split $a$ into one leaf on each incident edge and take their common face as the outer face. Draw the graph in a disk whose boundary meets it only at the leaves $b',d',e',c'$ on the four specified edges, in this order.

After deleting these four leaves, the graph consists of $H-a$ with the remaining leaves attached, so it is connected. Disjoint $b'$--$c'$ and $d'$--$e'$ paths would give the excluded paths in $H-a$. Thus Lemma~\ref{dis:four-disk-order} gives boundary order $b',d',c',e'$ while keeping all the new leaves on the outer face. Join $a$ to them in that face and suppress them. The four edges now alternate around $a$.
\end{proof}

\begin{lemma}\label{r1:twin-count}
Let $G$ be a three-connected obstruction with twin triples
$A_1=\{a,b,c\}$ and $A_2=\{a,d,e\}$ satisfying
$A_1\cap A_2=\{a\}$. Let $n_c$ and $n_e$ denote the numbers
of vertices of $G$ with neighborhoods $A_1$ and $A_2$,
respectively. Then at least one of these twin triples is
the neighborhood of exactly two vertices
\[
\min(n_c,n_e)=2.
\]
\end{lemma}
\begin{proof}
Otherwise, three vertices with neighborhood $\{a,b,c\}$ and three
with neighborhood $\{a,d,e\}$ give a subgraph $H$ of $G$ isomorphic
to $\mathfrak T_8$. Since $H$ has cut vertex $a$ and $G$ is
three-connected, $H\neq G$. Thus $H$ is a proper nontoroidal
subgraph of $G$, contradicting minimality.
\end{proof}

\begin{lemma}\label{r1:other-attachments}
Let $B$ be a $\Psi_3$-bridge whose attachment set contains $\{b,c,d,e\}$. Every other $\Psi_3$-bridge with nonempty interior has attachment set $\{a,b,c\}$, $\{a,d,e\}$ or $\{a,x,y\}$, where $x\in\{b,c\}$ and $y\in\{d,e\}$. The distinct pairs $xy$ for which $xy\in E(G)$ or a bridge has attachment set $\{a,x,y\}$ form a matching.
\end{lemma}
\begin{proof}
Contract $J_1,J_2$ to the triangles $abc,ade$ and the interior of $B$ to $w$. If another bridge has at least three attachments among $b,c,d,e$, exchange the two twin triples if necessary so that it meets both $b,c$, then interchange $d,e$ if necessary so that it also meets $d$. Contract its interior to $u$, then contract $ae$ into $a$ and $du$ into $d$. This gives a $K_5$ on $a,b,c,d,w$, a contradiction. By three-connectivity, every other bridge with nonempty interior
has at least three attachments. At most two are among $b,c,d,e$,
so its attachments are $a$ and two of $b,c,d,e$.

Suppose two distinct pairs share an end. After relabeling, we may
take them to be $bd,be$. If $bd$ is absent, we obtain it by
contracting the interior of a bridge on $\{a,b,d\}$ into $b$.
If $be$ is absent, we do the same with a bridge on $\{a,b,e\}$. With the triangles and $w$ as above, contracting $ac$ into $a$ gives a $K_5$ on $a,b,d,e,w$, again a contradiction. Hence the pairs form a matching.
\end{proof}

\subsection{A bridge with five attachments}\label{r1:five}

This subsection shows that a bridge with five attachments gives $\mathfrak T_{20}$ or $\mathfrak T_{23}$. Corollary~\ref{r1:five-criterion} gives the embedding test for such a bridge, and Proposition~\ref{r1:five-outcome} reduces the bridge and enumerates the possible forms.

Suppose a $\Psi_3$-bridge $B$ has all five attachments. Lemma~\ref{r1:exterior-faces} tells us it must lie in a disk with attachment sequence $a,b,a,d,a,c,a,e$. Put $H=G-\{p,q,r,z\}+ab+ac+ad+ae$. Contracting $ap,ar$ and deleting $q,z$ gives $H$, so $H$ has no $K_5$ minor.

\begin{corollary}\label{r1:five-criterion}
Let $F$ be a disk face of a torus embedding of $\Psi_3$ with attachment sequence $abadacae$, and let $U$ be a union of $\Psi_3$-bridges, including $a,b,c,d,e$. Then $U$ can be drawn in $F$ if and only if $U+ab+ac+ad+ae$ has a plane embedding in which $ab,ac$ alternate with $ad,ae$ at $a$. In particular, if a $\Psi_3$-bridge $B$ has all five attachments, then $G$ is toroidal if and only if
\[
H=G-\{p,q,r,z\}+ab+ac+ad+ae
\]
is planar and $H-a$ has no vertex-disjoint $b$--$c$ and $d$--$e$ paths.
\end{corollary}
\begin{proof}
Put $S=ab+ac+ad+ae$. Its plane embedding with cyclic order $ab,ad,ac,ae$ at $a$ has facial walk $abadacae$. The disks representing this face and $F$ have the same attachment appearances in the same cyclic order. A homeomorphism matching these appearances carries a drawing of $U$ in either disk to a drawing in the other.

If $B$ has all five attachments, Lemma~\ref{r1:exterior-faces} places it in the disk in every torus embedding of $G$, since the annulus does not meet $a$. Bridges with nonempty interior have at least three attachments, and every cross edge joins $\{b,c\}$ to $\{d,e\}$, so none can lie in the other two faces with attachment sequences $bc,de$. Thus, by the first claim and Figure~\ref{fig:psi3-torus}, $G$ is toroidal if and only if $H$ has the stated alternating plane embedding. Each bridge with $a$ deleted is connected and meets $B-a$, so $H-a$ is connected. Lemma~\ref{r1:alternation} now gives the stated test.
\end{proof}

\begin{proposition}\label{r1:five-outcome}
If some choice of the four twins gives a $\Psi_3$-bridge with five attachments, then $G$ is isomorphic to $\mathfrak T_{20}$ or $\mathfrak T_{23}$.
\end{proposition}
\begin{proof}
Let $B$ be such a bridge. If another component has neighborhood $\{a,b,c\}$, contract it and the interior of $B$ to single vertices. Keep $\Psi_3$ and the edges from these vertices to $\{a,b,c\}$ and $\{a,d,e\}$, respectively. This gives $\mathfrak T_8$~\supp{W9}, a proper nontoroidal minor since $a$ is a cut vertex. The case $\{a,d,e\}$ is symmetric.

By Lemma~\ref{r1:other-attachments}, the pairs supplied by cross edges or the remaining components form a matching. If there are two, label them $bd,ce$ and obtain either missing edge by contracting the corresponding component onto one of its attachments. Contract $J_1,J_2$ to the triangles $abc,ade$ and the interior of $B$ to $w$. Contracting $bd$ now gives a $K_5$, a contradiction. Thus we may assume $bd$ is the only possible cross edge and every other component has neighborhood~$\{a,b,d\}$.

Two such components can be contracted to two vertices with neighborhood $\{a,b,d\}$. A single bridge $B'$ with nonplanar triangle completion gives the same minor. Its completion is three-connected by Lemma~\ref{lem:triangle-completion}, and $apb$, $ard$ and a $b$--$d$ path through the interior of $B$ join its attachment pairs outside $B'$. Thus Lemma~\ref{endpoint:primitive} gives an allowed minor $T(s)$, and contracting its edge at $s$ gives the two vertices. In either case, keeping these vertices and $\Psi_3$ and contracting the interior of $B$ to one vertex gives $\mathfrak T_{22}$~\supp{W10}, after deleting any extra edges. This nontoroidal minor is proper since its twin triples $\{a,b,c\}$ and $\{a,b,d\}$ meet twice. Hence there is at most one other bridge with nonempty interior, and its triangle completion is planar.

Suppose $H-a$ has vertex-disjoint $b$--$c$ and $d$--$e$ paths $P,Q$. Any path using $bd$ or passing through the other bridge contains both $b$ and $d$. Since $P$ and $Q$ are disjoint, with $b$ on $P$ and $d$ on $Q$, neither path can do so. Thus both paths have their interiors in $B$, and both interiors are nonempty since $bc,de$ are absent. Contract the interiors of $P,Q$ to $u,v$, respectively, then contract the connected interior of $B$ to the edge $uv$, keeping $u,v$ distinct. This gives $ub,uc,vd,ve$ and at least one of $au,av$. Keep exactly one of $au,av$, together with $ub,uc,uv,vd,ve$ and $\Psi_3$. The resulting minor is $\mathfrak T_6$~\supp{W11}. It is proper since $a$ and its retained neighbor in $\{u,v\}$ form a two-cut. This contradicts minimality.

Thus no such paths exist, so $H$ is nonplanar by Corollary~\ref{r1:five-criterion}. If another bridge $B'$ on $\{a,b,d\}$ remains, the two bridges of $H$ meeting all of $\{a,b,d\}$ are $B'$ and $B+ac+ae$. The triangle completion of $B'$ is planar, so Lemma~\ref{core:assembly-general}(ii) makes the completion of the other bridge, $B+ab+ac+ad+ae+bd$, nonplanar. Contract the interior of $B'$ into $b$ and delete the resulting edge $ab$. This gives the proper minor $(\Psi_3\cup B)+bd$, which is nontoroidal by the same corollary. Hence $G=\Psi_3\cup B$, possibly with the edge $bd$.

The completion $H$ is nonplanar and has no $K_5$ minor. If $bd$ is present, $G-bd$ is toroidal, so $H-bd$ is planar by Corollary~\ref{r1:five-criterion}. Every $K_{3,3}$ subdivision of $H$ then uses $bd$. By minimality, replacing $B$ in $G$ by a proper allowed minor
gives a toroidal graph. The reduced bridge still has all five
attachments, so its completion is planar by
Corollary~\ref{r1:five-criterion}. Thus
Corollary~\ref{common:reduced-forms} applies with fixed edges
$ab,ac,ad,ae$, together with $bd$ if it is present in $G$. Of the $217{,}376$ candidates, all but eight have a $K_5$ minor, a cut of at most two vertices, or a proper nontoroidal minor~\supp{E20}. The eight are isomorphic to $\mathfrak T_{20},\mathfrak T_{23},\mathfrak T_{38},\mathfrak T_{40}$ or $\mathfrak T_{69}$. In each of the last three graphs, two twin triples are disjoint or meet in two vertices. This leaves $\mathfrak T_{20}$ and $\mathfrak T_{23}$.
\end{proof}

\subsection{A bridge with four attachments, not including \texorpdfstring{$a$}{a}}\label{r1:missing-a}
From now on, assume that no choice of the four twins gives a bridge
with five attachments.
Choose a $\Psi_3$-bridge $B$ with attachment set
$\{b,c,d,e\}$. In any torus embedding of $\Psi_3\cup B$,
Lemma~\ref{r1:exterior-faces} places $B$ in the disk or
annulus of Figure~\ref{fig:psi3-torus}. By Lemma~\ref{r1:other-attachments}, every other
$\Psi_3$-bridge with nonempty interior has three
attachments, including $a$. Let $B_1$ be their union,
and let $E$ be the set of cross edges. Put
\[
U=G-\{p,q,r,z\},\qquad H=U+ab+ac+ad+ae.
\]
Thus $U$ consists of $B$, $B_1$, the vertex $a$ and
the edges in $E$. The next lemma gives an embedding
test for each placement of $B$.

\begin{lemma}\label{r1:missing-a-criterion}
The graph $\Psi_3\cup U$ is toroidal if and only if at least one of the following holds.
\begin{enumerate}[label=\textup{(\arabic*)}]
\item $H$ has a plane embedding in which $ab,ac$ alternate with $ad,ae$ at $a$.
\item For some $x\in\{b,c\}$ and $y\in\{d,e\}$, every bridge in $B_1$ has attachment set $\{a,x,y\}$, and there is a set $E_0\subseteq E\cap\{xy\}$ for which both graphs
\[
B+bc+de+(E\setminus E_0),\qquad B_1+ax+ay+E_0
\]
are planar.
\end{enumerate}
\end{lemma}
\begin{proof}
Consider the two placements of $B$ given by Lemma~\ref{r1:exterior-faces}. If $B$ lies in the disk $abadacae$, the other faces meet only $b,c$ or only $d,e$. Thus all bridges must lie in this disk, and Corollary~\ref{r1:five-criterion} gives (1).

If $B$ lies in the annulus $bc\mid de$, every bridge in $B_1$ must lie in the disk $axay$, since each has three attachments including $a$. Hence each has attachment set $\{a,x,y\}$. The cross edges placed in this disk form a set $E_0\subseteq E\cap\{xy\}$, and the remaining cross edges lie in the annulus with $B$.

\looseness=-1
Lemma~\ref{lem:annulus-embedding} gives the first planarity test in (2). For the disk test, the plane path $xay$ has just one face, with boundary walk $axay$. A plane drawing of $B_1+ax+ay+E_0$ therefore draws $B_1+E_0$ in a face with the required attachment appearances. Conversely, placing a disk drawing of $B_1+E_0$ in this face gives a plane drawing of $B_1+ax+ay+E_0$. Every choice of $x,y$ occurs by symmetry in Figure~\ref{fig:psi3-torus}, so the two drawings together give a torus embedding.
\end{proof}

The next lemma shows that nonplanarity of $B+bc+de$
implies nonplanarity of $B+ab+ac+ad+ae$ when their union
has no $K_5$ minor.

\begin{lemma}\label{r1:matching-test}
Let $H$ be a planar graph with $H-a$ connected and
$N_H(a)=\{b,c,d,e\}$, where $b,c,d,e$ are distinct.
If $H+bc+de$ has no $K_5$ minor, then $(H-a)+bc+de$ is planar.
\end{lemma}
\begin{proof}
Suppose $(H-a)+bc+de$ is nonplanar. If $ab,ac$ do not alternate with $ad,ae$ in some plane embedding of $H$, draw $bc,de$ alongside $bac,dae$ and delete $a$. This gives a plane embedding of $(H-a)+bc+de$, a contradiction. Thus neither of the other two pairings alternates at $a$ in any plane embedding of $H$. Since $H-a$ is connected, Lemma~\ref{r1:alternation} gives disjoint $b$--$d$ and $c$--$e$ paths $P,Q$ in $H-a$, and disjoint $b$--$e$ and $c$--$d$ paths.

On each of the latter two paths, take a shortest subpath joining $P$ to $Q$. These subpaths are disjoint and meet $P,Q$ only at their ends. Together with the portions of $P,Q$ between their ends, they form a cycle $C$. Contract the remaining portions of $P,Q$ onto $C$. In the plane embedding of $H$, the vertices now occur on $C$ in the order of their edges at $a$, namely $b,d,c,e$ up to reversal. Shorten the four arcs of $C$ to single edges. Together with $bc,de$ and the four edges at $a$, this gives a $K_5$ minor of $H+bc+de$, a contradiction.
\end{proof}
We next restrict the other components and settle the case
where $B+bc+de$ is nonplanar.

\begin{lemma}\label{r1:missing-a-restrictions}
After relabeling within $\{b,c\}$ and $\{d,e\}$, every component of $G-V(\Psi_3)$
other than the interior of $B$ has neighborhood
$\{a,b,d\}$ or $\{a,c,e\}$, and every cross edge is
$bd$ or $ce$. If $B+bc+de$ is nonplanar, then
$G\cong\mathfrak T_{54}$.
\end{lemma}
\begin{proof}
If another component has neighborhood $\{a,b,c\}$,
contract it to $u$ and the interior of $B$ to $w$.
Keep only $\Psi_3$ and the edges at $u,w$, then delete
$bw$. This gives a proper nontoroidal minor
$\mathfrak T_3$~\supp{W12}, a contradiction.
The case $\{a,d,e\}$ is symmetric.
Lemma~\ref{r1:other-attachments} now gives the stated
neighborhoods and cross edges.

Put $H_2=B+bc+de$ and suppose it is nonplanar.
By Lemma~\ref{r1:triangle-minor},
$B+ab+ac+ad+ae+bc+de$ is a minor of $G$.
Thus $H_2$ has no $K_5$ minor. Let $B'$ be an
allowed minor of $B$ with $B'+bc+de$ nonplanar.
The graph $B'+ab+ac+ad+ae+bc+de$ is still a minor
of $G$, so Lemma~\ref{r1:matching-test} makes
$B'+ab+ac+ad+ae$ nonplanar. Thus $\Psi_3\cup B'$
fails both tests of Lemma~\ref{r1:missing-a-criterion}.
Taking $B'=B$, minimality gives $G=\Psi_3\cup B$, and
Corollary~\ref{common:reduced-forms}, with attachments
$b,c,d,e$ and added edges $bc,de$,
includes $B$ among $26{,}633$ forms. Of the resulting graphs,
all but four are toroidal or have a $K_5$ minor or a proper
nontoroidal minor. All four are isomorphic to $\mathfrak T_{54}$~\supp{E21}.
\end{proof}

We may now assume $H_2=B+bc+de$ is planar.
If $B_1$ is empty, condition (2) of
Lemma~\ref{r1:missing-a-criterion} forces $|E|=2$.
Let $E'$ consist of $E$ and the edges $xy$ supplied
by components with neighborhood $\{a,x,y\}$.
Thus $E'$ is nonempty, and
Lemma~\ref{r1:other-attachments} allows us to label
it as $\{bd\}$ or $\{bd,ce\}$. Put
\[
H_3=B+ab+ac+ad+ae+bc+de+E'.
\]

\begin{lemma}\label{r1:missing-a-cut}
The graph $H_3$ is three-connected and has no $K_5$
minor. If it is nonplanar, then the interior of $B$
contains twins $p',q'$ with neighborhood $\{x,y,w\}$,
where $w\in\operatorname{Int}(B)$ and $xy$ is $bd$
or $ce$. Moreover, $xy\notin E(G)$, and
$B[\operatorname{Int}(B)\setminus\{p',q'\}]$ is
connected and has neighbors at the other two
vertices of $\{b,c,d,e\}$.
\end{lemma}
\begin{proof}
By Lemma~\ref{r1:triangle-minor}, contracting each other component with neighborhood $\{a,x,y\}$ into $x$ gives $H_3$ as a minor of $G$, so it has no $K_5$ minor. Since the neighbors of each component of $G-V(H_3)$ form a triangle in $H_3$, any cut of $H_3$ also separates $G$. Thus $H_3$ is three-connected.

Suppose $H_3$ is nonplanar. The triangle $abc$ excludes $V_8$, so Halin's theorem gives a three-cut $X$ leaving at least three components, each with neighborhood $X$. If $a\in X$, each component contains a neighbor of $a$, but all four neighbors lie on the path $cbde$, which has at most two pieces after deleting the other two cut vertices. Hence $a\notin X$. Every attachment outside $X$ is adjacent to $a$, so all lie in the same component of $H_3-X$. At least two other components therefore lie entirely in the interior of $B$. These are components of $B-X$, so Lemma~\ref{lem:avoiding-components} gives exactly two single vertices $p',q'$ with neighborhood $X$ in $G$.

The twin triple $X$ meets each of $\{a,b,c\}$ and $\{a,d,e\}$ exactly once. Thus $X=\{x,y,w\}$ with $x\in\{b,c\}$, $y\in\{d,e\}$ and $w\in\operatorname{Int}(B)$. In the interior of $B$, $p',q'$ are leaves at $w$. Deleting them leaves the interior connected and preserves its neighbors at the other two attachments. Contract this remaining interior to $w$. If $xy=be$, contracting $bp',eq',cw$ gives a $K_5$ on $a,b,c,d,e$. Interchanging $b,d$ and $c,e$ excludes $cd$ as well. Hence $xy$ is $bd$ or $ce$, and $xy\notin E(G)$ since $N_G(p')$ is independent by Lemma~\ref{lem:independent-neighborhood}(a).
\end{proof}

We now complete the case where a bridge has attachment set
$\{b,c,d,e\}$.

\begin{proposition}\label{r1:missing-a-outcome}
If no choice of the four twins gives a bridge with five attachments, but some choice gives a bridge with attachment set $\{b,c,d,e\}$, then $G$ is isomorphic to one of $\mathfrak T_{29},\mathfrak T_{54},\mathfrak T_{55}$.
\end{proposition}
\begin{proof}
If $H_2$ is nonplanar, Lemma~\ref{r1:missing-a-restrictions} gives $G\cong\mathfrak T_{54}$. Assume $H_2$ is planar. We use the two tests of Lemma~\ref{r1:missing-a-criterion}.

Suppose $E'=\{bd,ce\}$ and $H_3$ is nonplanar. By Lemma~\ref{r1:missing-a-cut}, we may interchange $b,c$ and $d,e$ so that $B$ contains twins $p',q'$ with neighborhood $\{b,d,w\}$ and $bd\notin E(G)$. Since $bd\in E'$, some component has neighborhood $\{a,b,d\}$. Contract it to $u$. If $ce$ is absent, obtain it by contracting a component with neighborhood $\{a,c,e\}$ into $c$. Contract the remaining interior of $B$ apart from $p',q'$ to $w$, as allowed by the same lemma. Keep $\Psi_3$, the edges at $u,p',q'$, and $cw,ew,ce$. This gives $\mathfrak T_{29}$~\supp{W13}, so $G\cong\mathfrak T_{29}$. We may therefore assume that $H_3$ is planar whenever~$E'=\{bd,ce\}$.

If $B_1$ is empty, then $E=E'=\{bd,ce\}$. Thus $H_2+E$ is a subgraph of the planar graph $H_3$, and (2) holds with $E_0=\varnothing$, a contradiction.

For any remaining bridge on $\{a,x,y\}$, its triangle completion is three-connected by Lemma~\ref{lem:triangle-completion}. The paths $apx$, $ary$ and an $x$--$y$ path through the interior of $B$ join the three attachment pairs with disjoint interiors outside the bridge. Thus the completion is a minor of $G$. If it is nonplanar, Lemma~\ref{endpoint:primitive} gives an allowed minor $T(s)$, and contracting its edge at $s$ gives two vertices with neighborhood $\{a,x,y\}$.

Suppose $B_1$ consists of one bridge, with attachments $\{a,b,d\}$ by symmetry. Put $H_4=B_1+ab+ad$. If $H_4$ is nonplanar, deleting $E$ and contracting the interior of $B$ to one vertex still leaves $H_4$ in both tests. By minimality, $E=\varnothing$ and $B$ has one interior vertex. Likewise, no proper allowed reduction of $B_1$ keeps $H_4$ nonplanar. Corollary~\ref{common:reduced-forms}, with added edges $ab,ad$, gives the $349$ forms in~\supp{E22}. Each resulting graph has a $K_5$ minor, a cut of at most two vertices, or a proper nontoroidal minor.

Now suppose $H_4$ is planar. Taking $E_0=\varnothing$ in (2) shows that $H_2+E$ is nonplanar. Hence $E\ne\varnothing$ and $H_3$ is nonplanar, leaving $E=E'=\{bd\}$. Taking $E_0=\{bd\}$ shows that $H_4+bd$ is nonplanar, since $H_2$ is planar. This is the triangle completion of $B_1$, so $B_1$ has an allowed minor $T(s)$. Its triangle completion is still nonplanar by Lemma~\ref{core:three-graph-completions}. Thus the disk test and the annular choice $E_0=\{bd\}$ still fail, while $H_2+bd$ excludes $E_0=\varnothing$. Minimality gives $B_1=T(s)$. If $s=a$, its twins have adjacent neighbors $b,d$, contrary to Lemma~\ref{lem:independent-neighborhood}(a). If $s=b$ or $s=d$, their triple meets $\{a,d,e\}$ or $\{a,b,c\}$, respectively, twice.

Finally suppose $B_1$ contains at least two bridges. If one has nonplanar triangle completion, take its attachments to be $\{a,b,d\}$ and reduce it to two vertices with this neighborhood as above. In the disk test, the connected set consisting of $c$ and the interior of $B$ has neighbors at $a,b,d$, using the added edge $ac$. Together with the two new vertices this gives a $K_{3,3}$ minor. An annular placement requires every other bridge in $B_1$ to attach at $a,b,d$, in which case any one of them gives a third bridge on this triple in the disk $abad$. Thus both tests fail for a proper minor of $G$, a contradiction.

\looseness=-1
All bridges in $B_1$ therefore have planar triangle completion, so their interiors are single vertices by Lemma~\ref{lem:flat-replacement}. Two with the same neighborhood would give a twin triple meeting $\{a,b,c\}$ twice or $\{a,d,e\}$ twice. Hence there are exactly two vertices $u,v$, with neighborhoods $\{a,b,d\}$ and $\{a,c,e\}$. These neighborhoods are independent, so $E=\varnothing$.

Now $E'=\{bd,ce\}$, so $H_3$ is planar. In a plane drawing of $H_3$, each of the triangles $abd,ace$ bounds a face, since deleting its vertices leaves the connected interior of $B$ and the other two attachments. Insert $u,v$ in these faces and delete $bc,de,bd,ce$ to obtain a plane drawing of $H$. Since (1) fails and $H-a$ is connected, Lemma~\ref{r1:alternation} gives vertex-disjoint $b$--$c$ and $d$--$e$ paths in $H-a$. Neither path can use $u$ or $v$ without meeting an end of the other, so both have their interiors in $B$. These interiors are nonempty since $bc,de$ are absent. Contract them to $y,t$, respectively, then contract the connected interior of $B$ to the edge $yt$, keeping $by,cy,dt,et$. Together with $\Psi_3,u,v$ this gives $\mathfrak T_{55}$~\supp{W14}, so $G\cong\mathfrak T_{55}$.
\end{proof}

\subsection{A bridge with four attachments, including
\texorpdfstring{$a$}{a}}\label{r1:missing-b}

Assume that no choice of the four twins gives a bridge
with five attachments or four attachments excluding their
common vertex. Choose the four twins and a bridge $B$
with four attachments to maximize $|\operatorname{Int}(B)|$,
and call its attachments $\{a,c,d,e\}$.

\begin{lemma}\label{r1:missing-b-other-bridges}
After interchanging $d,e$ if necessary, every cross edge is $bd$ or $cd$, and every component of $G-V(\Psi_3)$ other than the interior of $B$ is a single vertex with neighborhood $\{a,d,e\}$, $\{a,b,d\}$ or $\{a,c,d\}$. Moreover, there is a $b$--$d$ path whose interior avoids $V(\Psi_3)\cup\operatorname{Int}(B)$.
\end{lemma}
\begin{proof}
Contract $J_1,J_2$ to the triangles $abc,ade$ and the interior of $B$ to $w$. Contracting $bc$ into $b$ gives all edges of $K_5$ on $a,b,d,e,w$ except possibly $bd,be$. If cross edges or other components join $\{b,c\}$ to both $d$ and $e$, first contract each such component onto a neighbor in $\{b,c\}$. This supplies both missing edges, a contradiction. Thus we may assume that all such connections go to $d$.

By three-connectivity, the only possible neighborhoods besides the three in the statement are $\{a,b,c\}$, $\{b,c,d\}$ and $\{a,b,c,d\}$. In these cases contract the component to $u$ and the interior of $B$ to $w$, retaining only $\Psi_3$ and the edges at $u,w$. Deleting $aw$, $cw$, or both $au,cw$, respectively, gives the proper nontoroidal minor $\mathfrak T_3$~\supp{W15}. Also, $b$ has only the neighbors $p,q$ in $\Psi_3$ and is not an attachment of $B$. Thus $bd$ is present or some component has neighborhood $\{a,b,d\}$, giving the required $b$--$d$ path.

\looseness=-1
Suppose some other bridge $B_1$ has nonplanar triangle completion. Choose a $c$--$d$ path through the interior of $B$. Together with $apbqc$ and $dreza$, it forms a cycle through $a,b,c,d,e$ that avoids the interior of $B_1$. The three arcs between the attachments of $B_1$ give its triangle completion as a minor of $G$, and Lemma~\ref{lem:triangle-completion} makes it three-connected.

If $B_1$ attaches at $\{a,d,e\}$ or $\{a,c,d\}$, the same three arcs give the outside paths in Lemma~\ref{endpoint:primitive}. Thus $B_1$ has an allowed minor $T(s)$, and contracting its edge at $s$ gives two vertices with the same attachment triple. The $b$--$d$ path above avoids the interior of $B_1$. Contract it into $d$. The original twin triples are now $\{a,c,d\}$ and $\{a,d,e\}$, and the two new vertices give four twins on one of these triples and two on the other. Contract the interior of $B$ to a vertex adjacent to $a,c,d,e$ and delete extra edges. This gives $\mathfrak T_2$~\supp{W16}, with $c,e$ interchanged if the four twins are on $\{a,c,d\}$. This minor is proper since $b,d$ were identified.

It remains to consider $B_1$ on $\{a,b,d\}$. Contract the interior of $B$ to $w$, delete $aw$, and retain only $\Psi_3$, $B_1$ and the edges $cw,dw,ew$. This proper minor is toroidal. Let $K$ consist of $\Psi_3$ and the three edges at $w$. It contains a $K_{3,3}$ subdivision with parts $\{a,d,e\}$ and $\{r,z,w\}$, using $apcw$ for the $a$--$w$ branch. Thus every face of $K$ in the torus embedding is a disk by Lemma~\ref{lem:torus-complement}. On any face containing $a,b,d$, the boundary sequence restricted to these three vertices is $abad$~\supp{R10}. Since $B_1$ lies in such a face, Lemma~\ref{cor:boundary-completions} gives a plane drawing of $B_1+ab+ad$.

Lemma~\ref{core:three-root-cut} now gives twins $p',q'$ in $B_1$ whose neighborhood contains $a$. Their triple meets each of $\{a,b,c\}$ and $\{a,d,e\}$ only at $a$, so it has the form $\{a,u,v\}$ with $u,v\in\operatorname{Int}(B_1)$. By part~(i) of the same lemma, we may label $u,v$ so that the component of $B_1-\{a,u,v\}$ containing $b$ has a neighbor at $u$.

Choose $p',q',r,z$ as the four twins. The component outside the new $\Psi_3$ containing $b$ also contains $p,q,c$ and the interior of $B$, since $p,q$ join $b$ to $c$ and $c$ has a neighbor in $B$. This component has neighbors at $a,d,e$ through $B$ and at $u$ through $B_1$. Five attachments are excluded, so its bridge has four attachments and at least $|\operatorname{Int}(B)|+4$ interior vertices, contrary to the choice of $B$. Hence every other bridge with nonempty interior has planar triangle completion, and Lemma~\ref{lem:flat-replacement} makes its interior a single vertex.
\end{proof}

The next lemma excludes the additional vertices and gives
the two nonplanar completions of $B$.

\begin{lemma}\label{r1:missing-b-exterior}
With $d,e$ labeled as in Lemma~\ref{r1:missing-b-other-bridges}, we have $G=(\Psi_3\cup B)+bd+E$ for some $E\subseteq\{cd\}$. Both~graphs
\begin{equation}\label{r1:c1-caps}
B+ac+E+de,\qquad B+ac+E+ad+ae
\end{equation}
are nonplanar.
\end{lemma}
\begin{proof}
By Lemma~\ref{r1:missing-b-other-bridges}, $n_c=2$. If $n_e\ge4$, contract the outside $b$--$d$ path into $d$ and the interior of $B$ to one vertex. Keeping four vertices with neighborhood $\{a,d,e\}$ gives the proper $\mathfrak T_2$ minor of~\supp{W16}. Thus $n_e=2$ or $3$.

There is at most one vertex with each neighborhood $\{a,b,d\}$ or $\{a,c,d\}$, since two with the same neighborhood would give a twin triple meeting $\{a,b,c\}$ twice. By the outside $b$--$d$ path and the independence of these neighborhoods, exactly one of $bd$ and a vertex with neighborhood $\{a,b,d\}$ is present. The remaining choice is $cd$, a vertex with neighborhood $\{a,c,d\}$, or neither.

For either choice at $b$, $G$ is toroidal exactly when one of the indicated sets of edges can be added to $B+ac$ to give a planar graph.
\begin{center}
\scaledtable{.9}{%
\begin{tabular}{@{}lll@{}}
\toprule
Remaining choice & $n_e=2$ & $n_e=3$\\
\midrule
Neither
  & $\{de\}$ or $\{ad,ae\}$ & $\{ad,ae,de\}$\\
$cd$
  & $\{cd,de\}$ or $\{cd,ad,ae\}$ & $\{cd,ae,de\}$\\
Vertex with neighborhood $\{a,c,d\}$
  & $cd$ and any two of $ad,ae,de$ & $\{cd,ae,de\}$\\
\bottomrule
\end{tabular}%
}
\end{center}

Indeed, the rotation tables in~\supp{R8} give the boundary sequences for all twelve choices of $G-\operatorname{Int}(B)$. For each disk face meeting all four attachments, Lemma~\ref{cor:boundary-completions} gives the test by joining consecutive distinct attachments in its boundary sequence. Only the first row with $n_e=2$ admits a genus-zero rotation. No single boundary walk then meets all four attachments, and the only annulus boundary walk that does is $ac\mid de$, giving the test $B+ac+de$ by Lemma~\ref{lem:annulus-embedding}. Keeping the inclusion-minimal edge sets gives the table.

If a vertex has neighborhood $\{a,b,d\}$, contract it into $b$ and delete the new edge $ab$. This replaces the vertex by $bd$ and leaves the tests unchanged, giving a proper nontoroidal minor. Thus $bd$ is present. If a vertex has neighborhood $\{a,c,d\}$, contract $bd$ into $d$. This leaves $B$ unchanged and gives three vertices with neighborhood $\{a,c,d\}$ and $n_e$ with neighborhood $\{a,d,e\}$. Apply Lemma~\ref{r2:four-tests} with $(a,d,c,e)$ in place of $(a,b,c,d)$. For $n_e=2$ its two tests are the choices containing $ad$ in the last row. For $n_e=3$ its only test is the one in the last column of that row. All are nonplanar, again giving a proper nontoroidal minor.

It remains to exclude $n_e=3$. In this case the applicable row gives a single nonplanar completion $B+E'$. It has no $K_5$ minor, since contracting $bd$ into $d$ in the minor of Lemma~\ref{r1:triangle-minor} with $U=B\cup bd$ gives $B+ac+ad+ae+cd+de$. The same row applies after any allowed reduction of $B$, so the resulting graph is nontoroidal whenever its completion is nonplanar. Corollary~\ref{common:reduced-forms} therefore applies with added edges $E'$. This gives $9{,}127$ candidates without $cd$ and $8{,}967$ with $cd$. In each, deleting or contracting one edge gives a nontoroidal minor~\supp{E24}. Thus $n_e=2$, and the first two rows give the claimed form of $G$ and the two nonplanar completions in~\eqref{r1:c1-caps}.
\end{proof}

We now reduce $B$ while keeping both completions in~\eqref{r1:c1-caps} nonplanar.

\begin{proposition}\label{r1:missing-b-outcome}
Suppose some choice of the four twins gives a bridge with four
attachments. If no choice gives a bridge with five attachments or
a bridge with four attachments missing the common vertex of the two
twin triples, then $G$ is isomorphic to one of
$\mathfrak T_{14},\mathfrak T_{37},\mathfrak T_{41},\mathfrak T_{59}$.
\end{proposition}
\begin{proof}
Put $F=B+ac+E$. As in the proof of Lemma~\ref{r1:missing-b-exterior},
$B+ac+ad+ae+cd+de$ is a minor of $G$, so it has no $K_5$ minor.
If $F$ is nonplanar, it therefore contains a $K_{3,3}$ subdivision
$H$. In this case put $E'=\{ac\}\cup E$.

Suppose $F$ is planar. Regard $d,e,a$ as its attachments. They are
independent, and $\operatorname{Int}(B)\cup\{c\}$ is connected and
adjacent to all three. Let $K=F+de+ad+ae$ be its triangle completion.
For any $S\subseteq V(B)$ with $|S|\le2$, every component of $B-S$
contains a surviving attachment of $B$, since $G-S$ is connected.
If $a\notin S$, the edges $ac,ad,ae$ join all surviving attachments
to $a$, so $K-S$ is connected. Also $K-a$ is connected, since $c,d,e$
each have a neighbor in the interior of $B$. Thus $K$ is two-connected
and every two-cut contains $a$. By~\eqref{r1:c1-caps}, both $F+de$
and $F+ad+ae$ are nonplanar. Lemma~\ref{r2:three-mark-entrance}
with $(x,y,z)=(d,e,a)$ now gives a $K_{3,3}$ subdivision
$H\subseteq F+de$ using $de$ such that $H-de+ad+ae$ is nonplanar.
In this case put $E'=\{ac,de\}\cup E$.

Apply Lemma~\ref{core:retained-subdivision} to $H\subseteq B+E'$
with attachments $a,c,d,e$. This gives an allowed minor $B'$ with
at most nine interior vertices and a $K_{3,3}$ subdivision
$H'\subseteq B'+E'$. If $F$ is nonplanar, $H'$ lies in both
completions of $B'$ in~\eqref{r1:c1-caps}. If $F$ is planar, $H'$
lies in the first completion. The reduction keeps $de$ and preserves
the order of the attachments on each branch path, so
$H'-de+ad+ae$ is homeomorphic to $H-de+ad+ae$ and lies in the
second completion. Thus both completions remain nonplanar. By the
tests in the proof of Lemma~\ref{r1:missing-b-exterior},
$(\Psi_3\cup B')+bd+E$ is a nontoroidal minor of $G$. By minimality
it is $G$, which therefore has at most $18$ vertices.

Construction~\ref{common:marked-construction} with added edges $E'$
now includes every possibility, with $de$ required on the subdivision
when $F$ is planar. The two choices of $E$ and the two cases for $F$
give $66{,}384$ candidates. After excluding toroidal graphs and those
with a $K_5$ minor or a proper nontoroidal minor, sixteen candidates
remain, isomorphic to thirteen of the $\mathfrak T_i$~\supp{E25}.
In nine of the thirteen, two twin triples are disjoint or meet in
two vertices, leaving
$\mathfrak T_{14},\mathfrak T_{37},\mathfrak T_{41},\mathfrak T_{59}$.
\end{proof}

\subsection{Every remaining bridge has three attachments}\label{r1:three-neighbors}

From now on, for every choice of the four twins, every $\Psi_3$-bridge with nonempty interior has exactly three attachments. Fix such a choice. We first restrict the possible neighborhoods and then show that every bridge with nonempty interior has a planar completion by two edges, allowing us to apply Proposition~\ref{prop:three-root-blocks}.

\begin{lemma}\label{r1:three-attachment-setup}
After relabeling within $\{b,c\}$ and $\{d,e\}$, the pairs $xy$
with $x\in\{b,c\}$ and $y\in\{d,e\}$ joined by an edge or a
component of $G-V(\Psi_3)$ form one of the sets
\begin{equation}\label{r1:availability}
\{bd,ce\}\quad\hbox{or}\quad\{bd,cd,ce\}.
\end{equation}
Also, $n_c+n_e\le6$, and the triangle completion of every bridge
with nonempty interior is a three-connected minor of $G$ and hence
has no $K_5$ minor.
\end{lemma}
\begin{proof}
Every vertex of $\{b,c,d,e\}$ occurs in one of these pairs.
Otherwise, by symmetry, $b$ occurs in none, and deleting $a,c$
separates $b$ from $d,e$. By Lemma~\ref{r1:triangle-minor}, we may
add the triangles $abc,ade$ while retaining all bridges. Each component
with neighbors in both $\{b,c\}$ and $\{d,e\}$ has only one neighbor
in at least one of these sets, since it has only three attachments.
Contract it into that neighbor to obtain all pairs it joins. If all
four pairs occurred, this would give a $K_5$ on $a,b,c,d,e$. Thus
the pairs form a matching of size two or a path of length three,
giving~\eqref{r1:availability} after relabeling. In the first case
the possible component neighborhoods are $\{a,b,c\}$, $\{a,d,e\}$,
$\{a,b,d\}$ and $\{a,c,e\}$. The second also allows $\{a,c,d\}$,
$\{b,c,d\}$ and $\{c,d,e\}$.

Choose a $b$--$d$ path and a $c$--$e$ path in $G$, each a cross
edge or a path with interior in a single component of $G-V(\Psi_3)$.
They are disjoint, since a component supplying both paths would
have four attachments. The paths also avoid the vertices counted
by $n_c,n_e$, whose neighborhoods are $\{a,b,c\}$ or $\{a,d,e\}$.
Contracting each path to a vertex gives a $K_{3,n_c+n_e}$ minor,
with $a$ and the two contracted vertices in one part. If
$n_c+n_e\ge7$, then $K_{3,7}$ is a proper nontoroidal minor of $G$,
contrary to minimality.

Now let $B$ be a bridge with nonempty interior. At least one of the
two chosen paths avoids the interior of $B$. In the minor with
triangles $abc,ade$, extend this path at its ends by $bc$ and $de$,
then join its new ends through $a$. This gives a cycle through
$a,b,c,d,e$ avoiding the interior of $B$. Shortening the three
subpaths between consecutive attachments of $B$ to single edges
gives its triangle completion as a minor of $G$. It is
three-connected by Lemma~\ref{lem:triangle-completion}.
\end{proof}

We first consider bridges with attachments $\{b,c,d\}$ or $\{c,d,e\}$, and then a bridge with attachments $\{a,b,d\}$.

\begin{lemma}\label{r1:bcd-two-edge}
Every bridge with nonempty interior and attachments $\{b,c,d\}$
or $\{c,d,e\}$ has a planar completion by two edges.
\end{lemma}
\begin{proof}
By symmetry, suppose $B$ has attachments $A=\{b,c,d\}$ and no
planar completion by two edges. Its triangle completion is nonplanar
and, by Lemma~\ref{r1:three-attachment-setup}, is a three-connected
minor of $G$. If $B+bc$ were planar, Lemma~\ref{core:three-root-setup}
would give a three-cut $X$ of the triangle completion that is a twin
triple. Since $X$ meets each of $\{a,b,c\}$ and $\{a,d,e\}$ exactly
once and $a,e\notin V(B)$, it contains $d$ and exactly one of $b,c$,
say $u$. Part~\textup{(vi)} of the same lemma then gives $B+bc+du$
planar, a contradiction. Thus $B+bc$ is nonplanar.

By~\eqref{r1:availability}, $ce$ is an edge or is supplied by a
component outside $\Psi_3$. Since $e$ is not an attachment of $B$,
shortening the corresponding $c$--$e$ path to $ce$ gives
$(\Psi_3\cup B)+ce$ as a minor of $G$. Consider a torus embedding
of this minor. In a genus-one rotation of $\Psi_3+ce$, the
sequences on $b,c,d$ along facial walks meeting all three are
$bcd$ and $bcdc$. For a planar rotation,
no boundary walk contains all three attachments, and the only possible
annuli are $bc\mid ad$ and $bc\mid de$~\supp{R9}. In each disk sequence
$b,c$ are consecutive, and in each annulus they form one boundary
sequence. Thus Lemmas~\ref{cor:boundary-completions}
and~\ref{lem:annulus-embedding} require $B+bc$ to be planar in every
case. The minor is therefore nontoroidal, so by minimality it is $G$.

The same argument applies to any allowed minor $B'$ of $B$ for which
$B'+bc$ is nonplanar. By Corollary~\ref{common:reduced-forms}, it suffices to consider the
$563$ candidates with attachments $b,c,d$ and added edge $bc$.
Of these, $562$ have a proper nontoroidal minor. The remaining graph
is isomorphic to $\mathfrak T_{41}$ and has $B+bd+cd$ planar~\supp{E26},
contrary to the choice of $B$.
\end{proof}

\begin{lemma}\label{r1:abd-others}
Suppose a bridge $B$ with attachments $\{a,b,d\}$ has no planar
completion by two edges. Then every other bridge with nonempty
interior has attachments $\{a,b,d\}$ or $\{a,c,e\}$, and every
cross edge is $bd$ or $ce$.
\end{lemma}
\begin{proof}
Suppose another bridge has attachments $\{a,b,c\}$, $\{a,d,e\}$,
$\{a,c,d\}$ or $\{c,d,e\}$. Contract its interior to a vertex $v$
and retain $\Psi_3$, $B$ and the edges at $v$. In the first three
cases also retain $ce$, shortening the $c$--$e$ path
from~\eqref{r1:availability} if necessary. Its interior avoids both
bridges, since neither has both $c,e$ as attachments.

In each of these four minors the graph outside $\operatorname{Int}(B)$
is nonplanar, so every face in a torus embedding of that graph is a
disk. On a facial walk meeting $a,b,d$, the sequence restricted to
these vertices is $abad$, with $abd$ also possible when
$N(v)=\{a,c,d\}$~\supp{R10}. Both require $B+ab+ad$ to
be planar by Lemma~\ref{cor:boundary-completions}. Thus the retained
minor is nontoroidal and by minimality is $G$. The same test applies
to every allowed minor of $B$, so Corollary~\ref{common:reduced-forms}
gives $349$ forms with added edges $ab,ad$ for each of the four
choices of $N(v)$. All but two of the resulting graphs have a proper
nontoroidal minor. Both exceptions have $N(v)=\{c,d,e\}$ and
$B+ad+bd$ planar~\supp{E27}, contrary to hypothesis.

It remains to rule out a bridge on $\{b,c,d\}$ and the edge $cd$.
In either case we may retain $(\Psi_3\cup B)+ce$ as above. In the
first case also retain the bridge with its interior contracted to a
vertex adjacent to $b,c,d$. In the second retain $cd$. Call the
resulting minor $H_1$. In both cases $H_1-\operatorname{Int}(B)$
is nonplanar. In its torus embeddings, the sequences on $a,b,d$
along facial walks meeting all three are $abad,adbd,abd$~\supp{R11}. The first
two give the completions $B+ab+ad$ and $B+ad+bd$, respectively.
The triangle completion for $abd$ contains both. Thus
Lemma~\ref{cor:boundary-completions} gives
\begin{equation}\label{r1:two-cap-test}
H_1\text{ is toroidal}\quad\Longleftrightarrow\quad
B+ab+ad\text{ is planar or }B+ad+bd\text{ is planar}.
\end{equation}
Both completions are nonplanar, so $G=H_1$ by minimality. The same
test holds after replacing $B$ by any allowed minor. Hence no proper
allowed minor keeps both completions nonplanar.

If $B+ad$ were nonplanar, Corollary~\ref{common:reduced-forms}
would give $563$ forms with added edge $ad$ for each choice of
$H_1$. All but four of the resulting graphs have a $K_5$ minor or a
proper nontoroidal minor~\supp{E28}. The four exceptions occur when
$cd$ is an edge. In $\mathfrak T_{52}$ and $\mathfrak T_{57}$, two
twin triples meet in two vertices. In $\mathfrak T_{58}$, two are
disjoint. In $\mathfrak T_{59}$, some choice of the four twins gives
a bridge with four attachments.
Thus $B+ad$ is planar.

By Lemma~\ref{r1:three-attachment-setup}, the triangle completion
of $B$ is three-connected and has no $K_5$ minor. We may now apply
Lemma~\ref{core:two-cap-critical} with $(x,y,z)=(a,d,b)$. Since
$B+ab+bd$ is nonplanar, only the two four-cycle forms remain.
In both, $b$ is adjacent to opposite vertices $u,u'$ of the cycle,
and the other two vertices $p',q'$ are twins with neighborhood
$\{d,u,u'\}$ or $\{a,u,u'\}$. The first triple is disjoint from
$\{a,b,c\}$. In the second case, choose $p',q'$ and $r,z$ as the
four twins. The path $bpc$ lies outside the new $\Psi_3$, and its
component has neighbors $u,u'$ through $b$, $a$ through $p$, and
$e$ through $ce$. Its bridge has at least four attachments, again
a contradiction. By Lemma~\ref{r1:three-attachment-setup}, this
leaves only the attachment triples $\{a,b,d\},\{a,c,e\}$ and the
cross edges $bd,ce$.
\end{proof}

The next lemma gives the planar completions needed for Proposition~\ref{prop:three-root-blocks}.

\begin{lemma}\label{r1:some-planar-center}
The attachments of every bridge $B$ with nonempty interior can be
labeled $x,y,t$ so that $B+xy+xt$ is planar.
\end{lemma}
\begin{proof}
Suppose $B$ has no planar completion by two edges. If its attachments
are $\{a,b,c\}$, $\{a,d,e\}$ or $\{a,c,d\}$, take the disjoint
$b$--$d$ and $c$--$e$ paths from the proof of
Lemma~\ref{r1:three-attachment-setup}. Each has just one end in $B$
and avoids its interior. Contract each path to its end in $B$ and
retain only $B$ and the edges at $p,q,r,z$. This gives a proper minor
consisting of $B$ and four twins whose neighborhood is the attachment
set of $B$. The minor is toroidal, so Lemma~\ref{lem:four-twins}
gives a planar completion of $B$ by two edges, a contradiction.

Lemma~\ref{r1:bcd-two-edge} excludes $\{b,c,d\}$ and $\{c,d,e\}$.
Interchanging $b,e$ and $c,d$ together if necessary, we may therefore
assume $B$ has attachments $\{a,b,d\}$. By Lemma~\ref{r1:abd-others},
every bridge with nonempty interior has attachments $\{a,b,d\}$ or
$\{a,c,e\}$, and the only possible cross edges are $bd,ce$. Let
$B_{bd}$ be the union of the bridges on $\{a,b,d\}$ together with
$bd$ if present, and define $B_{ce}$ similarly. If $B_{ce}+ac+ae$
were nonplanar, replacing $ab,ad$ in $B+ab+ad$ by $apb,ard$ and
$ac,ae$ in $B_{ce}+ac+ae$ by $aqc,aze$ would give two nonplanar
subgraphs meeting only at $a$. Their union has genus at least two
by block additivity and omits $pc$, contrary to minimality.
Thus $B_{ce}+ac+ae$ is planar.

By~\eqref{r1:availability}, $B_{ce}-a$ contains a $c$--$e$ path.
Contract it to a vertex $f$ and retain only $B_{bd}$ and the edges
at $p,q,r,z$. Fix a torus embedding of this proper minor. The edges
at the four twins form a copy $H$ of $\Psi_2$ with twin triples
$\{a,b,f\}$ and $\{a,d,f\}$. Since $B$ joins $b,d$ and has no
planar completion by two edges, Lemma~\ref{r2:cellular-center}
forces a closed regular neighborhood of $H$ to have genus zero.
Its boundary sequences are $ab,bf,ad,df,af,af$~\supp{R4}.
Every bridge in $B_{bd}$ and the edge $bd$ if present meet both
$b,d$. As no boundary walk contains both, all lie in the unique
annulus given by Lemma~\ref{lem:torus-complement}. Its two boundary
walks contain $b$ and $d$, respectively, so the two $af$ walks
bound~disks.

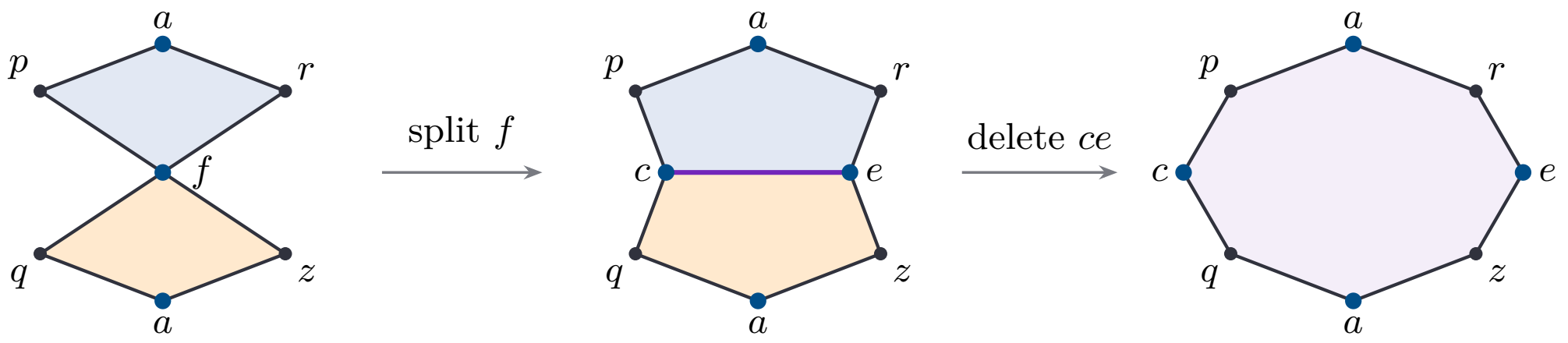
\begin{figure}[htbp]
\centering
\tikzsetnextfilename{psi3-split-f}%
\begin{tikzpicture}[
  x=1cm,y=1cm,font=\small,
  rim/.style={draw=figureInk,line width=.8pt,line join=round},
  vertex/.style={circle,fill=figureInk,inner sep=0pt,minimum size=3.2pt},
  attachment/.style={vertex,fill=figureNavy,minimum size=4pt},
  operation/.style={draw=figureInk!65,line width=.65pt,-{Stealth[length=4.5pt]}},
  every label/.append style={inner sep=1.4pt}
]
\foreach \stage/\offset in {0/0,1/5.05,2/10.1}{
  \begin{scope}[xshift=\offset cm]
    \coordinate (atop) at (0,1.09);
    \coordinate (abot) at (0,-1.09);
    \coordinate (p) at (-1.04,.69);
    \coordinate (r) at (1.04,.69);
    \coordinate (q) at (-1.04,-.69);
    \coordinate (z) at (1.04,-.69);
    \ifcase\stage
      \coordinate (c) at (0,0);
      \coordinate (e) at (0,0);
    \or
      \coordinate (c) at (-.78,0);
      \coordinate (e) at (.78,0);
    \or
      \coordinate (c) at (-1.44,0);
      \coordinate (e) at (1.44,0);
    \fi
    \ifnum\stage<2
      \fill[faceSimple] (c)--(p)--(atop)--(r)--(e)--cycle;
      \fill[faceRepeated] (c)--(q)--(abot)--(z)--(e)--cycle;
    \else
      \fill[figureBridge!8] (c)--(p)--(atop)--(r)--(e)--(z)--(abot)--(q)--cycle;
    \fi
    \draw[rim] (c)--(p)--(atop)--(r)--(e) (c)--(q)--(abot)--(z)--(e);
    \ifnum\stage=1
      \draw[figureBridge,line width=1.15pt] (c)--(e);
    \fi
    \foreach \v/\where in {p/above left,r/above right,q/below left,z/below right}
      \node[vertex,label=\where:{$\v$}] at (\v) {};
    \node[attachment,label=above:{$a$}] at (atop) {};
    \node[attachment,label=below:{$a$}] at (abot) {};
    \ifnum\stage=0
      \node[attachment,label={[label distance=3pt]right:$f$}] at (c) {};
    \else
      \node[attachment,label=left:{$c$}] at (c) {};
      \node[attachment,label=right:{$e$}] at (e) {};
    \fi
  \end{scope}
}
\draw[operation] (1.86,0)--(3.22,0) node[midway,above=2pt] {split $f$};
\draw[operation] (6.78,0)--(8.10,0) node[midway,above=2pt] {delete $ce$};
\end{tikzpicture}
\caption{The two $af$ disks in Lemma~\ref{r1:some-planar-center} become one disk with boundary sequence $caea$. The two appearances of $a$ in each drawing represent the same vertex.}
\label{fig:psi3-split-f}
\end{figure}
\looseness=-1
At $f$, the $bf$ and $df$ walks make $fp,fq$ and $fr,fz$ consecutive
pairs, with the two $af$ disks between them. Split $f$ into $c,e$
joined by an edge $ce$, keeping $p,q$ at $c$ and $r,z$ at $e$.
The two disks now share $ce$, so deleting it gives a disk with
boundary sequence $caea$ (Figure~\ref{fig:psi3-split-f}). This is the facial walk of the plane
path $cae$. Since $B_{ce}+ac+ae$ is planar, $B_{ce}$ can be drawn
in this disk as in the proof of Lemma~\ref{r1:missing-a-criterion}.
This gives a torus embedding of $G$, a contradiction.
\end{proof}

\subsection{Bounding the size and number of the remaining bridges}\label{r1:bounded-pieces}

By Lemma~\ref{r1:some-planar-center}, we may label the attachments of each bridge $B$ with nonempty interior as $x,y,t$ so that $B+xy+xt$ is planar. Then $B-x$ is a chain of edges and four-cycles by Proposition~\ref{prop:three-root-blocks}. We first consider $x\ne a$.

\begin{lemma}\label{r1:six-bridge-forms}
If $x\ne a$, then $B$ is the three-edge star or $B-x$ has one of
the following six forms, with end lengths as in
Lemma~\ref{r2:eight-bridge-forms}.
\begin{enumerate}[label=\textup{(\roman*)}]
\item One four-cycle with end lengths $(0,1),(0,2),(1,0)$ or $(2,0)$.
\item Two four-cycles with end lengths $(0,0)$, joined by a path
of length one or two.
\end{enumerate}
Each four-cycle contains exactly one end attachment, belonging to
$\{d,e\}$ if $x\in\{b,c\}$ and to $\{b,c\}$ if $x\in\{d,e\}$.
In each case $B$ has at most seven interior vertices.
\end{lemma}
\begin{proof}
If $B-x$ has no cycle, Proposition~\ref{prop:three-root-blocks}(v)
gives the three-edge star. On each four-cycle, the two twins have
neighborhood $\{x,u,v\}$, where $u,v$ are the other two vertices of
the cycle. One of $u,v$, say $u$, is interior to $B$, since otherwise both
twins would be isolated in the connected graph $B-\{x,y,t\}$.
Thus $\{x,u,v\}$ differs from both original twin triples and
meets each in exactly one vertex. If $x\in\{b,c\}$, the triple
already meets $\{a,b,c\}$ at $x$, so $v$ must be $d$ or $e$.
Similarly, if $x\in\{d,e\}$, then $v$ is $b$ or $c$. In either
case $v$ is the unique end attachment on the cycle.

Since $y,t$ each lie in only one block, there are at most two
four-cycles. With one four-cycle, one end length is zero and the
other is one or two, since no three consecutive blocks are edges.
With two four-cycles, both end lengths are zero. The cycles cannot
share a vertex, since their distinct twin triples would then share
both that vertex and $x$. Their joining path therefore has length
one or two. The largest form has two four-cycles joined by a path
of length two and has seven interior vertices.
\end{proof}

\begin{lemma}\label{r1:bridge-count}
At most two bridges have nonplanar triangle completion, and none has
attachment set $\{a,b,c\}$ or $\{a,d,e\}$. Among these bridges, the
choices $x=a$ and $x\ne a$ cannot occur together.
\begin{enumerate}[label=\textup{(\roman*)}]
\item If $x=a$ occurs, the available pairs in~\eqref{r1:availability}
are $bd,ce$. For each such bridge with attachments $a,y,t$, the edge
$yt$ is absent and no other component has both $y,t$ as neighbors.
\item If $x\ne a$ occurs, the chains $B-x$ of these bridges have at
most one four-cycle in total when the available pairs are $bd,ce$,
and at most two when they are $bd,cd,ce$.
\end{enumerate}
\end{lemma}
\begin{proof}
Let $B$ have nonplanar triangle completion. By
Proposition~\ref{prop:three-root-blocks}, $B-x$ contains a four-cycle,
since otherwise $B$ is the three-edge star.

First suppose $x=a$. On each four-cycle $up'vq'u$, the twins $p',q'$
have neighborhood $\{a,u,v\}$. At least one of $u,v$ is interior to
$B$, since otherwise both twins would be isolated in its interior.
Thus $\{a,u,v\}$ differs from both original twin triples and meets
each only at $a$. Hence both $u,v$ are interior. Every four-cycle
therefore avoids the end attachments $y,t$, so both have degree one
in $B$. Neither can serve as $x$, since any such choice must be
adjacent to the twins of a four-cycle.

Temporarily relabel so that $y=b$. Fix a four-cycle as above, with
the chain path from $b$ ending at $u$. Choose $p',q'$ and $r,z$ as
the four twins. The component outside the new $\Psi_3$ containing
$b,c,p,q$ has neighbors at $a$ through $p$ and at $u$ along the
chain path. Each of $d,e$ is joined to $b$ or $c$ by a cross edge
or a component by~\eqref{r1:availability}. If both were so joined
outside $B$, they would give two further attachments. Thus the
other end attachment of $B$ lies in $\{d,e\}$, say at $d$. Since
$e$ is joined outside $B$, no edge or component outside $B$ joins
$d$ to $b$ or $c$. Interchanging the two original twin pairs shows
likewise that no edge or component outside $B$ joins $b$ to $d$ or
$e$. Thus the available pairs form a matching, and $B$ alone supplies
the pair of its end attachments. This gives (i) and allows at most
two bridges with $x=a$.

Now suppose $x\ne a$. By Lemma~\ref{r1:six-bridge-forms}, each
four-cycle of $B-x$ gives a twin triple consisting of an available
pair and one interior vertex. Thus $B$ has attachments in both
$\{b,c\}$ and $\{d,e\}$, excluding both original triples as
attachment sets. Different cycles give different interior vertices,
since even the two cycles of a single chain are disjoint. Their twin
triples meet in exactly one vertex, so their available pairs also
meet in exactly one vertex. Hence no pair occurs twice, and at most
one of $bd,ce$ occurs. This gives (ii). Since each bridge under
consideration contributes a cycle, there are at most two.

Finally, the twins on a four-cycle of a bridge with $x=a$ have neighborhood
$\{a,u,v\}$, where $u$ and $v$ are interior to that bridge. A twin triple
from a four-cycle in another bridge with $x\ne a$ would contain none of
$a,u,v$. This contradicts the assumption that every two distinct twin
triples meet in exactly one vertex.
\end{proof}

It remains to determine the bridges with $x=a$.

\begin{lemma}\label{r1:a-chain-length}
Let $B$ be a bridge with attachments $a,y,t$. If its triangle
completion is nonplanar and $B+ay+at$ is planar, then $B=Q(a;y,t)$.
\end{lemma}
\begin{proof}
By Lemma~\ref{r1:bridge-count}, the available pairs are $bd,ce$.
Interchanging $b,c$ and $d,e$ simultaneously if necessary, suppose
$B$ has attachments $\{a,b,d\}$. Put $H_1=G-\operatorname{Int}(B)$.
There is no edge $bd$ and no other bridge has both $b,d$~as~attachments.

Suppose $H_1$ is planar. An additional bridge on $\{a,b,c\}$ or
$\{a,d,e\}$ would give a $K_{3,3}$ minor, since contracting its
interior gives a third vertex with the same neighborhood as the
corresponding twins. Thus the only additions to $\Psi_3$ in $H_1$
are bridges on $\{a,c,e\}$ and possibly $ce$. Let $B_{ce}$ be their
union. Shortening $apc$ and $are$ gives $B_{ce}+ac+ae$ as a planar
minor of $H_1$. Draw this graph and $B+ab+ad$ in disks meeting only
at $a$, with the four added edges in cyclic order $ab,ad,ac,ae$.
Corollary~\ref{r1:five-criterion} then gives a torus embedding of
$G$, a contradiction. Hence $H_1$ is nonplanar.

Suppose one of the shortenings in Lemma~\ref{r2:cycle-restoration}
gives $B'$. It preserves the three attachments and the connected
interior, so take a torus embedding of the proper minor $H_1\cup B'$.
Since $H_1$ is connected and nonplanar, every face of $H_1$ is a
disk by Lemma~\ref{lem:torus-complement}. The interior of $B'$ lies
in one such face. In the disk representing that face, the boundary
cycle meets $B'$ precisely in its attachment copies. Taking this
cycle as $H$ in Lemma~\ref{r2:cycle-restoration} restores $B$ with
the same copies, giving a torus embedding of $G$. Thus neither
shortening is possible.

By Proposition~\ref{prop:three-root-blocks} and the proof of
Lemma~\ref{r1:bridge-count}, the chain $B-a$ contains a four-cycle
and every four-cycle avoids $b,d$. Two four-cycles cannot share a
vertex, since their twin triples would then share both $a$ and that
vertex. As no three consecutive blocks are edges, two consecutive
four-cycles would be joined by a path of length one or two,
admitting shortening~\textup{(i)}. Thus there is just one four-cycle.
Each end path has length one or two, and length two would admit
shortening~\textup{(ii)}. Both end paths therefore have length one,
giving $B=Q(a;b,d)$.
\end{proof}

\subsection{Combining the remaining bridges}\label{r1:enumeration}

The preceding lemmas give the following construction.

\begin{construction}\label{r1:final-construction}
Start with $\Psi_3$ and choose $M=\{bd,ce\}$ or $M=\{bd,cd,ce\}$
as in~\eqref{r1:availability}. The additional attachment triples are
$\{a,b,d\}$ and $\{a,c,e\}$, together with $\{a,c,d\}$, $\{b,c,d\}$
and $\{c,d,e\}$ when $cd\in M$. Give the added bridges disjoint
interiors.
\begin{enumerate}[label=\textup{(\arabic*)}]
\item Choose at most two bridges with nonplanar triangle completion.
For each bridge $B$, choose an additional triple and label its
vertices $x,y,t$ in any order. If one chosen bridge has $x=a$, then
all do and $M=\{bd,ce\}$ by Lemma~\ref{r1:bridge-count}. In this
case use $B=Q(a;y,t)$ by Lemma~\ref{r1:a-chain-length}.

Otherwise, use one of the six forms for $B-x$ in
Lemma~\ref{r1:six-bridge-forms}, with end attachments $y,t$.
For each end attachment $u$ on a four-cycle, require $xu\in M$.
Join $x$ to each vertex of degree two in $B-x$ other than $y,t$,
as in Proposition~\ref{prop:three-root-blocks}. Allow at most one
four-cycle in total if $cd\notin M$, and at most two if $cd\in M$.

In either case, every two distinct triples among the two original
twin triples and those given by the four-cycles of the chosen
bridges must meet in exactly one vertex.

\item For each additional triple, choose whether to add a single
vertex with that neighborhood. Choose a subset of $M$ as cross
edges. These edges, together with all pairs $uv$ in the attachment
triples of the chosen bridges with $u\in\{b,c\}$ and $v\in\{d,e\}$,
must give exactly $M$. For each $Q(a;y,t)$ chosen in step~(1),
omit the edge $yt$ and allow no other bridge with both $y,t$ as
attachments.

\item Add zero, one or two further vertices, all with neighborhood
$\{a,b,c\}$ or all with neighborhood $\{a,d,e\}$.
\end{enumerate}
Keep one representative of each isomorphism class fixing
$a,b,c,d,e$ pointwise.
\end{construction}

\begin{proposition}\label{r1:three-outcome}
Suppose that for every choice of the four twins, every $\Psi_3$-bridge
with nonempty interior has exactly three attachments. Then $G$ is
isomorphic to one of $\mathfrak T_7,\mathfrak T_{33},\mathfrak T_{35},
\mathfrak T_{49},\mathfrak T_{70}$.
\end{proposition}
\begin{proof}
The preceding lemmas give all bridges with nonplanar triangle
completion in step~(1) of Construction~\ref{r1:final-construction}.
Every other bridge with nonempty interior is a three-edge star by
Lemma~\ref{lem:flat-replacement}. Two stars on the same additional
triple would give a twin triple meeting one of the original triples
in two vertices. Thus step~(2) includes all stars on the additional
triples and all cross edges allowed by~\eqref{r1:availability}.
Lemmas~\ref{r1:twin-count} and~\ref{r1:three-attachment-setup} give
$\min(n_c,n_e)=2$ and $n_c+n_e\le6$, exactly the five choices in
step~(3). Hence the construction includes $G$.

For $n_c=n_e=2$, the construction gives $31{,}063$ candidates.
Of these, $26{,}502$ have a degree-three vertex with two adjacent
neighbors, contrary to Lemma~\ref{lem:independent-neighborhood}(a).
In a further $4{,}368$, some choice of the four twins gives a bridge
with four or five attachments~\supp{E29}. In each exclusion the
recorded degree-three vertex or the four chosen twins lie outside
$\{a,b,c,d,e\}$, so adding vertices on the original triples preserves
their neighborhoods. In the second case, the component containing
the old component outside the same $\Psi_3$ still has all its
original attachments.

The $193$ remaining candidates give $965$ cases over the five choices
of $(n_c,n_e)$. Of these, $172$ are toroidal and $782$ have a proper
nontoroidal minor. The others are the five graphs in the
statement up to isomorphism~\supp{E29}.
\end{proof}

\begin{proof}[Proof of Theorem~\ref{thm:one-common}]
If some choice of the four twins gives a bridge with four or five
attachments, the conclusion follows from
Propositions~\ref{r1:five-outcome}, \ref{r1:missing-a-outcome}
and~\ref{r1:missing-b-outcome}. Otherwise, every choice of the four
twins gives only bridges with three attachments among those with
nonempty interior, so Proposition~\ref{r1:three-outcome} applies.
The converse follows from Proposition~\ref{prop:membership}.
\end{proof}

\section{Completion of the classification}\label{sec:computations}\label{sec:completion}

\begin{proposition}\label{prop:membership}
The graphs $\mathfrak T_1,\ldots,\mathfrak T_{74}$ are pairwise
nonisomorphic toroidal minor obstructions of genus two with no
$K_5$ minor. Each graph listed in Theorems~\ref{thm:disjoint},
\ref{thm:two-common} and~\ref{thm:one-common} is three-connected
and satisfies the twin-triple conditions of its theorem.
\end{proposition}
\begin{proof}
Both PAGE and \texttt{multi\_genus} computed $\gamma(\mathfrak T_i)=2$
for $1\le i\le74$, as described in Appendix~\ref{app:computations}.
The $74$ canonical forms in~\supp{O2} are distinct. For each graph
$G$, it also demonstrates that $G$ has no $K_5$ minor, gives its
connectivity and twin triples, and gives torus embeddings of $G-e$
and $G/e$ for every edge $e$.
\end{proof}

\begin{proof}[Proof of Theorem~\ref{thm:classification}]
Let $G$ be a toroidal minor obstruction with no $K_5$ minor.
As noted in the introduction, the disconnected case gives
$\mathfrak T_{25}$ and the cut-vertex case gives $\mathfrak T_8$.
The classification of Mohar and \v{S}koda~\cite{MoharSkoda2014}
and the comparison in~\supp{O3} give the $18$ listed graphs of
connectivity two.

Now suppose $G$ is three-connected. Unless $G\cong K_{3,7}=\mathfrak T_1$,
Theorem~\ref{thm:two-pages} gives two distinct twin triples. In this case, $G$ is one of the listed graphs by
Theorems~\ref{thm:disjoint}, \ref{thm:two-common} and~\ref{thm:one-common}. The converse follows from Proposition~\ref{prop:membership}.
\end{proof}

\section{Topological obstructions}\label{sec:topological}

A \emph{topological obstruction} for the torus is a nontoroidal graph
of minimum degree at least three for which every edge deletion is
toroidal. Every minor obstruction for the torus is a topological
obstruction.

A \emph{split} replaces a vertex by two adjacent vertices, dividing
its neighbors between them with at least two on each side. Every
topological obstruction can be obtained from a minor obstruction by
a sequence of splits~\mbox{\cite[Corollary~2.1.2]{BienstockLangston1995}}.

\begin{theorem}\label{thm:topological}
There are exactly $381$ topological obstructions for the torus with
no $K_5$ minor. They are listed in~\supp{TO1}, with the following
numbers of each order.
\[
\scaledtable{1.4}{$\displaystyle
\begin{array}{c|cccccccc}
n & 10 & 11 & 12 & 13 & 14 & 15 & 16 & 17\\ \hline
\text{number} & 1 & 9 & 43 & 111 & 151 & 58 & 7 & 1
\end{array}
$}
\]
\end{theorem}
\begin{proof}
Starting with the $74$ graphs of Theorem~\ref{thm:classification},
we considered all splits up to isomorphism, continuing from each
new topological obstruction with no $K_5$ minor. Every split of a
graph with a $K_5$ minor also has a $K_5$ minor. If $H-e$ is
nontoroidal, deleting the corresponding edge from a split of $H$
gives a graph that contracts to $H-e$ and is therefore also
nontoroidal.

The splits gave $7{,}006$ graphs up to isomorphism~\supp{TO2}. Of these, $4{,}884$
have a $K_5$ minor and a further $1{,}815$ have an edge deletion
containing some $\mathfrak T_i$ as a minor. The remaining $307$
have no $K_5$ minor and every edge deletion is toroidal. Contracting the split edge in each of these $307$ graphs
gives a nontoroidal proper minor. They are therefore distinct from the
$74$ minor obstructions, giving $381$ topological obstructions in all.
\end{proof}

\section{Further remarks}

Finding all the toroidal obstructions and proving that the list is complete seems difficult. Theorem~\ref{thm:classification} and the result of Gagarin, Myrvold and Chambers~\cite{GagarinMyrvoldChambers2009} settle the cases where one of the two forbidden minors for planarity is absent. By block additivity and the classification of the obstructions of connectivity two~\cite{MoharSkoda2014}, any toroidal obstruction still to be found must also be three-connected.

One class to consider next is the toroidal obstructions with no $K_5$ subdivision. Those with no $K_5$ minor are already on our list. For those that do contain a $K_5$ minor, we can no longer use Halin's theorem to find the three-cuts that gave us the twin triples.

One can also consider the problem of finding the obstructions for any closed connected surface $S$ other than the sphere. The proof of Theorem~\ref{thm:two-pages} shows that every three-connected minor obstruction for $S$ with no $K_5$ minor is either $K_{3,k}$ for some $k$ or has two distinct twin triples. The Wagner graph $V_8$ embeds in both the torus and the projective plane, so it embeds in $S$ by adding handles or crosscaps away from the graph, as even a noncellular embedding suffices to rule out $V_8$ as an obstruction when applying Halin's theorem. The disk replacement in Lemma~\ref{lem:flat-replacement} can be applied on any surface, and the neighbors of a degree-three vertex $v$ are still independent, since if two were adjacent, we could embed the graph with their edge deleted and then draw it back alongside the path through $v$.
    
\section*{Acknowledgements}
I would like to thank Dr.~Fran\c{c}ois Cl\'ement for advising me on the early work of this problem and for his helpful discussions, recommendations and feedback. I also thank Alexander Metzger for useful discussions and ideas, as well as his assistance with aspects of the computational work.

\section*{Funding}
This research received no funding.

\section*{Data Availability Statement}
The supplemental data cited in this paper, and the genus computations and programs described in Appendix~\ref{app:computations}, are available at \url{https://github.com/austinulrigg/toroidal-obstructions-no-k5}.

\section*{Conflict of Interest}
The author declares no conflicts of interest.
\appendix

\section{Computational details}\label{app:computations}

The supplemental data are in the folder \texttt{supplemental\_data} of the
\href{https://github.com/austinulrigg/toroidal-obstructions-no-k5}{repository}.
The labels in the text are also the names used in the repository. For example,
the rotation table cited as~\supp{R4} is \texttt{rotations/R4.json}.

O1 lists $\mathfrak T_1,\ldots,\mathfrak T_{74}$ in the order of
Table~\ref{tab:graph6}. The first character of each graph6 string gives
the number of vertices, from \texttt{I} for ten to \texttt{M} for fourteen.
Both PAGE~\cite{MetzgerUlrigg2026} and
\texttt{multi\_genus}~\cite{Brinkmann2022} computed
$\gamma(\mathfrak T_i)=2$ for $1\le i\le74$.
For each graph $G$, O2 gives a genus-two embedding of $G$, torus embeddings
of $G-e$ and $G/e$ for every edge $e$, and a certificate that $G$ has no
$K_5$ minor. It also gives the connectivity and twin triples of $G$.
The canonical forms computed with nauty~\cite{McKayPiperno2014} are
distinct.

\begin{table}[htbp]
\centering
\renewcommand{\encodingdefault}{T1}\small
\setlength{\tabcolsep}{4pt}
\begin{tabular}{@{}r@{\hspace{5pt}}l@{\hspace{14pt}}r@{\hspace{5pt}}l@{\hspace{14pt}}r@{\hspace{5pt}}l@{\hspace{14pt}}r@{\hspace{5pt}}l@{\hspace{14pt}}r@{\hspace{5pt}}l@{}}
\toprule
$i$ & graph6 & $i$ & graph6 & $i$ & graph6 & $i$ & graph6 & $i$ & graph6 \\
\midrule
1 & \verb!I???F~}~_! & 16 & \verb!K????Bort}Nw! & 31 & \verb!L?d_??B?v_QeYo! & 46 & \verb!L@G[C@?OE?C~E}! & 61 & \verb!L???EMeU@o?r}?! \\
2 & \verb!J?ACC@x^r}?! & 17 & \verb!KC`A@?_w?~b}! & 32 & \verb!L???WW?{DbMSrG! & 47 & \verb!L@GcC@?OE?B~F}! & 62 & \verb!Ms???????R_u^w^s?! \\
3 & \verb!J???wz_{FW?! & 18 & \verb!KW????^Huozw! & 33 & \verb!L???OCA{FWFktW! & 48 & \verb!LB?CC?KwB@@xFu! & 63 & \verb!MsG___K?q?G?BF?z?! \\
4 & \verb!JWr?__J@p`_! & 19 & \verb!K?@Hc?B@Vp~E! & 34 & \verb!L??`?_OqF_@|rw! & 49 & \verb!M?@HcA??_G^?Q]K]?! & 64 & \verb!Mr???????R_ufsZw?! \\
5 & \verb!J]o???fDo]?! & 20 & \verb!KAO`C?X@s{[E! & 35 & \verb!L????AL[cmLW~?! & 50 & \verb!M@G[CA?OA?G?H~E~?! & 65 & \verb!Mr????KB?E?WcNZo?! \\
6 & \verb!JBW?CMEw?{_! & 21 & \verb!KC`A@?X@tLNA! & 36 & \verb!L@G??B?@v?B|vo! & 51 & \verb!M?d_GO???B_U{B{{?! & 66 & \verb!ME?@@?_CWB@prArc?! \\
7 & \verb!J??EE?|Nvw?! & 22 & \verb!KSP???rBu@fx! & 37 & \verb!L?a??@bB?eno{S! & 52 & \verb!M?o????AWNEo|C{w?! & 67 & \verb!M??O`?_?]O]?Bwrc?! \\
8 & \verb!J??FF?[FFw?! & 23 & \verb!KSP???rBuEF`! & 38 & \verb!L?r@_??@WVEK{C! & 53 & \verb!M?d``a?_A?G?BN?~?! & 68 & \verb!M????cQo@W]?r_J[?! \\
9 & \verb!K????AN[vw^w! & 24 & \verb!L?o???r@oO}p|Q! & 39 & \verb!LI?GSA?o@o_|L{! & 54 & \verb!LsP@@?OC?]@e?{! & 69 & \verb!M?@b?oW_C??f?u}??! \\
10 & \verb!K????Bwfr]]w! & 25 & \verb!Ks\o?CA?_B_M! & 40 & \verb!L??E@_o?gYxa~?! & 55 & \verb!L????Fas`YBg~?! & 70 & \verb!MAP@@?__C?@fCufw?! \\
11 & \verb!K??@@B_v?~Zw! & 26 & \verb!L???@KY?^GzA~_! & 41 & \verb!L?KE?WC@F_EUo[! & 56 & \verb!LAO`?A@_V?BxFk! & 71 & \verb!MIK??KE_C?O?Q^Lf?! \\
12 & \verb!K???GGB}Fx^W! & 27 & \verb!L@J???B?v?BxuK! & 42 & \verb!MCaAA@?G@?C?F~F~?! & 57 & \verb!L??@e?K?W[]B~C! & 72 & \verb!M????B?Gow]?Bx~o?! \\
13 & \verb!K???GGC}Fv^[! & 28 & \verb!L?AR?OGAF_Q]X[! & 43 & \verb!L]o__O???J_u?{! & 58 & \verb!L@G?E?E@opZAxa! & 73 & \verb!M????B?Gow]?roNw?! \\
14 & \verb!K?CU@?BBVX]E! & 29 & \verb!L???B?RB^c\_f_! & 44 & \verb!L??EE?cE?Ff`~?! & 59 & \verb!L@G?E?E@op[Bv_! & 74 & \verb!M??????gYk^?fwZw?! \\
15 & \verb!K?K??BoppmXw! & 30 & \verb!L???GGww?wwX}o! & 45 & \verb!LI?CC?K?xwUARp! & 60 & \verb!L???C`_Bf_@|~_! &  &  \\
\bottomrule
\end{tabular}
\caption{The graphs $\mathfrak T_1,\ldots,\mathfrak T_{74}$ in graph6 format.}
\label{tab:graph6}
\end{table}

O4 compares the $74$ graphs with the $17{,}535$ minor obstructions of
Myrvold and Woodcock~\cite{MyrvoldWoodcock2018}. We used the file
\nolinkurl{all_minor_obs.txt.gz} on
\href{https://webhome.cs.uvic.ca/~wendym/torus/torus_obstructions.html}{Myrvold's
torus obstructions page}, revised August 16, 2017. O4 includes the
collection and gives an isomorphism from each of
$\mathfrak T_1,\ldots,\mathfrak T_{73}$ to one of its entries. No entry
has the canonical form of $\mathfrak T_{74}$.
O3 lists the $68$ pairwise nonisomorphic graphs of connectivity two in the
collection. By the classification of Mohar and \v{S}koda~\cite{MoharSkoda2014},
these are all the toroidal obstructions of connectivity two. For $50$ of
them, O3 gives branch sets for a $K_5$ minor. For each of the other $18$,
it gives an isomorphism to a graph in O1.

The tables R1--R12 give every rotation system of each graph considered,
the genus of the rotation, and the attachments in the order in which they
occur on each facial walk. An attachment is listed once for every appearance even when it may be repeated. R1--R3 also give the attachment sequences on the boundary of the
region along the specified path.

E1--E29 contain the graph and bridge enumerations. Each item gives the
construction, the choices used to form the candidates, the result for
each graph recorded, and the counts. E30 gives the genera
computed by PAGE for the eight choices of added edges on each of $T(z)$,
$S(x;y)$ and $Q(x;y,z)$ in Lemma~\ref{core:three-graph-completions}.

W1--W16 give individual embeddings, isomorphisms and minor models
used in the proofs. Here and in the
enumerations, an embedding is given by a rotation system. For an isomorphism we give the vertex
bijection, and for a minor we give the branch sets. For a cut we give the vertices whose deletion disconnects the
graph, and for nonplanarity we give a subdivision of $K_5$ or $K_{3,3}$.
To show that a graph has no $K_5$ minor, we exhibit it as a subgraph of a
graph formed from planar graphs and Wagner graphs by identifying
cliques of order at most three~\cite{Wagner1937}. 

When a graph is excluded by a hypothesis of the case, the enumeration
gives the reason, for example two twin triples meeting in two vertices.
For a proper nontoroidal minor, it gives an edge to delete or contract,
in the whole graph or in the bridge being reduced. The resulting graph contains some $\mathfrak T_i$ as a minor, fails the embedding tests for its case, or contains a subgraph of computed genus at least two. For the embedding tests,
the data give a subdivision showing that a required completion is
nonplanar, or a bridge that fits in none of the possible faces.
The enumerations use genus computations to show nontoroidality for
$311{,}739$ graphs up to isomorphism. 

In E3, E4 and E12, the graph for one combination is a subgraph of the
graph for another whenever the first uses at most as many of each bridge
form and edge as the second. We checked that every allowed
choice gives a subgraph of a recorded toroidal graph or contains a
recorded nontoroidal graph. In each case the propositional formula asks
whether there is an allowed choice with neither property.
The program \texttt{drat-trim}~\cite{WetzlerHeuleHunt2014} verified a
DRAT proof that the formula is unsatisfiable. E3, E4 and E12 record
$2{,}310$, $315{,}099$ and $9{,}878$ graphs, respectively.

TO1 lists the $381$ graphs of Theorem~\ref{thm:topological}. For each
of the $7{,}006$ graphs obtained by a split, TO2 gives the graph from
which it was split. For each of the $307$ further obstructions, it gives
a certificate that the graph has no $K_5$ minor and a torus embedding of
every edge deletion. For each of the other graphs, it gives either
branch sets for a $K_5$ minor or an edge whose deletion contains some
$\mathfrak T_i$ as a minor, together with the branch sets.

N1 lists the $120$ graphs mentioned in the introduction. None is
isomorphic to any of the $250{,}815$ graphs in Myrvold and Woodcock's
collection of topological obstructions, \texttt{all\_obs.txt.gz}.
For each graph $G$, N1 gives its
connectivity, a genus-two embedding, and a torus embedding of $G-e$ for
every edge $e$. For the $22$ minor obstructions, it also gives a torus
embedding of every $G/e$. For each of the other $98$, it gives an edge
whose contraction has genus two. The exact genus computations are
recorded for all $120$ graphs and for these $98$ contractions. N1 gives
branch sets for a $K_5$ minor in every graph except $\mathfrak T_{74}$.

The folder \texttt{genus} gives the results and running times for the
$74$ obstructions and the graphs from E1--E29.
The programs for the rotation tables, enumerations and splits are in
\texttt{programs}. The program in \texttt{checker} checks the witnesses
in the supplemental data.

\begin{samepage}
\section*{Declaration of AI Use}
The author used OpenAI's ChatGPT and Codex (Astra, Sol and earlier GPT models) and Anthropic's Claude (Opus 4.8 and later Opus models, and Fable) to write scripts, help analyze their outputs, organize the repository, and review related literature. Some exploratory computations helped the author develop ideas for definitions, conjectures and statements to prove. AI also sometimes suggested approaches to arguments and generated TikZ code from the author's instructions and, in some cases, the author's hand-drawn sketches. The author decided which suggestions to use, and developed all the final statements and proofs. The author takes full responsibility for the paper, including its arguments, computations, references and conclusions.
\par\end{samepage}

\end{document}